\documentclass[
a4paper,
11pt, 
DIV=14,
toc=bibliography, 
headsepline, 
parskip=half-,
abstract=true,
]{scrartcl}

\usepackage[english]{babel}
\usepackage[utf8]{inputenc}
\usepackage[T1]{fontenc}
\usepackage{lmodern}
\usepackage[babel=true]{microtype}
\usepackage{amsmath}
\usepackage{mathtools}
\usepackage{amssymb}
\usepackage{amsthm}
\usepackage{mathrsfs}
\usepackage{scrlayer-scrpage} 
\usepackage{enumitem}
\usepackage{tabularx}
\usepackage{booktabs} 
\usepackage{tikz} 
\usepackage{tikz-3dplot}   

\ihead{D.\;Platt, M.\;B.\;Schulz, A.\;Wood}
\ohead{\csname @title\endcsname}
\setkomafont{pageheadfoot}{\footnotesize}

\providecommand{\R}{\mathbb{R}}
\providecommand{\C}{\mathbb{C}}
\providecommand{\Z}{\mathbb{Z}}
\providecommand{\N}{\mathbb{N}}
\providecommand{\B}{\mathbb{B}} 
\providecommand{\Sp}{\mathbb{S}}  
\providecommand{\st}{\, :\ } 
\providecommand{\surf}[1][61]{\Upsilon_{#1}}
\providecommand{\Hopf}[1]{q_{#1}} 
\providecommand{\DN}{\Lambda}
\providecommand{\spec}{\mathrm{S}}
\providecommand{\ener}{\mathrm{E}}      
\providecommand{\area}{\mathrm{A}}

\DeclarePairedDelimiter\abs{\lvert}{\rvert}
\DeclarePairedDelimiter\nm{\lVert}{\rVert}
\DeclarePairedDelimiter\sk{\langle}{\rangle}
\DeclarePairedDelimiter\interval{]}{[}
\DeclarePairedDelimiter\Interval{[}{[}

\DeclareMathOperator{\divop}{div} 
\DeclareMathOperator{\Ogroup}{O} 
\DeclareMathOperator{\SOgroup}{SO}  
\DeclareMathOperator{\Conf}{Conf}
\DeclareMathOperator{\ind}{ind}
\DeclareMathOperator{\nul}{nul} 
\DeclareMathOperator{\II}{I\hspace*{-.5pt}I}

\usetikzlibrary{plotmarks} 
\pgfmathsetmacro{\picscale}{0.47}  
\pgfmathsetmacro{\figurescale}{\picscale*\textwidth/(2cm)}
\providecommand{\FBMS}[2][\picscale]{\draw(0,0,0)node[inner sep=0]{\includegraphics[width=#1\textwidth]{#2}};}
\pgfmathsetmacro{\thetaO}{72} 
\pgfmathsetmacro{\phiO}{90-45/2} 
\tdplotsetmaincoords{\thetaO}{\phiO}
\NewDocumentEnvironment{unitball}{O{\picscale} +b}{%
\begin{tikzpicture}[scale={#1*\textwidth/2cm},tdplot_main_coords,line cap=round,line join=round,semithick,baseline={(0,0,0)}]
\useasboundingbox[tdplot_screen_coords](-1.005,-1.005)rectangle(1.005,1.005);
\pgfmathsetmacro{\laengeA}{-135-45/2+\phiO}
\pgfmathsetmacro{\laengeB}{\laengeA+45}
\begin{scope}[black!30]
\path[tdplot_screen_coords,inner color=white,outer color=cyan!5!black!20](0,0)circle(1);
\foreach\laengengrad in {\laengeA,\laengeB,...,\phiO}{
\tdplotsetthetaplanecoords{\laengengrad}
\pgfmathsetmacro{\start}{-atan( 1/(tan(\thetaO)*( cos(\laengengrad)*sin(\phiO)-cos(\phiO)*sin(\laengengrad))))}
\tdplotdrawarc[tdplot_rotated_coords]{(0,0,0)}{1}{180+\start}{360+\start}{}{}
}
\foreach\breitengrad in {-60,-30,...,60}{
\pgfmathsetmacro{\breitenkreisradius}{cos(\breitengrad)}
\pgfmathsetmacro{\hoehe}{sin(\breitengrad)}\tdplotdrawarc{(0,0,\hoehe)}{\breitenkreisradius}{0}{360}{}{}
}
\end{scope}
#2%
\begin{scope}[black]
\foreach\laengengrad in {\laengeA,\laengeB,...,\phiO}{
\tdplotsetthetaplanecoords{\laengengrad}
\pgfmathsetmacro{\start}{-atan( 1/(tan(\thetaO)*( cos(\laengengrad)*sin(\phiO)-cos(\phiO)*sin(\laengengrad))))}
\tdplotdrawarc[tdplot_rotated_coords]{(0,0,0)}{1}{\start}{180+\start}{}{}
} 
\foreach\breitengrad in {-60,-30,...,60}{
\pgfmathsetmacro{\breitenkreisradius}{1*cos(\breitengrad)}
\pgfmathsetmacro{\hoehe}{1*sin(\breitengrad)}
\pgfmathsetmacro{\visiblerange}{acos(-sign(\breitengrad)*min(1,abs(tan(\breitengrad)*tan(90-\thetaO)))  )}
\ifdim\visiblerange pt>0pt
\tdplotdrawarc{(0,0,\hoehe)}{\breitenkreisradius}{\phiO-90+\visiblerange}{\phiO-90-\visiblerange}{}{}
\fi
}
\draw[tdplot_screen_coords] (0,0) circle (1);
\end{scope}%
\end{tikzpicture}%
}{}
\makeatletter
\tikzset{
scale plot marks/.is choice,
scale plot marks/true/.style={},	
scale plot marks/false/.code={
\def\pgfuseplotmark##1{\pgftransformresetnontranslations\csname pgf@plot@mark@##1\endcsname}
},
every mark/.append style={scale plot marks=false},
plus/.style={mark=+,mark size=2.25pt},
vdash/.style={mark=|,mark size=2.25pt},
hdash/.style={mark=-,mark size=2.25pt},
bullet/.style={mark=*,mark size=1.125pt},
}
\makeatother
\definecolor{FarbeA}{cmyk}{0,0.5,1,0}
\definecolor{FarbeB}{cmyk}{1,0.5,0,0}	

\usepackage[
pdftitle={Free boundary minimal surface in the unit ball with trivial symmetry group},
pdfauthor={Daniel Platt, Mario B. Schulz and Albert Wood}, 
pdfborder={0 0 0},
]{hyperref}

\theoremstyle{plain}
\newtheorem{theorem}{Theorem}[section]
\newtheorem{lemma}[theorem]{Lemma}
\newtheorem{corollary}[theorem]{Corollary} 
\newtheorem{proposition}[theorem]{Proposition}

\theoremstyle{definition}
\newtheorem{definition}[theorem]{Definition}
 
\theoremstyle{remark}
\newtheorem{remark}[theorem]{Remark}

\numberwithin{equation}{section} 

\title{Free boundary minimal surface in the unit ball with trivial symmetry group}
\author{Daniel Platt, Mario B. Schulz and Albert Wood}
\date{\vspace*{-2ex}}
 
\newcommand\printaddress{{
\setlength{\parindent}{17pt}
\medskip
\hfill November 2026
\par
{\scshape Daniel Platt}
\newline
Imperial College London, Department of Mathematics,
180 Queen's Gate, South Kensington, London SW7 2RH, United Kingdom
\newline
\textit{E-mail address:} \texttt{d.platt@imperial.ac.uk}
\par\medskip
{\scshape Mario B. Schulz}
\newline 
Università di Trento, 
Dipartimento di Matematica, 
via Sommarive 14, 
38123 Povo di Trento, 
Italy 
\newline
\textit{E-mail address:} 
\texttt{mario.schulz@unitn.it}
\par\medskip
{\scshape Albert Wood}
\newline 
King's College London,
Department of Mathematics,
Strand, London WC2R 2LS, United Kingdom
\newline
\textit{E-mail address:} 
\texttt{albert.wood@kcl.ac.uk}
}}

\begin{document}

\maketitle

\begin{abstract}
We develop a new method for constructing embedded free boundary minimal surfaces of genus zero in the Euclidean unit ball by rigorous numerical verification. 
The idea is to prove the existence of a free boundary minimal immersion near a numerically computed approximation. 
The method combines a boundary formulation of the variational problem with quantitative estimates for the Dirichlet-to-Neumann operator and a fix point argument.
Our construction does not rely on prescribed symmetries of the surface. 
As our main application, we obtain the first rigorous existence proof for an embedded free boundary minimal surface in the unit ball with trivial symmetry group.
It has genus zero, $61$ boundary components and Morse index equal to $180$, which we prove to be maximal among all embedded examples of this topological type. 
\end{abstract}

\tableofcontents

\clearpage
\section{Introduction}
 
Free boundary minimal surfaces in the three-dimensional Euclidean unit ball $\B^3$ are critical points of the area functional among all surfaces in $\B^3$ with boundary constrained to the unit sphere $\Sp^2=\partial\B^3$. 
Equivalently, they have vanishing mean curvature and meet $\partial\B^3$ orthogonally along their boundary. 
This setting has attracted intense interest since Fraser and Schoen \cite{FraserSchoen2011,FraserSchoen2016} discovered the connection with the optimisation problem for the first Steklov eigenvalue on surfaces with boundary. 
Examples of embedded free boundary minimal surfaces in $\B^3$ are challenging to construct because they are necessarily unstable, i.\,e.~saddle points of the area functional rather than minimisers. 
Several different methods have been developed to overcome this challenge, and we give a brief overview of three of them. 

\emph{Gluing methods} combine well-understood pieces, such as equatorial discs, catenoidal annuli and Karcher--Scherk towers, into an almost free boundary minimal surface. 
This initial surface is then perturbed to an exact solution, e.\,g.~by means of a fix point argument. 
Gluing methods have the advantage of providing direct control on the topology of the resulting surface. 
However, all gluing constructions in $\B^3$ have in common that they require the genus or the number of boundary components to be sufficiently large to ensure that the initial surface is close enough to an actual solution \cite{FolhaPacardZolotareva2017,KapouleasZou2021,KapouleasWiygul2023,CarlottoSchulzWiygul2025Stackings,KapouleasLi2021,CarlottoSchulzWiygul2025Nonuniqueness}; this is also used to control the kernel of the linearised operator. The symmetry group of the resulting surfaces typically grows with their genus or number of boundary components, resulting in examples with a high degree of symmetry.

\emph{Min-max theory} for the area functional, pioneered by Almgren \cite{Almgren1965} and Pitts \cite{Pitts1981}, can construct unstable embedded minimal surfaces as weak limits of min-max sequences for families of surfaces sweeping out the ambient manifold.  
In the Simon--Smith setting \cite{Smith1982,ColdingDeLellis2003}, adapted to the free boundary case by Li \cite{Li2015}, the topology of the limit surface can be controlled from above \cite{DeLellisPellandini2010,Ketover2019,FranzSchulz2026}.  
Controlling the topology from below typically relies on prescribing a finite group of ambient symmetries, see \cite{KetoverEquivariant}. 
This equivariant approach has proved the existence of free boundary minimal surfaces in $\B^3$ with connected boundary and arbitrary genus \cite{CarlottoFranzSchulz2022}, with genus zero and any number of boundary components \cite[§\,5]{FranzSchulz2026}, and with genus one and two boundary components \cite{FranzKetoverSchulz}, among others \cite{KetoverFBMS}.   

\emph{Steklov eigenvalue optimisation} is based on the fact that the ambient coordinate functions restricted to any free boundary minimal surface in $\B^3$ are Steklov eigenfunctions with eigenvalue one \cite{FraserSchoen2011}.
Conversely, Fraser and Schoen \cite{FraserSchoen2016} proved that if a smooth metric on a surface with boundary maximises the first (normalised) Steklov eigenvalue, then the corresponding first eigenfunctions define a free boundary minimal immersion into a Euclidean ball $\B^n$ of some dimension $n\in\N$. 
In the genus zero case, the resulting surfaces are actually embedded in $\B^3$ \cite[Proposition~8.1]{FraserSchoen2016}, and maximising metrics exist for any number $b\in\N$ of boundary components \cite[Corollary~1.11]{KarpukhinKusnerMcGrathStern}. 
Their behaviour as $b\to\infty$ has been studied in \cite{GirouardLagace2021,KarpukhinStern2024,McGrathZou2024}. 
However, the method provides no information on the symmetry group of the unconstrained maximisers. 
Karpukhin, Kusner, McGrath and Stern \cite{KarpukhinKusnerMcGrathStern} developed an equivariant version of this approach. 
By prescribing reflection symmetries, they proved the existence of embedded free boundary minimal surfaces in $\B^3$ of any orientable topological type and area less than $2\pi$.  

\begin{figure}
\tikz[baseline={(0,0,0)}]{\FBMS{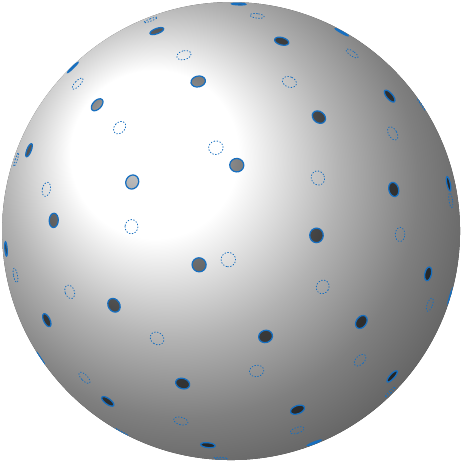}}
\hfill
\begin{unitball}
\FBMS{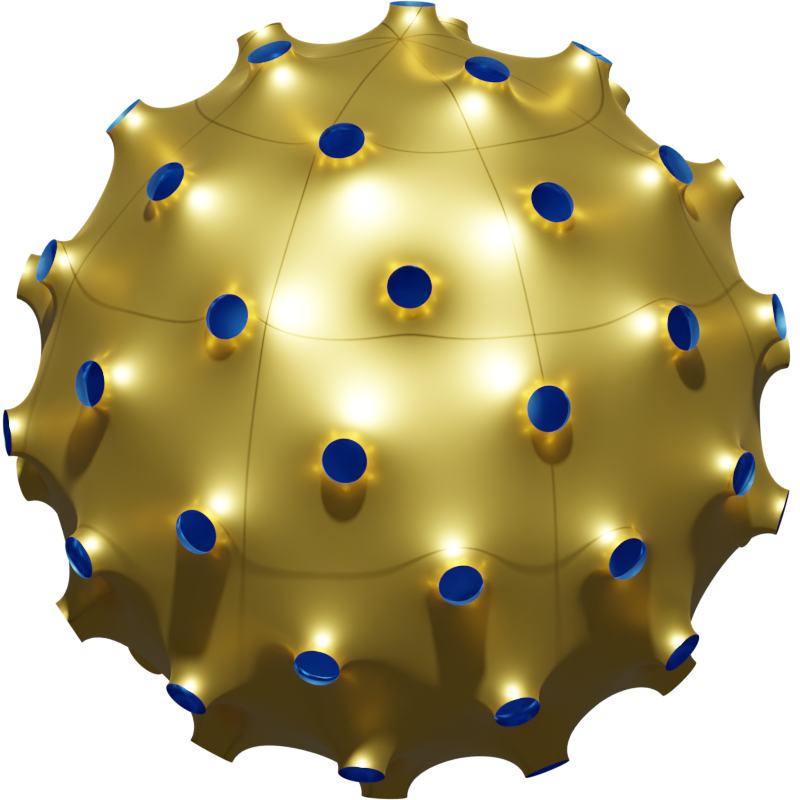}
\end{unitball}
\caption{The circle domain $\Omega_{c_0}\subset\Sp^2$ (left) and the corresponding numerical approximation of the free boundary minimal surface $\surf\subset\B^3$ (right).}
\label{fig:g0b61}
\end{figure}

To summarise, none of these methods yield a free boundary minimal surface in $\B^3$ which is known to have trivial symmetry group.  
In this article, we prove the existence of the first such example.
A numerical approximation of this surface is shown in Figure~\ref{fig:g0b61}.

\begin{theorem}
\label{thm:main}
The Euclidean unit ball $\B^3$ contains an embedded free boundary minimal surface $\surf$ with the following properties.
\begin{enumerate}[label={\normalfont(\roman*)},nosep]
\item\label{thm:main-topo} $\surf$ has genus zero and $61$ boundary components. 
\item\label{thm:main-sym} The symmetry group of $\surf$ is trivial. 
\item\label{thm:main-index} $\surf$ has Morse index equal to $180$, which is maximal among all embedded free boundary minimal surfaces in $\B^3$ satisfying \ref{thm:main-topo}, and nullity equal to $3$. 
\end{enumerate}
\end{theorem}

\begin{remark}
\label{rem:g0bN}
For comparison, combining \cite[Theorem 5.1]{FranzSchulz2026}, \cite[Remark~1.2]{FranzKetoverSchulz} 
and \cite[Theorem 1.1]{KarpukhinSchulz} or \cite{Tran} yields for each $3\leq b\in\N$ an embedded free boundary minimal surface $Z_b\subset\B^3$ of genus zero with $b$ boundary components, prismatic symmetry of order $4b$ and Morse index equal to $2b$. 
In particular, $Z_{61}$ and $\surf$ are non-congruent examples of the same topology. 
\end{remark}

Since the class of embedded free boundary minimal surfaces in $\B^3$ with genus zero and $b$ boundary components is nonempty for each $3\leq b\in\N$ by Remark~\ref{rem:g0bN}, the surfaces of largest area within this class are naturally distinguished objects of study \cite[Conjecture~7.10]{CarlottoSchulzWiygul2025Nonuniqueness}. 
By \cite[Theorem~1.1]{ChodoshGianocca}, every surface $\Sigma$ in this class has first Steklov eigenvalue equal to one, so that its first Steklov eigenvalue normalised by the boundary length equals $\operatorname{length}(\partial\Sigma)=2\operatorname{area}(\Sigma)$ (see \cite[Proposition~2.4]{Li2020}). 
Consequently, the area-maximising surfaces in this class are those induced by metrics which maximise the normalised first Steklov eigenvalue, and such metrics exist by \cite[Corollary~1.11]{KarpukhinKusnerMcGrathStern}, \cite[Proposition~8.1]{FraserSchoen2016}. 
We denote any such surface by $\Gamma_b^{\max}$.
For many values of $b$, numerical evidence \cite{OudetKaoOsting2021,SchulzGallery} (see also \cite[Open Question~5]{KarpukhinStern2024}) suggests that the boundary configuration of $\Gamma_b^{\max}$ resembles the minimum energy configuration of $b$ electrons constrained to the unit sphere -- a problem originally posed by Thomson \cite{Thomson1904} in 1904. 
For $b=61$, the best known numerical simulation of the Thomson problem yields a configuration with trivial symmetry group, and this inspired Theorem~\ref{thm:main}. 
In agreement with \cite[Conjecture~7.10,(iv)]{CarlottoSchulzWiygul2025Nonuniqueness}, we expect $\surf$ to be such a surface $\Gamma_{61}^{\max}$, although our construction involves neither Steklov eigenvalue nor area maximisation. 
In any case, Section~\ref{sec:asymmetry} provides rigorous bounds for the area of $\surf$, yielding 
$\operatorname{area}(\surf)\approx3.565685\,\pi$. 

We prove Theorem \ref{thm:main} using a boundary formulation of the variational problem.
For the classical Plateau problem of finding an area-minimising disc spanning a prescribed curve, such an approach goes back to the work of Douglas \cite{Douglas1931} and Struwe \cite{Struwe1988}, where the variational problem was reduced to finding a suitable parametrisation of the prescribed boundary curve.
In the free boundary setting, on the other hand, both the boundary and the conformal structure of the domain are unknown.
Building on the work of Da\;Lio and Pigati \cite{DaLioPigati2020}, we account for both types of variations by considering the Dirichlet energy of harmonic extensions of boundary maps. 
In particular, to vary the conformal structure we consider a \emph{circle domain} $\Omega_c\subset\Sp^2$ given by removing finitely many geodesic discs from the sphere and vary the disc configurations (see Figures~\ref{fig:g0b61}--\ref{fig:circle_variation}). 
Introducing local coordinates for the boundary maps and the circle configurations and eliminating the gauge freedoms, the problem of finding a free boundary minimal immersion $\Phi\colon\Omega_c\to\B^3$ is reduced to finding a zero of a smooth map between Banach spaces.
We prove the existence of such a zero near a numerically computed approximation by means of a quantitative fix point argument.
This \emph{rigorous numerical verification} requires bounds for the residual, the inverse of the linearised operator and the nonlinear remainder, which we derive analytically in terms of explicitly computable quantities.
  
There is an extensive literature on numerical methods for the computation of minimal surfaces. 
For an overview of this highly active research area, we refer to \cite{Polthier2005,DziukElliott2013,SapagovasEtAl2024}. 
While many of these methods rely on discretisations of surfaces  \cite{PinkallPolthier1993,Brakke1992}, others solve boundary formulations of the Plateau problem \cite{DziukHutchinson1996,Harbrecht2009,SakakibaraShimizu2024}. 
Numerical simulations of free boundary minimal surfaces based on Steklov eigenvalue optimisation are presented in \cite{OudetKaoOsting2021,ChenOsting2026}, and very recently in \cite{ChenOsting2026Zoo}, which focuses on the variety of examples of free boundary minimal surfaces of genus zero, without providing rigorous existence proofs for the new examples.
To compute our approximate solution, we start from an initial guess inspired by the Thomson problem and then use Newton iteration to optimise the circle configuration and the boundary map simultaneously (see Section~\ref{sec:initial}).

A central ingredient of our existence proof is the quantitative analysis of the Dirichlet-to-Neumann operator associated with a circle domain. 
By comparing it with the direct sum of the operators associated with the individual boundary circles and exploiting the regularising property of their difference, we derive explicit bounds for the operator and its derivatives with respect to the circle configuration.
For the numerical approximation of the Dirichlet-to-Neumann operator, we adapt the method of fundamental solutions used in \cite{HoangIzsakMaros2024}.
Rather than proving convergence of the numerical method, we derive explicit a posteriori error estimates 
(see \cite[Theorem~4.2]{BarrientosGaticaMaischak2002} for related work in the context of finite element approximations). 
These estimates also allow us to reduce the numerical verification of invertibility of the linearised operator to finite-dimensional computations.

Similarly to gluing methods, our existence proof relies on a perturbative argument near an approximate solution.
The use of computer assistance in such existence proofs in geometric analysis is relatively recent.  
An early example is \cite{ReitererTrubowitz2019}, followed by \cite{ButtsworthHodgkinson2026} and several follow-ups using ODE techniques, as well as more recent applications involving PDEs, including \cite{Platt2026}.  
We also highlight the related work \cite{CharltonHellerHellerTraizet2024}, which uses a numerically assisted fix point approach to prove existence of a deformation from a degenerate limit.
Our construction relies neither on a degenerate limit nor on prescribed symmetries of the surface. 
We emphasise that our variational formulation applies to free boundary minimal surfaces of genus zero irrespective of whether they arise from Steklov eigenvalue optimisation. 
In particular, our framework can be applied to investigate the existence of free boundary minimal surfaces of genus zero near other numerical approximations, including those obtained in \cite{SchulzGallery,ChenOsting2026Zoo}. 
The source code for the numerical verification in this article is available at \cite{PlattSchulzWoodCode}. 
 
Our argument for Theorem~\ref{thm:main}\;\ref{thm:main-index} also yields the following result, which is of independent interest as it provides general sharp upper bounds for the Morse index and nullity of free boundary minimal surfaces of genus zero.
The Morse index bound was previously noted in \cite[Remark~1.6 and Lemma~5.4]{Tran}.
We complement it with a bound for the sum of index and nullity and conclude that any example attaining the index bound is necessarily nondegenerate up to ambient rotations.
Here, $\ind_\area(\Sigma)$ denotes the Morse index of $\Sigma$, defined as the maximal dimension of a space of normal variations on which the second variation of \emph{area} is negative definite, and $\nul_\area(\Sigma)$ denotes the dimension of the kernel of the second variation. 
Our notation distinguishes it from the \emph{energy} index and nullity used in the proof in Section~\ref{sec:upper_index_bound}.  

\begin{theorem}
\label{thm:index-upper-bound}
Let $\Sigma\subset\B^3$ be a properly embedded free boundary minimal surface of genus zero with $3\leq b\in\N$ boundary components.
Then
\begin{align*}
\ind_\area(\Sigma)&\leq3b-3, \\
\ind_\area(\Sigma)+\nul_\area(\Sigma)&\leq3b.
\end{align*}
Moreover, if $\ind_\area(\Sigma)=3b-3$, then $\nul_\area(\Sigma)=3$.
\end{theorem}

\begin{remark}
Theorem~\ref{thm:main}\,\ref{thm:main-index} provides an example with $b=61$ for which Theorem~\ref{thm:index-upper-bound} is sharp. 
Another such example is the free boundary minimal trinoid with $b=3$, see Remark~\ref{rem:g0bN}. 
\end{remark}

\paragraph{Acknowledgments.}
This project has received funding from the European Research Council (ERC) under the European Union's Horizon 2020 research and innovation programme (grant agreement No.~947923). 
DP was supported by an Imperial College Research Fellowship, partly funded by Schmidt Futures.
AW was partially supported by research grants from the Research Grants Council of the Hong Kong Special Administrative Region, China [Project No.: CUHK 14304121, CUHK 14305122 and CUHK 14307123].
DP and MS thank the Chinese University of Hong Kong for their hospitality during a research visit.

\section{Conformal free boundary harmonic maps and the reduced energy}

Let $(M,g)$ be a compact, connected, oriented Riemannian surface with boundary, and let $\eta$ denote the outward unit conormal along $\partial M$. 
Given a smooth map $\Phi\colon M\to\R^3$, its \emph{Dirichlet energy} 
\begin{align}
\label{eqn:energy}
E_g(\Phi)\vcentcolon=\frac{1}{2}\int_{(M,g)}\abs{d\Phi}^2_g
\end{align}
depends only on the conformal class $[g]$ of $g$, since $M$ is two-dimensional. 
A map $\Phi\colon M\to\B^3$ satisfying $\Phi(\partial M)\subset\Sp^2$ is called \emph{free boundary harmonic} if it is a critical point of \eqref{eqn:energy} with respect to variations preserving the boundary constraint \cite[Definition~2.2]{KarpukhinSchulz}. 
The corresponding Euler--Lagrange equation is
\begin{align}\label{eqn:fbhm_eq}
\left\{\begin{aligned}
-\Delta_g \Phi &= 0 & &\text{ in $M$,} \\
\partial_\eta \Phi &= \abs{\partial_\eta\Phi}\Phi & &\text{ on $\partial M$.}
\end{aligned}\right.
\end{align}
If $\Phi\colon M\to\B^3$ is nonconstant, free boundary harmonic and weakly conformal in the sense that $\Phi^*g_{\R^3}=ug$ for some function $u\colon M\to\Interval{0,\infty}$, where $g_{\R^3}$ denotes the Euclidean metric, then $u=\tfrac{1}{2}\abs{d\Phi}_g^2$, and $\Phi$ is a branched free boundary minimal immersion. 
Conversely, the inclusion map $\Phi\colon\Sigma\to\B^3$ of any embedded free boundary minimal surface $\Sigma\subset\B^3$ is conformal free boundary harmonic, and satisfies $\abs{\partial_\eta\Phi}=1$, where $g=\Phi^*g_{\R^3}$ is the induced metric.

\paragraph{Two-variable Energy.}
In our construction, the conformal structure of $M$ is to be determined along with the map. 
Therefore, we follow \cite[equation~(8)]{Moore2007} and \cite[§\,3.1]{KarpukhinSchulz} in regarding the Dirichlet energy \eqref{eqn:energy} as a functional 
\begin{align}
\label{eqn:2variable_energy}
\widetilde{E}(\Phi,g)\vcentcolon=E_g(\Phi)
\end{align}
of the map and the metric. 
We write $\mathcal{S}(M)$ for the space of smooth symmetric bilinear forms on $M$ and $\mathcal{V}_\Phi$ for the space of maps $V\in C^\infty(M,\R^3)$ satisfying $\sk{V,\Phi}=0$ on $\partial M$, i.\,e.~variations of $\Phi$ preserving the boundary constraint. 
Differentiating the integrand in \eqref{eqn:energy} also with respect to $g$, we obtain for all $(V,h)\in\mathcal{V}_\Phi\times\mathcal{S}(M)$ 
\begin{align}\label{eqn:first-variation}
d\widetilde{E}(V,h)
&=dE_g(V)-\frac{1}{2}\int_{(M,g)}
\sk[\big]{\Hopf{\Phi},h}_g, & 
\Hopf{\Phi}&\vcentcolon=\Phi^*g_{\R^3}-\frac{1}{2}\abs{d\Phi}_g^2\,g.
\end{align}
The trace-free bilinear form $\Hopf{\Phi}$ defined in \eqref{eqn:first-variation} is the \emph{Hopf differential} of $\Phi$. 
By definition, $\Hopf{\Phi}\equiv0$ if and only if $\Phi$ is weakly conformal with respect to $g$. 
Hence, if $\Phi\colon M\to\B^3$ is nonconstant with $\Phi(\partial M)\subset\Sp^2$, then $(\Phi,g)$ is a critical point of $\widetilde{E}$ if and only if $\Phi$ is a branched free boundary minimal immersion which is weakly conformal with respect to $g$. 

If $\Phi\colon M\to\R^3$ is harmonic, then its Hopf differential is divergence-free. 
Indeed, in local coordinates and using the symmetry of the Hessian of $\Phi$, we have  
\begin{align*}
\bigl(\divop_g\Hopf{\Phi}\bigr)_j
&=\nabla^i\Bigl(\sk{\nabla_i\Phi,\nabla_j\Phi}-\tfrac{1}{2}\sk{\nabla_k\Phi,\nabla^k\Phi}\,g_{ij}\Bigr)
\\
&=\sk{\nabla^i\nabla_i\Phi,\nabla_j\Phi}+\sk{\nabla_i\Phi,\nabla^i\nabla_j\Phi}-\sk{\nabla_j\nabla_k\Phi,\nabla^k\Phi}
=\sk{\Delta_g\Phi,\nabla_j\Phi}=0. 
\end{align*}
Consequently, for any smooth vector field $Z$ on $M$, integration by parts yields
\begin{align}\label{eqn:hopf-boundary-general}
\frac{1}{2}\int_{(M,g)}\sk{\Hopf{\Phi},\mathcal{L}_Zg}_g
=\int_{\partial M}\Hopf{\Phi}(\eta,Z)-\int_{(M,g)}(\divop_g\Hopf{\Phi})(Z)
=\int_{\partial M}\Hopf{\Phi}(\eta,Z),
\end{align}
where $\mathcal{L}_Zg$ denotes the Lie derivative of $g$ along $Z$. 
By approximation, \eqref{eqn:hopf-boundary-general} remains valid for harmonic maps $\Phi\in H^2(M,\R^3)$ which are not necessarily smooth up to the boundary. 
If $\Phi\colon M\to\B^3$ is free boundary harmonic, then, given any vector $\tau$ tangent to $\partial M$, we have $g(\eta,\tau)=0$, and \eqref{eqn:fbhm_eq} together with $\abs{\Phi}=1$ on $\partial M$ implies 
\[
\Hopf{\Phi}(\eta,\tau)=\sk{\partial_\eta\Phi,\partial_\tau\Phi}=\abs{\partial_\eta\Phi}\sk{\Phi,\partial_\tau\Phi}=\abs{\partial_\eta\Phi}\tfrac{1}{2}\partial_\tau\abs{\Phi}^2=0.
\]
In this case, \eqref{eqn:hopf-boundary-general} simplifies to 
\begin{align}\label{eqn:hopf-boundary}
\frac{1}{2}\int_{(M,g)}\sk{\Hopf{\Phi},\mathcal{L}_Zg}_g
=\int_{\partial M}g(Z,\eta)\,\Hopf{\Phi}(\eta,\eta). 
\end{align}
In particular, the variation of $\widetilde{E}$ at $(\Phi,g)$ in the direction of $(0,\mathcal{L}_Zg)$ depends only on the normal component of $Z$ along $\partial M$, and it vanishes if $Z$ is tangent to $\partial M$.

\begin{figure}
\centering
\pgfmathsetmacro{\thetaO}{72}
\pgfmathsetmacro{\phiO}{90}
\tdplotsetmaincoords{\thetaO}{\phiO}
\def\caps{
70/25/35,
66/-82/40,
155/-10/33,
50/131/30,
115/200/28
}
\pgfmathsetmacro{\blickX}{sin(\thetaO)*sin(\phiO)}
\pgfmathsetmacro{\blickY}{-sin(\thetaO)*cos(\phiO)}
\pgfmathsetmacro{\blickZ}{cos(\thetaO)}
\begin{tikzpicture}[scale=\figurescale,line cap=round,tdplot_main_coords,
boundarystyle/.style={thick,FarbeB}
]
\def\setcapdata#1#2#3{%
\pgfmathsetmacro{\thetaK}{#1}%
\pgfmathsetmacro{\phiK}{#2}%
\pgfmathsetmacro{\rK}{#3}%
\tdplotsetrotatedcoords{\phiK}{\thetaK}{0}%
\tdplottransformmainrot{\blickX}{\blickY}{\blickZ}%
\pgfmathsetmacro{\nx}{\tdplotresx}%
\pgfmathsetmacro{\ny}{\tdplotresy}%
\pgfmathsetmacro{\nz}{\tdplotresz}%
\pgfmathsetmacro{\sichtB}{sqrt(\nx*\nx+\ny*\ny)}%
\pgfmathsetmacro{\azimut}{atan2(\ny,\nx)}%
\pgfmathsetmacro{\phiOrot}{atan2(\nx,-\ny)}%
\pgfmathsetmacro{\breitenkreisradius}{sin(\rK)}%
\pgfmathsetmacro{\hoehe}{cos(\rK)}%
\pgfmathsetmacro{\visiblerange}{acos(max(-1,min(1,-\nz*cos(\rK)/max(.001,\sichtB*sin(\rK)))))}%
\pgfmathsetmacro{\rimrange}{acos(max(-1,min(1,cos(\rK)/max(.001,\sichtB))))}%
\pgfmathsetmacro{\mindot}{\nz*cos(\rK)-\sichtB*sin(\rK)}%
\pgfmathsetmacro{\maxdot}{\nz*cos(\rK)+\sichtB*sin(\rK)}%
}%
\def\capfullcontour{%
\foreach \i in {0,...,120}{%
\pgfmathsetmacro{\t}{3*\i}%
\tdplottransformrotmain{\breitenkreisradius*cos(\t)}{\breitenkreisradius*sin(\t)}{\hoehe}%
\ifnum\i=0
\pgfpathmoveto{\pgfpointxyz{\tdplotresx}{\tdplotresy}{\tdplotresz}}%
\else
\pgfpathlineto{\pgfpointxyz{\tdplotresx}{\tdplotresy}{\tdplotresz}}%
\fi
}%
\pgfpathclose
}%
\def\capfrontcontour{%
\foreach \i in {0,...,80}{%
\pgfmathsetmacro{\t}{\azimut-\visiblerange+\visiblerange*\i/40}%
\tdplottransformrotmain
{\breitenkreisradius*cos(\t)}
{\breitenkreisradius*sin(\t)}
{\hoehe}%
\ifnum\i=0
\pgfpathmoveto{\pgfpointxyz{\tdplotresx}{\tdplotresy}{\tdplotresz}}%
\else
\pgfpathlineto{\pgfpointxyz{\tdplotresx}{\tdplotresy}{\tdplotresz}}%
\fi
}%
\capsilhouette
\pgfpathclose
}%
\def\capsilhouette{%
\foreach \i in {0,...,80}{%
\pgfmathsetmacro{\s}{-\rimrange+\rimrange*\i/40}%
\tdplottransformrotmain
{(-\nz*\nx*cos(\s)+\ny*sin(\s))/\sichtB}
{(-\nz*\ny*cos(\s)-\nx*sin(\s))/\sichtB}
{\sichtB*cos(\s)}%
\pgfpathlineto{\pgfpointxyz{1.005*\tdplotresx}{1.005*\tdplotresy}{1.005*\tdplotresz}}%
}%
}%
\def\cliphemisphere{%
\pgfpathmoveto{\pgfpointanchor{clipSW}{center}}%
\pgfpathlineto{\pgfpointanchor{clipSE}{center}}%
\pgfpathlineto{\pgfpointanchor{clipNE}{center}}%
\pgfpathlineto{\pgfpointanchor{clipNW}{center}}%
\pgfpathclose
\foreach \thetaCap/\phiCap/\rCap in \caps {%
\setcapdata{\thetaCap}{\phiCap}{\rCap}%
\ifdim\mindot pt>0pt
\capfullcontour
\else
\ifdim\maxdot pt>0pt
\capfrontcontour
\fi
\fi
}%
\pgfusepath{clip}%
} 
\def\shadehemisphere#1{%
\begin{scope}[tdplot_screen_coords]
\clip (0,0) circle (1);
\begin{scope}[even odd rule,tdplot_main_coords]
\cliphemisphere
\path[tdplot_screen_coords,ball color=#1] (0,0) circle (1.02);
\end{scope}
\end{scope}
}
\path[tdplot_screen_coords]
(-1.05,-1.05) coordinate (clipSW)
(1.05,-1.05) coordinate (clipSE)
(1.05,1.05) coordinate (clipNE)
(-1.05,1.05) coordinate (clipNW);
{%
\pgfmathsetmacro{\blickX}{-\blickX}%
\pgfmathsetmacro{\blickY}{-\blickY}%
\pgfmathsetmacro{\blickZ}{-\blickZ}%
\shadehemisphere{black!33}%
}
\begin{scope}[tdplot_main_coords,boundarystyle]
\foreach \thetaCap/\phiCap/\rCap in \caps {
\setcapdata{\thetaCap}{\phiCap}{\rCap}
\ifdim\mindot pt<0pt
\tdplotdrawarc[tdplot_rotated_coords]
{(0,0,\hoehe)}{\breitenkreisradius}
{\phiOrot-90-\visiblerange}
{\phiOrot-450+\visiblerange}{}{}
\fi
} 
\end{scope}
\shadehemisphere{white}
\begin{scope}[even odd rule,tdplot_main_coords]
\cliphemisphere
\begin{scope}[dotted,boundarystyle]
\foreach \thetaCap/\phiCap/\rCap in \caps {
\setcapdata{\thetaCap}{\phiCap}{\rCap}
\ifdim\mindot pt<0pt
\tdplotdrawarc[tdplot_rotated_coords]
{(0,0,\hoehe)}{\breitenkreisradius}
{\phiOrot-90-\visiblerange}
{\phiOrot-450+\visiblerange}{}{}
\fi
} 
\end{scope}
\end{scope}
\begin{scope}[tdplot_main_coords,boundarystyle]
\foreach \thetaCap/\phiCap/\rCap in \caps {
\setcapdata{\thetaCap}{\phiCap}{\rCap}
\ifdim\maxdot pt>0pt
\tdplotdrawarc[tdplot_rotated_coords]
{(0,0,\hoehe)}{\breitenkreisradius}
{\phiOrot-90+\visiblerange}
{\phiOrot-90-\visiblerange}{}{}
\fi
}
\end{scope}
\foreach \thetaCap/\phiCap/\rCap [count=\k] in \caps {
\tdplotsetcoord{p\k}{1}{\thetaCap}{\phiCap}
\expandafter\xdef\csname capdata\k\endcsname{%
\noexpand\setcapdata{\thetaCap}{\phiCap}{\rCap}%
} 
}
\def\usecap#1{\csname capdata#1\endcsname}
\begin{scope}[FarbeA,semithick]
\usecap{1}
\path[tdplot_rotated_coords,FarbeB]
({\breitenkreisradius*cos(203)},{\breitenkreisradius*sin(203)},\hoehe)
node[above,rotate=30,inner sep=1.5pt] {$\partial_1\Omega_c$};
\tdplotsetrotatedthetaplanecoords{0}
\tdplotdrawarc[tdplot_rotated_coords,latex-latex]{(0,0,0)}{1}{0}{\rK}{anchor=west}{$r_1$} 
\usecap{2}
\pgfmathsetmacro{\winkelT}{180+60}
\tdplotdrawarc[tdplot_rotated_coords,black,very thin]{(p2)}{0.1}{\winkelT}{\winkelT+90}{}{} 
\draw[tdplot_rotated_coords,black,thin,-stealth](0,0,1)--++(0+\winkelT:1)node[below right]{$e_{1,2}$};
\draw[tdplot_rotated_coords,black,thin,-stealth](0,0,1)--++(90+\winkelT:1)node[left]{$e_{2,2}$};
\tdplotsetrotatedthetaplanecoords{\winkelT-90}
\tdplotdrawarc[tdplot_rotated_coords,FarbeA,latex-latex]{(0,0,0)}{1}{0}{\rK}{anchor=120}{$r_2$}
\usecap{3}
\path[tdplot_rotated_coords,FarbeB]
({\breitenkreisradius*cos(150)},{\breitenkreisradius*sin(150)},\hoehe)
node[above,rotate=-8,inner sep=1.5pt] {$\partial_3\Omega_c$};
\tdplotsetrotatedthetaplanecoords{-90}
\tdplotdrawarc[tdplot_rotated_coords,FarbeA,latex-latex]{(0,0,0)}{1}{0}{\rK}{anchor=north}{$r_3$}
\end{scope} 
\draw plot[bullet] coordinates {(p1)};
\draw plot[bullet] coordinates {(p2)};
\draw plot[bullet] coordinates {(p3)};
\node[above] at (p1) {$p_1$};
\node[above left] at (p2) {$p_2$};
\node[below] at (p3) {$p_3$};
\node[tdplot_screen_coords] at (-135:0.5) {$\Omega_c$};
\end{tikzpicture}
\caption{A circle domain $\Omega_c\subset\Sp^2$ with configuration $c=(p_1,\ldots,p_b;r_1,\ldots,r_b)$ in the case $b=5$.}
\label{fig:circle_domain}
\end{figure}
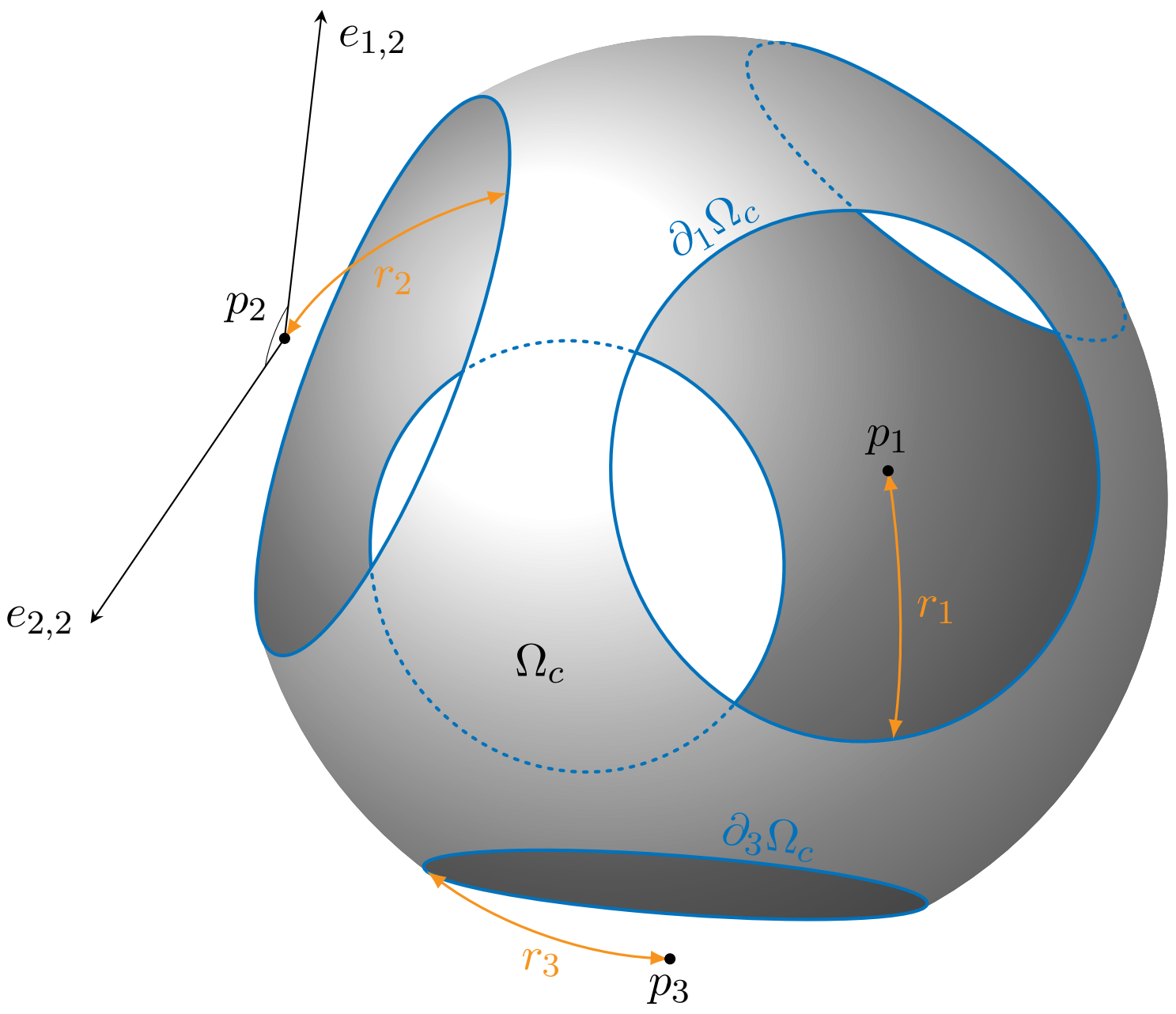

\paragraph{Circle domains.}
In the following, we assume that the Riemannian surface $(M,g)$ has genus zero and $b\in\N$ boundary components. 
Given points $p_1,\ldots,p_b\in\Sp^2$ and geodesic radii $r_1,\ldots,r_b\in\interval{0,\pi}$ such that the closures of the open geodesic discs $B_{r_k}(p_k)\subset\Sp^2$ for $k\in\{1,\ldots,b\}$ are pairwise disjoint, let
\begin{align*} 
\Omega_c\vcentcolon=\Sp^2\setminus\bigcup_{k=1}^bB_{r_k}(p_k)
\end{align*}
be the \emph{circle domain} with configuration $c=(p_1,\ldots,p_b;r_1,\ldots,r_b)$ as visualised in Figure~\ref{fig:circle_domain}.  
The set $\mathcal{C}$ of all such configurations is an open subset of $(\Sp^2)^b\times\interval{0,\pi}^b$ of dimension $3b$. 
Given $k\in\{1,\ldots,b\}$ we denote the $k$-th boundary circle of $\Omega_c$ by 
\begin{align}
\label{eqn:k-th_boundary_circle}
\partial_k\Omega_c\vcentcolon=\partial B_{r_k}(p_k)
=\{q\in\Sp^2\mid \sk{q,p_k}=\cos r_k\}.
\end{align}
Given an orthonormal basis $(e_{1,k},e_{2,k})$ of $T_{p_k}\Sp^2$, we parametrise $\partial_k\Omega_c$ by 
\begin{align}
\label{eqn:20261009b}
\begin{aligned}
\psi_{c,k}\colon 
\Interval{0,2\pi}&\to\partial_k\Omega_c \\ 
\theta&\mapsto p_k\cos r_k+(e_{1,k}\cos\theta+e_{2,k}\sin\theta)\sin r_k, 
\end{aligned}
\end{align}
which allows us to identify functions on $\partial_k\Omega_c$ with $2\pi$-periodic functions of $\theta$.

By Koebe's theorem (see \cite[§\,V.6]{Goluzin1969}), there exists $c\in\mathcal{C}$ such that $(M,g)$ is conformally equivalent to $(\Omega_c,g_{\Sp^2})$, where $g_{\Sp^2}$ denotes the restriction of the round metric to $\Omega_c$. 
Moreover, any conformal equivalence between two circle domains is the restriction of a Möbius transformation of $\Sp^2$. 
The group $\Conf(\Sp^2)$ of orientation-preserving Möbius transformations acts on $\mathcal C$ by $\phi\cdot c=c'$, where $c'=((p'_k,r'_k))_{k=1}^b$ is uniquely determined by
\begin{align}
\label{eqn:Conf-action}
\phi\bigl(B_{r_k}(p_k)\bigr)=B_{r'_k}(p'_k)
\end{align}
for all $k\in\{1,\ldots,b\}$. 
This action is free provided that $b\geq3$. 
Indeed, if $\phi\in\Conf(\Sp^2)$ fixes $c\in\mathcal{C}$, then $\phi$ maps each of the $b$ disjoint closed discs to itself and thus has a fix point in each of them by Brouwer's theorem \cite[Chapter~1, Theorem~19.3]{Taylor2011}, while a Möbius transformation with three distinct fix points is the identity. 
Conformally identifying $(M,g)$ with $(\Omega_c,g_{\Sp^2})$, we choose for all configurations $c'\in\mathcal{C}$ near $c$ diffeomorphisms 
\begin{align}
\label{eqn:diffeos_f_c}
f_{c'}\colon\Omega_c&\to\Omega_{c'}, 
\intertext{depending smoothly on $c'$ and with $f_c=\mathrm{id}$, such that}
\label{eqn:diffeos_compatibility}
f_{c'}\circ\psi_{c,k}&=\psi_{c',k}
\end{align}
for each $k\in\{1,\ldots,b\}$, where the bases $(e_{1,k},e_{2,k})$ of $T_{p_k'}\Sp^2$ from \eqref{eqn:20261009b} 
are chosen to depend smoothly on the centres $p_k'$ of $c'$.
Pulling back the round metric yields a family of metrics $f_{c'}^*g_{\Sp^2}$ on $M$, and these are the only variations of the metric we will consider.

\paragraph{Circle variations.} 
Given a circle domain $(\Omega_c,g_{\Sp^2})$ and some $\epsilon>0$, let 
\begin{align*}
\interval{-\epsilon,\epsilon}&\to\mathcal{C} \\
t&\mapsto c(t)=\bigl(p_1(t),\ldots,p_b(t);r_1(t),\ldots,r_b(t)\bigr)
\end{align*}
be a smooth curve and let 
\[
\dot c\vcentcolon=\frac{dc}{dt}(0)=(\dot p_1,\ldots,\dot p_b;\dot r_1,\ldots,\dot r_b)\in T_c\mathcal{C}
\]
denote its velocity at $t=0$. 
Then, recalling \eqref{eqn:diffeos_f_c}, the maps $f_{c(t)}\colon\Omega_c\to\Omega_{c(t)}$ define a smooth family of diffeomorphisms with $f_{c(0)}=\mathrm{id}$, and $Z\vcentcolon=\frac{d}{dt}\big\vert_{t=0}f_{c(t)}$ is a smooth vector field on $\Omega_c$ tangent to $\Sp^2$ such that 
\[
\frac{d}{dt}\Big\vert_{t=0}f_{c(t)}^*g_{\Sp^2}=\mathcal{L}_Zg_{\Sp^2}.
\]%
\begin{figure}\centering
\pgfmathsetmacro{\thetaO}{72}
\pgfmathsetmacro{\phiO}{70}
\tdplotsetmaincoords{\thetaO}{\phiO}
\def\caps{
\thetaCap/\phiCap/\rCap,
125/45/30,
135/180/20,
115/-60/15
}
\pgfmathsetmacro{\thetaCap}{5}
\pgfmathsetmacro{\phiCap}{0}
\pgfmathsetmacro{\rCap}{55}
\pgfmathsetmacro{\thetaNew}{\thetaCap+15}
\pgfmathsetmacro{\phiNew}{\phiCap+15}
\pgfmathsetmacro{\rNew}{\rCap-12}
\pgfmathsetmacro{\variationTime}{0.4}
\pgfmathsetmacro{\blickX}{sin(\thetaO)*sin(\phiO)}
\pgfmathsetmacro{\blickY}{-sin(\thetaO)*cos(\phiO)}
\pgfmathsetmacro{\blickZ}{cos(\thetaO)}
\begin{tikzpicture}[scale=\figurescale,line cap=round,tdplot_main_coords,
boundarystyle/.style={thick,FarbeB}]
\def\setcapdata#1#2#3{%
\pgfmathsetmacro{\thetaK}{#1}%
\pgfmathsetmacro{\phiK}{#2}%
\pgfmathsetmacro{\rK}{#3}%
\tdplotsetrotatedcoords{\phiK}{\thetaK}{0}%
\tdplottransformmainrot{\blickX}{\blickY}{\blickZ}%
\pgfmathsetmacro{\nx}{\tdplotresx}%
\pgfmathsetmacro{\ny}{\tdplotresy}%
\pgfmathsetmacro{\nz}{\tdplotresz}%
\pgfmathsetmacro{\sichtB}{sqrt(\nx*\nx+\ny*\ny)}%
\pgfmathsetmacro{\azimut}{atan2(\ny,\nx)}%
\pgfmathsetmacro{\phiOrot}{atan2(\nx,-\ny)}%
\pgfmathsetmacro{\breitenkreisradius}{sin(\rK)}%
\pgfmathsetmacro{\hoehe}{cos(\rK)}%
\pgfmathsetmacro{\visiblerange}{acos(max(-1,min(1,-\nz*cos(\rK)/max(.001,\sichtB*sin(\rK)))))}%
\pgfmathsetmacro{\rimrange}{acos(max(-1,min(1,cos(\rK)/max(.001,\sichtB))))}%
\pgfmathsetmacro{\mindot}{\nz*cos(\rK)-\sichtB*sin(\rK)}%
\pgfmathsetmacro{\maxdot}{\nz*cos(\rK)+\sichtB*sin(\rK)}%
}%
\def\capfullcontour{%
\foreach \i in {0,2,...,120}{%
\pgfmathsetmacro{\t}{3*\i}%
\tdplottransformrotmain{\breitenkreisradius*cos(\t)}{\breitenkreisradius*sin(\t)}{\hoehe}%
\ifnum\i=0
\pgfpathmoveto{\pgfpointxyz{\tdplotresx}{\tdplotresy}{\tdplotresz}}%
\else
\pgfpathlineto{\pgfpointxyz{\tdplotresx}{\tdplotresy}{\tdplotresz}}%
\fi
}%
\pgfpathclose
}%
\def\capfrontcontour{%
\foreach \i in {0,2,...,80}{%
\pgfmathsetmacro{\t}{\azimut-\visiblerange+\visiblerange*\i/40}%
\tdplottransformrotmain
{\breitenkreisradius*cos(\t)}
{\breitenkreisradius*sin(\t)}
{\hoehe}%
\ifnum\i=0
\pgfpathmoveto{\pgfpointxyz{\tdplotresx}{\tdplotresy}{\tdplotresz}}%
\else
\pgfpathlineto{\pgfpointxyz{\tdplotresx}{\tdplotresy}{\tdplotresz}}%
\fi
}%
\capsilhouette
\pgfpathclose
}%
\def\capsilhouette{%
\foreach \i in {0,2,...,80}{%
\pgfmathsetmacro{\s}{-\rimrange+\rimrange*\i/40}%
\tdplottransformrotmain
{(-\nz*\nx*cos(\s)+\ny*sin(\s))/\sichtB}
{(-\nz*\ny*cos(\s)-\nx*sin(\s))/\sichtB}
{\sichtB*cos(\s)}%
\pgfpathlineto{\pgfpointxyz{1.005*\tdplotresx}{1.005*\tdplotresy}{1.005*\tdplotresz}}%
}%
}%
\def\cliphemisphere{%
\pgfpathmoveto{\pgfpointanchor{clipSW}{center}}%
\pgfpathlineto{\pgfpointanchor{clipSE}{center}}%
\pgfpathlineto{\pgfpointanchor{clipNE}{center}}%
\pgfpathlineto{\pgfpointanchor{clipNW}{center}}%
\pgfpathclose
\foreach \thetaCircle/\phiCircle/\rCircle in \caps {%
\setcapdata{\thetaCircle}{\phiCircle}{\rCircle}%
\ifdim\mindot pt>0pt
\capfullcontour
\else
\ifdim\maxdot pt>0pt
\capfrontcontour
\fi
\fi
}%
\pgfusepath{clip}%
}
\def\shadehemisphere#1{%
\begin{scope}[tdplot_screen_coords]
\clip(0,0)circle(1);
\begin{scope}[even odd rule,tdplot_main_coords]
\cliphemisphere
\path[tdplot_screen_coords,ball color=#1] (0,0) circle (1.02);
\end{scope}
\end{scope}
}
\path[tdplot_screen_coords]
(-1.05,-1.05) coordinate (clipSW)
(1.05,-1.05) coordinate (clipSE)
(1.05,1.05) coordinate (clipNE)
(-1.05,1.05) coordinate (clipNW);
{%
\pgfmathsetmacro{\blickX}{-\blickX}%
\pgfmathsetmacro{\blickY}{-\blickY}%
\pgfmathsetmacro{\blickZ}{-\blickZ}%
\shadehemisphere{cyan!10!black!25}%
}
\begin{scope}[tdplot_main_coords,boundarystyle]
\foreach \thetaCircle/\phiCircle/\rCircle in \caps {
\setcapdata{\thetaCircle}{\phiCircle}{\rCircle}
\ifdim\mindot pt<0pt
\tdplotdrawarc[tdplot_rotated_coords]
{(0,0,\hoehe)}{\breitenkreisradius}
{\phiOrot-90-\visiblerange}
{\phiOrot-450+\visiblerange}{}{}
\fi
}
\end{scope}
\shadehemisphere{white}
\begin{scope}[even odd rule,tdplot_main_coords]
\cliphemisphere
\begin{scope}[dotted,boundarystyle]
\foreach \thetaCircle/\phiCircle/\rCircle in \caps {
\setcapdata{\thetaCircle}{\phiCircle}{\rCircle}
\ifdim\mindot pt<0pt
\tdplotdrawarc[tdplot_rotated_coords]
{(0,0,\hoehe)}{\breitenkreisradius}
{\phiOrot-90-\visiblerange}
{\phiOrot-450+\visiblerange}{}{}
\fi
}
\end{scope}
\end{scope}
\begin{scope}[tdplot_main_coords,boundarystyle]
\foreach[count=\k] \thetaCircle/\phiCircle/\rCircle in \caps {
\tdplotsetcoord{p\k}{1}{\thetaCircle}{\phiCircle}
\setcapdata{\thetaCircle}{\phiCircle}{\rCircle}
\ifdim\maxdot pt>0pt
\tdplotdrawarc[tdplot_rotated_coords]
{(0,0,\hoehe)}{\breitenkreisradius}
{\phiOrot-90+\visiblerange}
{\phiOrot-90-\visiblerange}{}{}
\fi
}
\end{scope} 
\tdplotsetrotatedcoords{\phiNew}{\thetaNew}{0}
\pgfmathsetmacro{\hoeheNew}{cos(\rNew)}
\pgfmathsetmacro{\breitenkreisradiusNew}{sin(\rNew)}
\tdplotdrawarc[tdplot_rotated_coords,FarbeB,thick,densely dashed]
{(0,0,\hoeheNew)}{\breitenkreisradiusNew}{0}{360}{}{} 
\tdplotsetcoord{pk}{1}{\thetaCap}{\phiCap}
\tdplotsetcoord{pkNew}{1}{\thetaNew}{\phiNew}
\pgfmathsetmacro{\centreSeparation}{acos(max(-1,min(1,cos(\thetaCap)*cos(\thetaNew)+sin(\thetaCap)*sin(\thetaNew)*cos(\phiNew-\phiCap))))}
\pgfmathsetmacro{\centreAzimuth}{atan2(sin(\thetaNew)*sin(\phiNew-\phiCap),cos(\thetaCap)*sin(\thetaNew)*cos(\phiNew-\phiCap)-sin(\thetaCap)*cos(\thetaNew))}
\pgfmathsetmacro{\velocityEOne}{rad(\centreSeparation)*cos(\centreAzimuth)/\variationTime}
\pgfmathsetmacro{\velocityETwo}{rad(\centreSeparation)*sin(\centreAzimuth)/\variationTime}
\pgfmathsetmacro{\rDot}{rad(\rNew-\rCap)/\variationTime}
\tdplotsetrotatedcoords{\phiCap}{\thetaCap}{0}
\tdplotsetrotatedthetaplanecoords{\centreAzimuth}
\tdplotdrawarc[tdplot_rotated_coords,thin,densely dotted]
{(0,0,0)}{1}{0}{\centreSeparation}{}{}
\tdplotsetrotatedcoords{\phiCap}{\thetaCap}{0}
\draw[tdplot_rotated_coords,FarbeA,thick,-stealth]
(0,0,1) -- (\velocityEOne,\velocityETwo,1)
node[below] {$\dot p_k$};
\draw plot[bullet] coordinates {(pk)};
\draw plot[bullet] coordinates {(pkNew)};
\node[left] at (pk) {$p_k$};
\node[below] at (pkNew) {$p_k(t)$}; 
\pgfmathsetmacro{\tPoint}{-70}
\tdplotsetrotatedcoords{\phiCap}{\thetaCap}{0}
\pgfmathsetmacro{\qX}{sin(\rCap)*cos(\tPoint)}
\pgfmathsetmacro{\qY}{sin(\rCap)*sin(\tPoint)}
\pgfmathsetmacro{\qZ}{cos(\rCap)}
\pgfmathsetmacro{\etaX}{-cos(\rCap)*cos(\tPoint)}
\pgfmathsetmacro{\etaY}{-cos(\rCap)*sin(\tPoint)}
\pgfmathsetmacro{\etaZ}{sin(\rCap)}
\pgfmathsetmacro{\tauX}{-sin(\tPoint)}
\pgfmathsetmacro{\tauY}{cos(\tPoint)}
\pgfmathsetmacro{\tauZ}{0}
\pgfmathsetmacro{\alphaDot}{-\rDot-\velocityEOne*cos(\tPoint)-\velocityETwo*sin(\tPoint)}
\pgfmathsetmacro{\betaDot}{cos(\rCap)*(\velocityETwo*cos(\tPoint)-\velocityEOne*sin(\tPoint))}
\pgfmathsetmacro{\zX}{\alphaDot*\etaX+\betaDot*\tauX}
\pgfmathsetmacro{\zY}{\alphaDot*\etaY+\betaDot*\tauY}
\pgfmathsetmacro{\zZ}{\alphaDot*\etaZ+\betaDot*\tauZ}
\path[tdplot_rotated_coords] (\qX,\qY,\qZ) coordinate (q);
\path[tdplot_rotated_coords] ({\qX+\alphaDot*\etaX},{\qY+\alphaDot*\etaY},{\qZ+\alphaDot*\etaZ}) coordinate (qAlpha);
\path[tdplot_rotated_coords] ({\qX+\betaDot*\tauX},{\qY+\betaDot*\tauY},{\qZ+\betaDot*\tauZ}) coordinate (qBeta);
\path[tdplot_rotated_coords] ({\qX+\zX},{\qY+\zY},{\qZ+\zZ}) coordinate (qZ); 
\tdplotsetrotatedcoords{\phiCap}{\thetaCap}{0}
\begin{scope}[tdplot_rotated_coords,thin,densely dotted,FarbeB]
\foreach \i in {0,2,...,80}{%
\pgfmathsetmacro{\s}{\i/80}%
\pgfmathsetmacro{\ds}{\s*\centreSeparation}%
\pgfmathsetmacro{\rs}{\rCap+\s*(\rNew-\rCap)}%
\pgfmathsetmacro{\qsU}{cos(\rs)*sin(\ds)+sin(\rs)*cos(\tPoint-\centreAzimuth)*cos(\ds)}%
\pgfmathsetmacro{\qsV}{sin(\rs)*sin(\tPoint-\centreAzimuth)}%
\pgfmathsetmacro{\qsZ}{cos(\rs)*cos(\ds)-sin(\rs)*cos(\tPoint-\centreAzimuth)*sin(\ds)}%
\pgfmathsetmacro{\qsX}{\qsU*cos(\centreAzimuth)-\qsV*sin(\centreAzimuth)}%
\pgfmathsetmacro{\qsY}{\qsU*sin(\centreAzimuth)+\qsV*cos(\centreAzimuth)}%
\ifnum\i=80
\xdef\qNewX{\qsX}%
\xdef\qNewY{\qsY}%
\xdef\qNewZ{\qsZ}%
\fi
\ifnum\i=0
\pgfpathmoveto{\pgfpointxyz{\qsX}{\qsY}{\qsZ}}%
\else
\pgfpathlineto{\pgfpointxyz{\qsX}{\qsY}{\qsZ}}%
\fi
}%
\pgfusepath{stroke}
\draw plot[bullet](\qNewX,\qNewY,\qNewZ)node[anchor=200] {$q(t)$};
\end{scope}
\draw[tdplot_rotated_coords,FarbeA,thin,densely dotted] 
(qZ)--(qAlpha) 
(qZ)--(qBeta) 
;
\draw[tdplot_rotated_coords,black,semithick,-stealth](q)--({\qX+\etaX},{\qY+\etaY},{\qZ+\etaZ})node[below left,inner sep=0]{$\eta~~$};
\draw[tdplot_rotated_coords,black,semithick,-stealth](q)--({\qX+\tauX},{\qY+\tauY},{\qZ+\tauZ})node[right]{$\tau$};
\draw[tdplot_rotated_coords,FarbeA,thick,-stealth](q)--(qZ)node[right] {$Z$};
\draw[tdplot_rotated_coords,thick,magenta] 
(q)--(qAlpha)node[left,midway]{$\alpha_{\dot c}$}
(q)--(qBeta)node[below,near end]{$\beta_{\dot c}$}
;
\draw[FarbeB] plot[bullet] (q)node[below left] {$q$};
\draw[black] plot[bullet] (p4);
\draw[black] plot[bullet] (p2);
\end{tikzpicture}
\caption{Variation of the boundary circle $\partial_k\Omega_c\subset\Sp^2$ with velocity $Z=\alpha_{\dot c}\eta+\beta_{\dot c}\tau$.}
\label{fig:circle_variation}
\end{figure}
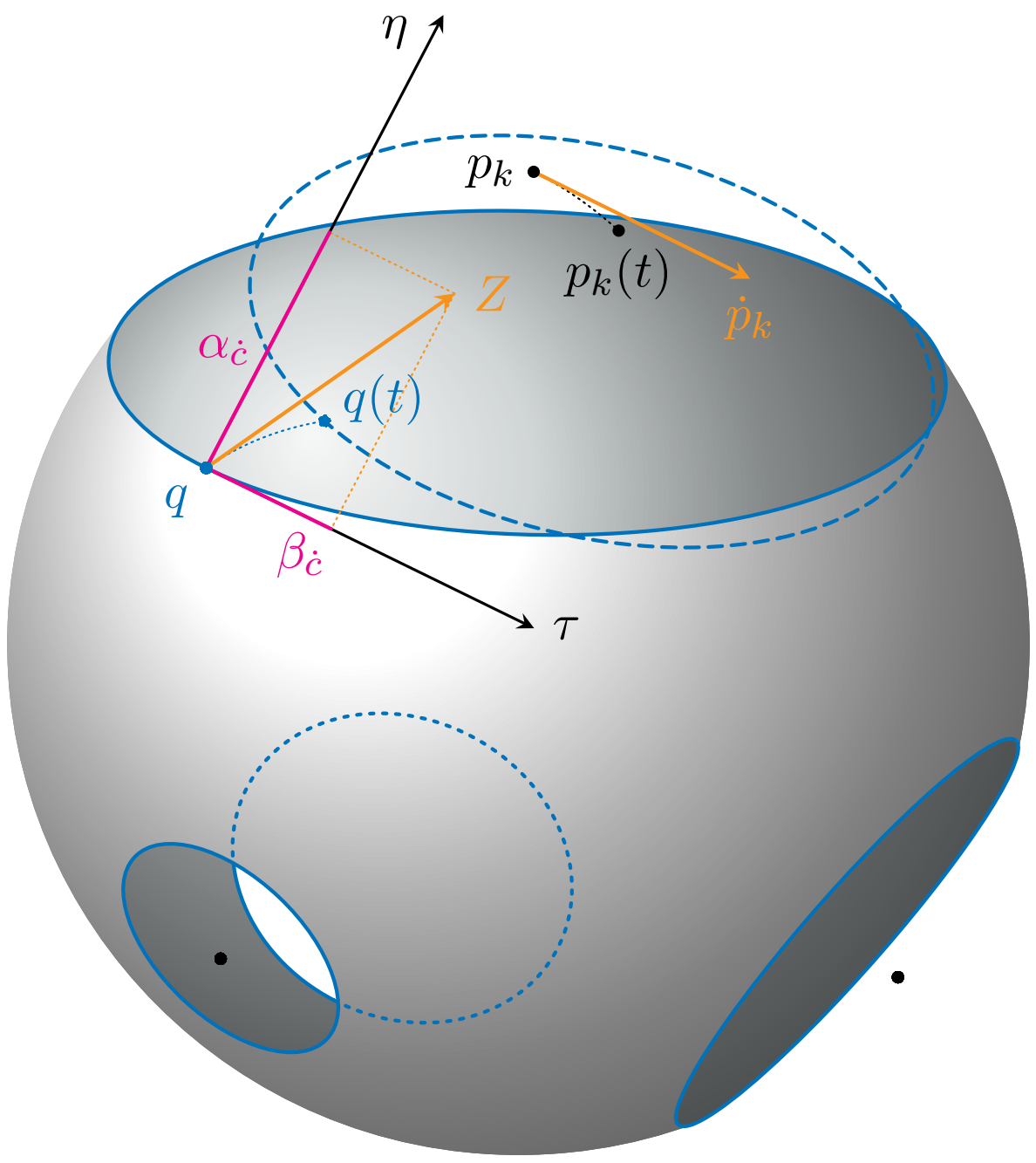%
Given $k\in\{1,\ldots,b\}$, the outward unit conormal of $\Omega_c$ at $q\in\partial_k\Omega_c$ as in \eqref{eqn:k-th_boundary_circle} is given by
\[
\eta(q)=\frac{p_k-q\cos r_k}{\sin r_k}.
\]
By \eqref{eqn:diffeos_compatibility} the diffeomorphism $f_{c(t)}$ maps $\partial_k\Omega_c$ onto $\partial_k\Omega_{c(t)}$, so that
\begin{align}
\label{eqn:20261009a}
\sk{f_{c(t)}(q),p_k(t)}=\cos r_k(t)
\end{align}
for all $t\in\interval{-\epsilon,\epsilon}$ and all $q\in\partial_k\Omega_c$.
Differentiating \eqref{eqn:20261009a} at $t=0$ yields $\sk{Z(q),p_k}+\sk{q,\dot p_k}=-\dot r_k\sin r_k$. 
Since $Z(q)$ is tangent to $\Sp^2$, we have $\sk{Z(q),q}=0$, and therefore 
\begin{align}
\label{eqn:normal-velocity}
\alpha_{\dot c}(q)\vcentcolon=g_{\Sp^2}\bigl(Z(q),\eta(q)\bigr)
=\frac{\sk{Z(q),p_k}}{\sin r_k}
=-\dot r_k-\frac{\sk{\dot p_k,q}}{\sin r_k}
\end{align}
as visualised in Figure~\ref{fig:circle_variation}. 
In particular, $\alpha_{\dot c}$ depends only on $\dot c$ and not on the choice of the diffeomorphisms $f_{c'}$. 
To make this more explicit, we recall the parametrisation $\psi_{c,k}\colon 
\Interval{0,2\pi}\to\partial_k\Omega_c$ from \eqref{eqn:20261009b}. 
Since $\abs{p_k(t)}=1$ for all $t$, we have $\sk{\dot p_k,p_k}=0$, and \eqref{eqn:normal-velocity} becomes 
\begin{align}\label{eqn:normal-velocity-fourier}
\alpha_{\dot c}\bigl(\psi_{c,k}(\theta)\bigr)
&=-\dot r_k-\sk{\dot p_k,e_{1,k}}\cos\theta-\sk{\dot p_k,e_{2,k}}\sin\theta.
\end{align}
Thus, the radius velocity $\dot r_k$ yields the constant mode, and the two components of the centre velocity $\dot p_k$ yield the two first-order Fourier modes on $\partial_k\Omega_c$.
Conversely, every function on $\partial\Omega_c$ which on each $\partial_k\Omega_c$ is a linear combination of $1$, $\cos\theta$ and $\sin\theta$ is of the form $\alpha_{\dot c}$ for a unique $\dot c\in T_c\mathcal{C}$, since $1$, $\cos\theta$ and $\sin\theta$ are linearly independent and $\dot p_k\in T_{p_k}\Sp^2$ is determined by its components $\sk{\dot p_k,e_{1,k}}$ and $\sk{\dot p_k,e_{2,k}}$.

Let $\tau$ denote the unit tangent vector field along $\partial\Omega_c$ satisfying $\psi_{c,k}'=\sin(r_k)\,\tau\circ\psi_{c,k}$ for each $k\in\{1,\ldots,b\}$, such that $(\eta,\tau)$ is an orthonormal frame of $T\Sp^2$ along $\partial\Omega_c$. 
By \eqref{eqn:diffeos_compatibility}, we have $Z\circ\psi_{c,k}=\frac{d}{dt}\big\vert_{t=0}\psi_{c(t),k}$. 
Hence, also the tangential component 
\[
\beta_{\dot c}\vcentcolon=g_{\Sp^2}(Z,\tau)
\]
of $Z$ along $\partial\Omega_c$ is independent of the choice of the diffeomorphisms $f_{c'}$. 
Indeed, taking the inner product of $\tau\circ\psi_{c,k}(\theta)=-e_{1,k}\sin\theta+e_{2,k}\cos\theta$ with $\frac{d}{dt}\big\vert_{t=0}\psi_{c(t),k}(\theta)$ yields, analogously to \eqref{eqn:normal-velocity-fourier}, 
\begin{align}\label{eqn:tangential-velocity-fourier}
\beta_{\dot c}\bigl(\psi_{c,k}(\theta)\bigr)
=\sk{\dot p_k,e_{2,k}}\cos(r_k)\cos\theta-\sk{\dot p_k,e_{1,k}}\cos(r_k)\sin\theta+\sk{\dot e_{1,k},e_{2,k}}\sin(r_k), 
\end{align}
where $\dot e_{1,k}$ denotes the velocity at $t=0$ of the basis vector $e_{1,k}$ of $T_{p_k(t)}\Sp^2$ from \eqref{eqn:diffeos_compatibility}. 
For the choice of bases in Section~\ref{sec:gauge}, the last term in \eqref{eqn:tangential-velocity-fourier} vanishes for curves through a fixed reference configuration (see Remark~\ref{rem:parallel-transport}).

\paragraph{Conformality.}
For a free boundary harmonic map $\Phi\colon(\Omega_c,g_{\Sp^2})\to\B^3$, the identity \eqref{eqn:hopf-boundary} with the vector fields $Z$ associated with circle variations as in \eqref{eqn:normal-velocity} tests the Hopf differential $\Hopf{\Phi}$ only against metric variations induced by motions of the boundary circles.
Following \cite[§\,5.3]{Tran}, a \emph{Teichmüller direction} of $\Omega_c$ is a smooth trace-free symmetric bilinear form on $\Omega_c$ modulo trace-free parts of $\mathcal{L}_Yg_{\Sp^2}$, where $Y$ is any smooth vector field on $\Omega_c$ which is tangent to $\partial\Omega_c$ along $\partial\Omega_c$. 
By \cite[Theorem~5.9]{Tran}, every Teichmüller direction of $\Omega_c$ is induced by a motion of the boundary circles. 
We only need this surjectivity, which reads as follows in our notation, recalling $\alpha_{\dot c}$ from \eqref{eqn:normal-velocity}.

\begin{lemma}[{\cite[Theorem~5.9]{Tran}}]
\label{lem:Tran}
For every smooth trace-free symmetric bilinear form $h$ on $\Omega_c$ there exist $\dot c\in T_c\mathcal{C}$ and a smooth vector field $W$ on $\Omega_c$ tangent to $\Sp^2$ such that $g_{\Sp^2}(W,\eta)=\alpha_{\dot c}$ on $\partial\Omega_c$, and such that $h$ is the trace-free part of $\mathcal{L}_Wg_{\Sp^2}$.  
\end{lemma}

\begin{proof}
The configuration space in \cite[§\,5.4]{Tran} coincides with $\mathcal{C}$, with the caps given by the closed discs $\overline{B}_{r_k}(p_k)$. 
In the notation therein, the normal $\mathbf{n}$ points into $\Omega_c$, i.\,e.~$\mathbf{n}=-\eta$, and the normal speed $\theta_{\dot\omega}$ for $\dot\omega=\dot c$ equals $-\alpha_{\dot c}$ by \eqref{eqn:normal-velocity}. 
Hence, \cite[Theorem~5.9]{Tran} yields $\dot c\in T_c\mathcal{C}$ and smooth vector fields $Z$, $Y$ on $\Omega_c$ such that $g_{\Sp^2}(Z,\eta)=\alpha_{\dot c}$ and $g_{\Sp^2}(Y,\eta)=0$ on $\partial\Omega_c$, and such that $h$ is the trace-free part of $\mathcal{L}_Zg_{\Sp^2}+\mathcal{L}_Yg_{\Sp^2}$. 
The claim follows with $W=Z+Y$. 
\end{proof}

The following statement combines \eqref{eqn:hopf-boundary} and Lemma~\ref{lem:Tran} to reduce weak conformality of a free boundary harmonic map to finitely many conditions on the boundary, which amount to criticality with respect to circle variations by Lemma~\ref{lem:critical-c} below.

\begin{proposition}
\label{prop:critical-conformal}
Let $\Phi\colon(\Omega_c,g_{\Sp^2})\to\B^3$ be a nonconstant free boundary harmonic map. 
If for each $k\in\{1,\ldots,b\}$, the restriction of $\Hopf{\Phi}(\eta,\eta)$ to $\partial_k\Omega_c$ is $L^2$-orthogonal to $1$, $\cos\theta$ and $\sin\theta$ from \eqref{eqn:20261009b},  
then $\Phi$ is weakly conformal and thus a branched free boundary minimal immersion. 
\end{proposition}

\begin{proof}
Applying Lemma~\ref{lem:Tran} to the Hopf differential $h=\Hopf{\Phi}$ yields $\dot c\in T_c\mathcal{C}$ and a vector field $W$ as stated therein. 
Since $\Hopf{\Phi}$ is trace-free, it is pointwise orthogonal to the trace part of $\mathcal{L}_Wg_{\Sp^2}$. 
Hence, \eqref{eqn:hopf-boundary} applied with $W$ in place of $Z$ yields 
\[
\int_{\Omega_c}\abs{\Hopf{\Phi}}^2
=\int_{\Omega_c}\sk{\Hopf{\Phi},\mathcal{L}_Wg_{\Sp^2}}
=2\int_{\partial\Omega_c}\alpha_{\dot c}\,\Hopf{\Phi}(\eta,\eta).
\]
By \eqref{eqn:normal-velocity-fourier}, the restriction of $\alpha_{\dot c}$ to each $\partial_k\Omega_c$ is a linear combination of $1$, $\cos\theta$ and $\sin\theta$, and hence the right-hand side vanishes by assumption. 
Therefore, $\Hopf{\Phi}\equiv0$, i.\,e.~$\Phi$ is weakly conformal, and the claim follows from the discussion following \eqref{eqn:fbhm_eq}. 
\end{proof}

\paragraph{Reduced energy.}
We reduce the problem of finding critical points of the two-variable energy \eqref{eqn:2variable_energy} by restricting the metric to the round metric on circle domains $\Omega_c$ with $c\in\mathcal{C}$, and by determining the map via harmonic extension of its boundary values.
The resulting variational problem depends only on the boundary values and on the $3b$ parameters of $c$. 
These boundary values are given by maps on the disjoint union 
\begin{align}\label{eqn:circle_union}
\Xi\vcentcolon=\bigsqcup_{k=1}^b\Sp^1
\end{align}
of $b$ circles, where we identify $\Sp^1$ with $\R/2\pi\Z$. 
The parametrisations \eqref{eqn:20261009b} combine to a diffeomorphism 
\begin{align*}
\psi_c\colon\Xi\to\partial\Omega_c
\end{align*}
which restricts to $\psi_{c,k}$ on the $k$-th copy of $\Sp^1$. 
Given $\gamma\in L^2(\Xi,\R^3)$, we denote its restriction to the $k$-th copy of $\Sp^1$ by $\gamma_k$ and its Fourier coefficients by 
\begin{align}\label{eqn:Fourier-coefficients}
\hat\gamma_k(n)\vcentcolon=\frac{1}{2\pi}\int_0^{2\pi}\gamma_k(\theta)e^{-in\theta}\,d\theta.
\end{align}
for $n\in\Z$. 
Given $s\geq0$ and $d \geq 1$, the Sobolev space $H^s(\Xi,\R^d)$ consists of all $\gamma\in L^2(\Xi,\R^d)$ with finite norm
\begin{align}\label{eqn:Sobolev-norm}
\nm{\gamma}_{H^s}^2\vcentcolon=2\pi\sum_{k=1}^b\sum_{n\in\Z}{(1+n^2)}^s\abs{\hat\gamma_k(n)}^2. 
\end{align}
The pointwise norm $\abs{\cdot}$ denotes the Euclidean norm, and analogously we define the Sobolev norm for matrix-valued functions by taking the Frobenius norm.
If $s\in\N\cup\{0\}$, then Parseval's identity yields $\nm{\gamma}_{H^s}^2=\sum_{j=0}^s\binom{s}{j}\nm{\gamma^{(j)}}_{L^2}^2$, where $\gamma^{(j)}$ denotes the $j$-th weak derivative with respect to $\theta$ and $\nm{\cdot}_{L^2}$ the $L^2$-norm with respect to $d\theta$. 
Moreover, the Cauchy--Schwarz inequality and the identity $\sum_{n\in\Z}(1+n^2)^{-1}=\pi\coth(\pi)$ imply 
\begin{align}
\label{eqn:Sobolev-constant}
\sum_{n\in\Z}\abs{\hat\gamma_k(n)}
&\leq\Bigl(\sum_{n\in\Z}\frac{1}{1+n^2}\Bigr)^{\frac{1}{2}}\Bigl(\sum_{n\in\Z}(1+n^2)\abs{\hat\gamma_k(n)}^2\Bigr)^{\frac{1}{2}}
\leq c_\infty\nm{\gamma}_{H^1}, 
&
c_\infty&\vcentcolon=\sqrt{\frac{\coth(\pi)}{2}}, 
\end{align}
for all $k\in\{1,\ldots,b\}$. 
Hence, if $\gamma\in H^1(\Xi,\R^3)$, then the Fourier series $\sum_{n\in\Z}\hat\gamma_k(n)e^{in\theta}$ converges absolutely and uniformly to $\gamma_k$, which is therefore continuous with $\abs{\gamma_k(\theta)}\leq c_\infty\nm{\gamma}_{H^1}$ for all $\theta\in\Sp^1$. 
Applying this also to $\gamma'$ and using $\nm{\gamma'}_{H^1}\leq\nm{\gamma}_{H^2}$, we obtain 
\begin{align}\label{eqn:Sobolev-embedding}
\max\bigl\{\nm{\gamma}_{L^\infty},\nm{\gamma'}_{L^\infty}\bigr\}\leq c_\infty\nm{\gamma}_{H^2} 
\end{align}
and in particular $H^2(\Xi,\R^3)\subset C^1(\Xi,\R^3)$. 
This allows us to define  
\[
H^2(\Xi,\Sp^2)\vcentcolon=\{\gamma\in H^2(\Xi,\R^3)\mid\abs{\gamma}=1\},
\]
where the constraint is understood pointwise.

\begin{definition}[Reduced energy]
\label{defn:reduced-energy}
Given a configuration $c\in\mathcal{C}$ and boundary data $\gamma\in H^2(\Xi,\R^3)$, the \emph{harmonic extension} $\Phi_{\gamma,c}\colon\Omega_c\to\R^3$ of $\gamma$ is the unique solution of 
\begin{align*}
\left\{\begin{aligned}
\Delta_{g_{\Sp^2}}\Phi&=0 & &\text{ in $\Omega_c$,} \\
\Phi&=\gamma\circ\psi_c^{-1} & &\text{ on $\partial\Omega_c$.}
\end{aligned}\right.
\end{align*}
We define the associated \emph{reduced energy} as
\begin{align*}
\begin{aligned}
E\colon H^2(\Xi,\R^3)\times\mathcal{C}&\to\R,
\\[-.5ex]
(\gamma,c)&\mapsto E_{g_{\Sp^2}}(\Phi_{\gamma,c})=\frac{1}{2}\int_{\Omega_c}\abs{d\Phi_{\gamma,c}}^2. 
\end{aligned}
\end{align*}
\end{definition}

\begin{remark}
The harmonic extension $\Phi_{\gamma,c}$ minimises the Dirichlet energy among all maps $\Phi\colon\Omega_c\to\R^3$ of class $H^1$ which coincide with $\gamma\circ\psi_c^{-1}$ on $\partial\Omega_c$,
see \cite[§\,2.2.5]{Evans2010PDE}.
In other words, $E(\gamma,c)$ is the least energy needed to span the parametrised boundary curves $\gamma\circ\psi_c^{-1}$ by a map on the circle domain $\Omega_c$.
\end{remark}

\begin{remark}
\label{rem:bases_dependence}
The reduced energy depends on the bases $(e_{1,k},e_{2,k})$ in \eqref{eqn:20261009b} only through the parametrisations $\psi_{c,k}$. 
Replacing $(e_{1,k},e_{2,k})$ by another orthonormal basis of $T_{p_k}\Sp^2$ amounts to replacing $\psi_{c,k}(\theta)$ by $\psi_{c,k}(a\pm\theta)$ for some $a\in\R$, depending on whether the orientation of the basis is preserved. 
Hence, the reduced energy of $(\gamma,c)$ with respect to the new bases coincides with the reduced energy of $(\tilde\gamma,c)$ with respect to the original bases, where $\tilde\gamma_k(a\pm\theta)=\gamma_k(\theta)$, and the corresponding harmonic extensions coincide as maps on $\Omega_c$. 
In particular, the choice of the bases does not affect the geometric conclusions of our construction, and we fix it in Section~\ref{sec:gauge}. 
\end{remark}

\paragraph{Dirichlet-to-Neumann operator.} 
The free boundary condition in \eqref{eqn:fbhm_eq} involves the normal derivative of the map along the boundary. 
For harmonic extensions, this normal derivative is determined by the boundary values via the linear, but non-local, \emph{Dirichlet-to-Neumann operator}. 
For a fixed configuration $c\in\mathcal{C}$, and in terms of the parametrisations $\psi_{c,k}$ from \eqref{eqn:20261009b}, this operator is given by
\begin{align}
\label{eqn:DtN}
\DN_c\colon H^2(\Xi,\R^3)&\to H^1(\Xi,\R^3), &
(\DN_c\gamma)_k&\vcentcolon=\sin(r_k)\,(\partial_\eta\Phi_{\gamma,c})\circ\psi_{c,k}
\end{align}
for each $k\in\{1,\ldots,b\}$, where $\eta$ denotes the outward unit conormal along $\partial\Omega_c$ as before. 
The factor $\sin(r_k)=\abs{\psi_{c,k}'}$ in \eqref{eqn:DtN} accounts for the arclength element $\sin(r_k)\,d\theta$ of $\partial_k\Omega_c$ in the parametrisation $\psi_{c,k}$. 
In other words, $\DN_c\gamma\,d\theta$ is the pullback of $\partial_\eta\Phi_{\gamma,c}\,ds$ via $\psi_c$, which ensures that Green's formula takes the form 
\begin{align}\label{eqn:DtN-Green}
\int_{\Omega_c}\sk{d\Phi_{\beta,c},d\Phi_{\gamma,c}}
=\int_{\partial\Omega_c}\sk{\Phi_{\beta,c},\partial_\eta\Phi_{\gamma,c}}
=\int_\Xi\sk{\beta,\DN_c\gamma}\,d\theta 
\end{align}
for all $\beta,\gamma\in H^2(\Xi,\R^3)$. 
In particular, $\DN_c$ is symmetric with respect to the $L^2$-inner product on $\Xi$, and for all $(\gamma,c)\in H^2(\Xi,\R^3)\times\mathcal{C}$ we conclude
\begin{align}\label{eqn:reduced-energy-DtN}
E(\gamma,c)=\frac{1}{2}\int_\Xi\sk{\gamma,\DN_c\gamma}\,d\theta.
\end{align}

\paragraph{Critical points.}
We call $(\gamma,c)\in H^2(\Xi,\Sp^2)\times\mathcal{C}$ a \emph{critical point} of the reduced energy \eqref{eqn:reduced-energy-DtN} if 
\begin{align}\label{eqn:20261013b}
\frac{d}{dt}\Big\vert_{t=0}E\bigl(\gamma(t),c\bigr)&=0, 
&
\frac{d}{dt}\Big\vert_{t=0}E\bigl(\gamma,c(t)\bigr)&=0
\end{align}
for every $C^1$-curve $t\mapsto\gamma(t)$ in $H^2(\Xi,\Sp^2)$ with $\gamma(0)=\gamma$ and every smooth curve $t\mapsto c(t)$ in $\mathcal{C}$ with $c(0)=c$, where both derivatives in \eqref{eqn:20261013b} exist by Lemmata~\ref{lem:critical-gamma} and~\ref{lem:critical-c} stated below. 
Given $\gamma\in H^2(\Xi,\Sp^2)$, the space of variations of $\gamma$ preserving the constraint $\abs{\gamma}=1$ is 
\begin{align}\label{eqn:variations-gamma}
\mathcal{V}_\gamma\vcentcolon=\{V\in H^2(\Xi,\R^3)\mid\sk{V,\gamma}=0\}
\end{align}  
and we denote the pointwise orthogonal projection onto the tangent planes $T_{\gamma(\theta)}\Sp^2$ by 
\begin{align}
\label{eqn:Pgamma}
\begin{aligned}
P_\gamma\colon H^1(\Xi,\R^3)&\to H^1(\Xi,\R^3) \\
w&\mapsto w-\sk{w,\gamma}\gamma.
\end{aligned}
\end{align}
Note that $P_\gamma$ maps $H^2(\Xi,\R^3)$ onto $\mathcal{V}_\gamma$. 

\begin{lemma}
\label{lem:critical-gamma}
Given $(\gamma,c)\in H^2(\Xi,\Sp^2)\times\mathcal{C}$, let $\Phi=\Phi_{\gamma,c}$ be the harmonic extension of $\gamma$ as in Definition~\ref{defn:reduced-energy}, 
and let $t\mapsto\gamma(t)$ be a $C^1$-curve in $H^2(\Xi,\Sp^2)$ with $\gamma(0)=\gamma$. 
Then $V=\frac{d}{dt}\big\vert_{t=0}\gamma(t)\in\mathcal{V}_\gamma$, and 
\begin{align}\label{eqn:variation-gamma}
\frac{d}{dt}\Big\vert_{t=0}E\bigl(\gamma(t),c\bigr)=\int_\Xi\sk{V,\DN_c\gamma}\,d\theta. 
\end{align} 
Moreover, the right-hand side of \eqref{eqn:variation-gamma} vanishes for all $V\in\mathcal{V}_\gamma$ if and only if $P_\gamma\DN_c\gamma=0$, which is the case if and only if $\Phi$ is free boundary harmonic. 
In this case, $\Phi$ is smooth up to the boundary. 
\end{lemma}

\begin{proof}
Differentiating $\abs{\gamma(t)}^2=1$ pointwise, which is justified by \eqref{eqn:Sobolev-embedding}, yields $\sk{V,\gamma}=0$. 
Since $\DN_c\colon H^2(\Xi,\R^3)\to H^1(\Xi,\R^3)$ is linear, bounded and symmetric, \eqref{eqn:reduced-energy-DtN} implies 
\[
\frac{d}{dt}\Big\vert_{t=0}E\bigl(\gamma(t),c\bigr)
=\frac{1}{2}\int_\Xi\sk{V,\DN_c\gamma}\,d\theta
 +\frac{1}{2}\int_\Xi\sk{\gamma,\DN_cV}\,d\theta
=\int_\Xi\sk{V,\DN_c\gamma}\,d\theta. 
\]
Since $P_\gamma$ maps $H^2(\Xi,\R^3)$ onto $\mathcal{V}_\gamma$ and $\int_\Xi\sk{P_\gamma w,\DN_c\gamma}\,d\theta=\int_\Xi\sk{w,P_\gamma\DN_c\gamma}\,d\theta$ for all smooth maps $w\colon\Xi\to\R^3$, the right-hand side of \eqref{eqn:variation-gamma} vanishes for all $V\in\mathcal{V}_\gamma$ if and only if $P_\gamma\DN_c\gamma=0$. 
By \eqref{eqn:DtN}, this is the case if and only if the component of $\partial_\eta\Phi$ tangent to $\Sp^2$ vanishes along $\partial\Omega_c$, i.\,e.~if and only if $\partial_\eta\Phi$ is parallel to $\Phi$ along $\partial\Omega_c$. 
In particular, if $\Phi$ is free boundary harmonic, then $P_\gamma\DN_c\gamma=0$ by \eqref{eqn:fbhm_eq}. 
Conversely, if $P_\gamma\DN_c\gamma=0$, then $\gamma\circ\psi_c^{-1}\colon\partial\Omega_c\to\Sp^2$ is $\frac{1}{2}$-harmonic in the sense of \cite[Definition~1.1]{DaLioPigati2020} by \cite[(1.5)]{DaLioPigati2020}, and hence smooth by \cite[Theorem~1.4]{DaLioPigati2020}.
Consequently, $\Phi$ is smooth up to the boundary. 
Since $\abs{\Phi}^2$ is subharmonic and $\abs{\Phi}=1$ on $\partial\Omega_c$, the maximum principle implies $\Phi(\Omega_c)\subset\B^3$ and $\sk{\partial_\eta\Phi,\Phi}=\frac{1}{2}\partial_\eta\abs{\Phi}^2\geq0$ along $\partial\Omega_c$. 
Therefore, $\partial_\eta\Phi=\abs{\partial_\eta\Phi}\Phi$, and $\Phi$ solves \eqref{eqn:fbhm_eq}. 
\end{proof}

\begin{lemma}
\label{lem:critical-c}
Given $(\gamma,c)\in H^2(\Xi,\R^3)\times\mathcal{C}$, let $\Phi=\Phi_{\gamma,c}$ be the harmonic extension of $\gamma$ as in Definition~\ref{defn:reduced-energy}, 
and let $t\mapsto c(t)$ be a smooth curve in $\mathcal{C}$ with $c(0)=c$. 
Let $\alpha_{\dot c}$ and $\beta_{\dot c}$ be the normal and tangential components along $\partial\Omega_c$ of the associated vector field $Z$ as in \eqref{eqn:normal-velocity} and \eqref{eqn:tangential-velocity-fourier}. 
Then 
\begin{align}\label{eqn:variation-c}
\frac{d}{dt}\Big\vert_{t=0}E\bigl(\gamma,c(t)\bigr)
=-\int_{\partial\Omega_c}\Bigl(\alpha_{\dot c}\,\Hopf{\Phi}(\eta,\eta)+\beta_{\dot c}\,\Hopf{\Phi}(\eta,\tau)\Bigr). 
\end{align}
Moreover, if $\Phi$ is free boundary harmonic, then 
\begin{align}\label{eqn:circle-variation}
\frac{d}{dt}\Big\vert_{t=0}E\bigl(\gamma,c(t)\bigr)=-\int_{\partial\Omega_c}\alpha_{\dot c}\,\Hopf{\Phi}(\eta,\eta), 
\end{align}
and this vanishes for all such curves if and only if, for each $k\in\{1,\ldots,b\}$, the restriction of $\Hopf{\Phi}(\eta,\eta)$ to $\partial_k\Omega_c$ is $L^2$-orthogonal to $1$, $\cos\theta$ and $\sin\theta$. 
\end{lemma}

\begin{proof}
Let $g(t)=f_{c(t)}^*g_{\Sp^2}$. 
Since $f_{c(t)}\colon(\Omega_c,g(t))\to(\Omega_{c(t)},g_{\Sp^2})$ is an isometry, the map 
\[
\Phi(t)\vcentcolon=\Phi_{\gamma,c(t)}\circ f_{c(t)}
\]
is harmonic with respect to $g(t)$ and satisfies $E_{g(t)}(\Phi(t))=E(\gamma,c(t))$. 
Moreover, \eqref{eqn:diffeos_compatibility} implies 
\[
\Phi(t)\circ\psi_{c,k}=\Phi_{\gamma,c(t)}\circ\psi_{c(t),k}=\gamma_k
\] 
for each $k\in\{1,\ldots,b\}$. 
Hence, $W(t)\vcentcolon=\Phi-\Phi(t)$ vanishes on $\partial\Omega_c$, and integration by parts yields 
\[
E_{g(t)}(\Phi)=E_{g(t)}\bigl(\Phi(t)\bigr)+E_{g(t)}\bigl(W(t)\bigr). 
\]
Since $\Phi(t)$ is $g(t)$-harmonic, and $\Phi$ is $g(0)$-harmonic, we also have 
\begin{align}
\label{eqn:20261013a}
\int_{(\Omega_c,g(t))}\abs[\big]{dW(t)}_{g(t)}^2
=\int_{(\Omega_c,g(t))}\sk[\big]{d\Phi,dW(t)}_{g(t)}-\int_{(\Omega_c,g(0))}\sk[\big]{d\Phi,dW(t)}_{g(0)}. 
\end{align}
The metrics $g(t)$ depend smoothly on $t$, so the right-hand side 
of equation \eqref{eqn:20261013a} is bounded by $C\abs{t}\,\nm{d\Phi}_{L^2}\nm{dW(t)}_{L^2}$ for some constant $C$ and all sufficiently small $\abs{t}$. 
Therefore, $E_{g(t)}(W(t))=O(t^2)$ as $t\to0$, and 
\[
\frac{d}{dt}\Big\vert_{t=0}E\bigl(\gamma,c(t)\bigr)
=\frac{d}{dt}\Big\vert_{t=0}E_{g(t)}(\Phi)
=-\frac{1}{2}\int_{\Omega_c}\sk{\Hopf{\Phi},\mathcal{L}_Zg_{\Sp^2}}
\]
by \eqref{eqn:first-variation}, recalling $\frac{d}{dt}\big\vert_{t=0}g(t)=\mathcal{L}_Zg_{\Sp^2}$. 
Since $\gamma\in H^2(\Xi,\R^3)$, elliptic regularity implies $\Phi\in H^2(\Omega_c,\R^3)$. 
Hence, \eqref{eqn:hopf-boundary-general} applies, and \eqref{eqn:variation-c} follows from $Z=\alpha_{\dot c}\eta+\beta_{\dot c}\tau$ along $\partial\Omega_c$. 
If $\Phi$ is free boundary harmonic, then $\Hopf{\Phi}(\eta,\tau)=0$ as shown before \eqref{eqn:hopf-boundary}, which proves \eqref{eqn:circle-variation}. 
The final claim follows since, by \eqref{eqn:normal-velocity-fourier} and the subsequent discussion, the restrictions of $\alpha_{\dot c}$ to $\partial_k\Omega_c$ range over all linear combinations of $1$, $\cos\theta$ and $\sin\theta$.
\end{proof}

\begin{corollary}
\label{cor:variation-c-theta}
In the setting of Lemma~\ref{lem:critical-c}, let $\alpha_k=\alpha_{\dot c}\circ\psi_{c,k}$ and $\beta_k=\beta_{\dot c}\circ\psi_{c,k}$ for each $k\in\{1,\ldots,b\}$. 
Then 
\begin{align*} 
\frac{d}{dt}\Big\vert_{t=0}E\bigl(\gamma,c(t)\bigr)
=\sum_{k=1}^b\int_0^{2\pi}\biggl(\frac{\alpha_k}{2\sin r_k}\Bigl(\abs{\gamma_k'}^2-\abs{(\DN_c\gamma)_k}^2\Bigr)-\frac{\beta_k}{\sin r_k}\sk{\gamma_k',(\DN_c\gamma)_k}\biggr)\,d\theta. 
\end{align*}
\end{corollary}

\begin{proof}
Since $(\eta,\tau)$ is an orthonormal frame along $\partial\Omega_c$, we have $\abs{d\Phi}^2=\abs{\partial_\eta\Phi}^2+\abs{\partial_\tau\Phi}^2$ there, and the definition of $\Hopf{\Phi}$ in \eqref{eqn:first-variation} yields 
\begin{align*}
\Hopf{\Phi}(\eta,\eta)&=\tfrac{1}{2}\bigl(\abs{\partial_\eta\Phi}^2-\abs{\partial_\tau\Phi}^2\bigr), 
&
\Hopf{\Phi}(\eta,\tau)&=\sk{\partial_\eta\Phi,\partial_\tau\Phi}. 
\intertext{Differentiating $\Phi\circ\psi_{c,k}=\gamma_k$ and recalling \eqref{eqn:DtN}, we obtain} 
(\partial_\tau\Phi)\circ\psi_{c,k}&=\frac{\gamma_k'}{\sin r_k}, 
&
(\partial_\eta\Phi)\circ\psi_{c,k}&=\frac{(\DN_c\gamma)_k}{\sin r_k}. 
\end{align*}
Inserting these identities into \eqref{eqn:variation-c} and using that the arclength element of $\partial_k\Omega_c$ is $\sin(r_k)\,d\theta$ proves the claim. 
\end{proof}

\begin{proposition}
\label{prop:critical-points}
Let $(\gamma,c)\in H^2(\Xi,\Sp^2)\times\mathcal{C}$ with $\gamma$ nonconstant, and let $\Phi=\Phi_{\gamma,c}$ be its harmonic extension. 
Then $(\gamma,c)$ is a critical point of the reduced energy \eqref{eqn:reduced-energy-DtN} if and only if $\Phi$ is conformal free boundary harmonic. 
In this case, $\Phi$ is a branched free boundary minimal immersion which is smooth up to the boundary and has no branch points on $\partial\Omega_c$. 
\end{proposition}

\begin{proof}
Suppose that $(\gamma,c)$ is a critical point of the reduced energy. 
Given $V\in\mathcal{V}_\gamma$, the curve $t\mapsto(\gamma+tV)/\abs{\gamma+tV}$ in $H^2(\Xi,\Sp^2)$ is of class $C^1$ with velocity $V$ at $t=0$, since $\abs{\gamma+tV}^2=1+t^2\abs{V}^2$. 
Hence, the right-hand side of \eqref{eqn:variation-gamma} vanishes for all $V\in\mathcal{V}_\gamma$, and $\Phi$ is free boundary harmonic by Lemma~\ref{lem:critical-gamma}. 
Lemma~\ref{lem:critical-c} then implies that, for each $k\in\{1,\ldots,b\}$, the restriction of $\Hopf{\Phi}(\eta,\eta)$ to $\partial_k\Omega_c$ is $L^2$-orthogonal to $1$, $\cos\theta$ and $\sin\theta$. 
Since $\gamma$ is nonconstant, so is $\Phi$, and Proposition~\ref{prop:critical-conformal} implies that $\Phi$ is weakly conformal. 
Conversely, if $\Phi$ is conformal free boundary harmonic, then both derivatives in \eqref{eqn:20261013b} vanish by Lemma~\ref{lem:critical-gamma} and by \eqref{eqn:circle-variation}, since $\Hopf{\Phi}\equiv0$. 

In this case, $\Phi$ is smooth up to the boundary by Lemma~\ref{lem:critical-gamma}, and a branched free boundary minimal immersion by the discussion following \eqref{eqn:fbhm_eq}. 
Since $\Phi$ is nonconstant, $\abs{\Phi}^2$ is a nonconstant subharmonic function attaining its maximum $1$ on $\partial\Omega_c$. 
Hence, the strong maximum principle implies $\abs{\Phi}<1$ in the interior of $\Omega_c$, and Hopf's boundary point lemma yields $\sk{\partial_\eta\Phi,\Phi}=\frac{1}{2}\partial_\eta\abs{\Phi}^2>0$ along $\partial\Omega_c$. 
By conformality, $\abs{\partial_\tau\Phi}=\abs{\partial_\eta\Phi}>0$, i.\,e.~$d\Phi$ does not vanish on $\partial\Omega_c$ (cf.~\cite[Corollary~1.7]{DaLioPigati2020}). 
\end{proof}

\section{Gauge fixing and the perturbation problem}
\label{sec:gauge}

In the following, we assume that $b\geq3$. 
In this section, we reduce the problem of finding a critical point of the reduced energy \eqref{eqn:reduced-energy-DtN} near an approximate solution $(\gamma_0,c_0)\in H^2(\Xi,\Sp^2)\times\mathcal{C}$ to a fix point problem. 
Given $\phi\in\Conf(\Sp^2)$ and $c\in\mathcal{C}$, the restriction of $\phi$ to $\Omega_c$ is a conformal diffeomorphism onto $\Omega_{\phi\cdot c}$ by \eqref{eqn:Conf-action}. 
Since harmonicity and the Dirichlet energy of maps on surfaces are conformally invariant, the reduced energy satisfies
\begin{align}\label{eqn:energy-symmetries}
E(R\circ\gamma,c)&=E(\gamma,c), &
E(\gamma,\phi\cdot c)&=E(\gamma\circ\psi_{\phi\cdot c}^{-1}\circ\phi\circ\psi_c,c),
\end{align}
for all $R\in\SOgroup(3)$ and $\phi\in\Conf(\Sp^2)$.
This induces two distinct gauge freedoms: 
the $3$-dimensional group of rotations of the \emph{target} and the $6$-dimensional group of Möbius transformations of the \emph{domain}. 
We fix the former by imposing three orthogonality conditions on perturbations of the boundary map, and the latter by restricting the configurations to a local Möbius slice. 
The following lemma shows that criticality with respect to the boundary map $\gamma$ in the sense of Lemma~\ref{lem:critical-gamma} automatically implies criticality with respect to Möbius variations of the circle configuration $c$. 

\begin{lemma}
\label{lem:critical-Mobius}
Let $(\gamma,c)\in H^2(\Xi,\Sp^2)\times\mathcal{C}$ satisfy $P_\gamma\DN_c\gamma=0$.
Then, for every smooth curve $t\mapsto\phi(t)$ in $\Conf(\Sp^2)$ with $\phi(0)=\mathrm{id}$, we have
\begin{align*}
\frac{d}{dt}\Big\vert_{t=0}E\bigl(\gamma,\phi(t)\cdot c\bigr)=0.
\end{align*}
\end{lemma}

\begin{proof}
The curve $t\mapsto c(t)\vcentcolon=\phi(t)\cdot c$ in $\mathcal{C}$ is smooth with $c(0)=c$, and $X\vcentcolon=\frac{d}{dt}\big\vert_{t=0}\phi(t)$ is a conformal vector field on $\Sp^2$, in the sense that $\mathcal{L}_Xg_{\Sp^2}=(\divop X)\,g_{\Sp^2}$. 
By \eqref{eqn:Conf-action}, $\phi(t)$ maps $\partial_k\Omega_c$ onto $\partial_k\Omega_{c(t)}$ for each $k\in\{1,\ldots,b\}$, so \eqref{eqn:20261009a} holds with $\phi(t)$ in place of $f_{c(t)}$. 
Hence, the derivation of \eqref{eqn:normal-velocity} yields $g_{\Sp^2}(X,\eta)=\alpha_{\dot c}$ along $\partial\Omega_c$. 
Moreover, $\Phi=\Phi_{\gamma,c}$ is free boundary harmonic by Lemma~\ref{lem:critical-gamma}. 
Therefore, \eqref{eqn:circle-variation} and \eqref{eqn:hopf-boundary} applied with $X$ in place of $Z$ yield 
\[
\frac{d}{dt}\Big\vert_{t=0}E\bigl(\gamma,\phi(t)\cdot c\bigr)
=-\int_{\partial\Omega_c}\alpha_{\dot c}\,\Hopf{\Phi}(\eta,\eta)
=-\frac{1}{2}\int_{\Omega_c}(\divop X)\sk{\Hopf{\Phi},g_{\Sp^2}}
=0, 
\]
since $\Hopf{\Phi}$ is trace-free. 
\end{proof}

\paragraph{Möbius slice.}
Given $c\in\mathcal{C}$ and recalling \eqref{eqn:Conf-action}, the orbit map 
\begin{align}
\label{eqn:Conf-orbit_map}
\begin{aligned}
\Conf(\Sp^2)&\to\mathcal{C},\\ 
\phi&\mapsto\phi\cdot c
\end{aligned} 
\end{align} 
is injective since the action is free, and of constant rank since $(\phi_1\circ\phi_2)\cdot c=\phi_1\cdot(\phi_2\cdot c)$ for all $\phi_1,\phi_2\in\Conf(\Sp^2)$, where composition with $\phi_1$ and the action of $\phi_1$ on $\mathcal{C}$ are diffeomorphisms. 
Hence, \eqref{eqn:Conf-orbit_map} is an immersion by the global rank theorem (cf.~\cite[Theorem~7.15]{Lee2003}), and the image $T_c(\Conf(\Sp^2)\cdot c)$ of the differential of \eqref{eqn:Conf-orbit_map} at the identity has dimension $\dim\Conf(\Sp^2)=6$. 
Since \eqref{eqn:Conf-orbit_map} is an immersion, we may choose a smooth submanifold $U\subset\mathcal{C}$ of dimension $3b-6$ containing $c_0$ such that, after shrinking $U$ if necessary, the direct sum decomposition $T_c\mathcal{C}=T_cU\oplus T_c\bigl(\Conf(\Sp^2)\cdot c\bigr)$ holds for every $c\in U$.
We specify the choice of $U$ and the corresponding local coordinates in Section~\ref{sec:residual}.

\begin{corollary}
\label{cor:critical-slice}
Let $(\gamma,c)\in H^2(\Xi,\Sp^2)\times U$.
If $P_\gamma\DN_c\gamma=0$ and
\begin{align*}
\frac{d}{dt}\Big\vert_{t=0}E\bigl(\gamma,c(t)\bigr)=0
\end{align*}
for every smooth curve $t\mapsto c(t)$ in $U$ with $c(0)=c$, then $(\gamma,c)$ is a critical point of the reduced energy \eqref{eqn:reduced-energy-DtN}.
\end{corollary}

\begin{proof}
By Lemma~\ref{lem:critical-gamma}, $\Phi_{\gamma,c}$ is free boundary harmonic, and the first derivative in \eqref{eqn:20261013b} vanishes for every $C^1$-curve $t\mapsto\gamma(t)$ in $H^2(\Xi,\Sp^2)$ with $\gamma(0)=\gamma$.
Hence, by \eqref{eqn:circle-variation} and \eqref{eqn:normal-velocity}, the derivative of $t\mapsto E(\gamma,c(t))$ at $t=0$ depends only and linearly on $\dot c\in T_c\mathcal{C}$. 
It vanishes for $\dot c\in T_cU$ by assumption, and for $\dot c\in T_c(\Conf(\Sp^2)\cdot c)$ by Lemma~\ref{lem:critical-Mobius}. 
Therefore, it vanishes for all $\dot c\in T_c\mathcal{C}$, and $(\gamma,c)$ is a critical point by \eqref{eqn:20261013b}. 
\end{proof}

\paragraph{Perturbations of the boundary map.} 
To search for a critical point near a given boundary map $\gamma_0\colon\Xi\to\Sp^2$, it is convenient to describe nearby maps 
$\gamma\colon\Xi\to\Sp^2$ by unconstrained maps $v\colon\Xi\to\R^2$. 
The idea is to perturb $\gamma_0$ within its tangent planes, using a frame along $\gamma_0$, and to project back onto $\Sp^2$. 
Let $\gamma_0\in H^2(\Xi,\Sp^2)$ be smooth with nowhere vanishing derivative $\gamma_0'$. 
Then the columns of the matrix-valued map 
\begin{align}
\label{eqn:frame}
\Gamma\vcentcolon=\biggl(\frac{\gamma_0\times\gamma_0'}{\abs{\gamma_0'}},\ \frac{\gamma_0'}{\abs{\gamma_0'}}\biggr)\colon\Xi\to\R^{3\times2}
\end{align}
form an orthonormal basis of $T_{\gamma_0(\theta)}\Sp^2$ for every $\theta\in\Xi$, in the sense that 
\begin{align*}
\Gamma(\theta)^{\intercal}\,\Gamma(\theta)
&=\begin{pmatrix} 1&0\\0&1\end{pmatrix},
&
\Gamma(\theta)\,\Gamma(\theta)^{\intercal}&=\begin{pmatrix}1&0&0\\0&1&0\\0&0&1\end{pmatrix}-\gamma_0(\theta)\,\gamma_0(\theta)^{\intercal}.
\end{align*}
In particular, $\Gamma\Gamma^{\intercal}w=P_{\gamma_0}w$ for all $w\in H^1(\Xi,\R^3)$, recalling \eqref{eqn:Pgamma}. 
Given $v\in H^2(\Xi,\R^2)$, we then have $\abs{\gamma_0+\Gamma v}^2=1+\abs{v}^2$, and 
\begin{align}\label{eqn:perturbation}
\gamma_v\vcentcolon=\frac{\gamma_0+\Gamma v}{\sqrt{1+\abs{v}^2}}
\colon\Xi\to\Sp^2
\end{align}
is the central projection of the perturbation $\gamma_0+\Gamma v$ onto $\Sp^2$. 
By \eqref{eqn:Sobolev-embedding}, products of functions in $H^2(\Xi,\R)$, as well as their compositions with smooth functions on open intervals containing their range, belong to $H^2(\Xi,\R)$. 
Hence, $\gamma_v\in H^2(\Xi,\Sp^2)$, and $\gamma_v$ depends smoothly on $v$.
Given $c\in\mathcal{C}$, the reduced energy $E(\gamma_v,c)$ becomes a function of $v$ in the linear space $H^2(\Xi,\R^2)$. 
The following lemma computes its gradient with respect to $v$ and shows that this gradient vanishes if and only if $\gamma_v$ is critical with respect to variations of the boundary map as in Lemma~\ref{lem:critical-gamma}.
 
\begin{lemma}
\label{lem:perturbation-gradient}
Given $(v,c)\in H^2(\Xi,\R^2)\times\mathcal{C}$, let 
$\gamma_v$ be as in \eqref{eqn:perturbation}.
Then 
\begin{align*}
G_c(v)\vcentcolon=\frac{1}{\sqrt{1+\abs{v}^2}}\biggl(\Gamma^{\intercal}\DN_c\gamma_v-\frac{\sk{\gamma_v,\DN_c\gamma_v}}{\sqrt{1+\abs{v}^2}}\,v\biggr)\in H^1(\Xi,\R^2)
\end{align*}
satisfies 
\begin{align*}
\frac{d}{dt}\Big\vert_{t=0}E(\gamma_{v+tw},c)=\int_\Xi\sk{w,G_c(v)}\,d\theta
\end{align*}
for all $w\in H^2(\Xi,\R^2)$. 
Moreover, $G_c(v)=0$ if and only if $P_{\gamma_v}\DN_c\gamma_v=0$. 
\end{lemma} 

\begin{proof}
Differentiating \eqref{eqn:perturbation} yields 
\begin{align}
\label{eqn:20261017a}
\frac{d}{dt}\Big\vert_{t=0}\gamma_{v+tw}=\frac{1}{\sqrt{1+\abs{v}^2}}\biggl(\Gamma w-\frac{\sk{v,w}}{\sqrt{1+\abs{v}^2}}\,\gamma_v\biggr), 
\end{align}
and the definition of $G_c(v)$ implies the pointwise identity 
\begin{align}
\label{eqn:20261017b}
\sk{w,G_c(v)}=\Bigl\langle\frac{d}{dt}\Big\vert_{t=0}\gamma_{v+tw},\,\DN_c\gamma_v\Bigr\rangle. 
\end{align}
Since $t\mapsto\gamma_{v+tw}$ is a smooth curve in $H^2(\Xi,\Sp^2)$, the first claim follows from Lemma~\ref{lem:critical-gamma}. 

For the second claim, let $\theta\in\Xi$. 
The value of the right-hand side of \eqref{eqn:20261017a} at $\theta$ depends only on $w(\theta)$, and is in $T_{\gamma_v(\theta)}\Sp^2$ since $\gamma_{v+tw}$ maps to $\Sp^2$. 
Since $\sk{\Gamma w,\gamma_0}=0$ and $\sk{\gamma_v,\gamma_0}=(1+\abs{v}^2)^{-\frac{1}{2}}$, the inner product of \eqref{eqn:20261017a} with $\gamma_0$ is 
\begin{align*}
\Bigl\langle\frac{d}{dt}\Big\vert_{t=0}\gamma_{v+tw},\,\gamma_0\Bigr\rangle=-\frac{\sk{v,w}}{(1+\abs{v}^2)^{\frac{3}{2}}}. 
\end{align*}
Hence, if the right-hand side of \eqref{eqn:20261017a} vanishes at $\theta$, then $\sk{v,w}=0$ at $\theta$, so the vanishing of \eqref{eqn:20261017a} reduces to $\Gamma w=0$ at $\theta$, which implies $w(\theta)=\Gamma(\theta)^{\intercal}\Gamma(\theta)w(\theta)=0$. 
Thus, the linear map $\R^2\to T_{\gamma_v(\theta)}\Sp^2$ mapping $w(\theta)$ to the value of \eqref{eqn:20261017a} at $\theta$ is injective, and hence bijective since both spaces are two-dimensional.
By \eqref{eqn:20261017b}, $G_c(v)(\theta)$ vanishes if and only if $(\DN_c\gamma_v)(\theta)$ is orthogonal to $T_{\gamma_v(\theta)}\Sp^2$, which is equivalent to $(P_{\gamma_v}\DN_c\gamma_v)(\theta)=0$ by \eqref{eqn:Pgamma}. 
\end{proof}

\paragraph{Perturbations of the configuration.}
Analogously, we describe configurations in the Möbius slice $U$ near $c_0\in U$ by unconstrained parameters by choosing an open neighbourhood $B\subset \R^{3b-6}$ of the origin and a smooth map 
\begin{align}
\label{eqn:perturbation-configuration}
\begin{aligned}
B&\to U,\\
\mu&\mapsto c_\mu=\bigl(p_1(\mu),\ldots,p_b(\mu);r_1(\mu),\ldots,r_b(\mu)\bigr)
\end{aligned}
\end{align}
which is a diffeomorphism onto an open neighbourhood of $c_0=c_\mu\vert_{\mu=0}$ in $U$. 
Recalling Remark~\ref{rem:bases_dependence}, we fix the bases in \eqref{eqn:20261009b} as follows. 
For each $k\in\{1,\ldots,b\}$, we choose an orthonormal basis $(e_{1,k}(0),e_{2,k}(0))$ of $T_{p_k(0)}\Sp^2$, and for $\mu\in B$ we let $(e_{1,k}(\mu),e_{2,k}(\mu))$ be its parallel transport along the shortest geodesic in $\Sp^2$ from $p_k(0)$ to $p_k(\mu)$, where we shrink $B$ if necessary to ensure $p_k(\mu)\neq-p_k(0)$. 
Since parallel transport along a great circle is the restriction of the rotation of $\R^3$ which fixes the axis orthogonal to the great circle, the bases depend smoothly on $\mu\in B$, and so do the parametrisations $\psi_{c_\mu,k}$. 

\begin{remark}
\label{rem:parallel-transport}
With this choice of bases, the last term in \eqref{eqn:tangential-velocity-fourier} vanishes at $c_0$. 
Indeed, given a smooth curve $t\mapsto\mu(t)$ in $B$ with $\mu(0)=0$, differentiating the rotations at $t=0$ yields 
\begin{align*}
\dot e_{1,k}=(p_k\times\dot p_k)\times e_{1,k}=-\sk{\dot p_k,e_{1,k}}\,p_k, 
\end{align*}
where $p_k=p_k(0)$, and hence $\sk{\dot e_{1,k},e_{2,k}}=0$. 
\end{remark}

Given $\gamma\in H^2(\Xi,\Sp^2)$, the reduced energy $E(\gamma,c_\mu)$ becomes a function of $\mu\in B$, and its partial derivatives with respect to $\mu_1,\ldots,\mu_{3b-6}$ can be computed explicitly by Corollary~\ref{cor:variation-c-theta}, applied with $\dot c=\partial c_\mu/\partial\mu_j$. 
The following lemma relates these partial derivatives to the derivatives along curves in $U$ appearing in Corollary~\ref{cor:critical-slice}. 
It is a direct consequence of Lemma~\ref{lem:critical-c}, since the right-hand side of \eqref{eqn:variation-c} depends linearly on $\dot c$.

\begin{lemma}
\label{lem:slice-gradient}
Let $\gamma\in H^2(\Xi,\Sp^2)$ and $\mu\in B$. 
Then, for every smooth curve $t\mapsto\mu(t)$ in $B$ with $\mu(0)=\mu$, 
\begin{align*}
\frac{d}{dt}\Big\vert_{t=0}E\bigl(\gamma,c_{\mu(t)}\bigr)=\sum_{j=1}^{3b-6}\dot\mu_j\,\frac{\partial}{\partial\mu_j}E(\gamma,c_\mu), 
\end{align*}
where $\dot\mu$ denotes the velocity of $\mu(t)$ at $t=0$. 
In particular, all partial derivatives of $\mu\mapsto E(\gamma,c_\mu)$ vanish at $\mu$ if and only if the derivative of $t\mapsto E(\gamma,c(t))$ at $t=0$ vanishes for every smooth curve $t\mapsto c(t)$ in $U$ with $c(0)=c_\mu$. 
\end{lemma}

\begin{proof}
Given a smooth curve $t\mapsto\mu(t)$ in $B$ with $\mu(0)=\mu$, we have 
\begin{align*}
\dot{c}&\vcentcolon=\frac{d}{dt}\Big\vert_{t=0} c_{\mu(t)}
=\sum_{j=1}^{3b-6}\dot\mu_j\frac{\partial c_\mu}{\partial\mu_j}, 
&
\dot e_{1,k}&\vcentcolon=\frac{d}{dt}\Big\vert_{t=0} e_{1,k}(\mu(t))=\sum_{j=1}^{3b-6}\dot\mu_j\frac{\partial e_{1,k}}{\partial\mu_j}
\end{align*}
for each $k\in\{1,\ldots,b\}$. 
Hence, $\alpha_{\dot c}$ and $\beta_{\dot c}$ depend linearly on $\dot\mu$ by \eqref{eqn:normal-velocity} and \eqref{eqn:tangential-velocity-fourier}, and Lemma~\ref{lem:critical-c} yields $\ell_1,\ldots,\ell_{3b-6}\in\R$, depending on $\gamma$ and $\mu$ but not on the curve, such that 
\begin{align}
\label{eqn:20261017c}
\frac{d}{dt}\Big\vert_{t=0}E\bigl(\gamma,c_{\mu(t)}\bigr)=\sum_{j=1}^{3b-6}\dot\mu_j\,\ell_j. 
\end{align}
For the straight line $t\mapsto\mu(t)$ in the $j$-th unit direction, the left-hand side of \eqref{eqn:20261017c} is by definition the partial derivative of $\mu\mapsto E(\gamma,c_\mu)$ with respect to $\mu_j$, while the right-hand side equals $\ell_j$. 
Hence, this partial derivative exists and equals $\ell_j$, which proves the first claim.
  
The second claim follows since \eqref{eqn:perturbation-configuration} is a diffeomorphism onto an open subset of $U$, so that every smooth curve in $U$ through $c_\mu$ is, for small $\abs{t}$, of the form $t\mapsto c_{\mu(t)}$, and every $\dot\mu\in\R^{3b-6}$ occurs.
\end{proof}

\paragraph{Rotation gauge.}
By \eqref{eqn:energy-symmetries}, the reduced energy is invariant under rotations of the boundary map. 
Consequently, if $(\gamma,c)$ is a critical point of \eqref{eqn:reduced-energy-DtN}, then so is $(R\circ\gamma,c)$ for every $R\in\SOgroup(3)$, which means that critical points of the reduced energy are not isolated. 
In order to obtain isolated solutions near the given boundary map $\gamma_0\colon\Xi\to\Sp^2$, we fix a basis $\xi_1,\xi_2,\xi_3$ of $\R^3$ and, recalling the frame $\Gamma\colon\Xi\to\R^{3\times2}$ along $\gamma_0$ from \eqref{eqn:frame}, define
\begin{align}
\label{eqn:rotation-fields}
\rho_i\vcentcolon=\Gamma^{\intercal}(\xi_i\times\gamma_0)\in H^2(\Xi,\R^2)
\end{align}
for $i\in\{1,2,3\}$. 
Since $\xi_i\times\gamma_0$ is tangent to $\Sp^2$ along $\gamma_0$, we have $\Gamma\rho_i=\xi_i\times\gamma_0$. 
Thus, the components of $\rho_i$ 
are the components of the velocity of $\gamma_0$ under the rotations about the axis through $\xi_i$ with respect to the frame $\Gamma$. 

\begin{definition}[Gradient map]
\label{defn:gradient-map}
Recalling $B\subset\R^{3b-6}$ from \eqref{eqn:perturbation-configuration}, $G_c(v)$ from Lemma~\ref{lem:perturbation-gradient} and $\rho_1,\rho_2,\rho_3$ from \eqref{eqn:rotation-fields}, we define 
\begin{align*}
F\colon H^2(\Xi,\R^2)\times B\times\R^3&\to H^1(\Xi,\R^2)\times\R^{3b-6}\times\R^3
\\
(v,\mu,\lambda)&\mapsto\biggl(G_{c_\mu}(v)+\sum_{i=1}^3\lambda_i\rho_i,\;
\Bigl(\frac{\partial}{\partial\mu_j}E(\gamma_v,c_\mu)\Bigr)_{j=1,\ldots,3b-6},\;
\Bigl(\int_\Xi\sk{v,\rho_i}\,d\theta\Bigr)_{i=1,\ldots,3}\biggr). 
\end{align*}
\end{definition}

\begin{remark}[Lagrange multipliers]
\label{rem:Lagrange_multipliers}
The partial derivatives in Definition~\ref{defn:gradient-map} exist by Lemma~\ref{lem:slice-gradient}. 
The components of the gradient map $F$ are the partial derivatives of the Lagrangian 
\begin{align*}
L(v,\mu,\lambda)\vcentcolon=E(\gamma_v,c_\mu)+\sum_{i=1}^3\lambda_i\int_\Xi\sk{v,\rho_i}\,d\theta 
\end{align*}
with respect to $v$, $\mu$ and $\lambda$, where the derivative with respect to $v$ is given by Lemma~\ref{lem:perturbation-gradient} via the $L^2$-inner product. 
Formally, the zeros of $F$ are therefore the critical points of $(v,\mu)\mapsto E(\gamma_v,c_\mu)$ subject to the constraints
\begin{align}
\label{eqn:rho_constraints}
\int_\Xi\sk{v,\rho_1}\,d\theta= 
\int_\Xi\sk{v,\rho_2}\,d\theta= 
\int_\Xi\sk{v,\rho_3}\,d\theta=0, 
\end{align}
and $\lambda_1,\lambda_2,\lambda_3$ are the corresponding Lagrange multipliers. 
Heuristically, since $\rho_1,\rho_2,\rho_3$ describe the rotations of $\gamma_0$, the constraints \eqref{eqn:rho_constraints} select among rotated boundary maps near $\gamma_0$ those whose coordinates in \eqref{eqn:perturbation} are $L^2$-orthogonal to the directions of rotation.
\end{remark}

The following lemma shows that the Lagrange multipliers in Remark~\ref{rem:Lagrange_multipliers} vanish at every zero of the gradient map $F$, provided a certain $3\times3$ matrix is invertible, and hence that zeros of $F$ yield critical points of the reduced energy without constraints.
Its proof requires a generalisation of \eqref{eqn:rotation-fields} of the form 
\begin{align}
\label{eqn:rotation-fields-v}
\rho_i^v\vcentcolon=\rho_i+\Gamma^{\intercal}(\xi_i\times\Gamma v)-\sk{\gamma_0,\xi_i\times\Gamma v}\,v\in H^2(\Xi,\R^2)
\end{align}
for $i\in\{1,2,3\}$ and $v\in H^2(\Xi,\R^2)$, such that $\rho_i^0=\rho_i$. 
In the proof of Lemma~\ref{lem:zeros-F}, we show that $\rho_i^v$ describes, in terms of the coordinate $v$, the velocity of the perturbation $\gamma_v$ from \eqref{eqn:perturbation} under the rotations about the axis through $\xi_i$.

\begin{lemma}
\label{lem:zeros-F}
Let $(v,\mu,\lambda)\in H^2(\Xi,\R^2)\times B\times\R^3$ with $F(v,\mu,\lambda)=0$, and suppose that the matrix $\mathcal{R}^v\in\R^{3\times3}$ given by 
\begin{align*}
\mathcal{R}^v_{ij}=\int_\Xi\sk{\rho_j,\rho_i^v}\,d\theta
\end{align*}
for $i,j\in\{1,2,3\}$ is invertible. 
Then $\lambda=0$, the boundary map $\gamma_v$ is nonconstant, and $(\gamma_v,c_\mu)$ is a critical point of the reduced energy \eqref{eqn:reduced-energy-DtN}. 
\end{lemma}

\begin{proof} 
Let $c=c_\mu$ as in \eqref{eqn:perturbation-configuration}. 
Given $w\in H^2(\Xi,\R^2)$, we have $\sk{\Gamma w,\gamma_0}=0$ and $\sk{\Gamma w,\Gamma v}=\sk{w,v}$, and hence $\sk{\Gamma w,\gamma_v}=\sk{w,v}/\sqrt{1+\abs{v}^2}$ by \eqref{eqn:perturbation}. 
Therefore, recalling \eqref{eqn:20261017a} and the projection $P_{\gamma_v}$ from \eqref{eqn:Pgamma}, we have 
\begin{align}
\label{eqn:20261019c}
\frac{d}{dt}\Big\vert_{t=0}\gamma_{v+tw}=\frac{1}{\sqrt{1+\abs{v}^2}}\biggl(\Gamma w-\frac{\sk{v,w}}{\sqrt{1+\abs{v}^2}}\,\gamma_v\biggr)
=\frac{1}{\sqrt{1+\abs{v}^2}}\,P_{\gamma_v}(\Gamma w).
\end{align}
Recalling $\rho_i^v$ from \eqref{eqn:rotation-fields-v} as well as $\Gamma\rho_i=\xi_i\times\gamma_0$ and $\Gamma\Gamma^{\intercal}y=y-\sk{y,\gamma_0}\gamma_0$ for any $y\in\R^3$, in particular for $y=\xi_i\times\Gamma v$, we obtain, for each $i\in\{1,2,3\}$,
\begin{align*} 
\Gamma\rho_i^v
&=\xi_i\times\gamma_0+\Bigl(\xi_i\times\Gamma v-\sk{\gamma_0,\xi_i\times\Gamma v}\gamma_0\Bigr)-\sk{\gamma_0,\xi_i\times\Gamma v}\Gamma v
\\ 
&=\xi_i\times(\gamma_0+\Gamma v)-\sk{\gamma_0,\xi_i\times\Gamma v}\,(\gamma_0+\Gamma v).
\end{align*} 
The term $\gamma_0+\Gamma v=\sqrt{1+\abs{v}^2}\,\gamma_v$ is radial, so that $P_{\gamma_v}(\gamma_0+\Gamma v)=0$,
while $\xi_i\times\gamma_v$ is tangent to $\Sp^2$ along $\gamma_v$. 
Hence, inserting $w=\rho_i^v$ in \eqref{eqn:20261019c} yields 
\begin{align}
\label{eqn:20261019a}
\frac{d}{dt}\Big\vert_{t=0}\gamma_{v+t\rho_i^v}=P_{\gamma_v}(\xi_i\times\gamma_v)=\xi_i\times\gamma_v. 
\end{align}
Let $t\mapsto R(t)$ be the rotations about the axis through $\xi_i$ such that $\frac{d}{dt}\big\vert_{t=0}R(t)y=\xi_i\times y$ for all $y\in\R^3$. 
Then $t\mapsto R(t)\circ\gamma_v$ is a smooth curve in $H^2(\Xi,\Sp^2)$ with velocity $\xi_i\times\gamma_v$ at $t=0$, along which the reduced energy is constant by \eqref{eqn:energy-symmetries}. 
Hence, Lemma~\ref{lem:critical-gamma}, \eqref{eqn:20261017b} and \eqref{eqn:20261019a} imply 
\begin{align}
\label{eqn:20261019b}
\int_\Xi\sk{\rho_i^v,G_c(v)}\,d\theta=\int_\Xi\sk{\xi_i\times\gamma_v,\DN_c\gamma_v}\,d\theta=0. 
\end{align}
By assumption, the first component of $F(v,\mu,\lambda)$ vanishes, which 
by Definition~\ref{defn:gradient-map} means that 
\begin{align}
\label{eqn:20261019e}
G_c(v)&=-\lambda_1\rho_1-\lambda_2\rho_2-\lambda_3\rho_3.
\end{align}
Inserting \eqref{eqn:20261019e} in \eqref{eqn:20261019b} yields $\sum_{j=1}^3\mathcal{R}^v_{ij}\lambda_j=0$ for each $i\in\{1,2,3\}$, and therefore $\lambda=0$ as claimed since $\mathcal{R}^v$ is invertible by assumption. 
Consequently, $G_c(v)=0$, and Lemma~\ref{lem:perturbation-gradient} implies $P_{\gamma_v}\DN_c\gamma_v=0$. 
Moreover, the second component of $F(v,\mu,\lambda)$ vanishes, and Lemma~\ref{lem:slice-gradient} implies that the derivative of $t\mapsto E(\gamma_v,c(t))$ at $t=0$ vanishes for every smooth curve $t\mapsto c(t)$ in $U$ with $c(0)=c$. 
Hence, $(\gamma_v,c)$ is a critical point of the reduced energy by Corollary~\ref{cor:critical-slice}. 

It remains to show that $\gamma_v$ is nonconstant. 
Suppose to the contrary that the map $\gamma_v$ is constant with value 
$y=\sum_{i=1}^3 a_i\xi_i\in \Sp^2$ for some $a\in\R^3\setminus\{0\}$. 
Since the right-hand side of \eqref{eqn:20261019c} is linear in $w$, equation \eqref{eqn:20261019a} implies that $w=\sum_{i=1}^3 a_i\rho_i^v$ satisfies 
\begin{align*}
P_{\gamma_v}(\Gamma w)=\sqrt{1+\abs{v}^2}\;y\times\gamma_v=\sqrt{1+\abs{v}^2}\;y\times y=0, 
\end{align*}
which means that $\Gamma w=\sk{\Gamma w,\gamma_v}\gamma_v$. 
Taking the inner product with $\gamma_0$ and using $\sk{\Gamma w,\gamma_0}=0$ and $\sk{\gamma_v,\gamma_0}=(1+\abs{v}^2)^{-1/2}>0$, we obtain $\sk{\Gamma w,\gamma_v}=0$, and hence $\Gamma w=0$ and $w=\Gamma^{\intercal}\Gamma w=0$.
Therefore, 
\[
\sum_{i=1}^3 a_i\mathcal{R}^v_{ij}=\int_\Xi\sk{\rho_j,w}\,d\theta=0
\]
for each $j\in\{1,2,3\}$, which contradicts the invertibility of $\mathcal{R}^v$. 
\end{proof}

Our fix point argument uses quantitative estimates for the Fréchet derivative of $F$, whereas Lemmata~\ref{lem:critical-c}, \ref{lem:perturbation-gradient} and~\ref{lem:slice-gradient} only provide derivatives along curves. 
The following lemma provides the required regularity. 
Its key ingredient is that the normal derivative of the harmonic extension depends smoothly on the positions of the boundary circles. 
The idea of the proof is to pull the harmonic extensions back to the fixed circle domain $\Omega_{c_0}$. 
Moving the boundary circles then amounts to varying the metric smoothly, while the boundary data remain fixed. 
After subtracting a fixed extension of the boundary data, the harmonic extension is obtained by inverting the Laplacian on maps vanishing along the boundary. 

\begin{lemma}
\label{lem:smoothness}
The Dirichlet-to-Neumann operator $\DN_{c_\mu}\colon H^2(\Xi,\R^3)\to H^1(\Xi,\R^3)$ from \eqref{eqn:DtN} depends smoothly on $\mu\in B$ with respect to the operator norm. 
Moreover, the gradient map $F$ from Definition~\ref{defn:gradient-map} is smooth. 
\end{lemma}

\begin{proof}
Smoothness is local in $\mu$.
We give the proof near $0$ but the same proof applies near any point in $B$.
For $\mu$ sufficiently close to $0$, choose diffeomorphisms $f_\mu\colon\Omega_{c_0}\to\Omega_{c_\mu}$ as in \eqref{eqn:diffeos_f_c}--\eqref{eqn:diffeos_compatibility}, depending smoothly on $\mu\in B$, with $f_0=\mathrm{id}$ and $f_\mu\circ\psi_{c_0,k}=\psi_{c_\mu,k}$ for each $k\in\{1,\ldots,b\}$. 
Then the pulled back metrics $g_\mu=f_\mu^*g_{\Sp^2}$ on $\Omega_{c_0}$ depend smoothly on $\mu$, in the sense that their coefficients in local coordinates are smooth functions of $(x,\mu)$ for $\mu$ sufficiently close to $0$. 
Given $\gamma\in H^2(\Xi,\R^3)$, let 
\[
h=\gamma\circ\psi_{c_0}^{-1}\in H^2(\partial\Omega_{c_0},\R^3).
\]
Since $f_\mu\colon(\Omega_{c_0},g_\mu)\to(\Omega_{c_\mu},g_{\Sp^2})$ is an isometry, the map $u_\mu=\Phi_{\gamma,c_\mu}\circ f_\mu\in H^1(\Omega_{c_0},\R^3)$ is the $g_\mu$-harmonic extension of $h$, and \eqref{eqn:DtN} yields 
\begin{align}
\label{eqn:20261020a}
(\DN_{c_\mu}\gamma)_k=\sin(r_k(\mu))\,(\partial_{\eta_\mu}u_\mu)\circ\psi_{c_0,k}, 
\end{align}
where $\eta_\mu$ denotes the outward unit conormal of $\partial\Omega_{c_0}$ with respect to $g_\mu$. 
In order to express $u_\mu$ in terms of operators between fixed function spaces, we recall from \cite[Chapter~5, Proposition~1.7]{Taylor2011} that the $g_0$-harmonic extension 
\begin{align*}
H^2(\partial\Omega_{c_0},\R^3)&\to H^{\frac{5}{2}}(\Omega_{c_0},\R^3), \\
h&\mapsto u_0
\end{align*}
is a bounded operator. 
By \cite[Chapter~5, Theorem~1.3]{Taylor2011} and interpolation as in the proof of \cite[Chapter~5, Proposition~1.7]{Taylor2011}, the restriction
\begin{align*}
L_\mu\colon H^{\frac{5}{2}}(\Omega_{c_0},\R^3)\cap H^1_0(\Omega_{c_0},\R^3)\to H^{\frac{1}{2}}(\Omega_{c_0},\R^3)
\end{align*}
of the Laplace--Beltrami operator $\Delta_{g_\mu}\colon H^{\frac{5}{2}}(\Omega_{c_0},\R^3)\to H^{\frac{1}{2}}(\Omega_{c_0},\R^3)$ is an isomorphism. 
Therefore, $w_\mu=L_\mu^{-1}\Delta_{g_\mu}u_0$ is the unique solution in $H^{\frac{5}{2}}(\Omega_{c_0},\R^3)$ of 
\begin{align*}
\left\{
\begin{aligned}
\Delta_{g_\mu}w_\mu&=\Delta_{g_\mu}u_0 
&&\text{ in $\Omega_{c_0}$,} \\
w_\mu&=0 
&&\text{ on $\partial\Omega_{c_0}$.}
\end{aligned}\right.
\end{align*}
Since $u_0$ coincides with $h$ on $\partial\Omega_{c_0}$, the difference $u_0-w_\mu$ is $g_\mu$-harmonic with boundary values $h$, and uniqueness of the harmonic extension yields 
\begin{align}
\label{eqn:20261020b}
u_\mu=u_0-L_\mu^{-1}\Delta_{g_\mu}u_0. 
\end{align}
The coefficients of $\Delta_{g_\mu}$ are smooth functions of $(x,\mu)\in\Omega_{c_0}\times B$. 
Hence, $\Delta_{g_\mu}$ and $L_\mu$ depend smoothly on $\mu$ with respect to the operator norm, and so does $L_\mu^{-1}$, since inversion is smooth on the open set of invertible operators. 
Similarly, by \cite[Chapter~4, Proposition~4.5]{Taylor2011}, the conormal derivative $u\mapsto\partial_{\eta_\mu}u$ is a bounded operator $H^{\frac{5}{2}}(\Omega_{c_0},\R^3)\to H^1(\partial\Omega_{c_0},\R^3)$ depending smoothly on $\mu$. 
Since $\gamma\mapsto h\mapsto u_0$ does not depend on $\mu$, the first claim follows from \eqref{eqn:20261020a} and \eqref{eqn:20261020b}.

For the second claim, we recall from \eqref{eqn:perturbation} that $v\mapsto\gamma_v$ is a smooth map $H^2(\Xi,\R^2)\to H^2(\Xi,\R^3)$, and that pointwise products are bounded bilinear maps $H^2\times H^1\to H^1$ by \eqref{eqn:Sobolev-embedding}. 
Hence, the first component of $F$, given explicitly in Lemma~\ref{lem:perturbation-gradient}, is smooth by the first claim. 
By \eqref{eqn:reduced-energy-DtN}, the second component of $F$ is given by 
\begin{align*}
\frac{\partial}{\partial\mu_j}E(\gamma_v,c_\mu)=\frac{1}{2}\int_\Xi\Bigl\langle\gamma_v,\frac{\partial\DN_{c_\mu}}{\partial\mu_j}\gamma_v\Bigr\rangle\,d\theta, 
\end{align*}
which is smooth by the first claim. 
The third component of $F$ is linear and bounded. 
\end{proof}

\paragraph{Linearisation.}
Lemma~\ref{lem:zeros-F} relates a zero of the gradient map $F$ near the origin with a critical point of the reduced energy. 
Our fix point argument relies on the nondegeneracy of the linearisation of $F$, which involves derivatives of the Dirichlet-to-Neumann operator with respect to the configuration.

\begin{lemma}
\label{lem:DtN-derivative}
Given a smooth curve $t\mapsto\mu(t)$ in $B\subset\R^{3b-6}$ from \eqref{eqn:perturbation-configuration}, let $c(t)=c_{\mu(t)}$ and $c=c(0)$.
Recalling the boundary parametrisations $\psi_{c,k}$ from \eqref{eqn:20261009b} and the normal and tangential velocities $\alpha_{\dot c}$ and $\beta_{\dot c}$ from \eqref{eqn:normal-velocity-fourier}--\eqref{eqn:tangential-velocity-fourier}, let $\tilde\alpha,\tilde\beta\colon\Xi\to\R$ be defined on the $k$-th circle by
\begin{align*}
\tilde\alpha&=\frac{\alpha_{\dot c}\circ\psi_{c,k}}{\sin r_k},
&
\tilde\beta&=\frac{\beta_{\dot c}\circ\psi_{c,k}}{\sin r_k}
\end{align*}
for each $k\in\{1,\ldots,b\}$.
Then, for all $\zeta,\gamma\in H^2(\Xi,\R^3)$,
\begin{align*}
\int_\Xi\sk[\Big]{\zeta,
\Bigl(\frac{d}{dt}\Big\vert_{t=0}\DN_{c(t)}\Bigr)\gamma}\,d\theta
&=\int_\Xi
\tilde\alpha\bigl(\sk{\zeta',\gamma'}-\sk{\DN_c\zeta,\DN_c\gamma}\bigr)
-\tilde\beta\bigl(\sk{\zeta',\DN_c\gamma}+\sk{\gamma',\DN_c\zeta}\bigr)
\,d\theta.
\end{align*}
\end{lemma}

\begin{proof}
By Lemma~\ref{lem:smoothness}, the derivative
\begin{align*}
\dot{\DN}\vcentcolon=\frac{d}{dt}\Big\vert_{t=0}\DN_{c(t)}
\colon H^2(\Xi,\R^3)\to H^1(\Xi,\R^3)
\end{align*}
is a bounded linear operator.
Moreover, $\dot{\DN}$ is symmetric with respect to the $L^2$-inner product by \eqref{eqn:DtN-Green}.
Consequently, \eqref{eqn:reduced-energy-DtN} and Corollary~\ref{cor:variation-c-theta} imply
\begin{align*}
\frac{1}{2}\int_\Xi\sk{\gamma,\dot{\DN}\gamma}\,d\theta
&=\frac{d}{dt}\Big\vert_{t=0}E(\gamma,c(t))
=\frac{1}{2}\int_\Xi
\tilde\alpha\bigl(\abs{\gamma'}^2-\abs{\DN_c\gamma}^2\bigr)
-2\tilde\beta\sk{\gamma',\DN_c\gamma}\,d\theta
\end{align*}
for every $\gamma\in H^2(\Xi,\R^3)$.
The claim follows by polarisation, i.\,e.~by applying this identity with $\gamma$ replaced by $\zeta+\gamma$ and subtracting the identities for $\zeta$ and $\gamma$.
\end{proof}

By Remark~\ref{rem:Lagrange_multipliers} and Lemma~\ref{lem:smoothness}, the components of the gradient map $F$ are the partial derivatives of the smooth Lagrangian $L$ from Remark~\ref{rem:Lagrange_multipliers}. 
The following statement obtains $DF\vert_0$ by computing the second derivatives of $L$ at the origin.

\begin{proposition}
\label{prop:linearisation}
The derivative of the gradient map $F$ from Definition~\ref{defn:gradient-map} at the origin is 
\begin{align*}
DF\vert_0
\begin{pmatrix} w\\ h\\ \ell\end{pmatrix}
=\begin{pmatrix}
S_{vv} & S_{v\mu} & \rho\\
S_{v\mu}^* & S_{\mu\mu} & 0\\
\rho^* & 0 & 0
\end{pmatrix}
\begin{pmatrix} w\\ h\\ \ell\end{pmatrix}
=\begin{pmatrix}
S_{vv}w+S_{v\mu}h+\rho\,\ell\\
S_{v\mu}^*w+S_{\mu\mu}h\\
\rho^*w
\end{pmatrix}
\in 
H^1(\Xi,\R^2)\times\R^{3b-6}\times\R^3
\end{align*}
for $(w,h,\ell)\in H^2(\Xi,\R^2)\times\R^{3b-6}\times\R^3$, where 
\begin{align*}
S_{vv}\colon H^2(\Xi,\R^2)&\to H^1(\Xi,\R^2), 
&
S_{vv}w&=\Gamma^{\intercal}\DN_{c_0}(\Gamma w)-\sk{\gamma_0,\DN_{c_0}\gamma_0}\,w,
\\
S_{v\mu}\colon\R^{3b-6}&\to H^1(\Xi,\R^2), 
&
S_{v\mu}h&=\sum_{j=1}^{3b-6}h_j\,\Gamma^{\intercal}\Bigl(\frac{\partial}{\partial\mu_j}\Big\vert_{\mu=0}\DN_{c_\mu}\Bigr)\gamma_0,
\\
S_{v\mu}^*\colon H^2(\Xi,\R^2)&\to\R^{3b-6}, 
&
(S_{v\mu}^*w)_j&=\int_\Xi\sk[\Big]{\Gamma w,\Bigl(\frac{\partial}{\partial\mu_j}\Big\vert_{\mu=0}\DN_{c_\mu}\Bigr)\gamma_0}\,d\theta,
\\
S_{\mu\mu}\colon\R^{3b-6}&\to\R^{3b-6}, 
&
(S_{\mu\mu}h)_i&=\frac12\sum_{j=1}^{3b-6}h_j\int_\Xi\sk[\Big]{\gamma_0,\Bigl(\frac{\partial^2}{\partial\mu_i\partial\mu_j}\Big\vert_{\mu=0}\DN_{c_\mu}\Bigr)\gamma_0}\,d\theta,
\\
\rho\colon\R^3&\to H^1(\Xi,\R^2), 
&
\rho\,\ell&=\sum_{i=1}^3\ell_i\rho_i,
\\
\rho^*\colon H^2(\Xi,\R^2)&\to\R^3, 
&
(\rho^*w)_i&=\int_\Xi\sk{w,\rho_i}\,d\theta.
\end{align*} 
The upper left $2\times2$ block of $DF\vert_0$ represents the Hessian of $(v,\mu)\mapsto E(\gamma_v,c_\mu)$ at the origin with respect to the $L^2$-inner product on $H^2(\Xi,\R^2)$ and the Euclidean inner product on $\R^{3b-6}$.
Moreover, $DF\vert_0$ is symmetric with respect to the product of these inner products and the Euclidean inner product on $\R^3$.
\end{proposition}

\begin{proof}
Recalling $\gamma_v$ from \eqref{eqn:perturbation}, we have for all $w,z\in H^2(\Xi,\R^2)$ 
\begin{align*}
{D\gamma_v}\vert_{v=0}(w)&=\Gamma w,
&
{D^2\gamma_v}\vert_{v=0}(w,z)&=-\sk{w,z}\gamma_0.
\end{align*} 
Hence, recalling the Lagrangian $L$ from Remark~\ref{rem:Lagrange_multipliers}, we obtain from \eqref{eqn:reduced-energy-DtN} that
\begin{align*}
D^2_{vv}L\vert_0(w,w)
&=\int_\Xi
\sk{\Gamma w,\DN_{c_0}(\Gamma w)}
-\sk{\gamma_0,\DN_{c_0}\gamma_0}\abs{w}^2
\,d\theta
=\int_\Xi\sk{w,S_{vv}w}\,d\theta.
\end{align*}
By Lemmata~\ref{lem:perturbation-gradient} and \ref{lem:smoothness}, for every $h\in\R^{3b-6}$ we have
\begin{align*}
D^2_{v\mu}L\vert_0(w,h)
&=\sum_{j=1}^{3b-6}h_j
\int_\Xi\sk[\Big]{\Gamma w,
\Bigl(\frac{\partial}{\partial\mu_j}\Big\vert_{\mu=0}\DN_{c_\mu}\Bigr)\gamma_0}\,d\theta
=\int_\Xi\sk{w,S_{v\mu}h}\,d\theta.
\end{align*}
Since $L$ is smooth, its Hessian is symmetric, and therefore
\begin{align*}
D^2_{\mu v}L\vert_0(h,w)
&=D^2_{v\mu}L\vert_0(w,h)
=\sum_{j=1}^{3b-6}h_j(S_{v\mu}^*w)_j.
\end{align*}
Since the Lagrange multiplier term in $L$ is independent of $\mu$, \eqref{eqn:reduced-energy-DtN} implies
\begin{align*}
D^2_{\mu\mu}L\vert_0(h,h)
&=\frac12\sum_{i,j=1}^{3b-6}h_i h_j
\int_\Xi\sk[\Big]{\gamma_0,
\Bigl(\frac{\partial^2}{\partial\mu_i\partial\mu_j}
\Big\vert_{\mu=0}\DN_{c_\mu}\Bigr)\gamma_0}\,d\theta
=\sum_{i=1}^{3b-6}h_i(S_{\mu\mu}h)_i
\end{align*}
for all $h\in\R^{3b-6}$.
Since $S_{vv}$ and $S_{\mu\mu}$ are symmetric, the equalities of quadratic forms determine the $vv$- and $\mu\mu$-derivatives of $L$. 
Finally, we have
\begin{align*}
D^2_{v\lambda}L\vert_0(w,\ell)
&=\int_\Xi\sk{w,\rho\,\ell}\,d\theta,
&
D^2_{\mu\lambda}L\vert_0&=0,
&
D^2_{\lambda\lambda}L\vert_0&=0.
\end{align*}
By symmetry of the Hessian, these identities yield the remaining blocks of $DF\vert_0$.
\end{proof}
 
\paragraph{Fix point argument.}
To prove that $F$ has a zero near the origin, we use the following quantitative version of Banach's fix point theorem.

\begin{theorem}
\label{thm:contraction}
Let $X$ and $Y$ be Banach spaces.
Given $r>0$, let
$X_r\vcentcolon=\{x\in X\mid\nm{x}_X\leq r\}$, let $F\colon X_r\to Y$ be Fréchet differentiable at the origin, and define
\begin{align*}
\mathscr{L}&\vcentcolon=DF\vert_0,
&
\mathscr{N}(x)&\vcentcolon=F(x)-F(0)-\mathscr{L}x.
\end{align*}
Suppose that $\mathscr{L}\colon X\to Y$ is an isomorphism, and that $\varepsilon,q,K\geq0$ satisfy
\begin{align}
\label{eqn:assumptions}
{\nm{F(0)}}_Y&\leq\varepsilon,
&
\nm{\mathscr{L}^{-1}}_{Y\to X}&\leq q,
&
\nm{\mathscr N(x)-\mathscr N(y)}_Y&\leq K\nm{x-y}_X 
\end{align}
for all $x,y\in X_r$. 
If
\begin{align}
\label{eqn:radii-conditions}
q\varepsilon+qKr&\leq r,
&
qK&<1,
\end{align}
then $F$ has a unique zero in $X_r$.
\end{theorem}

\begin{proof}
Let $T(x)\vcentcolon=-\mathscr{L}^{-1}\bigl(F(0)+\mathscr{N}(x)\bigr)$. 
Since $\mathscr{N}(0)=0$, assumptions \eqref{eqn:assumptions} and \eqref{eqn:radii-conditions} imply
\begin{align*}
\nm{T(x)}_X
&\leq q\varepsilon+qK\nm{x}_X
\leq q\varepsilon+qKr
\leq r,
\\[1ex]
\nm{T(x)-T(y)}_X
&\leq q\nm{\mathscr{N}(x)-\mathscr{N}(y)}_Y
\leq qK\nm{x-y}_X
\end{align*}
for all $x,y\in X_r$.
Since $qK<1$, the map $T$ is a contraction of the complete metric space $X_r$.
By the contraction mapping principle \cite[Chapter~1, Theorem~2.2]{Taylor2011}, $T$ has a unique fix point in $X_r$. 
By definition, $x=T(x)$ if and only if $F(x)=0$.  
\end{proof}

Combining Theorem~\ref{thm:contraction} with Lemma~\ref{lem:zeros-F}, we obtain the following criterion for the existence of a free boundary minimal immersion near the given data $(\gamma_0,c_0)$. 

\begin{proposition}
\label{prop:existence}
Recalling $B\subset\R^{3b-6}$ from \eqref{eqn:perturbation-configuration}, let $r>0$ such that $\abs{\mu}\leq r$ implies $\mu\in B$. 
Let 
\begin{align*}
X&=H^2(\Xi,\R^2)\times\R^{3b-6}\times\R^3, 
&
Y&=H^1(\Xi,\R^2)\times\R^{3b-6}\times\R^3 
\end{align*}
be equipped with the respective norms $\nm{(v,\mu,\lambda)}^2=\nm{v}^2+\abs{\mu}^2+\abs{\lambda}^2$, where $\nm{v}$ denotes the $H^2$-norm in $X$ and the $H^1$-norm in $Y$, as defined in \eqref{eqn:Sobolev-norm}, and let
$F\colon H^2(\Xi,\R^2)\times B\times\R^3\to Y$ be as in Definition~\ref{defn:gradient-map}.
Suppose that 
\begin{enumerate}[label={\normalfont(\alph*)}]
\item\label{prop:existence-a}
the derivative $DF\vert_0\colon X\to Y$ is an isomorphism and there exist $\varepsilon,q,K\geq0$ satisfying \eqref{eqn:assumptions}--\eqref{eqn:radii-conditions},
\item\label{prop:existence-b} 
the matrix $\mathcal{R}^v$ from Lemma~\ref{lem:zeros-F} is invertible for every $v\in H^2(\Xi,\R^2)$ with $\nm{v}_{H^2}\leq r$.
\end{enumerate} 
Then there exist $v\in H^2(\Xi,\R^2)$ and $\mu\in B$ with $F(v,\mu,0)=0$ and $\nm{v}_{H^2}^2+\abs{\mu}^2\leq r^2$.
Moreover, $\gamma_v$ is nonconstant and $(\gamma_v,c_\mu)$ is a critical point of the reduced energy \eqref{eqn:reduced-energy-DtN}. 
In particular, the harmonic extension $\Phi_{\gamma_v,c_\mu}\colon\Omega_{c_\mu}\to\B^3$ is a branched free boundary minimal immersion which is smooth up to the boundary and has no branch points on $\partial\Omega_{c_\mu}$. 
\end{proposition}

\begin{proof}
Since $B$ is open and contains the closed ball of radius $r$ in $\R^{3b-6}$, the gradient map $F$ is defined on an open neighbourhood of $X_r$ in $X$, and smooth there by Lemma~\ref{lem:smoothness}. 
Hence, Theorem~\ref{thm:contraction} yields a zero $(v,\mu,\lambda)\in X_r$ of $F$. 
Since $\nm{v}_{H^2}\leq r$, the matrix $\mathcal{R}^v$ is invertible, and Lemma~\ref{lem:zeros-F} implies that $\lambda=0$, that $\gamma_v$ is nonconstant, and that $(\gamma_v,c_\mu)$ is a critical point of the reduced energy. 
The last claim follows from Proposition~\ref{prop:critical-points}. 
\end{proof}


\section{Quantitative estimates for the Dirichlet-to-Neumann operator}
\label{sec:Dirichlet-to-Neumann}
  
Verifying the assumptions of Proposition~\ref{prop:existence} requires estimates for the Dirichlet-to-Neumann operator $\DN_c$ defined in \eqref{eqn:DtN}. 
The main idea is to compare $\DN_c$ with the direct sum of the Dirichlet-to-Neumann operators associated with the individual boundary circles.
This decoupled operator has a simple Fourier multiplier representation, and its difference from $\DN_c$ is regularising. 
We use this comparison to derive bounds for $\DN_c$ and its derivatives with respect to the configuration, and to approximate $\DN_c$ by the method of fundamental solutions with rigorous error bounds.

\begin{definition}[Decoupled Dirichlet-to-Neumann operator]
Given $c=(p_1,\ldots,p_b;r_1,\ldots,r_b)\in\mathcal{C}$ and $\Xi$ as in \eqref{eqn:circle_union}, we define  
\[
\DN_{\mathrm{dec}}\colon H^2(\Xi,\R^3)\to H^1(\Xi,\R^3)
\]
by $(\DN_{\mathrm{dec}}f)_k\vcentcolon=\DN_{(p_k,r_k)}f_k$ for each $k\in\{1,\ldots,b\}$, where $\DN_{(p_k,r_k)}$ is the Dirichlet-to-Neumann operator \eqref{eqn:DtN} for the single-circle configuration $(p_k,r_k)$ in place of $c$. 
\end{definition}

\begin{lemma}
\label{lem:decoupled-Fourier}
For every $f\in H^2(\Xi,\R^3)$, the decoupled Dirichlet-to-Neumann operator satisfies
\begin{align*}
\widehat{(\DN_{\mathrm{dec}}f)_k}(n)=\abs{n}\hat f_k(n)
\end{align*}
for all $k\in\{1,\ldots,b\}$ and $n\in\Z$, where we recall the notation for Fourier coefficients from \eqref{eqn:Fourier-coefficients}. 
In particular, $\DN_{\mathrm{dec}}$ is independent of the configuration $c\in\mathcal{C}$ and the choice of bases in \eqref{eqn:20261009b}.
\end{lemma}

\begin{proof}
Let $c=(p_1,\ldots,p_b;r_1,\ldots,r_b)\in\mathcal{C}$ and $k\in\{1,\ldots,b\}$ be fixed.
The composition of stereographic projection from $p_k\in\Sp^2$ with a suitable dilation and rotation maps the single-circle domain $\Omega_{(p_k,r_k)}=\Sp^2\setminus B_{r_k}(p_k)$ conformally onto the Euclidean unit disc such that $\psi_{c,k}(\theta)\mapsto e^{i\theta}$.
In polar coordinates $(\rho,\theta)$ on the unit disc, the pull-back of the spherical metric is of the form 
\begin{align*}
g_{\Sp^2}=u^2(d\rho^2+\rho^2d\theta^2)
\end{align*}
for some smooth function $u>0$.
Since $\abs{\psi_{c,k}'}=\sin r_k$, we have $u(1,\theta)=\sin r_k$
and the outward unit conormal $\eta$ along the boundary $\{\rho=1\}$ 
satisfies  
\begin{align*}
\sin(r_k)\partial_\eta=\partial_\rho. 
\end{align*} 
By conformal invariance of harmonicity in dimension $2$, the harmonic extension of $e^{in\theta}$ is $\rho^{\abs{n}}e^{in\theta}$ for every $n\in\Z$.
Hence, by definition \eqref{eqn:DtN}, 
\begin{align*}
\DN_{(p_k,r_k)}e^{in\theta}
&=\sin(r_k) \partial_\eta \bigl(\rho^{\abs{n}}e^{in\theta}\bigr)
=\partial_\rho\bigl(\rho^{\abs{n}}e^{in\theta}\bigr)\Big\vert_{\rho=1}
=\abs{n}e^{in\theta}.
\end{align*}
The claim then follows by linearity and density of Fourier polynomials in $H^2(\Sp^1,\R^3)$. 
\end{proof}

\paragraph{Conformal collar estimates.} 
The harmonic extension problems defining $\DN_c$ and $\DN_{\mathrm{dec}}$ have the same local geometry near each boundary circle, while only $\DN_c$ accounts for interactions between different circles.
We use conformal collars to quantify the difference between these operators.
Given $c=(p_1,\ldots,p_b;r_1,\ldots,r_b)\in\mathcal{C}$, we choose $R_k\in\interval{r_k,\pi}$ for each $k\in\{1,\ldots,b\}$ such that the closures of the spherical caps $B_{R_k}(p_k)\subset\Sp^2$ are pairwise disjoint, and set
\begin{align}
\label{eqn:conformal_collar_width}
\omega\vcentcolon=\min_{k\in\{1,\ldots,b\}}
\log\frac{\tan(R_k/2)}{\tan(r_k/2)}>0.
\end{align}
Da\;Lio and Pigati \cite[Lemma~2.6]{DaLioPigati2020} proved that the difference between the Dirichlet-to-Neumann operator and the corresponding Fourier multiplier on each boundary circle is regularising, in the sense that it maps distributions to smooth functions.
The following lemma provides explicit Sobolev bounds for this difference in terms of the conformal collar width \eqref{eqn:conformal_collar_width}.

\begin{lemma}
\label{lem:decoupled-collar}
For every $s\geq0$, the difference $\DN_c-\DN_{\mathrm{dec}}$ extends (from smooth functions) to a bounded operator $L^2(\Xi,\R^3)\to H^s(\Xi,\R^3)$ with operator norm 
\begin{align*}
{\nm{\DN_c-\DN_{\mathrm{dec}}}}_{L^2\to H^s}
\leq\frac{1}{\omega}
\max_{n\in\N_0}\bigl(1+n^2\bigr)^{s/2}e^{-\omega n/2}.
\end{align*}
\end{lemma}

\begin{proof}
For each $k\in\{1,\ldots,b\}$, let $(\varrho,\theta)$ denote geodesic polar coordinates centred at $p_k$, with respect to the orthonormal basis $e_{1,k},e_{2,k}$ from \eqref{eqn:20261009b}.
Changing the radial coordinate to
\begin{align*}
\ell\vcentcolon=\log\frac{\tan(\varrho/2)}{\tan(r_k/2)},
\end{align*}
the spherical metric reads
\(
g_{\Sp^2}=(\sin\varrho)^2(d\ell^2+d\theta^2).
\)
By \eqref{eqn:conformal_collar_width}, the regions $\{0\leq\ell\leq\omega\}$ are pairwise disjoint annuli, each conformally equivalent to the flat cylinder $[0,\omega]\times\Sp^1$.
By the conformal invariance of the Dirichlet energy \eqref{eqn:energy}, we may indeed work with the flat metric on these cylinders. 
Given $f\in H^2(\Xi,\R^3)$, let $u_k^-$ and $u_k^+$ be harmonic on $[0,\omega]\times\Sp^1$, satisfying the boundary conditions  
\begin{align}
\label{eqn:20261029a}
u_k^\pm(0,\theta)&=f_k(\theta),
&
\frac{\partial u_k^-}{\partial\ell}(\omega,\theta)&=0,
&
u_k^+(\omega,\theta)&=0.
\end{align}
We denote by $A_-$ and $A_+$ the corresponding Dirichlet-to-Neumann operators at $\ell=0$, defined by
\begin{align*}
(A_\pm f)_k\vcentcolon=-\frac{\partial u_k^\pm}{\partial\ell}(0,\cdot).
\end{align*}
Taking Fourier coefficients with respect to $\theta$, the harmonic equation for $u_k^\pm$ becomes
\begin{align}\label{eqn:20261029b}
\frac{\partial^2}{\partial\ell^2}\widehat{u}_k^\pm(n)
=n^2\widehat{u}_k^\pm(n).
\end{align}
Solving equation \eqref{eqn:20261029b} subject to the boundary conditions \eqref{eqn:20261029a}, we obtain
\begin{align}
\label{eqn:20261029h}
\widehat{(A_-f)_k}(n)
&=\abs{n}\tanh(\omega\abs{n})\hat f_k(n),
&
\widehat{(A_+f)_k}(n)
&=\abs{n}\coth(\omega\abs{n})\hat f_k(n)
\intertext{for $n\in\Z\setminus\{0\}$.
For $n=0$, the solutions are affine in $\ell$, and hence}
\label{eqn:20261029c}
\widehat{(A_-f)_k}(0)&=0,
&
\widehat{(A_+f)_k}(0)&=\frac{1}{\omega}\hat f_k(0).
\end{align}
Let $\Phi=\Phi_{f,c}$ be the harmonic extension from Definition~\ref{defn:reduced-energy}.
Since $u_k^-$ minimises the Dirichlet energy on the $k$-th cylinder among all maps with boundary values $f_k$ at $\ell=0$, the energy of $\Phi$ is bounded below by the sum of the energies of $u_k^-$.
Conversely, extending each $u_k^+$ by zero to the complement of the corresponding annulus in $\Omega_c$ yields an admissible competitor for $\Phi$.
Hence, by the Dirichlet principle, integration by parts on the cylinders using \eqref{eqn:20261029a}, and Green's formula \eqref{eqn:DtN-Green}, we obtain
\begin{align}
\label{eqn:20261029d}
\int_\Xi\sk{f,A_-f}\,d\theta
&\leq\int_{\Omega_c}\abs{d\Phi}^2
=\int_\Xi\sk{f,\DN_cf}\,d\theta
\leq\int_\Xi\sk{f,A_+f}\,d\theta.
\end{align} 
For every $x>0$, we have
\begin{align}
\label{eqn:20261029f}
0\leq x\bigl(1-\tanh(\omega x)\bigr)
&\leq x\bigl(\coth(\omega x)-1\bigr)
=\frac{xe^{-\omega x}}{\sinh(\omega x)}
\leq\frac{e^{-\omega x}}{\omega}. 
\end{align} 
For every $n\in\Z\setminus\{0\}$, inserting $x=\abs{n}$ into \eqref{eqn:20261029f} yields
\begin{align*}
\abs{n}-\frac{1}{\omega}e^{-\omega\abs{n}}
&\leq\abs{n}\tanh(\omega\abs{n}),
&
\abs{n}\coth(\omega\abs{n})
&\leq\abs{n}+\frac{1}{\omega}e^{-\omega\abs{n}}.
\end{align*}
Recalling the Fourier multiplier formulas in Lemma~\ref{lem:decoupled-Fourier} and \eqref{eqn:20261029h}--\eqref{eqn:20261029c}, we conclude  
\begin{align}
\label{eqn:20261029g}
\DN_{\mathrm{dec}}-\frac{1}{\omega}e^{-\omega\DN_{\mathrm{dec}}}
\leq A_-
\leq A_+
\leq\DN_{\mathrm{dec}}+\frac{1}{\omega}e^{-\omega\DN_{\mathrm{dec}}}
\end{align}
as quadratic forms by Parseval's identity, where $e^{-\omega\DN_{\mathrm{dec}}}$ denotes the bounded operator on $L^2(\Xi,\R^3)$ defined by
\begin{align*}
\widehat{(e^{-\omega\DN_{\mathrm{dec}}}f)_k}(n)
=e^{-\omega\abs{n}}\hat f_k(n).
\end{align*}
Let $T\vcentcolon=e^{-\omega\DN_{\mathrm{dec}}/2}$ be defined analogously, so that $T$ is self-adjoint and $T^2=e^{-\omega\DN_{\mathrm{dec}}}$. 
Then, \eqref{eqn:20261029g} and \eqref{eqn:20261029d} imply 
\begin{align*}
\abs[\Big]{\omega\int_\Xi\sk{f,(\DN_c-\DN_{\mathrm{dec}})f}\,d\theta}
\leq\int_\Xi\abs{Tf}^2\,d\theta
=\nm{Tf}_{L^2}^2
\end{align*}
for every Fourier polynomial $f$.
Since $T$ maps Fourier polynomials bijectively onto Fourier polynomials, substituting $f=T^{-1}g$ yields
\begin{align}
\label{eqn:20261029i}
\abs[\Big]{\omega\int_\Xi\sk[\Big]{T^{-1}g,(\DN_c-\DN_{\mathrm{dec}})T^{-1}g}\,d\theta}
\leq\nm{g}_{L^2}^2.
\end{align}
By \eqref{eqn:DtN-Green} and Lemma~\ref{lem:decoupled-Fourier}, the bilinear form
\begin{align*}
(g,h)\mapsto\omega\int_\Xi\sk[\Big]{T^{-1}g,(\DN_c-\DN_{\mathrm{dec}})T^{-1}h}\,d\theta
\end{align*}
defined on Fourier polynomials is symmetric with norm 
at most $1$ by polarisation of \eqref{eqn:20261029i}. 
By density, it extends uniquely to a bounded symmetric bilinear form on $L^2(\Xi,\R^3)$.
By the Riesz representation theorem, there exists a bounded self-adjoint operator $S_c$ on $L^2(\Xi,\R^3)$ with $\nm{S_c}_{L^2\to L^2}\leq1$ such that
\begin{align}
\label{eqn:decoupled-factorization}
\DN_c-\DN_{\mathrm{dec}}=\frac{1}{\omega}TS_cT
\end{align}
on Fourier polynomials.
Since both sides of \eqref{eqn:decoupled-factorization} define bounded operators $H^2(\Xi,\R^3)\to H^1(\Xi,\R^3)$, the identity extends to $H^2(\Xi,\R^3)$ by density.
By Lemma~\ref{lem:decoupled-Fourier} and the definition of the Sobolev norms in \eqref{eqn:Sobolev-norm}, we have 
\begin{align*}
{\nm{T}}_{L^2\to L^2}&=1,
&
{\nm{T}}_{L^2\to H^s}
&=\max_{n\in\N_0}(1+n^2)^{s/2}e^{-\omega n/2}.
\end{align*}
Consequently, \eqref{eqn:decoupled-factorization} implies
\begin{align}\label{eqn:20261029e}
\begin{aligned}
\nm{(\DN_c-\DN_{\mathrm{dec}})f}_{H^s}
&\leq\frac{1}{\omega}
{\nm{T}}_{L^2\to H^s}
{\nm{S_c}}_{L^2\to L^2}
{\nm{T}}_{L^2\to L^2}
{\nm{f}}_{L^2}
\\
&\leq\frac{1}{\omega}
\max_{n\in\N_0}(1+n^2)^{s/2}e^{-\omega n/2}
{\nm{f}}_{L^2}
\end{aligned}
\end{align}
for every Fourier polynomial $f$. 
Since Fourier polynomials are dense in $L^2(\Xi,\R^3)$, estimate \eqref{eqn:20261029e} proves the claim.
\end{proof} 

\begin{corollary}
\label{cor:decoupled-explicit-bounds}
For every $j\in\N\cup\{0\}$, we have
\begin{align*} 
{\nm{\DN_c-\DN_{\mathrm{dec}}}}_{L^2\to H^j}
\leq C_j(\omega)\vcentcolon=
\frac{1}{\omega}
\sqrt{1+\sum_{m=1}^j\binom{j}{m}
\left(\frac{2m}{e\omega}\right)^{2m}}.
\end{align*} 
\end{corollary}

\begin{proof}
For every $m\in\N$ and every $j,n\in\N\cup\{0\}$, we have 
\begin{align*}
\max_{x\geq0}x^{2m}e^{-\omega x}
&=\left(\frac{2m}{e\omega}\right)^{2m}, 
&
(1+n^2)^je^{-\omega n}
&\leq 1+\sum_{m=1}^j\binom{j}{m}
\left(\frac{2m}{e\omega}\right)^{2m}
\end{align*}
by the binomial theorem. 
The claim then follows from Lemma~\ref{lem:decoupled-collar}.
\end{proof}

\paragraph{Configuration derivatives.} 
To estimate the dependence of the Dirichlet-to-Neumann operator $\DN_c$ on the configuration, we recall the parametrisation $\mu\mapsto c_\mu$ from \eqref{eqn:perturbation-configuration}.
Given $\mu\in B\subset\R^{3b-6}$ and $h\in\R^{3b-6}$, let $\tilde\alpha_{\mu,h},\tilde\beta_{\mu,h}\colon\Xi\to\R$ denote the normalised boundary velocities (as in Lemma~\ref{lem:DtN-derivative}) associated with the curve $t\mapsto c_{\mu+th}$ at $t=0$, and let 
\begin{align}\label{eqn:dotDN}
\dot\DN_{\mu,h}\vcentcolon=\frac{d}{dt}\Big\vert_{t=0}\DN_{c_{\mu+th}}.
\end{align}
The following proposition provides uniform estimates in terms of a common conformal collar width $\omega$ and a constant $C_c$ controlling the boundary velocities and their Lipschitz dependence on $\mu$.

\begin{proposition}
\label{prop:DtN-differences}
Suppose that $B\subset\R^{3b-6}$ is convex and that there exist constants $\omega,C_c>0$ such that every configuration $c_\mu$ with $\mu\in B$ admits pairwise disjoint conformal collars of width at least $\omega$, and such that
\begin{align*}
\nm{\tilde\alpha_{\mu,h}}_{L^\infty}
+\nm{\tilde\beta_{\mu,h}}_{L^\infty}
&\leq C_c\abs{h},
\\
\nm{\tilde\alpha_{\mu,h}-\tilde\alpha_{\nu,h}}_{L^\infty}
+\nm{\tilde\beta_{\mu,h}-\tilde\beta_{\nu,h}}_{L^\infty}
&\leq C_c\abs{\mu-\nu}\abs{h}
\end{align*}
for all $\mu,\nu\in B$ and $h\in\R^{3b-6}$.
Set
$C_\Lambda\vcentcolon=6C_cC_3(\omega)\bigl(2+C_3(\omega)\bigr)$, 
where $C_3(\omega)$ is defined in Corollary~\ref{cor:decoupled-explicit-bounds}.
Then, for all $\mu,\nu\in B$ and $h\in\R^{3b-6}$, we have
\begin{align}
\label{eqn:DtN-differences-1}
\nm{\DN_{c_\mu}}_{H^2\to H^1}
&\leq1+C_1(\omega),
\\\label{eqn:DtN-differences-2}
\nm{\dot\DN_{\mu,h}}_{H^1\to H^2}
&\leq C_\Lambda\abs{h},
\\\label{eqn:DtN-differences-3}
\nm{\DN_{c_\mu}-\DN_{c_\nu}}_{H^1\to H^2}
&\leq C_\Lambda\abs{\mu-\nu},
\\\label{eqn:DtN-differences-4}
\nm{\dot\DN_{\mu,h}-\dot\DN_{\nu,h}}_{H^2\to H^1}
&\leq C_\Lambda\bigl(1+12C_c(1+C_3(\omega))\bigr)
\abs{\mu-\nu}\abs{h}.
\end{align}
\end{proposition}

\begin{remark}
Claims \eqref{eqn:DtN-differences-2} and \eqref{eqn:DtN-differences-3} indicate that the configuration derivatives and differences of Dirichlet-to-Neumann operators are regularising, even though $\DN_c$ itself is an operator of order one.    
\end{remark}

\begin{proof}
By Lemma~\ref{lem:decoupled-Fourier} and Corollary~\ref{cor:decoupled-explicit-bounds}, we have
\begin{align*}
\nm{\DN_{c_\mu}}_{H^2\to H^1}
&\leq\nm{\DN_{\mathrm{dec}}}_{H^2\to H^1}
+\nm{\DN_{c_\mu}-\DN_{\mathrm{dec}}}_{H^2\to H^1}
\leq1+C_1(\omega)
\end{align*}
as claimed in \eqref{eqn:DtN-differences-1}.
For the remaining estimates \eqref{eqn:DtN-differences-2}--\eqref{eqn:DtN-differences-4}, let $\mu\in B$ and $h\in\R^{3b-6}$ be fixed, and set $\Theta_\mu\vcentcolon=\DN_{c_\mu}-\DN_{\mathrm{dec}}$.
For smooth $f,\zeta\colon\Xi\to\R^3$, Lemma~\ref{lem:DtN-derivative} yields
\begin{align*}
\int_\Xi\sk{\zeta,\dot\DN_{\mu,h}f}\,d\theta
&=\int_\Xi
\sk[\big]{\zeta',\;\tilde\alpha_{\mu,h}f'-\tilde\beta_{\mu,h}\DN_{c_\mu}f}
-\sk[\big]{\DN_{c_\mu}\zeta,\;\tilde\alpha_{\mu,h}\DN_{c_\mu}f+\tilde\beta_{\mu,h}f'}
\,d\theta
\\
&=\int_\Xi\sk[\Big]{\zeta,\,
-(\tilde\alpha_{\mu,h}f')'
+(\tilde\beta_{\mu,h}\DN_{c_\mu}f)'
-\DN_{c_\mu}(\tilde\alpha_{\mu,h}\DN_{c_\mu}f+\tilde\beta_{\mu,h}f')}
\,d\theta,
\end{align*}
where the second equality follows by integration by parts and the symmetry of $\DN_{c_\mu}$ from \eqref{eqn:DtN-Green}.
Since $\zeta$ is arbitrary, we conclude
\begin{align}
\label{eqn:DtN-derivative-operator}
\dot\DN_{\mu,h}f
&=-(\tilde\alpha_{\mu,h}f')'
+(\tilde\beta_{\mu,h}\DN_{c_\mu}f)'
-\DN_{c_\mu}(\tilde\alpha_{\mu,h}\DN_{c_\mu}f+\tilde\beta_{\mu,h}f').
\end{align}
By \eqref{eqn:normal-velocity-fourier} and \eqref{eqn:tangential-velocity-fourier}, the restrictions of $\tilde\alpha_{\mu,h}$ and $\tilde\beta_{\mu,h}$ to each component of $\Xi$ are Fourier polynomials of degree at most one.
Given $n\in\Z$, $m\in\{-1,0,1\}$ and $v\in\C^3$, let $a(\theta)=e^{im\theta}$ and $f(\theta)=e^{in\theta}v$.
Multiplication by $a$ shifts the Fourier frequency from $n$ to $n+m$.
Since $n$ and $n+m$ cannot have opposite signs, Lemma~\ref{lem:decoupled-Fourier} implies
\begin{align*}
-(af')'
&=n(n+m)af
=\abs{n}\abs{n+m}af
=\DN_{\mathrm{dec}}(a\DN_{\mathrm{dec}}f),
\\
(a\DN_{\mathrm{dec}}f)'
&=i(n+m)\abs{n}af
=in\abs{n+m}af
=\DN_{\mathrm{dec}}(af').
\end{align*}
By density of Fourier polynomials, we conclude
\begin{align}
\label{eqn:decoupled-product-identities}
-(af')'&=\DN_{\mathrm{dec}}(a\DN_{\mathrm{dec}}f),
&
(a\DN_{\mathrm{dec}}f)'&=\DN_{\mathrm{dec}}(af')
\end{align}
for every smooth $f\colon\Xi\to\R^3$ and every scalar Fourier polynomial $a$ of degree at most one on each component of $\Xi$.
Substituting $\DN_{c_\mu}=\DN_{\mathrm{dec}}+\Theta_\mu$ into \eqref{eqn:DtN-derivative-operator}, the terms involving only $\DN_{\mathrm{dec}}$ cancel by \eqref{eqn:decoupled-product-identities} and we obtain
\begin{align}
\label{eqn:DtN-derivative-remainder}
\dot\DN_{\mu,h}f
&=(\tilde\beta_{\mu,h}\Theta_\mu f)'
-\DN_{\mathrm{dec}}(\tilde\alpha_{\mu,h}\Theta_\mu f)
-\Theta_\mu(\tilde\alpha_{\mu,h}\DN_{c_\mu}f+\tilde\beta_{\mu,h}f').
\end{align}
For every real-valued scalar Fourier polynomial $a$ of degree at most one on each component of $\Xi$, we have
\begin{align}
\label{eqn:degree-one-product}
\nm{ag}_{H^j}\leq2^j\nm{a}_{L^\infty}\nm{g}_{H^j}
\end{align}
for all $j\in\N\cup\{0\}$ and $g\in H^j(\Xi,\R^3)$.
Indeed, the Fourier coefficients of $a$ satisfy
\begin{align*}
\sum_{m=-1}^1\abs{\hat a_k(m)}
&=\nm{a_k}_{L^\infty},
&
1+(n+m)^2&\leq4(1+n^2)
\end{align*}
for every $k\in\{1,\ldots,b\}$, $n\in\Z$ and $m\in\{-1,0,1\}$, so \eqref{eqn:degree-one-product} follows from Parseval's identity and the triangle inequality.
By Lemma~\ref{lem:decoupled-Fourier} and Corollary~\ref{cor:decoupled-explicit-bounds}, we have
\begin{align*}
\nm{\DN_{c_\mu}f}_{L^2}
&\leq\nm{\DN_{\mathrm{dec}}f}_{L^2}
+\nm{\Theta_\mu f}_{L^2}
\leq\nm{f'}_{L^2}+C_3(\omega)\nm{f}_{L^2}.
\end{align*}
Consequently, \eqref{eqn:DtN-derivative-remainder}, \eqref{eqn:degree-one-product} and Corollary~\ref{cor:decoupled-explicit-bounds} imply
\begin{align*}
\nm{\dot\DN_{\mu,h}f}_{H^2}
&\leq C_3(\omega)\Bigl(
8\bigl(\nm{\tilde\alpha_{\mu,h}}_{L^\infty}
+\nm{\tilde\beta_{\mu,h}}_{L^\infty}\bigr)\nm{f}_{L^2}
+\nm{\tilde\alpha_{\mu,h}\DN_{c_\mu}f
+\tilde\beta_{\mu,h}f'}_{L^2}
\Bigr)
\\
&\leq C_3(\omega)\bigl(9+C_3(\omega)\bigr)
\bigl(\nm{\tilde\alpha_{\mu,h}}_{L^\infty}
+\nm{\tilde\beta_{\mu,h}}_{L^\infty}\bigr)\nm{f}_{H^1}
\\[1ex]
&\leq C_cC_3(\omega)\bigl(9+C_3(\omega)\bigr)\abs{h}\nm{f}_{H^1}
\leq C_\Lambda\abs{h}\nm{f}_{H^1}.
\end{align*} 
By density, $\dot\DN_{\mu,h}$ extends to a bounded operator $H^1(\Xi,\R^3)\to H^2(\Xi,\R^3)$.
In particular, the second claim \eqref{eqn:DtN-differences-2} follows.

Given $\mu,\nu\in B$ let $\mu(t)\vcentcolon=\nu+t(\mu-\nu)$ for $t\in\Interval{0,1}$.
Then $\mu(t)\in B$ for every $t\in\Interval{0,1}$ since $B$ is convex by assumption. 
By Lemma~\ref{lem:smoothness}, the fundamental theorem of calculus yields
\begin{align*}
(\DN_{c_\mu}-\DN_{c_\nu})f
&=\int_0^1\dot\DN_{\mu(t),\mu-\nu}f\,dt
\end{align*}
in $H^1(\Xi,\R^3)$ for every $f\in H^2(\Xi,\R^3)$.
Taking Fourier coefficients of this identity and applying Minkowski's inequality to finite Fourier sums, we obtain
\begin{align*}
\biggl(
2\pi\sum_{k=1}^b\sum_{\abs{n}\leq N}
(1+n^2)^2
\abs[\big]{\widehat{((\DN_{c_\mu}-\DN_{c_\nu})f)_k}(n)}^2
\biggr)^{\frac{1}{2}}
&\leq\int_0^1\nm{\dot\DN_{\mu(t),\mu-\nu}f}_{H^2}\,dt
\leq C_\Lambda\abs{\mu-\nu}\nm{f}_{H^1}
\end{align*}
for every $N\in\N$, where we used \eqref{eqn:DtN-differences-2}.
Letting $N\to\infty$ and recalling \eqref{eqn:Sobolev-norm}, we conclude
\begin{align*}
\nm{(\DN_{c_\mu}-\DN_{c_\nu})f}_{H^2}
\leq C_\Lambda\abs{\mu-\nu}\nm{f}_{H^1}.
\end{align*}
By density of Fourier polynomials in $H^1(\Xi,\R^3)$, the difference $\DN_{c_\mu}-\DN_{c_\nu}$ extends to a bounded operator $H^1(\Xi,\R^3)\to H^2(\Xi,\R^3)$, proving \eqref{eqn:DtN-differences-3}.
It remains to prove \eqref{eqn:DtN-differences-4}.
Subtracting \eqref{eqn:DtN-derivative-remainder} at $\mu$ and $\nu$, and using $\Theta_\mu-\Theta_\nu=\DN_{c_\mu}-\DN_{c_\nu}$, we obtain
\begin{align}
\label{eqn:DtN-derivative-difference}
\begin{aligned}
(\dot\DN_{\mu,h}-\dot\DN_{\nu,h})f
&=\bigl((\tilde\beta_{\mu,h}-\tilde\beta_{\nu,h})\Theta_\mu f\bigr)'
-\DN_{\mathrm{dec}}\bigl((\tilde\alpha_{\mu,h}-\tilde\alpha_{\nu,h})\Theta_\mu f\bigr)
\\
&\hphantom{{}={}}-\Theta_\nu\bigl((\tilde\alpha_{\mu,h}-\tilde\alpha_{\nu,h})\DN_{c_\mu}f
+(\tilde\beta_{\mu,h}-\tilde\beta_{\nu,h})f'\bigr)
\\
&\hphantom{{}={}}+\bigl(\tilde\beta_{\nu,h}(\Theta_\mu-\Theta_\nu)f\bigr)'
-\DN_{\mathrm{dec}}\bigl(\tilde\alpha_{\nu,h}(\Theta_\mu-\Theta_\nu)f\bigr)
\\
&\hphantom{{}={}}-(\Theta_\mu-\Theta_\nu)
\bigl(\tilde\alpha_{\mu,h}\DN_{c_\mu}f+\tilde\beta_{\mu,h}f'\bigr)
\\
&\hphantom{{}={}}-\Theta_\nu\bigl(\tilde\alpha_{\nu,h}(\Theta_\mu-\Theta_\nu)f\bigr).
\end{aligned}
\end{align}
By \eqref{eqn:DtN-differences-1}, \eqref{eqn:degree-one-product} and the assumptions on the boundary velocities, we have
\begin{align*}
\nm{\tilde\alpha_{\mu,h}\DN_{c_\mu}f+\tilde\beta_{\mu,h}f'}_{H^1}
&\leq2C_c(1+C_3(\omega))\abs{h}\nm{f}_{H^2},
\\
\nm{(\tilde\alpha_{\mu,h}-\tilde\alpha_{\nu,h})\DN_{c_\mu}f
+(\tilde\beta_{\mu,h}-\tilde\beta_{\nu,h})f'}_{L^2}
&\leq C_c(1+C_3(\omega))
\abs{\mu-\nu}\abs{h}\nm{f}_{H^2}.
\end{align*}
Consequently, applying \eqref{eqn:degree-one-product}, \eqref{eqn:DtN-differences-3} and Corollary~\ref{cor:decoupled-explicit-bounds} to \eqref{eqn:DtN-derivative-difference} yields
\begin{align*}
\nm{(\dot\DN_{\mu,h}-\dot\DN_{\nu,h})f}_{H^1}
&\leq C_c\Bigl(
C_3(\omega)(5+C_3(\omega))
+C_\Lambda(6+3C_3(\omega))
\Bigr)
\abs{\mu-\nu}\abs{h}\nm{f}_{H^2}
\\
&\leq C_\Lambda\bigl(1+12C_c(1+C_3(\omega))\bigr)
\abs{\mu-\nu}\abs{h}\nm{f}_{H^2},
\end{align*}
where the last inequality follows from the definition of $C_\Lambda$ and the elementary inequalities
\begin{align*}
5+C_3(\omega)&\leq6(2+C_3(\omega)),
&
6+3C_3(\omega)&\leq12(1+C_3(\omega)).
\end{align*}
By density, this proves \eqref{eqn:DtN-differences-4}.
\end{proof} 

\paragraph{Approximation of the Dirichlet-to-Neumann operator.}
To approximate $\DN_cf$, we use the decomposition $\DN_c=\DN_{\mathrm{dec}}+(\DN_c-\DN_{\mathrm{dec}})$ and choose a Fourier polynomial $\varphi$ approximating $f$ in $H^2$.
The decoupled part $\DN_{\mathrm{dec}}\varphi$ can then be computed from Fourier coefficients. 
The regularising property of $\DN_c-\DN_{\mathrm{dec}}$ allows us to replace $\varphi$ by auxiliary boundary data $\tilde{f}$ in this difference, with an error controlled just in $L^2$, such that the Fourier coefficients of $(\DN_c-\DN_{\mathrm{dec}})\tilde{f}$ can be computed explicitly.  

\begin{definition}[Fourier truncation]
\label{def:fourier-truncation}
Given $N\in\N\cup\{0\}$ and $f\in L^2(\Xi)$, we define the \emph{degree $N$ Fourier truncation} of $f$ by
\begin{align*}
{(\Pi_N f)}_k(\theta)\vcentcolon=\sum_{\abs{n}\leq N}\hat f_k(n)e^{in\theta}
\end{align*}
for each $k\in\{1,\ldots,b\}$.
For vector- and matrix-valued functions, $\Pi_N$ acts componentwise.
\end{definition}

\begin{proposition}
\label{prop:DtN-error}
Let $c\in\mathcal{C}$ and $N\in\N_0$.
Given $f,\tilde{f}\in H^2(\Xi,\R^3)$ and a Fourier polynomial $\varphi\colon\Xi\to\R^3$, let
\begin{align*}
\tilde{h}\vcentcolon=\DN_{\mathrm{dec}}\varphi
+\Pi_N(\DN_c-\DN_{\mathrm{dec}})\tilde{f}.
\end{align*}
Then, recalling the conformal collar width $\omega$ from \eqref{eqn:conformal_collar_width}, 
\begin{align*} 
\nm{\DN_cf-\tilde{h}}_{H^1}
&\leq\bigl(1+C_1(\omega)\bigr)\nm{f-\varphi}_{H^2}
+C_1(\omega)\nm[\big]{\varphi-\tilde{f}}_{L^2}
+\nm[\Big]{(\mathrm{Id}-\Pi_N)(\DN_c-\DN_{\mathrm{dec}})\tilde{f}}_{H^1}.
\end{align*}
\end{proposition}

\begin{proof}
By the definition of $\tilde{h}$ and the decomposition $\DN_c=\DN_{\mathrm{dec}}+(\DN_c-\DN_{\mathrm{dec}})$, we have
\begin{align*}
\DN_cf-\tilde{h}
&=\DN_c(f-\varphi)
+(\DN_c-\DN_{\mathrm{dec}})(\varphi-\tilde{f})
+(\mathrm{Id}-\Pi_N)(\DN_c-\DN_{\mathrm{dec}})\tilde{f}.
\end{align*}
By Lemma~\ref{lem:decoupled-Fourier} and Corollary~\ref{cor:decoupled-explicit-bounds}, we have
\begin{align*}
\nm{\DN_c}_{H^2\to H^1}
&\leq1+C_1(\omega),
&
\nm{\DN_c-\DN_{\mathrm{dec}}}_{L^2\to H^1}
&\leq C_1(\omega).
\end{align*}
The claim then follows from the triangle inequality.
\end{proof}

\paragraph{Method of fundamental solutions.}
We construct auxiliary boundary data $\tilde{f}$ for Proposition~\ref{prop:DtN-error} using the method of fundamental solutions, placing singularities inside the spherical caps removed from $\Omega_c$.
Given $c=(p_1,\ldots,p_b;r_1,\ldots,r_b)\in\mathcal{C}$ and $L\in\N$, we choose points $s_1,\ldots,s_L\in\bigcup_{k=1}^b B_{r_k}(p_k)$ and coefficients $a_0,\ldots,a_L\in\R^3$ to define $\Phi\colon\overline{\Omega}_c\to\R^3$ and its boundary traces $\tilde{f}_k$ by
\begin{align}
\label{eqn:mfs-field}
\Phi(x)&\vcentcolon=a_0-\frac{1}{4\pi}\sum_{\ell=1}^L a_\ell\log\bigl(1-\sk{x,s_\ell}\bigr),
&
\tilde{f}_k&\vcentcolon=\Phi\circ\psi_{c,k}.
\end{align}
Setting $u=\sk{x,s_\ell}$, we have $\Delta_{g_{\Sp^2}}u=-2u$ and $\abs{du}^2=1-u^2$, so that
\begin{align*}
\Delta_{g_{\Sp^2}}\log(1-u)
&=-\frac{\Delta_{g_{\Sp^2}}u}{1-u}-\frac{\abs{du}^2}{(1-u)^2}
=\frac{2u}{1-u}-\frac{1+u}{1-u}
=-1.
\end{align*} 
Hence, $\Phi$ is harmonic in a neighbourhood of $\overline{\Omega}_c\subset\Sp^2$, provided that $\sum_{\ell=1}^L a_\ell=0$.  
The Fourier coefficients of the logarithmic terms in \eqref{eqn:mfs-field} can be expressed using the complex parameter  
\begin{align}
\label{eqn:mfs-fourier-ratios}
q_{k\ell}
\vcentcolon=\frac{\sin(r_k)\bigl(\sk{e_{1,k},s_\ell}-i\sk{e_{2,k},s_\ell}\bigr)}
{1-\cos(r_k)\sk{p_k,s_\ell}+\abs{\sk{p_k,s_\ell}-\cos(r_k)}} 
\in\C
\end{align}
for the boundary circle with label $k\in\{1,\ldots,b\}$ and the source with label $\ell\in\{1,\ldots,L\}$, 
where $(e_{1,k},e_{2,k})$ is the orthonormal basis used in \eqref{eqn:20261009b}.

\begin{proposition}
\label{prop:mfs-fourier-coefficients}
Suppose that $\sum_{\ell=1}^L a_\ell=0$.
Then $\abs{q_{k\ell}}<1$ for all $k\in\{1,\ldots,b\}$ and $\ell\in\{1,\ldots,L\}$.
Moreover, the Fourier coefficients of the boundary traces $\tilde{f}_k$ from \eqref{eqn:mfs-field} satisfy 
\begin{align*}
\widehat{\tilde{f}_k}(0)
&=a_0-\frac{1}{4\pi}\sum_{\ell=1}^L a_\ell
\log\biggl(
\frac{1-\cos r_k\,\sk{p_k,s_\ell}+\abs{\sk{p_k,s_\ell}-\cos r_k}}{2}
\biggr),
\\
\widehat{\tilde{f}_k}(n)
&=\frac{1}{4\pi n}\sum_{\ell=1}^L a_\ell q_{k\ell}^{\,n},
&&n\geq1,
\\
\widehat{\bigl((\DN_c-\DN_{\mathrm{dec}})\tilde{f}\bigr)_k}(n)
&=-\frac{1}{2\pi}
\sum_{\substack{1\leq\ell\leq L\\s_\ell\notin B_{r_k}(p_k)}}
a_\ell q_{k\ell}^{\,n},
&&n\geq0. 
\end{align*}
The Fourier coefficients for $n<0$ are determined by complex conjugation. 
\end{proposition}

\begin{remark}
The formula for $(\DN_c-\DN_{\mathrm{dec}})\tilde{f}$ stated in Proposition~\ref{prop:mfs-fourier-coefficients} shows that its restriction to the $k$-th boundary component depends only on singularities outside $B_{r_k}(p_k)$.
In particular, singularities inside this cap do not contribute to the correction on its boundary.
\end{remark}

\begin{proof}
Given $k\in\{1,\ldots,b\}$ and $\ell\in\{1,\ldots,L\}$, let 
$u\vcentcolon=\sk{p_k,s_\ell}$ 
and 
$z\vcentcolon=\sk{e_{1,k},s_\ell}-i\sk{e_{2,k},s_\ell}$. 
For $r$ in a neighbourhood of $r_k$, let
\begin{align*}
A(r)&\vcentcolon=\frac{1-\cos(r)u+\abs{u-\cos(r)}}{2},
&
q(r)&\vcentcolon=\frac{\sin(r)z}{2A(r)}.
\end{align*}
Since $\abs{z}^2=1-u^2$, we have
$\bigl(1-\cos(r)u\bigr)^2-\sin^2(r)\abs{z}^2
=\bigl(u-\cos(r)\bigr)^2$. 
As $s_\ell\notin\partial B_{r_k}(p_k)$, this implies $A(r)>0$ and $\abs{q(r)}<1$ for $r$ sufficiently close to $r_k$.
In particular, $q(r_k)=q_{k\ell}$.

Let $\psi_k(r,\theta)$ denote the parametrisation \eqref{eqn:20261009b} with $r_k$ replaced by $r$.
A direct calculation yields 
\begin{align*}
1-\sk{\psi_k(r,\theta),s_\ell}
&=A(r)\abs[\big]{1-q(r)e^{i\theta}}^2.
\shortintertext{Consequently,}
\log\bigl(1-\sk{\psi_k(r,\theta),s_\ell}\bigr)
&=\log A(r)
-\sum_{n=1}^\infty
\frac{q(r)^n e^{in\theta}+\overline{q(r)}^{\,n}e^{-in\theta}}{n}.
\end{align*}
Evaluating at $r=r_k$ and recalling \eqref{eqn:mfs-field} proves the first two Fourier formulas.
For the last claim we observe that $u\neq\cos(r_k)$ since the source points lie inside the pairwise disjoint spherical caps.
Consequently, the sign of $u-\cos(r)$ is constant for $r$ sufficiently close to $r_k$, and
\begin{align*}
2A(r)
&=\begin{cases}
(1+u)(1-\cos(r)) , &u>\cos(r_k),\\
(1-u)(1+\cos(r)), &u<\cos(r_k).
\end{cases}
\end{align*}
It follows that
\begin{align*}
\sin(r)q'(r)
&=\begin{cases}
-q(r), &u>\cos(r_k),\\
q(r), &u<\cos(r_k),
\end{cases}
&
\sin(r)\frac{A'(r)}{A(r)}
&=\begin{cases}
1+\cos(r), &u>\cos(r_k),\\
\cos(r)-1, &u<\cos(r_k).
\end{cases}
\end{align*}
Since the outward unit conormal to $\Omega_c$ points towards the removed cap, \eqref{eqn:DtN} implies
\begin{align*}
(\DN_c\tilde{f})_k(\theta)
&=-\sin(r_k)\frac{\partial}{\partial r}
\Phi\bigl(\psi_k(r,\theta)\bigr)\Big\vert_{r=r_k}.
\end{align*}
Since $\abs{q(r)}<1$ uniformly for $r$ sufficiently close to $r_k$, we may differentiate the Fourier expansion termwise.
Using the formulas for $q'(r)$ and $A'(r)$, we obtain, for every $n\geq1$,
\begin{align*}
\widehat{(\DN_c\tilde{f})_k}(n)
&=\frac{1}{4\pi}
\biggl(
\sum_{s_\ell\in B_{r_k}(p_k)}a_\ell q_{k\ell}^{n}
-\sum_{s_\ell\notin B_{r_k}(p_k)}a_\ell q_{k\ell}^{n}
\biggr). 
\shortintertext{Since} 
\widehat{(\DN_{\mathrm{dec}}\tilde{f})_k}(n)
&=n\widehat{\tilde{f}_k}(n)
=\frac{1}{4\pi}\sum_{\ell=1}^L a_\ell q_{k\ell}^{\,n} 
\end{align*}
this proves the third claim for $1\leq n\in\N$. 
For $n=0$, we have $\widehat{(\DN_{\mathrm{dec}}\tilde{f})_k}(0)=0$, while
\begin{align*}
\widehat{(\DN_c\tilde{f})_k}(0)
&=\frac{1}{4\pi}
\biggl(
(1+\cos(r_k))\sum_{s_\ell\in B_{r_k}(p_k)}a_\ell
-(1-\cos(r_k))\sum_{s_\ell\notin B_{r_k}(p_k)}a_\ell
\biggr)
=-\frac{1}{2\pi}
\sum_{s_\ell\notin B_{r_k}(p_k)}a_\ell,
\end{align*}
where the last equality uses $\sum_{\ell=1}^L a_\ell=0$.
Finally, the formulas for negative $n$ follow from the fact that the components of $\tilde{f}$ and $(\DN_c-\DN_{\mathrm{dec}})\tilde{f}$ are real-valued.
\end{proof}

Proposition~\ref{prop:mfs-fourier-coefficients} expresses the Fourier coefficients of $\tilde{f}$ and $(\DN_c-\DN_{\mathrm{dec}})\tilde{f}$ explicitly in terms of the source points and coefficients.
Since $\abs{q_{k\ell}}<1$, these coefficients decay geometrically as  $n$ increases.
This allows us to bound the Fourier truncation error $(\mathrm{Id}-\Pi_N)(\DN_c-\DN_{\mathrm{dec}})\tilde{f}$ and the difference between $\tilde{f}$ and a Fourier polynomial $\varphi$ using only finitely many Fourier coefficients and explicit remainder bounds.
The following proposition thus provides computable bounds for the last two terms in Proposition~\ref{prop:DtN-error}.

\begin{proposition}
\label{prop:mfs-truncation-bounds}
Suppose that $\sum_{\ell=1}^L a_\ell=0$, and let $N\in\N\cup\{0\}$.
Then, for each $k\in\{1,\ldots,b\}$,
\begin{align*}
&\nm[\Big]{\bigl((\mathrm{Id}-\Pi_N)(\DN_c-\DN_{\mathrm{dec}})\tilde{f}\bigr)_k}_{H^1}
\\
&\leq\frac{1}{\sqrt{\pi}}
\sum_{s_\ell\notin B_{r_k}(p_k)}
\abs{a_\ell}\abs{q_{k\ell}}^{N+1}
\biggl(
\frac{1+(N+1)^2}{1-\abs{q_{k\ell}}^2}
+\frac{2(N+1)\abs{q_{k\ell}}^2}{(1-\abs{q_{k\ell}}^2)^2}
+\frac{\abs{q_{k\ell}}^2(1+\abs{q_{k\ell}}^2)}
{(1-\abs{q_{k\ell}}^2)^3}
\biggr)^{\frac{1}{2}}.
\end{align*}
Moreover, for every Fourier polynomial $\varphi\colon\Xi\to\R^3$ of degree at most $N$,
\begin{align*}
\nm{\varphi-\tilde{f}}_{L^2}^2
&\leq2\pi\sum_{k=1}^b\sum_{\abs{n}\leq N}
\abs{\widehat{\varphi_k}(n)-\widehat{\tilde{f}_k}(n)}^2
+\frac{1}{4\pi(N+1)^2}
\sum_{k=1}^b
\biggl(
\sum_{\ell=1}^L
\frac{\abs{a_\ell}\abs{q_{k\ell}}^{N+1}}
{\sqrt{1-\abs{q_{k\ell}}^2}}
\biggr)^2.
\end{align*}
\end{proposition}

\begin{proof}
Let $v\vcentcolon=(\DN_c-\DN_{\mathrm{dec}})\tilde{f}$.
By Proposition~\ref{prop:mfs-fourier-coefficients}, Parseval's identity and the triangle inequality, we have
\begin{align*}
\nm[\Big]{\bigl((\mathrm{Id}-\Pi_N)v\bigr)_k}_{H^1}
&=\biggl(
4\pi\sum_{n=N+1}^\infty
(1+n^2)\abs{\widehat{v_k}(n)}^2
\biggr)^{\frac{1}{2}}
\leq\frac{1}{\sqrt{\pi}}
\sum_{s_\ell\notin B_{r_k}(p_k)}
\abs{a_\ell}
\biggl(
\sum_{n=N+1}^\infty
(1+n^2)\abs{q_{k\ell}}^{2n}
\biggr)^{\frac{1}{2}}.
\end{align*}
The first claim then follows from the fact that $t=\abs{q_{k\ell}}^2\in\Interval{0,1}$ satisfies  
\begin{align*}
\sum_{n=N+1}^\infty(1+n^2)t^n
&=t^{N+1}
\biggl(
\frac{1+(N+1)^2}{1-t}
+\frac{2(N+1)t}{(1-t)^2}
+\frac{t(1+t)}{(1-t)^3}
\biggr).
\end{align*} 
Since $\varphi$ is a Fourier polynomial of degree at most $N$, Parseval's identity yields
\begin{align}
\label{eqn:mfs-truncation-parseval}
\nm[\big]{\varphi-\tilde{f}}_{L^2}^2
&=2\pi\sum_{k=1}^b\sum_{\abs{n}\leq N}
\abs[\Big]{\widehat{\varphi_k}(n)-\widehat{\tilde{f}_k}(n)}^2
+4\pi\sum_{k=1}^b\sum_{n=N+1}^\infty
\abs[\Big]{\widehat{\tilde{f}_k}(n)}^2.
\end{align}
Using Proposition~\ref{prop:mfs-fourier-coefficients} again, we obtain
\begin{align*}
\biggl(
\sum_{n=N+1}^\infty\abs{\widehat{\tilde{f}_k}(n)}^2
\biggr)^{\frac{1}{2}}
&\leq\frac{1}{4\pi}
\sum_{\ell=1}^L\abs{a_\ell}
\biggl(
\sum_{n=N+1}^\infty
\frac{\abs{q_{k\ell}}^{2n}}{n^2}
\biggr)^{\frac{1}{2}}
\leq\frac{1}{4\pi(N+1)}
\sum_{\ell=1}^L
\frac{\abs{a_\ell}\abs{q_{k\ell}}^{N+1}}
{\sqrt{1-\abs{q_{k\ell}}^2}}.
\end{align*}
Inserting this bound into \eqref{eqn:mfs-truncation-parseval} proves the second estimate.
\end{proof}

\paragraph{Normalisation of Fourier polynomials.}
In applications of Proposition~\ref{prop:DtN-error}, the boundary data may be obtained by normalising a Fourier polynomial $\varphi$.
Choosing $f=\varphi/\abs{\varphi}$, we need to control the normalisation error $f-\varphi$ in $H^2$.
The following lemma provides an explicit bound in terms of the Fourier coefficients of $\abs{\varphi}^2-1$.

\begin{lemma}
\label{lem:normalization-trace-errors}
Given $M\in\N$, let $\varphi\colon\Xi\to\R^3$ be a Fourier polynomial of degree at most $M$. 
If
\begin{align*}
\delta\vcentcolon=\max_{1\leq k\leq b}
\sum_{\abs{n}\leq 2M}
\abs[\Big]{\widehat{\bigl(\abs{\varphi_k}^2-1\bigr)}(n)}<1
\end{align*}
then $\gamma\vcentcolon=\varphi/\abs{\varphi}$ is well-defined and
\begin{align*}
\nm{\gamma-\varphi}_{H^2}
&\leq
\frac{\delta\bigl((1-\delta)^2+9M^2\bigr)}
{2(1-\delta)^3}
\nm{\varphi}_{L^2}.
\end{align*} 
Moreover,
\begin{align}
\label{eqn:reciprocal-bounds}
\nm[\big]{\abs{\varphi}^{-1}}_{L^\infty}
&\leq\frac{1}{\sqrt{1-\delta}},
&
\nm[\big]{(\abs{\varphi}^{-1})'}_{L^\infty}
&\leq\frac{M\delta}{(1-\delta)^{3/2}}.
\end{align}
\end{lemma}

\begin{proof}
Let $w\vcentcolon=\abs{\varphi}^2-1$.
Since $w$ is a Fourier polynomial of degree at most $2M$, the definition of $\delta$ implies $\nm{w}_{L^\infty}\leq\delta<1$ and $\nm{w'}_{L^\infty}\leq2M\delta$.
Consequently, $\abs{\varphi}^2\geq1-\delta>0$, so $\gamma$ is well-defined.
This proves the first part of \eqref{eqn:reciprocal-bounds}.
The second follows from
$(\abs{\varphi}^{-1})'
=-\frac{w'}{2(1+w)^{3/2}}$.

By the binomial series, we have
\begin{align}
\label{eqn:20261103a}
\gamma-\varphi
&=\varphi\bigl((1+w)^{-\frac12}-1\bigr)
=\sum_{m=1}^\infty(-1)^m c_m\varphi w^m,
&
c_m&\vcentcolon=\frac{1}{4^m}\binom{2m}{m}.
\end{align}
For every $m\in\N$, the product $\varphi w^m$ is a Fourier polynomial of degree at most $(2m+1)M$.
Hence, Parseval's identity and $\nm{w}_{L^\infty}\leq\delta$ imply
\begin{align*}
\nm[\Big]{\frac{\partial^j(\varphi w^m)}{\partial\theta^j}}_{L^2}
&\leq\bigl((2m+1)M\bigr)^j\nm{\varphi w^m}_{L^2}
\leq\bigl((2m+1)M\bigr)^j\delta^m\nm{\varphi}_{L^2}
\end{align*}
for $j\in\{0,1,2\}$.
The coefficients in \eqref{eqn:20261103a} satisfy 
$c_1=\frac{1}{2}$ and 
$(2m+2)c_{m+1}=(2m+1)c_m$, and hence, by induction, 
\begin{align*}
c_m\leq\frac{m+1}{4m}\leq\frac12
\end{align*}
for every $m\in\N$. 
Since $\varphi w^m$ has Fourier degree at most $(2m+1)M$, we also have
\begin{align*}
\nm{\varphi w^m}_{H^2}
&\leq\bigl(1+(2m+1)^2M^2\bigr)\nm{\varphi w^m}_{L^2}
\leq\bigl(1+(2m+1)^2M^2\bigr)\delta^m\nm{\varphi}_{L^2}.
\end{align*}
Since $c_m\leq\frac12$ and $\delta<1$, the series in \eqref{eqn:20261103a} converges absolutely in $H^2$.
Consequently, it can be differentiated termwise up to order two, with convergence in $L^2$.
Using $2m+1\leq3m$ and the geometric series identities
\begin{align*}
\sum_{m=1}^\infty\delta^m
&=\frac{\delta}{1-\delta},
&
\sum_{m=1}^\infty m\delta^m
&=\frac{\delta}{(1-\delta)^2},
&
\sum_{m=1}^\infty m(m+1)\delta^m
&=\frac{2\delta}{(1-\delta)^3},
\end{align*}
we conclude from the triangle inequality that
\begin{align*}
\nm{\gamma-\varphi}_{L^2}
&\leq\frac{\delta\nm{\varphi}_{L^2}}{2(1-\delta)},
&
\nm{(\gamma-\varphi)'}_{L^2}
&\leq\frac{3M\delta\nm{\varphi}_{L^2}}{2(1-\delta)^2},
&
\nm{(\gamma-\varphi)''}_{L^2}
&\leq\frac{9M^2\delta\nm{\varphi}_{L^2}}{2(1-\delta)^3},
\end{align*}
where for the last inequality we also used
$c_m(2m+1)^2\leq\frac94m(m+1)$.
Finally, recalling the Sobolev norm from \eqref{eqn:Sobolev-norm}, we obtain
\begin{align*}
\nm{\gamma-\varphi}_{H^2}
&\leq
\frac{\delta\nm{\varphi}_{L^2}}{2(1-\delta)^3}
\sqrt{(1-\delta)^4+18M^2(1-\delta)^2+81M^4}
=\frac{\delta\bigl((1-\delta)^2+9M^2\bigr)}
{2(1-\delta)^3}
\nm{\varphi}_{L^2},
\end{align*}
which proves the claim.
\end{proof}


\section{Existence of a free boundary minimal immersion}

In this section we prove the existence of a branched free boundary minimal immersion by applying Proposition~\ref{prop:existence} to the gradient map $F$ from Definition~\ref{defn:gradient-map}.
This requires quantitative estimates for the residual, the inverse of the linearisation and the nonlinear remainder.

\subsection{Approximate solution}
\label{sec:initial}
 
Let $\mathcal P_M$ denote the space of $\R^3$-valued Fourier polynomials on $\Xi$ with Fourier modes $\abs{n}\leq M$, where our choice of $M$ is recorded in Appendix \ref{section-numerical-values}. 
The approximate solution $(\gamma_0,c_0)$ is given by a circle configuration $c_0 \in \mathcal{C}$ and a Fourier polynomial $P \in \mathcal P_M$ defining $\gamma_0 \vcentcolon= P/\abs{P}$.
The configuration $c_0$ is stored as a list of $b=61$ vectors in $\R^3$, whose normalisations are the centres $p_1,\ldots,p_b\in\Sp^2$, and a list of the corresponding geodesic radii $r_1,\ldots,r_b\in\interval{0,\pi}$.  
For each $k\in\{1,\ldots,b\}$, we also store a vector in $\R^3$ whose normalised orthogonal projection onto $T_{p_k}\Sp^2$ yields $e_{1,k}$, and set $e_{2,k}=p_k\times e_{1,k}$.
The resulting circle domain $\Omega_{c_0}$ is visualised in Figure~\ref{fig:g0b61}.
The Fourier polynomial $P_k$ on the $k$-th boundary component is stored as a list of $2M+1$ coefficients $a_{k,0},a_{k,1},b_{k,1},\ldots,a_{k,M},b_{k,M}\in\R^3$ with respect to the real orthonormal Fourier basis $f_0,\ldots,f_{2M}\colon\Sp^1\to\R$ of degree at most $M$ defined by 
\begin{align}\label{eq:realFourierBasis}
f_0(\theta)&=\frac{1}{\sqrt{2\pi}},
&
f_{2n-1}(\theta)&=\frac{\cos(n\theta)}{\sqrt{\pi}},
&
f_{2n}(\theta)&=\frac{\sin(n\theta)}{\sqrt{\pi}}
\end{align}
for $n\in\{1,\ldots,M\}$, so that
\begin{align}\label{eq:PkFourierCoeffs}
P_k(\theta)&=\frac{a_{k,0}}{\sqrt{2\pi}}+\frac{1}{\sqrt{\pi}}\sum_{n=1}^M\bigl(a_{k,n}\cos(n\theta)+b_{k,n}\sin(n\theta)\bigr).
\end{align}
The complex Fourier coefficients of $P_k$ as defined in \eqref{eqn:Fourier-coefficients} are thus given by
\begin{align*}
\widehat P_k(0)&=\frac{a_{k,0}}{\sqrt{2\pi}},
&
\widehat P_k(n)&=\frac{a_{k,n}-ib_{k,n}}{2\sqrt{\pi}},
&
\widehat P_k(-n)&=\frac{a_{k,n}+ib_{k,n}}{2\sqrt{\pi}}
\end{align*}
for $n\in\{1,\ldots,M\}$.

We use non-rigorous numerics to compute the circle configuration $c_0 \in \mathcal{C}$ and Fourier coefficients of $P$ that give an approximate critical point of the reduced Dirichlet energy from equation~\eqref{eqn:reduced-energy-DtN}.
We carry this out in two stages:
\begin{itemize}[wide]
\item
\emph{Initialisation:}
we initialise $b$ points following a randomly perturbed Fibonacci configuration on $\Sp^2$ and then move them so as to locally minimise their Thomson energy $\sum_{i<j}\abs{p_i-p_j}^{-1}$.
We use the resulting points as the centres of both the domain holes and the initial target circles.
We note that it is not expected in general that the positions that minimise this energy are equal to the centres of mass of the boundary components for free boundary minimal surfaces, but this initial guess is sufficient for our optimisation scheme.
We then predict the total boundary length and distribute it among the components based on heuristic extrapolation formulae from examples with few boundary components.
This determines the target radii.
We determine each domain radius $r_*$ from the target radius $R_*$ by the heuristic $\tan (r_*/2)=1/(b-1) \tan ((\pi/2-R_*)/2)$, which recovers the data of the critical catenoid when $b=2$.
Each initial target circle is thus parametrised by an $\R^3$-valued Fourier polynomial of degree one.

\item
\emph{Fourier polishing:}
We then use Newton iteration to correct the domain moduli and boundary parametrisation simultaneously.
We express each boundary correction as two Fourier polynomials in a tangent frame along the current boundary map.
We take several Newton steps at successive degree cutoffs of $4,8,12,22,36,44$ for these components.
After each correction we normalise the boundary maps and compute the Fourier expansion of the result.
We store a final approximate solution for $P$ of degree (at most) $64$.
\end{itemize}

\subsection{Residual bound}
\label{sec:residual}

After obtaining the approximate solution $(\gamma_0,c_0)\in H^2(\Xi,\Sp^2)\times\mathcal{C}$ visualised in Figure~\ref{fig:g0b61}, we need to verify that the residual $F(0,0,0)$ is sufficiently small.
Definition \ref{defn:gradient-map} of $F$ depends on the choice of a slice for the Möbius action.
We begin by specifying this slice, thereby completing the definition of $F$.
We then derive a rigorous bound for $F(0,0,0)$.

\paragraph{Choice of a Möbius slice.}
In the following, we choose the Möbius slice $U$ and the local coordinates $c_\mu$ introduced in Section~\ref{sec:gauge}. 
Recalling the action \eqref{eqn:Conf-action} of the Möbius group $\Conf(\Sp^2)$ on configurations, we identify the tangent space of $\Conf(\Sp^2)$ at the identity with the space of conformal vector fields
\begin{align*}
X_{u,w}(x)=u\times x+w-\sk{w,x}x\in T_x\Sp^2
\end{align*}
defined for $u,w\in\R^3$ and $x\in\Sp^2$. 
Here $u\times x$ generates rotations of the domain sphere, while $w-\sk{w,x}x$ describes the remaining three conformal directions. 
Let $t\mapsto\phi(t)\in\Conf(\Sp^2)$ be the flow of $X_{u,w}$ with $\phi(0)=\mathrm{id}$, and consider the curve
\begin{align}
\label{eqn:20261027a}
c(t)\vcentcolon=\phi(t)\cdot c_0=\bigl(p_1(t),\ldots,p_b(t);r_1(t),\ldots,r_b(t)\bigr)
\end{align}
in $\mathcal{C}$. 
By \eqref{eqn:Conf-action} and \eqref{eqn:k-th_boundary_circle}, we have
$\sk{\phi(t)(x),p_k(t)}=\cos r_k(t)$ 
for every $x\in\partial_k\Omega_{c_0}$ and $k\in\{1,\ldots,b\}$.
Differentiating at $t=0$ and using $\sk{\dot p_k,p_k}=0$, we obtain
\begin{align}
\label{eqn:20261027b}
\dot p_k&=u\times p_k+\cos r_k\bigl(w-\sk{w,p_k}p_k\bigr),
&
\dot r_k&=-\sin r_k\sk{w,p_k}. 
\end{align}
As before, the dot denotes derivatives evaluated at $t=0$ and $p_k,r_k$ are the centre and radius of the $k$-th circle in the reference configuration $c_0$. 
For each $k\in\{1,\ldots,b\}$, recalling the orthonormal basis $(e_{1,k},e_{2,k})$ of $T_{p_k}\Sp^2$ from \eqref{eqn:20261009b}, we introduce local coordinates $(a_1,a_2,s)$ on $\Sp^2\times\interval{0,\pi}$ near $(p_k,r_k)$ through the parametrisation
\begin{align}\notag
\R^3&\to\Sp^2\times\interval{0,\pi},
\\\label{eqn:configuration-local-coordinates}
(a_1,a_2,s)&\mapsto
\biggl(
\frac{p_k+a_1e_{1,k}+a_2e_{2,k}}{\sqrt{1+a_1^2+a_2^2}},~
2\arctan\Bigl(e^s\tan\frac{r_k}{2}\Bigr)
\biggr).
\end{align}
Differentiating \eqref{eqn:20261027a} in these coordinates at $t=0$ and recalling \eqref{eqn:20261027b}, we obtain
\begin{align}
\label{eqn:20261027c}
\begin{pmatrix}
\dot a_1\\
\dot a_2\\
\dot s
\end{pmatrix}
&=\begin{pmatrix}
\sk{\dot p_k,e_{1,k}}\\
\sk{\dot p_k,e_{2,k}}\\
\dot r_k/\sin r_k
\end{pmatrix}
=\begin{pmatrix}
(p_k\times e_{1,k})^{\intercal} & \cos r_k\,e_{1,k}^{\intercal}\\
(p_k\times e_{2,k})^{\intercal} & \cos r_k\,e_{2,k}^{\intercal}\\
0 & -p_k^{\intercal}
\end{pmatrix}
\begin{pmatrix}
u\\
w
\end{pmatrix}.
\end{align}
We select the two centre coordinates of three circles to remain fixed during the search for a zero of $F$, and denote by
\begin{align}
\label{eqn:infinitesimal-mobius-on-fixed-parameters}
\zeta_\mu\in\R^{6\times6}
\end{align}
the matrix consisting of the corresponding six rows of the matrices in \eqref{eqn:20261027c}, for $k\in\{1,\ldots,b\}$, with $r_k$ replaced by $r_k(\mu)$. 
The concrete choice of \eqref{eqn:infinitesimal-mobius-on-fixed-parameters} is specified in the companion code. 
Fixing these $6$ coordinates at their values at $c_0$ and denoting the remaining $3b-6$ coordinates by $\mu\in B\subset\R^{3b-6}$, we obtain the parametrisation $\mu\mapsto c_\mu$ from \eqref{eqn:perturbation-configuration}. 
For concreteness, we set
\begin{align*}
B&\vcentcolon=\{\mu\in\R^{3b-6}\mid \abs{\mu}<r_*\},
&
U\vcentcolon=\{c_\mu\in\mathcal{C}\mid \mu\in B\},
\end{align*}
where $r_*$ is recorded in Appendix~\ref{section-numerical-values}. 
The submanifold $U\subset\mathcal{C}$ is a Möbius slice as introduced in Section~\ref{sec:gauge} if 
\begin{align*}
T_{c_\mu}\mathcal{C}
=T_{c_\mu}U\oplus T_{c_\mu}\bigl(\Conf(\Sp^2)\cdot c_\mu\bigr)
\end{align*}
for every $\mu\in B$, which by \eqref{eqn:20261027c} holds if and only if $\zeta_\mu$ is invertible for every $\mu\in B$.
The following lemma verifies this condition quantitatively.
 
\begin{lemma}
\label{lem:Mobius-transversality}
There exist constants $r_*>0$ and $\epsilon_2>\sqrt{2}r_*$, recorded in Appendix~\ref{section-numerical-values}, such that
for all $z\in\R^6$ and $\mu\in B$ 
\begin{align*}
\abs{\zeta_\mu z}&\geq\bigl(\epsilon_2-\sqrt2\abs{\mu}\bigr)\abs{z}.
\end{align*}
\end{lemma}

\begin{proof}
The companion code computes $\zeta_0$ as an interval matrix.
We compute an interval enclosure of $\zeta_0^{-1}$ and obtain the bound\begin{align}
\label{eqn:zeta-invertible-1}
\abs{\zeta_0 z}
&\geq\frac{\abs{z}}{\nm{\zeta_0^{-1}}_2}
\geq\frac{\abs{z}}{\nm{\zeta_0^{-1}}_{\mathrm{Frob}}}
\geq\epsilon_2\abs{z}.
\end{align}
By the definition of $\zeta_\mu$ in \eqref{eqn:infinitesimal-mobius-on-fixed-parameters}, the selected centres and their orthonormal frames remain fixed as $\mu$ varies.
Hence, by \eqref{eqn:20261027c}, the only entries of $\zeta_\mu$ depending on $\mu$ are those involving $\cos r_k(\mu)$.
For the logarithmic radius coordinate $s$ defined in \eqref{eqn:configuration-local-coordinates}, we have $\frac{\partial}{\partial s}\cos r_k=-\sin^2r_k$, so these factors are Lipschitz continuous with constant $1$.
Each such factor occurs in two rows with coefficients given by the orthonormal vectors $e_{1,k}$ and $e_{2,k}$, and different radii affect distinct rows.
Therefore,
\begin{align}
\label{eqn:zeta-invertible-2}
\nm{\zeta_\mu-\zeta_0}_2
\leq\nm{\zeta_\mu-\zeta_0}_{\mathrm{Frob}}
\leq\sqrt{2}\abs{\mu}.
\end{align}
Combining \eqref{eqn:zeta-invertible-1} and \eqref{eqn:zeta-invertible-2} proves the claim.
\end{proof}

To estimate the residual, we use the following lemma to bound the projection onto the frame $\Gamma$ in terms of a cross product with the Fourier polynomial $P$.

\begin{lemma}
\label{lem:boundary-projection}
Let $P\colon\Xi\to\R^3$ be a Fourier polynomial such that $P\times P'$ is nowhere zero, and let $\Gamma$ be the frame from \eqref{eqn:frame} along $\gamma_0=P/\abs{P}$.
Then, for every $w\in H^1(\Xi,\R^3)$, we have
\begin{align*}
\nm{\Gamma^{\intercal}w}_{H^1(\Xi)}
&\leq
\bigl(1+\nm{\Gamma'}_{L^\infty(\Xi)}\bigr)
\Bigl(
\nm{\abs{P}^{-1}}_{L^\infty(\Xi)}
+\nm{(\abs{P}^{-1})'}_{L^\infty(\Xi)}
\Bigr)
\nm{P\times w}_{H^1(\Xi)}.
\end{align*}
\end{lemma}

\begin{proof}
By definition \eqref{eqn:frame}, the columns of $\Gamma$ are orthonormal and satisfy
\begin{align*}
\Gamma^{\intercal}(\gamma_0\times w)
&=\begin{pmatrix}0&1\\-1&0\end{pmatrix}\Gamma^{\intercal}w.
\end{align*}
Since the matrix on the right is constant and orthogonal, we obtain
\begin{align}
\label{eqn:boundary-projection-rotation}
\nm{\Gamma^{\intercal}w}_{H^1(\Xi)}
&=\nm{\Gamma^{\intercal}(\gamma_0\times w)}_{H^1(\Xi)}.
\end{align}
Moreover, the product rule and the pointwise orthonormality of the columns of $\Gamma$ imply
\begin{align}  
\label{eqn:frame-multiplication-bound}
\nm{\Gamma^{\intercal}v}_{H^1(\Xi)}
&\leq\bigl(1+\nm{\Gamma'}_{L^\infty(\Xi)}\bigr)\nm{v}_{H^1(\Xi)}
\end{align}
for every $v\in H^1(\Xi,\R^3)$, in particular for $v=\gamma_0\times w$. 
Similarly, the product rule implies
\begin{align}\label{eqn:20261107c}
\nm{\gamma_0\times w}_{H^1(\Xi)}
&=\nm{\abs{P}^{-1}(P\times w)}_{H^1(\Xi)}
\leq
\Bigl(
\nm{\abs{P}^{-1}}_{L^\infty(\Xi)}
+\nm{(\abs{P}^{-1})'}_{L^\infty(\Xi)}
\Bigr)
\nm{P\times w}_{H^1(\Xi)}.
\end{align}
Combining \eqref{eqn:boundary-projection-rotation}--\eqref{eqn:20261107c} proves the claim.
\end{proof}

\paragraph{Bounding the residual.}
The first estimate needed to apply Proposition~\ref{prop:existence} is an upper bound $\varepsilon$ for the size of the residual $F(0,0,0)$, which we subsequently abbreviate by $F(0)$.

\begin{proposition}
\label{prop:residual-bound}
Recalling the gradient map $F$ from Definition~\ref{defn:gradient-map} and the Banach space $Y$ from Proposition~\ref{prop:existence}, we have ${\nm{F(0)}}_Y\leq\varepsilon$, where $\varepsilon$ is recorded in Appendix~\ref{section-numerical-values}.
\end{proposition}

\begin{proof}
Recalling the expression for $G_c(v)$ from Lemma~\ref{lem:perturbation-gradient}, we have by definition 
\begin{align*}
F(0)
&=\biggl(
\Gamma^{\intercal}\DN_{c_0}\gamma_0,\;
\Bigl(\frac{\partial}{\partial\mu_j}\Big\vert_{\mu=0}E(\gamma_0,c_\mu)\Bigr)_{j=1,\ldots,3b-6},\;
0
\biggr), & 
\gamma_0&=\frac{P}{\abs{P}}
\end{align*}
where the Fourier polynomial $P$ of degree $M$ was introduced in Section~\ref{sec:initial}.
Given $N\geq M$, let $\tilde{f}$ be the boundary trace of the harmonic map $\Phi$ defined in \eqref{eqn:mfs-field}, with coefficients satisfying $a_1+\ldots+a_L=0$.
Applying Proposition~\ref{prop:DtN-error} with $c=c_0$, $f=\gamma_0$ and $\varphi=P$, we obtain the Fourier polynomial
\begin{align}
\label{eqn:tilde_h}
\tilde{h}&=\DN_{\mathrm{dec}}P+\Pi_N(\DN_{c_0}-\DN_{\mathrm{dec}})\tilde{f},
\end{align}
which approximates $\DN_{c_0}\gamma_0$. 
For each $j\in\{1,\ldots,3b-6\}$, let $\tilde{\alpha}_j,\tilde{\beta}_j\colon\Xi\to\R$ denote the normalised boundary velocities from Lemma~\ref{lem:DtN-derivative} corresponding to $\dot c=\frac{\partial c_\mu}{\partial\mu_j}\big\vert_{\mu=0}$ and let 
\begin{align}\label{eqn:20261107b}
z_j\vcentcolon=\int_\Xi\Bigl(
\frac{\tilde{\alpha}_j}{2}\bigl(\abs{P'}^2-\abs{\tilde{h}}^2\bigr)
-\tilde{\beta}_j\sk{P',\tilde{h}}
\Bigr)\,d\theta.
\end{align} 
By Corollary~\ref{cor:variation-c-theta} and the definition of $z_j$, polarisation yields
\begin{align}
\label{eqn:20261106a}
\frac{\partial}{\partial\mu_j}\Big\vert_{\mu=0}E(\gamma_0,c_\mu)-z_j
&=\frac{1}{2}\int_\Xi\biggl(
\tilde{\alpha}_j\sk{\gamma_0'-P',\gamma_0'+P'}
-\tilde{\alpha}_j\sk{\DN_{c_0}\gamma_0-\tilde{h},\DN_{c_0}\gamma_0+\tilde{h}}
\\\notag
&\hphantom{{}=\frac{1}{2}\int_\Xi\biggl(}
-\tilde{\beta}_j\sk{\gamma_0'-P',\DN_{c_0}\gamma_0+\tilde{h}}
-\tilde{\beta}_j\sk{\DN_{c_0}\gamma_0-\tilde{h},\gamma_0'+P'}
\biggr)\,d\theta.
\end{align}
For all $a,b,c,d\in\R^3$ and $\alpha,\beta\in\R$, the Cauchy--Schwarz inequality implies
\begin{align}
\notag
\abs[\big]{\alpha(\sk{a,c}-\sk{b,d})-\beta(\sk{a,d}+\sk{b,c})}^2
&\leq
\bigl(\abs{\alpha a-\beta b}^2+\abs{\beta a+\alpha b}^2\bigr)
\bigl(\abs{c}^2+\abs{d}^2\bigr)
\\\label{eqn:20261106b}
&=(\alpha^2+\beta^2)
\bigl(\abs{a}^2+\abs{b}^2\bigr)
\bigl(\abs{c}^2+\abs{d}^2\bigr).
\end{align}
Applying \eqref{eqn:20261106b} to the integrand in \eqref{eqn:20261106a}, integrating over $\Xi$ and using the Cauchy--Schwarz inequality, we obtain
\begin{align}
\notag
&\abs[\bigg]{\frac{\partial}{\partial\mu_j}\Big\vert_{\mu=0}E(\gamma_0,c_\mu)-z_j}^2
\\\notag
&\leq\frac{1}{4}
\nm[\Big]{\tilde{\alpha}_j^2+\tilde{\beta}_j^2}_{L^\infty(\Xi)}
\Bigl(
\nm{\gamma_0'-P'}_{L^2(\Xi)}^2
+\nm{\DN_{c_0}\gamma_0-\tilde{h}}_{L^2(\Xi)}^2
\Bigr)
\Bigl(
\nm{\gamma_0'+P'}_{L^2(\Xi)}^2
+\nm{\DN_{c_0}\gamma_0+\tilde{h}}_{L^2(\Xi)}^2
\Bigr)
\\\label{eqn:residual-proof-E-part}
&=\vcentcolon\varepsilon_j^2.
\end{align} 
Applying inequality \eqref{eqn:frame-multiplication-bound} with $v=\DN_{c_0}\gamma_0-\tilde{h}$, the triangle inequality yields  
\begin{align}
\notag
\nm[\big]{\Gamma^{\intercal}\DN_{c_0}\gamma_0}_{H^1(\Xi)}
&\leq
\nm[\big]{\Gamma^{\intercal}\tilde{h}}_{H^1(\Xi)}
+\nm[\big]{\Gamma^{\intercal}(\DN_{c_0}\gamma_0-\tilde{h})}_{H^1(\Xi)}
\\\label{eqn:residual-proof-Gamma-part}
&\leq
\nm[\big]{\Gamma^{\intercal}\tilde{h}}_{H^1(\Xi)}
+\bigl(1+\nm{\Gamma'}_{L^\infty(\Xi)}\bigr)
\nm[\big]{\DN_{c_0}\gamma_0-\tilde{h}}_{H^1(\Xi)}.
\end{align}
Combining \eqref{eqn:residual-proof-E-part} and \eqref{eqn:residual-proof-Gamma-part} yields
\begin{align}
\label{eqn:residual-total-bound}
\nm{F(0)}_Y^2
&\leq
\Bigl(
\nm[\big]{\Gamma^{\intercal}\tilde{h}}_{H^1(\Xi)}
+\bigl(1+\nm{\Gamma'}_{L^\infty(\Xi)}\bigr)
\nm[\big]{\DN_{c_0}\gamma_0-\tilde{h}}_{H^1(\Xi)}
\Bigr)^2
+\sum_{j=1}^{3b-6}\bigl(\abs{z_j}+\varepsilon_j\bigr)^2.
\end{align} 
It remains to bound the terms on the right-hand side of \eqref{eqn:residual-total-bound} by explicitly computable quantities for which the companion code provides rigorous upper bounds, verifying that $\nm{F(0)}_Y\leq\varepsilon$. 

Applying Proposition~\ref{prop:DtN-error} with $c=c_0$, $f=\gamma_0$ and $\varphi=P$, we obtain
\begin{align*}
&\nm{\DN_{c_0}\gamma_0-\tilde{h}}_{H^1(\Xi)}
\\
&\leq
\bigl(1+C_1(\omega)\bigr)\nm{\gamma_0-P}_{H^2(\Xi)}
+C_1(\omega)\nm{P-\tilde{f}}_{L^2(\Xi)}
+\nm[\big]{(\mathrm{Id}-\Pi_N)(\DN_{c_0}-\DN_{\mathrm{dec}})\tilde{f}}_{H^1(\Xi)},
\end{align*}
where $\omega$ is the conformal collar width of $c_0$ from \eqref{eqn:conformal_collar_width}.
The first term is bounded by Lemma~\ref{lem:normalization-trace-errors}, while the remaining two terms are bounded by Proposition~\ref{prop:mfs-truncation-bounds}.
Since $\tilde{h}$ is a Fourier polynomial, Lemma~\ref{lem:boundary-projection} yields
\begin{align}
\label{eqn:20261107a}
\nm{\Gamma^{\intercal}\tilde{h}}_{H^1(\Xi)}
&\leq
\bigl(1+\nm{\Gamma'}_{L^\infty(\Xi)}\bigr)
\Bigl(
\nm{\abs{P}^{-1}}_{L^\infty(\Xi)}
+\nm{(\abs{P}^{-1})'}_{L^\infty(\Xi)}
\Bigr)
\nm{P\times\tilde{h}}_{H^1(\Xi)}.
\end{align} 
Since $P\times\tilde{h}$ is a Fourier polynomial, its $H^1(\Xi)$-norm can be computed using finite Fourier arithmetic.
The remaining factors in \eqref{eqn:20261107a} are bounded by interval arithmetic using the explicit Fourier polynomial $P$ and lower bounds for $\abs{P}$ and $\abs{P\times P'}$ on $\Xi$.
For the factors $\nm{\abs{P}^{-1}}_{L^\infty(\Xi)}$ and $\nm{(\abs{P}^{-1})'}_{L^\infty(\Xi)}$, we use \eqref{eqn:reciprocal-bounds} with $\varphi=P$.

By \eqref{eqn:normal-velocity-fourier} and \eqref{eqn:tangential-velocity-fourier}, the normalised boundary velocities $\tilde{\alpha}_j$ and $\tilde{\beta}_j$ are Fourier polynomials of degree at most one.
Hence, $z_j$ defined in \eqref{eqn:20261107b} can be computed using finite Fourier arithmetic.
To bound $\varepsilon_j$ in \eqref{eqn:residual-proof-E-part}, we use Lemma~\ref{lem:normalization-trace-errors} and Proposition~\ref{prop:DtN-error}, together with Proposition~\ref{prop:mfs-truncation-bounds}, and the following elementary inequalities.
\begin{align*}
\nm{\gamma_0'-P'}_{L^2(\Xi)}
&\leq\nm{\gamma_0-P}_{H^2(\Xi)},
\\
\nm{\gamma_0'+P'}_{L^2(\Xi)}
&\leq2\nm{P'}_{L^2(\Xi)}+\nm{\gamma_0'-P'}_{L^2(\Xi)},
\\
\nm{\DN_{c_0}\gamma_0+\tilde{h}}_{L^2(\Xi)}
&\leq2\nm{\tilde{h}}_{L^2(\Xi)}+\nm{\DN_{c_0}\gamma_0-\tilde{h}}_{L^2(\Xi)}.
\end{align*}
The norms of $P'$ and $\tilde{h}$ are computed using finite Fourier arithmetic, while $\nm{\sqrt{\tilde{\alpha}_j^2+\tilde{\beta}_j^2}}_{L^\infty(\Xi)}$ is bounded by interval arithmetic.  
\end{proof}

\subsection{Linear estimate}

The second estimate required by Proposition~\ref{prop:existence} controls the linearised operator $\mathscr{L}\vcentcolon=DF\vert_0$ from Proposition~\ref{prop:linearisation}. 
In this subsection, we prove that $\mathscr{L}\colon X\to Y$ is an isomorphism and obtain a quantitative bound for its inverse.
 
\paragraph{High-mode/low-mode splitting.}
Let $N\in\N$ be fixed. 
(Our choice is recorded in Appendix~\ref{section-numerical-values}.)
We take the basis $\xi_1,\xi_2,\xi_3$ from \eqref{eqn:rotation-fields} to be the standard basis of $\R^3$.
We introduce the Hilbert space
\begin{align*}
\mathcal{H}\vcentcolon=L^2(\Xi,\R^2)\times\R^{3b-6}\times\R^3
\end{align*}
with the canonical inner product, and extend the Fourier truncation $\Pi_N$ from Definition~\ref{def:fourier-truncation} to $\mathcal{H}$ by $\Pi_N(v,\mu,\lambda)=(\Pi_Nv,\mu,\lambda)$.
This allows the orthogonal decomposition
\begin{align}
\label{eqn:orthogonal_decomposition}
\mathcal{H}&=\mathcal{H}_{\leq N}\oplus\mathcal{H}_{>N},
&
\mathcal{H}_{\leq N}&\vcentcolon=\Pi_N\mathcal{H},
&
\mathcal{H}_{>N}&\vcentcolon=(\mathrm{Id}-\Pi_N)\mathcal{H}.
\end{align}
In particular, all configuration and rotation variables belong to the  vector space $\mathcal{H}_{\leq N}$ which has finite dimension 
\begin{align*}
\dim \mathcal{H}_{\leq N}&=2b(2N+1)+3b-3.
\end{align*}

\begin{lemma}
\label{lem:linearisation-self-adjoint}
The linearised operator $\mathscr{L}=DF\vert_0$ extends to an unbounded self-adjoint operator
\[
\mathscr{L}\colon Y=H^1(\Xi,\R^2)\times\R^{3b-6}\times\R^3\subset\mathcal{H}\to\mathcal{H}.
\]
\end{lemma}
 
\begin{proof}
By Lemma~\ref{lem:decoupled-Fourier}, the operator $\DN_{\mathrm{dec}}$ acting componentwise on vector-valued $L^2$-functions is self-adjoint with domain $H^1(\Xi,\R^2)$.
Consequently, the operator $\DN_{\mathrm{dec}}\oplus0\oplus0$ is self-adjoint on $\mathcal{H}$ with domain $Y$.
We claim that $\mathscr{L}-(\DN_{\mathrm{dec}}\oplus0\oplus0)$ extends to a bounded symmetric operator on $\mathcal{H}$.
By Proposition~\ref{prop:linearisation} and the identity $\Gamma^{\intercal}\Gamma=\mathrm{Id}$, we have
\begin{align*}
S_{vv}w-\DN_{\mathrm{dec}}w
&=\Gamma^{\intercal}
\Bigl(\DN_{\mathrm{dec}}(\Gamma w)-\Gamma\DN_{\mathrm{dec}}w\Bigr)
+\Gamma^{\intercal}
\bigl(\DN_{c_0}-\DN_{\mathrm{dec}}\bigr)(\Gamma w)
-\sk{\gamma_0,\DN_{c_0}\gamma_0}\,w
\end{align*}
for all $w\in H^2(\Xi,\R^2)$.
By Lemma~\ref{lem:decoupled-Fourier} and the inequality
$\abs{\abs{n+m}-\abs{n}}\leq\abs{m}$, we have
\begin{align*}
\nm[\big]{\DN_{\mathrm{dec}}(\Gamma w)-\Gamma\DN_{\mathrm{dec}}w}_{L^2}
&\leq
\biggl(
\max_{k\in\{1,\ldots,b\}}
\sum_{m\in\Z}\abs{m}\abs[\big]{\widehat{\Gamma_k}(m)}
\biggr)\nm{w}_{L^2}.
\end{align*}
The sum is finite since $\Gamma$ is smooth.
By Lemma~\ref{lem:decoupled-collar} and the pointwise orthonormality of $\Gamma$, the operator $w\mapsto\Gamma^{\intercal}(\DN_{c_0}-\DN_{\mathrm{dec}})(\Gamma w)$ is bounded on $L^2(\Xi,\R^2)$.
Moreover, multiplication by $\sk{\gamma_0,\DN_{c_0}\gamma_0}$ is bounded on $L^2(\Xi,\R^2)$ since $\gamma_0\in H^2(\Xi,\R^3)$ and $\DN_{c_0}\gamma_0\in H^1(\Xi,\R^3)\subset L^\infty(\Xi,\R^3)$.
Consequently, $S_{vv}-\DN_{\mathrm{dec}}$ extends to a bounded operator $L^2(\Xi,\R^2)\to L^2(\Xi,\R^2)$.

By Proposition~\ref{prop:linearisation}, the operators $S_{v\mu}\colon\R^{3b-6}\to L^2(\Xi,\R^2)$ and $\rho\colon\R^3\to L^2(\Xi,\R^2)$ are bounded. 
Moreover, $S_{\mu\mu}$ is a symmetric matrix.
Consequently, $\mathscr{L}-(\DN_{\mathrm{dec}}\oplus0\oplus0)$ extends to a bounded symmetric operator on $\mathcal{H}$, which is therefore self-adjoint.
Since $\DN_{\mathrm{dec}}\oplus0\oplus0$ is self-adjoint on $\mathcal{H}$ with domain $Y$, the operator $\mathscr{L}$ extends to a self-adjoint operator on $\mathcal{H}$ with the same domain $Y$.
\end{proof} 

With respect to the orthogonal decomposition \eqref{eqn:orthogonal_decomposition} of $\mathcal{H}=\mathcal{H}_{\leq N}\oplus\mathcal{H}_{>N}$, we write
\begin{align}
\label{eqn:LABC}
\mathscr{L}&=
\begin{pmatrix}
\mathscr{A}&\mathscr{B}\\
\mathscr{B}^*&\mathscr{C}
\end{pmatrix},
&
\begin{aligned}
\mathscr{A}&\vcentcolon=\Pi_N\mathscr{L}\Pi_N,
\\
\mathscr{B}&\vcentcolon=\Pi_N\mathscr{L}(\mathrm{Id}-\Pi_N),
\\
\mathscr{C}&\vcentcolon=(\mathrm{Id}-\Pi_N)\mathscr{L}(\mathrm{Id}-\Pi_N).
\end{aligned}
\end{align}
The block decomposition holds on the domain $Y\subset\mathcal{H}$,  while $\mathscr{B}$ extends to a bounded operator $\mathscr{B}\colon\mathcal{H}_{>N}\to\mathcal{H}_{\leq N}$. 
The operators $\mathscr{A}\colon\mathcal{H}_{\leq N}\to\mathcal{H}_{\leq N}$ and $\mathscr{C}\colon Y\cap\mathcal{H}_{>N}\to\mathcal{H}_{>N}$ are self-adjoint.
Indeed, since $\Pi_N$ commutes with $\DN_{\mathrm{dec}}\oplus0\oplus0$, the boundedness of $\mathscr{B}$ and the self-adjointness of $\mathscr{C}$ follow from the proof of Lemma~\ref{lem:linearisation-self-adjoint}. 
 
\begin{lemma}
\label{lem:configuration-linearisation}
Given $\mu\in B$ and $i\in\{1,\ldots,3b-6\}$ let $\tilde{\alpha}_i,\tilde{\beta}_i$ denote the normalised boundary velocities from Lemma~\ref{lem:DtN-derivative} corresponding to $\partial c_\mu/\partial\mu_i$, and set
\begin{align*}
h_\mu\vcentcolon=\DN_{c_\mu}\gamma_0.
\end{align*}
Then, for all $i,j\in\{1,\ldots,3b-6\}$, the matrix $S_{\mu\mu}$ from Proposition~\ref{prop:linearisation} satisfies
\begin{align*}
{(S_{\mu\mu})}_{ij}
&=\int_\Xi
\biggl(
\frac{1}{2}\frac{\partial\tilde{\alpha}_i}{\partial\mu_j}
\Bigl(\abs{\gamma_0'}^2-\abs{h_\mu}^2\Bigr)
-\frac{\partial\tilde{\beta}_i}{\partial\mu_j}\sk{\gamma_0',h_\mu}
-\sk[\bigg]{\tilde{\alpha}_i h_\mu+\tilde{\beta}_i\gamma_0',
\frac{\partial h_\mu}{\partial\mu_j}}
\biggr)\bigg\vert_{\mu=0}\,d\theta.
\end{align*}
\end{lemma}
 
\begin{proof}
By Corollary~\ref{cor:variation-c-theta}, we have
\begin{align*}
\frac{\partial}{\partial\mu_i}E(\gamma_0,c_\mu)
&=\int_\Xi
\biggl(
\frac{1}{2}\tilde{\alpha}_i
\Bigl(\abs{\gamma_0'}^2-\abs{h_\mu}^2\Bigr)
-\tilde{\beta}_i\sk{\gamma_0',h_\mu}
\biggr)\,d\theta.
\end{align*}
By Lemma~\ref{lem:smoothness} and the smooth dependence of the boundary velocities on $\mu$, we may differentiate under the integral sign with respect to $\mu_j$. 
The claim then follows from the product rule, since $(S_{\mu\mu})_{ij}$ is the second partial derivative of $\mu\mapsto E(\gamma_0,c_\mu)$ at the origin by Proposition~\ref{prop:linearisation}.
\end{proof}

\paragraph{Fourier approximations.}
Recalling $\gamma_0=P/\abs{P}$ from Section~\ref{sec:initial}, the frame $\Gamma$ from \eqref{eqn:frame} satisfies
\begin{align}
\label{eqn:frame-Fourier}
\Gamma
&=\biggl(
\frac{P\times P'}{\abs{P\times P'}},
\frac{(P\times P')\times P}{\abs{P\times P'}\abs{P}}
\biggr).
\end{align}
We use \eqref{eqn:frame-Fourier} to construct a Fourier polynomial approximation $G\colon\Xi\to\R^{3\times2}$ of $\Gamma$ with errors bounded in $H^2(\Xi,\R^{3\times2})$. 
The error bounds are computed as follows: we subdivide each boundary circle into intervals. 
For $j=0,1,2$, evaluate $(\Gamma-G)^{(j)}$ at the interval midpoints and use a coarse bound for $(\Gamma-G)^{(j+1)}$ to bound it on the whole interval by the mean value inequality. The coarse derivative bounds are obtained by Fourier differentiation and interval evaluation applied to Eqn.~\ref{eqn:frame-Fourier}.
We also recall the Fourier polynomial $\tilde{h}$ from \eqref{eqn:tilde_h}, which approximates $\DN_{c_0}\gamma_0$ with an error bounded by Proposition~\ref{prop:DtN-error}.

\paragraph{Low-mode inverse bound.}
We estimate the norm of the inverse of $\mathscr{A}\colon\mathcal{H}_{\leq N}\to\mathcal{H}_{\leq N}$ by representing it as a symmetric matrix with respect to an orthonormal basis of the finite-dimensional space $\mathcal{H}_{\leq N}$.

For $k\in\{1,\ldots,b\}$, $s\in\{0,\ldots,2N\}$ and $j\in\{1,2\}$, let
\begin{align}
\label{eqn:Fourier-basis}
\varphi_{(k-1)(4N+2)+2s+j}\in\Pi_NL^2(\Xi,\R^2)
\end{align}
be the function equal to $f_s e_j$ on the $k$-th boundary component and zero elsewhere, where $e_1,e_2$ denote the standard basis of $\R^2$ and $\{f_0, \ldots, f_{2N}\}$ is the real orthonormal Fourier basis of degree at most $N$ as in \eqref{eq:realFourierBasis}. Then $\varphi_1,\ldots,\varphi_{2b(2N+1)}$ form an orthonormal basis of $\Pi_NL^2(\Xi,\R^2)$, which we combine with the standard Euclidean bases of $\R^{3b-6}$ and $\R^3$.

For the approximations of the Dirichlet-to-Neumann operator we fix a truncation degree $\widetilde N\in\N$, recorded in Appendix~\ref{section-numerical-values}.

\emph{Fourier--Fourier entries:}
For $i,j\in\{1,\ldots,(4N+2)b\}$, identifying $\varphi_j$ with  $(\varphi_j,0,0)\in\mathcal{H}_{\leq N}$, we have 
\begin{align}
\label{eqn:low-mode-Fourier-entries}
\sk{\mathscr{A}\varphi_i,\varphi_j}_{\mathcal{H}}
&=\int_\Xi
\sk[\big]{\DN_{c_0}(\Gamma\varphi_i),\Gamma\varphi_j}
-\sk[\big]{\gamma_0,\DN_{c_0}\gamma_0}\sk{\varphi_i,\varphi_j}
\,d\theta 
\end{align}
by Proposition~\ref{prop:linearisation}. 
For each $i\in\{1,\ldots,(4N+2)b\}$, let $\tilde{f}_i$ be the boundary trace of a harmonic map $\Phi_i\colon\overline{\Omega}_{c_0}\to\R^3$ of the form \eqref{eqn:mfs-field}, with common source points, chosen to approximate $G\varphi_i$. 
As in Proposition~\ref{prop:DtN-error}, we define the Fourier polynomial
\begin{align}
\label{eqn:Fourier-Fourier-DtN-approximation}
\tilde{h}_i
&\vcentcolon=\DN_{\mathrm{dec}}(G\varphi_i)
+\Pi_{\widetilde N}(\DN_{c_0}-\DN_{\mathrm{dec}})\tilde{f}_i
\end{align}
approximating $\DN_{c_0}(\Gamma\varphi_i)$.
For $i\leq j$, we approximate \eqref{eqn:low-mode-Fourier-entries} by 
\begin{align}
\label{eqn:Fourier-Fourier-approximation}
\tilde{a}_{ij}
&\vcentcolon=\int_\Xi
\biggl(
\sk{\tilde{h}_i,\Gamma\varphi_j}
-\sk{P,\tilde{h}}\sk{\varphi_i,\varphi_j}
\biggr)\,d\theta.
\end{align}
For $i<j$, we set $\tilde{a}_{ji}=\tilde{a}_{ij}$.
Since $\tilde{h}_i$ and $\varphi_j$ are Fourier polynomials, \eqref{eqn:Fourier-Fourier-approximation} can be computed by enclosing finitely many Fourier coefficients of $\Gamma$ and then using finite Fourier arithmetic.
By the pointwise orthonormality of $\Gamma$ and the Cauchy--Schwarz inequality, we obtain for $i \leq j$
\begin{align}
\label{eqn:Fourier-Fourier-error}
\abs[\Big]{\sk[\big]{\mathscr{A}\varphi_i,\varphi_j}_{\mathcal{H}}-\tilde{a}_{ij}}
&\leq
\nm[\big]{\DN_{c_0}(\Gamma\varphi_i)-\tilde{h}_i}_{L^2(\Xi)}
\\\notag
&\hphantom{{}={}}+\nm[\big]{\sk[\big]{\varphi_i,\varphi_j}}_{L^2(\Xi)}
\biggl(
\nm[\big]{\DN_{c_0}\gamma_0-\tilde{h}}_{L^2(\Xi)}
+\nm{\tilde{h}}_{L^\infty(\Xi)}
\nm{\gamma_0-P}_{L^2(\Xi)}
\biggr).
\end{align}
It remains to bound the Dirichlet-to-Neumann approximation errors in \eqref{eqn:Fourier-Fourier-error}.
Applying Proposition~\ref{prop:DtN-error} with $f=\Gamma\varphi_i$, $\tilde{f}=\tilde{f}_i$ and $\varphi=G\varphi_i$ in the notation therein, we obtain
\begin{align}
&\nm[\big]{\DN_{c_0}(\Gamma\varphi_i)-\tilde{h}_i}_{L^2(\Xi)}
\notag\\\label{eqn:Fourier-Fourier-DtN-error}
&\leq\nm[\big]{\DN_{c_0}(\Gamma\varphi_i)-\tilde{h}_i}_{H^1(\Xi)}
\\\notag
&\leq
\bigl(1+C_1(\omega)\bigr)
\nm{(\Gamma-G)\varphi_i}_{H^2(\Xi)}
+
C_1(\omega)\nm{G\varphi_i-\tilde{f}_i}_{L^2(\Xi)}
+\nm[\big]{(\mathrm{Id}-\Pi_{\widetilde N})(\DN_{c_0}-\DN_{\mathrm{dec}})\tilde{f}_i}_{H^1(\Xi)}.
\end{align}
By the product rule and the definition of the Sobolev norms in \eqref{eqn:Sobolev-norm}, we have
\begin{align}
\label{eqn:Fourier-Fourier-frame-error}
\nm[\big]{(\Gamma-G)\varphi_i}_{H^2(\Xi)}
&\leq
\Bigl(
\nm{\varphi_i}_{L^\infty(\Xi)}
+\sqrt{2}\nm{\varphi_i'}_{L^\infty(\Xi)}
+\nm{\varphi_i''}_{L^\infty(\Xi)}
\Bigr)
\nm{\Gamma-G}_{H^2(\Xi)}.
\end{align}
The norms of $\varphi_i$ appearing in \eqref{eqn:Fourier-Fourier-frame-error} are explicitly determined by the orthonormal Fourier basis.
The term $\nm{G\varphi_i-\tilde{f}_i}_{L^2(\Xi)}$ in \eqref{eqn:Fourier-Fourier-DtN-error} is bounded by the second estimate of Proposition~\ref{prop:mfs-truncation-bounds}.
The remaining term is bounded by the first estimate of Proposition~\ref{prop:mfs-truncation-bounds}.

\emph{Fourier--configuration entries:}
For $i\in\{1,\ldots,(4N+2)b\}$ and $j\in\{1,\ldots,3b-6\}$, let $e_j$ denote the $j$-th standard basis vector of $\R^{3b-6}$, identified with $(0,e_j,0)\in\mathcal{H}_{\leq N}$.
Let $\tilde{\alpha}_j,\tilde{\beta}_j$ be the normalised boundary velocities from Lemma~\ref{lem:DtN-derivative} corresponding to ${\partial c_\mu}/{\partial\mu_j}\vert_{\mu=0}$.
By Proposition~\ref{prop:linearisation} and Lemma~\ref{lem:DtN-derivative}, we have
\begin{align}
\label{eqn:low-mode-Fourier-configuration-entries}
\sk{\mathscr{A}\varphi_i,e_j}_{\mathcal{H}}
&=\int_\Xi
\Bigl(
\sk[\big]{(\Gamma\varphi_i)',\;\tilde{\alpha}_j\gamma_0'-\tilde{\beta}_j\DN_{c_0}\gamma_0}
-\sk[\big]{\DN_{c_0}(\Gamma\varphi_i),\;
\tilde{\alpha}_j\DN_{c_0}\gamma_0+\tilde{\beta}_j\gamma_0'}
\Bigr)\,d\theta.
\end{align}
Extending \eqref{eqn:Fourier-Fourier-approximation}, we approximate \eqref{eqn:low-mode-Fourier-configuration-entries} by
\begin{align}
\label{eqn:Fourier-configuration-approximation}
\tilde{a}_{i,j+(4N+2)b}
&\vcentcolon=\int_\Xi
\Bigl(
\sk[\big]{(\Gamma\varphi_i)',\;\tilde{\alpha}_jP'-\tilde{\beta}_j\tilde{h}}
-\sk[\big]{\tilde{h}_i,\;\tilde{\alpha}_j\tilde{h}+\tilde{\beta}_jP'}
\Bigr)\,d\theta.
\end{align}  
By \eqref{eqn:low-mode-Fourier-configuration-entries}, \eqref{eqn:Fourier-configuration-approximation} and the Cauchy--Schwarz inequality \eqref{eqn:20261106b}, we obtain
\begin{align}
\notag
&\abs[\Big]{\sk{\mathscr{A}\varphi_i,e_j}_{\mathcal{H}}-\tilde{a}_{i,j+(4N+2)b}}^2
\\\label{eqn:Fourier-configuration-error}
&\leq
\nm[\Big]{\tilde{\alpha}_j^2+\tilde{\beta}_j^2}_{L^\infty(\Xi)}
\biggl(
\nm[\big]{((\Gamma-G)\varphi_i)'}_{L^2(\Xi)}^2
+\nm[\big]{\DN_{c_0}(\Gamma\varphi_i)-\tilde{h}_i}_{L^2(\Xi)}^2
+\nm{\gamma_0'-P'}_{L^2(\Xi)}^2
\\\notag
&\hphantom{{}\leq{}}
+\nm[\big]{\DN_{c_0}\gamma_0-\tilde{h}}_{L^2(\Xi)}^2
\biggr)
\Bigl(
\nm{\gamma_0'}_{L^2(\Xi)}^2
+\nm{\DN_{c_0}\gamma_0}_{L^2(\Xi)}^2
+\nm[\big]{(\Gamma \varphi_i)'}_{L^2(\Xi)}^2
+\nm{\tilde{h}_i}_{L^2(\Xi)}^2
\Bigr).
\end{align}
By \eqref{eqn:Fourier-Fourier-frame-error}, we have
\begin{align*}
\nm[\big]{((\Gamma-G)\varphi_i)'}_{L^2(\Xi)}
&\leq\nm{(\Gamma-G)\varphi_i}_{H^2(\Xi)}
\\
&\leq\Bigl(
\nm{\varphi_i}_{L^\infty(\Xi)}
+\sqrt{2}\nm{\varphi_i'}_{L^\infty(\Xi)}
+\nm{\varphi_i''}_{L^\infty(\Xi)}
\Bigr)
\nm{\Gamma-G}_{H^2(\Xi)}.
\end{align*}
Moreover, the triangle inequality yields
\begin{align*}
\nm[\big]{(\Gamma\varphi_i)'}_{L^2(\Xi)}
&\leq\nm[\big]{(G\varphi_i)'}_{L^2(\Xi)}
+\nm[\big]{((\Gamma-G)\varphi_i)'}_{L^2(\Xi)},
\\
\nm{\gamma_0'}_{L^2(\Xi)}
&\leq\nm{P'}_{L^2(\Xi)}
+\nm{\gamma_0-P}_{H^2(\Xi)},
\\
\nm{\DN_{c_0}\gamma_0}_{L^2(\Xi)}
&\leq\nm{\tilde{h}}_{L^2(\Xi)}
+\nm[\big]{\DN_{c_0}\gamma_0-\tilde{h}}_{H^1(\Xi)}.
\end{align*}
The errors involving $\gamma_0-P$ and $\DN_{c_0}\gamma_0-\tilde{h}$ are bounded by Lemma~\ref{lem:normalization-trace-errors} and Proposition~\ref{prop:DtN-error}, respectively, as in the proof of Proposition~\ref{prop:residual-bound}.
The approximation error involving $\tilde{h}_i$ is bounded by \eqref{eqn:Fourier-Fourier-DtN-error}.
The norms of $(G\varphi_i)'$, $\tilde{h}_i$ and $P'$ are computed using finite Fourier arithmetic, while the remaining $L^\infty$ norms are bounded by interval arithmetic.
 
\emph{Fourier--rotation entries:}
For $i\in\{1,\ldots,(4N+2)b\}$ and $j\in\{1,2,3\}$, we identify $\xi_j$ with $(0,0,\xi_j)\in\mathcal{H}_{\leq N}$.
By Proposition~\ref{prop:linearisation} and \eqref{eqn:rotation-fields}, we have
\begin{align}
\label{eqn:low-mode-Fourier-rotation-entries}
\sk{\mathscr{A}\varphi_i,\xi_j}_{\mathcal{H}}
&=\int_\Xi\sk{\varphi_i,\rho_j}\,d\theta.
\end{align}
Denoting the columns of $\Gamma$ and $G$ by $\Gamma_1,\Gamma_2$ and $G_1,G_2$, respectively, we have
\begin{align*}
\rho_j
&=\Gamma^{\intercal}(\xi_j\times\gamma_0)
=\begin{pmatrix}
-\sk{\xi_j,\Gamma_2}\\
\hphantom{-}\sk{\xi_j,\Gamma_1}
\end{pmatrix}
\end{align*}
by \eqref{eqn:frame}.
Extending \eqref{eqn:Fourier-configuration-approximation}, we approximate \eqref{eqn:low-mode-Fourier-rotation-entries} by
\begin{align}
\label{eqn:Fourier-rotation-approximation}
\tilde{a}_{i,(4N+2)b+3b-6+j}
&\vcentcolon=\int_\Xi
\sk[\bigg]{\varphi_i,
\begin{pmatrix}
-\sk{\xi_j,\Gamma_2}\\
\hphantom{-}\sk{\xi_j,\Gamma_1}
\end{pmatrix}}
\,d\theta.
\end{align}
Since $\varphi_i$ is a Fourier basis element, the integral can be computed using finite Fourier arithmetic after enclosing one Fourier mode of $\Gamma_1$ and $\Gamma_2$.
Moreover, the Cauchy--Schwarz inequality and the orthonormality of $\varphi_i$ imply
\begin{align}
\label{eqn:Fourier-rotation-error}
\abs[\Big]{\sk{\mathscr{A}\varphi_i,\xi_j}_{\mathcal{H}}
-\tilde{a}_{i,(4N+2)b+3b-6+j}}
&\leq\nm{\Gamma-G}_{L^2(\Xi)}.
\end{align}

\emph{Configuration--configuration entries:}
For $i,j\in\{1,\ldots,3b-6\}$, identifying $e_i$ and $e_j$ with $(0,e_i,0)$ and $(0,e_j,0)$ in $\mathcal{H}_{\leq N}$, we have
\begin{align}
\label{eqn:low-mode-configuration-configuration-entries}
\sk{\mathscr{A}e_j,e_i}_{\mathcal{H}}
&=\int_\Xi
\biggl(
\frac{1}{2}\frac{\partial\tilde{\alpha}_i}{\partial\mu_j}
\bigl(\abs{\gamma_0'}^2-\abs{h_\mu}^2\bigr)
-\frac{\partial\tilde{\beta}_i}{\partial\mu_j}
\sk{\gamma_0',h_\mu}
\\\notag
&\hphantom{{}={}}
+\sk[\Big]{\tilde{\alpha}_i h_\mu+\tilde{\beta}_i\gamma_0',\;
\bigl(\tilde{\alpha}_j\gamma_0'-\tilde{\beta}_j h_\mu\bigr)'
+
\DN_{c_\mu}\bigl(\tilde{\alpha}_j h_\mu+\tilde{\beta}_j\gamma_0'\bigr)}
\biggr)\Big\vert_{\mu=0}\,d\theta
\end{align}
by Lemma~\ref{lem:configuration-linearisation} and \eqref{eqn:DtN-derivative-operator}. 
For each $j\in\{1,\ldots,3b-6\}$, let $\tilde{f}_j^{\mathrm{c}}$ be the boundary trace of a harmonic map of the form \eqref{eqn:mfs-field}, chosen to approximate $\tilde{\alpha}_j\tilde{h}+\tilde{\beta}_jP'$ at $\mu=0$.
In analogy with Proposition~\ref{prop:DtN-error}, we define the Fourier polynomial
\begin{align}
\label{eqn:configuration-DtN-approximation}
\tilde{h}_j^{\mathrm{c}}
&\vcentcolon=
\DN_{\mathrm{dec}}(\tilde{\alpha}_j\tilde{h}+\tilde{\beta}_jP')
+\Pi_{\widetilde N}(\DN_{c_0}-\DN_{\mathrm{dec}})\tilde{f}_j^{\mathrm{c}},
\end{align}
approximating 
$\DN_{c_0}(\tilde{\alpha}_j\DN_{c_0}\gamma_0+\tilde{\beta}_j\gamma_0')$, where the boundary velocities are evaluated at $\mu=0$. 
(The superscript $\mathrm{c}$ distinguishes \eqref{eqn:configuration-DtN-approximation} from \eqref{eqn:Fourier-Fourier-DtN-approximation} notationwise.)
For $i,j\in\{1,\ldots,3b-6\}$, we approximate \eqref{eqn:low-mode-configuration-configuration-entries} by
\begin{align}
\label{eqn:configuration-configuration-approximation}
\tilde{s}_{ij}
&\vcentcolon=\int_\Xi
\biggl(
\frac{1}{2}\frac{\partial\tilde{\alpha}_i}{\partial\mu_j}
\bigl(\abs{P'}^2-\abs{\tilde{h}}^2\bigr)
-\frac{\partial\tilde{\beta}_i}{\partial\mu_j}
\sk{P',\tilde{h}}
+\sk[\Big]{\tilde{\alpha}_i\tilde{h}+\tilde{\beta}_iP',\;
\bigl(\tilde{\alpha}_jP'-\tilde{\beta}_j\tilde{h}\bigr)'
+\tilde{h}_j^{\mathrm{c}}}
\biggr)\,d\theta,
\end{align}
where the boundary velocities and their configuration derivatives are evaluated at $\mu=0$.
We set
\begin{align}
\label{eqn:configuration-configuration-symmetrisation}
\tilde{a}_{(4N+2)b+i,(4N+2)b+j}
&\vcentcolon=\frac{1}{2}\bigl(\tilde{s}_{ij}+\tilde{s}_{ji}\bigr).
\end{align}
By construction, \eqref{eqn:configuration-configuration-symmetrisation} is symmetric in $i,j$. 
Since all functions in \eqref{eqn:configuration-configuration-approximation} are Fourier polynomials, its entries can be computed using finite Fourier arithmetic.

By Lemma~\ref{lem:decoupled-Fourier} and Corollary~\ref{cor:decoupled-explicit-bounds},  $\DN_{c_0}$ extends to a bounded operator $\DN_{c_0}\colon H^1(\Xi,\R^3)\to L^2(\Xi,\R^3)$ such that for all $f\in H^1(\Xi,\R^3)$
\begin{align*}
\nm{\DN_{c_0}f}_{L^2(\Xi)}
&\leq\bigl(1+C_0(\omega)\bigr)\nm{f}_{H^1(\Xi)}.
\end{align*} 
Applying this with $f=\tilde{\alpha}_jv+\tilde{\beta}_ju$ 
and using the product rule, we obtain for all $u,v\in H^1(\Xi,\R^3)$
\begin{align*}
&\nm[\Big]{(\tilde{\alpha}_ju-\tilde{\beta}_jv)'
+\DN_{c_0}(\tilde{\alpha}_jv+\tilde{\beta}_ju)}_{L^2(\Xi)}^2
\\
&\leq
\biggl(
\nm{(\tilde{\alpha}_ju-\tilde{\beta}_jv)'}_{L^2(\Xi)}
+\bigl(1+C_0(\omega)\bigr)
\nm{\tilde{\alpha}_jv+\tilde{\beta}_ju}_{H^1(\Xi)}
\biggr)^2
\\
&\leq
\bigl(2+C_0(\omega)\bigr)^2
\biggl(
\nm[\Big]{\sqrt{\tilde{\alpha}_j^2+\tilde{\beta}_j^2}}_{L^\infty(\Xi)}
+\nm[\Big]{\sqrt{(\tilde{\alpha}_j')^2+(\tilde{\beta}_j')^2}}_{L^\infty(\Xi)}
\biggr)^2
\Bigl(
\nm{u}_{H^1(\Xi)}^2+\nm{v}_{H^1(\Xi)}^2
\Bigr)
\\
&\leq
2\bigl(2+C_0(\omega)\bigr)^2
\biggl(
\nm[\big]{\tilde{\alpha}_j^2+\tilde{\beta}_j^2}_{L^\infty(\Xi)}
+\nm[\big]{(\tilde{\alpha}_j')^2+(\tilde{\beta}_j')^2}_{L^\infty(\Xi)}
\biggr)
\Bigl(
\nm{u}_{H^1(\Xi)}^2+\nm{v}_{H^1(\Xi)}^2
\Bigr).
\end{align*}
Therefore, 
\begin{align}
\label{eqn:configuration-configuration-quadratic-bound}
&\abs[\Bigg]{\int_\Xi
\biggl(
\frac{1}{2}\frac{\partial\tilde{\alpha}_i}{\partial\mu_j}
\bigl(\abs{u}^2-\abs{v}^2\bigr)
-\frac{\partial\tilde{\beta}_i}{\partial\mu_j}\sk{u,v}
+\sk[\Big]{\tilde{\alpha}_i v+\tilde{\beta}_i u,\;
(\tilde{\alpha}_j u-\tilde{\beta}_j v)'
+\DN_{c_0}(\tilde{\alpha}_j v+\tilde{\beta}_j u)}
\biggr)\,d\theta}
\\\notag
&\leq K_{ij}
\Bigl(
\nm{u}_{H^1(\Xi)}^2+\nm{v}_{H^1(\Xi)}^2
\Bigr),
\end{align}
where $K_{ij}\geq0$ is defined by
\begin{align}
\label{eqn:configuration-configuration-K}
K_{ij}^2
&\vcentcolon=
\frac{1}{2}
\nm[\Big]{
\Bigl(\frac{\partial\tilde{\alpha}_i}{\partial\mu_j}\Bigr)^2
+\Bigl(\frac{\partial\tilde{\beta}_i}{\partial\mu_j}\Bigr)^2
}_{L^\infty(\Xi)}
\\\notag&\hphantom{{}\vcentcolon={}}
+4\bigl(2+C_0(\omega)\bigr)^2
\nm[\big]{\tilde{\alpha}_i^2+\tilde{\beta}_i^2}_{L^\infty(\Xi)}
\biggl(
\nm[\big]{\tilde{\alpha}_j^2+\tilde{\beta}_j^2}_{L^\infty(\Xi)}
+\nm[\big]{(\tilde{\alpha}_j')^2+(\tilde{\beta}_j')^2}_{L^\infty(\Xi)}
\biggr),
\end{align}
and all boundary velocities and their configuration derivatives are evaluated at $\mu=0$.
The integral in \eqref{eqn:configuration-configuration-quadratic-bound} defines a quadratic form on $H^1(\Xi,\R^3)\times H^1(\Xi,\R^3)$, whose associated symmetric bilinear form has operator norm at most $K_{ij}$. 
By polarisation of \eqref{eqn:configuration-configuration-quadratic-bound} and the Cauchy--Schwarz inequality, we obtain
\begin{align}
\label{eqn:configuration-configuration-error}
\abs[\Big]{\sk{\mathscr{A}e_j,e_i}_{\mathcal{H}}-\tilde{s}_{ij}}^2
&\leq
\biggl(
K_{ij}^2
\Bigl(
\nm{\gamma_0'-P'}_{H^1(\Xi)}^2
+\nm{\DN_{c_0}\gamma_0-\tilde{h}}_{H^1(\Xi)}^2
\Bigr)
+\nm{\tilde{\alpha}_i\tilde{h}+\tilde{\beta}_iP'}_{L^2(\Xi)}^2
\biggr)
\\\notag
&\quad\cdot
\biggl(
\nm{\gamma_0'+P'}_{H^1(\Xi)}^2
+\nm{\DN_{c_0}\gamma_0+\tilde{h}}_{H^1(\Xi)}^2
+\nm[\Big]{\DN_{c_0}(\tilde{\alpha}_j\tilde{h}+\tilde{\beta}_jP')
-\tilde{h}_j^{\mathrm{c}}}_{L^2(\Xi)}^2
\biggr).
\end{align}
Applying Proposition~\ref{prop:DtN-error} with $f=\varphi=\tilde{\alpha}_j\tilde{h}+\tilde{\beta}_jP'$ and $\tilde{f}=\tilde{f}_j^{\mathrm{c}}$ in the notation therein, we obtain
\begin{align}
\notag
&\nm[\Big]{\DN_{c_0}(\tilde{\alpha}_j\tilde{h}+\tilde{\beta}_jP')-\tilde{h}_j^{\mathrm{c}}}_{L^2(\Xi)}
\\\label{eqn:configuration-DtN-error}
&\leq
C_1(\omega)
\nm[\big]{\tilde{\alpha}_j\tilde{h}+\tilde{\beta}_jP'-\tilde{f}_j^{\mathrm{c}}}_{L^2(\Xi)}
+\nm[\Big]{(\mathrm{Id}-\Pi_{\widetilde N})(\DN_{c_0}-\DN_{\mathrm{dec}})\tilde{f}_j^{\mathrm{c}}}_{H^1(\Xi)}.
\end{align}
The first term in \eqref{eqn:configuration-DtN-error} is bounded by the second estimate of Proposition~\ref{prop:mfs-truncation-bounds}, with truncation degree at least that of $\tilde{\alpha}_j\tilde{h}+\tilde{\beta}_jP'$.
The second term is bounded by the first estimate of Proposition~\ref{prop:mfs-truncation-bounds}, using the fixed degree $\widetilde N$.
By self-adjointness of $\mathscr{A}$ and the definition \eqref{eqn:configuration-configuration-symmetrisation}, we have
\begin{align}
\label{eqn:configuration-configuration-symmetric-error}
\abs[\Big]{\sk{\mathscr{A}e_j,e_i}_{\mathcal{H}}
-\tilde{a}_{(4N+2)b+i,(4N+2)b+j}}
&\leq\frac{1}{2}
\Bigl(
\abs[\big]{\sk{\mathscr{A}e_j,e_i}_{\mathcal{H}}-\tilde{s}_{ij}}
+\abs[\big]{\sk{\mathscr{A}e_i,e_j}_{\mathcal{H}}-\tilde{s}_{ji}}
\Bigr).
\end{align}
By Proposition~\ref{prop:linearisation}, the configuration--rotation and rotation--rotation entries of $\mathscr{A}$ vanish.
We set the corresponding entries $\tilde{a}_{ij}$ to zero and define the configuration--Fourier and rotation--Fourier entries by symmetry, so that $\tilde{a}_{ji}=\tilde{a}_{ij}$.
Thus, the matrix
\(
\bigl(\tilde{a}_{ij}\bigr)_{1\leq i,j\leq (4N+2)b+3b-3}
\)
is symmetric.

The estimates \eqref{eqn:Fourier-Fourier-error}, \eqref{eqn:Fourier-configuration-error}, \eqref{eqn:Fourier-rotation-error} and \eqref{eqn:configuration-configuration-symmetric-error} yield intervals containing each entry of the symmetric matrix representing $\mathscr{A}$ in the chosen orthonormal basis.
We form the symmetric matrix whose entries are the midpoints of these intervals and numerically compute an approximate inverse $J$, which we choose to be symmetric.
We then verify that
\begin{align*}
{\nm{\mathrm{Id}-J\mathscr{A}}}_2&\leq\delta_{\mathscr{A}}<1
\end{align*}
using interval arithmetic (see e.\,g.~\cite{OishiTanabeOgitaRump2007}). 
This implies that $\mathscr{A}$ is invertible and
\begin{align}
\label{eqn:linearisation-A-bound}
{\nm{\mathscr{A}^{-1}}}_2
&\leq\frac{\nm{J}_2}{1-\delta_{\mathscr{A}}},
&
{\nm{\mathscr{A}^{-1}-J}}_2
&\leq\frac{\delta_{\mathscr{A}}\nm{J}_2}{1-\delta_{\mathscr{A}}}.
\end{align}
For a symmetric matrix $M$, let $\lambda_{\max}(M)$ denote its largest eigenvalue.
Since $\mathscr{A}^{-1}$ and $J$ are symmetric, the second estimate in \eqref{eqn:linearisation-A-bound} implies
\begin{align}
\label{eqn:linearisation-A-eigenvalue-bound}
\abs[\big]{\lambda_{\max}(\mathscr{A}^{-1})-\lambda_{\max}(J)}
&\leq\frac{\delta_{\mathscr{A}}\nm{J}_2}{1-\delta_{\mathscr{A}}}.
\end{align} 

\paragraph{Off-diagonal bound.}
By Proposition~\ref{prop:linearisation},
the adjoint of $\mathscr{B}$ defined in \eqref{eqn:LABC} is given by 
\begin{align*}
\mathscr{B}^*
&=(\mathrm{Id}-\Pi_N)
\begin{pmatrix}
S_{vv}\Pi_N, &
S_{v\mu}, &
\rho
\end{pmatrix}.
\end{align*}
We estimate $\mathscr{B}$ by constructing Fourier polynomial approximations of the columns of $\mathscr{B}^*$ and bounding their errors. 

\emph{Fourier columns:}
For $i\in\{1,\ldots,(4N+2)b\}$, recalling the Fourier basis $\varphi_i$ from \eqref{eqn:Fourier-basis}, Proposition~\ref{prop:linearisation} and \eqref{eqn:LABC} imply
\begin{align*}
\mathscr{B}^*\varphi_i
&=(\mathrm{Id}-\Pi_N)
\Bigl(
\Gamma^{\intercal}\DN_{c_0}(\Gamma\varphi_i)
-\sk{\gamma_0,\DN_{c_0}\gamma_0}\varphi_i
\Bigr).
\end{align*}
We approximate this by the Fourier polynomial
$(\mathrm{Id}-\Pi_N)(G^{\intercal}\tilde{h}_i-\sk{P,\tilde{h}}\varphi_i)$, recalling $\tilde{h}$ and $\tilde{h}_i$ from \eqref{eqn:tilde_h} and \eqref{eqn:Fourier-Fourier-DtN-approximation}.
By the pointwise orthonormality of $\Gamma$, we obtain
\begin{align}
&\nm[\Big]{\mathscr{B}^*\varphi_i
-(\mathrm{Id}-\Pi_N)\bigl(G^{\intercal}\tilde{h}_i-\sk{P,\tilde{h}}\varphi_i\bigr)}_{L^2(\Xi)}
\notag\\\label{eqn:B-Fourier-columns}
&\leq
\nm[\big]{\DN_{c_0}(\Gamma\varphi_i)-\tilde{h}_i}_{L^2(\Xi)}
+\nm{\tilde{h}_i}_{L^\infty(\Xi)}\nm{\Gamma-G}_{L^2(\Xi)}
\\\notag&\hphantom{{}\leq{}}
+
\nm{\varphi_i}_{L^\infty(\Xi)}
\Bigl(
\nm{\DN_{c_0}\gamma_0-\tilde{h}}_{L^2(\Xi)}
+\nm{\tilde{h}}_{L^\infty(\Xi)}
\nm{\gamma_0-P}_{L^2(\Xi)}
\Bigr).
\end{align}
Here we also used that $\mathrm{Id}-\Pi_N$ is an orthogonal projection.
The approximation errors are bounded by \eqref{eqn:Fourier-Fourier-DtN-error}, Lemma~\ref{lem:normalization-trace-errors} and Proposition~\ref{prop:DtN-error}, together with the certified approximation of $\Gamma$ by $G$.
The remaining norms involve Fourier polynomials and can be computed using finite Fourier arithmetic.

\emph{Configuration columns:}
For $j\in\{1,\ldots,3b-6\}$, Proposition~\ref{prop:linearisation}, definition \eqref{eqn:LABC} and equation \eqref{eqn:DtN-derivative-operator} imply
\begin{align*}
\mathscr{B}^*e_j
&=-(\mathrm{Id}-\Pi_N)\Gamma^{\intercal}
\Bigl(
(\tilde{\alpha}_j\gamma_0'-\tilde{\beta}_j\DN_{c_0}\gamma_0)'
+\DN_{c_0}(\tilde{\alpha}_j\DN_{c_0}\gamma_0+\tilde{\beta}_j\gamma_0')
\Bigr),
\end{align*}
where $e_j$ and the boundary velocities $\tilde{\alpha}_j,\tilde{\beta}_j$ are as in \eqref{eqn:low-mode-Fourier-configuration-entries}.
Recalling $\tilde{h}_j^{\mathrm{c}}$ from \eqref{eqn:configuration-DtN-approximation}, we approximate this column by the Fourier polynomial
\begin{align*}
-(\mathrm{Id}-\Pi_N)G^{\intercal}
\Bigl((\tilde{\alpha}_jP'-\tilde{\beta}_j\tilde{h})'
+\tilde{h}_j^{\mathrm{c}}\Bigr).
\end{align*}
Using $\gamma_0'=P'+(\gamma_0'-P')$ and $\DN_{c_0}\gamma_0=\tilde{h}+(\DN_{c_0}\gamma_0-\tilde{h})$, we obtain
\begin{align}
\label{eqn:B-configuration-difference}
&(\tilde{\alpha}_j\gamma_0'-\tilde{\beta}_j\DN_{c_0}\gamma_0)'
+\DN_{c_0}(\tilde{\alpha}_j\DN_{c_0}\gamma_0+\tilde{\beta}_j\gamma_0')
-(\tilde{\alpha}_jP'-\tilde{\beta}_j\tilde{h})'
-\tilde{h}_j^{\mathrm{c}}
\\\notag
&=\bigl(\tilde{\alpha}_j(\gamma_0'-P')-\tilde{\beta}_j(\DN_{c_0}\gamma_0-\tilde{h})\bigr)'
+\DN_{c_0}\bigl(\tilde{\alpha}_j(\DN_{c_0}\gamma_0-\tilde{h})+\tilde{\beta}_j(\gamma_0'-P')\bigr)
+\DN_{c_0}(\tilde{\alpha}_j\tilde{h}+\tilde{\beta}_jP')-\tilde{h}_j^{\mathrm{c}}.
\end{align}
By \eqref{eqn:B-configuration-difference}, the product-rule estimate used in proving \eqref{eqn:configuration-configuration-quadratic-bound}, and the fact that $\mathrm{Id}-\Pi_N$ is an orthogonal projection, we obtain
\begin{align}
\label{eqn:B-moduli-columns}
&\nm[\Big]{\mathscr{B}^*e_j
+(\mathrm{Id}-\Pi_N)G^{\intercal}
\Bigl((\tilde{\alpha}_jP'-\tilde{\beta}_j\tilde{h})'
+\tilde{h}_j^{\mathrm{c}}\Bigr)}_{L^2(\Xi)}
\\\notag
&\leq
\nm{\Gamma-G}_{L^2(\Xi)}
\nm[\Big]{(\tilde{\alpha}_jP'-\tilde{\beta}_j\tilde{h})'
+\tilde{h}_j^{\mathrm{c}}}_{L^\infty(\Xi)}
+\nm[\Big]{\DN_{c_0}(\tilde{\alpha}_j\tilde{h}+\tilde{\beta}_jP')
-\tilde{h}_j^{\mathrm{c}}}_{L^2(\Xi)}
\\\notag&\hphantom{{}={}}
+\bigl(2+C_0(\omega)\bigr)
\biggl(
\nm[\Big]{\sqrt{\tilde{\alpha}_j^2+\tilde{\beta}_j^2}}_{L^\infty(\Xi)}
+\nm[\Big]{\sqrt{(\tilde{\alpha}_j')^2+(\tilde{\beta}_j')^2}}_{L^\infty(\Xi)}
\biggr)
\\\notag
&\hphantom{{}={}}
\cdot
\Bigl(
\nm[\big]{\gamma_0'-P'}_{H^1(\Xi)}^2
+\nm[\big]{\DN_{c_0}\gamma_0-\tilde{h}}_{H^1(\Xi)}^2
\Bigr)^{\frac{1}{2}}.
\end{align}
The norm of $\DN_{c_0}(\tilde{\alpha}_j\tilde{h}+\tilde{\beta}_jP')-\tilde{h}_j^{\mathrm{c}}$ in \eqref{eqn:B-moduli-columns} is bounded by \eqref{eqn:configuration-DtN-error}. 
The remaining approximation errors are bounded by Lemma~\ref{lem:normalization-trace-errors}, Proposition~\ref{prop:DtN-error} and the approximation of $\Gamma$ by $G$.

\emph{Rotation columns:}
For $j\in\{1,2,3\}$, recalling $\rho_j$ from \eqref{eqn:rotation-fields}, Proposition~\ref{prop:linearisation} and \eqref{eqn:LABC} imply
\begin{align*}
\mathscr{B}^*(0,0,e_j)
&=(\mathrm{Id}-\Pi_N)\rho_j
=(\mathrm{Id}-\Pi_N)
\begin{pmatrix}
-\sk{\xi_j,\Gamma_2}\\
\hphantom{-}\sk{\xi_j,\Gamma_1}
\end{pmatrix},
\end{align*}
where $e_j$ denotes the $j$th standard basis vector of $\R^3$.
Replacing $\Gamma$ by its Fourier polynomial approximation $G$, we obtain
\begin{align}
\label{eqn:B-rotation-columns}
\nm[\Big]{\mathscr{B}^*(0,0,e_j)-(\mathrm{Id}-\Pi_N)
\begin{pmatrix}
-\sk{\xi_j,G_2}\\
\hphantom{-}\sk{\xi_j,G_1}
\end{pmatrix}}_{L^2(\Xi)}
&\leq\abs{\xi_j}\nm{\Gamma-G}_{L^2(\Xi)},
\end{align}
since $\mathrm{Id}-\Pi_N$ is an orthogonal projection.

\emph{Combining the column bounds:}
We define a bounded operator $\tilde{\mathscr{B}}^*\colon\mathcal{H}_{\leq N}\to\mathcal{H}_{>N}$ by its action on the chosen orthonormal basis:
\begin{align*}
\tilde{\mathscr{B}}^*\varphi_i
&\vcentcolon=(\mathrm{Id}-\Pi_N)
\bigl(G^{\intercal}\tilde{h}_i-\sk{P,\tilde{h}}\varphi_i\bigr),
\\
\tilde{\mathscr{B}}^*e_j
&\vcentcolon=-(\mathrm{Id}-\Pi_N)G^{\intercal}
\Bigl((\tilde{\alpha}_jP'-\tilde{\beta}_j\tilde{h})'
+\tilde{h}_j^{\mathrm{c}}\Bigr),
\\
\tilde{\mathscr{B}}^*(0,0,e_j)
&\vcentcolon=(\mathrm{Id}-\Pi_N)
\begin{pmatrix}
-\sk{\xi_j,G_2}\\
\hphantom{-}\sk{\xi_j,G_1}
\end{pmatrix},
\end{align*}
where the indices range as in \eqref{eqn:B-Fourier-columns}, \eqref{eqn:B-moduli-columns} and \eqref{eqn:B-rotation-columns}.
Since these columns are Fourier polynomials, $\tilde{\mathscr{B}}^*$ is represented by a finite matrix with respect to the chosen orthonormal basis of $\mathcal{H}_{\leq N}$ and finitely many orthonormal Fourier modes in $\mathcal{H}_{>N}$. 
By the Cauchy--Schwarz inequality applied to the columns, we have
\begin{align*}
\nm{\mathscr{B}^*-\tilde{\mathscr{B}}^*}^2
&\leq\sum_{i=1}^{\mathclap{(4N+2)b}}
\nm{(\mathscr{B}^*-\tilde{\mathscr{B}}^*)\varphi_i}_{L^2(\Xi)}^2
+\sum_{j=1}^{3b-6}
\nm{(\mathscr{B}^*-\tilde{\mathscr{B}}^*)e_j}_{L^2(\Xi)}^2
+\sum_{j=1}^3
\nm{(\mathscr{B}^*-\tilde{\mathscr{B}}^*)(0,0,e_j)}_{L^2(\Xi)}^2,
\end{align*}
which is bounded by \eqref{eqn:B-Fourier-columns}, \eqref{eqn:B-moduli-columns} and \eqref{eqn:B-rotation-columns}. 
Moreover, 
\begin{align}
\label{eqn:linearisation-B-bound}
\nm{\mathscr{B}}
&\leq\nm{\tilde{\mathscr{B}}^*}_2
+\nm{\mathscr{B}^*-\tilde{\mathscr{B}}^*}
\leq\sqrt{\nm{\tilde{\mathscr{B}}^*}_1\nm{\tilde{\mathscr{B}}^*}_\infty}
+\nm{\mathscr{B}^*-\tilde{\mathscr{B}}^*},
\end{align}
where the matrix norms of $\tilde{\mathscr{B}}^*$ can be computed using finite Fourier arithmetic.

\paragraph{High-mode coercivity.} 
Let $v\in H^1(\Xi,\R^2)$ such that $\Pi_Nv=0$. 
Then, by Lemma~\ref{lem:decoupled-Fourier},  
\begin{align*}
\widehat{\bigl(\DN_{\mathrm{dec}}(\Gamma v)-\Gamma\DN_{\mathrm{dec}}v\bigr)}_k(m)
&=\sum_{n\in\Z}\bigl(\abs{m}-\abs{m-n}\bigr)\hat\Gamma_k(n)\hat v_k(m-n)
\end{align*}
for $k\in\{1,\ldots,b\}$ and $m\in\Z$.
Since $\abs{\abs{m}-\abs{m-n}}\leq\abs{n}$ and 
\begin{align*}
\sum_{n\neq0}\abs{n}\abs{\hat\Gamma_k(n)}
&\leq
\Bigl(\sum_{n\neq0}\frac{1}{n^2}\Bigr)^{\frac{1}{2}}
\Bigl(\sum_{n\neq0}n^4\abs{\hat\Gamma_k(n)}^2\Bigr)^{\frac{1}{2}}
=\sqrt{\frac{\pi}{6}}\nm{\Gamma_k''}_{L^2}
\end{align*}
we obtain 
\begin{align}
\label{eqn:C-commutator}
\nm[\big]{\DN_{\mathrm{dec}}(\Gamma v)-\Gamma\DN_{\mathrm{dec}}v}_{L^2(\Xi)}
&\leq
\sqrt{\frac{\pi}{6}}\,
\max_{k\in\{1,\ldots,b\}}\nm{\Gamma_k''}_{L^2}
\nm{v}_{L^2(\Xi)}.
\end{align}
Recalling $\mathscr{C}$ from \eqref{eqn:LABC}, Proposition~\ref{prop:linearisation} and the pointwise identity $\Gamma^{\intercal}\Gamma=\mathrm{Id}$ imply
\begin{align}
\notag
\int_\Xi\sk{\mathscr{C}v,v}\,d\theta
&=\int_\Xi\sk{\DN_{\mathrm{dec}}v,v}\,d\theta
+\int_\Xi\sk{\DN_{\mathrm{dec}}(\Gamma v)-\Gamma\DN_{\mathrm{dec}}v,\Gamma v}\,d\theta
\\\notag
&\quad+\int_\Xi\sk{(\DN_{c_0}-\DN_{\mathrm{dec}})(\Gamma v),\Gamma v}\,d\theta
-\int_\Xi\sk{\gamma_0,\DN_{c_0}\gamma_0}\abs{v}^2\,d\theta
\\ 
&\geq
\Bigl(
N+1
-\sqrt{\frac{\pi}{6}}\max_{1\leq k\leq b}\nm{\Gamma_k''}_{L^2}
-C_0(\omega)
-\nm{\DN_{c_0}\gamma_0}_{L^\infty(\Xi)}
\Bigr)\nm{v}_{L^2(\Xi)}^2
\label{eqn:linearisation-C-bound}
\end{align}
by \eqref{eqn:C-commutator}, Lemma~\ref{lem:decoupled-Fourier} and Corollary~\ref{cor:decoupled-explicit-bounds}. 
The norms of $\Gamma_k''$ are bounded by
\begin{align*}
\nm{\Gamma_k''}_{L^2}
&\leq\nm{G_k''}_{L^2}+\nm{\Gamma-G}_{H^2(\Xi)}.
\end{align*}
Using $\abs{\gamma_0}=1$ and the Sobolev embedding $H^1(\Xi)\subset L^\infty(\Xi)$ with constant $c_\infty=\sqrt{\frac{1}{2}\coth(\pi)}$ as in \eqref{eqn:Sobolev-constant}, 
we may bound 
\[
\nm{\DN_{c_0}\gamma_0}_{L^\infty(\Xi)}
\leq
\nm{\tilde{h}}_{L^\infty(\Xi)}
+c_\infty\nm{\DN_{c_0}\gamma_0-\tilde{h}}_{H^1(\Xi)}.
\]
The norm of $\DN_{c_0}\gamma_0-\tilde{h}$ is bounded by Proposition~\ref{prop:DtN-error}, as in the proof of Proposition~\ref{prop:residual-bound}. 
The remaining terms are bounded using finite Fourier arithmetic and interval arithmetic.
Using these bounds, we verify that the coefficient of $\nm{v}_{L^2(\Xi)}^2$ in \eqref{eqn:linearisation-C-bound} is positive.
Since $\mathscr{C}$ is self-adjoint, it follows that $\mathscr{C}$ is invertible and $\nm{\mathscr{C}^{-1}}$ is bounded by the reciprocal of this coefficient. 
Here and in the following, unlabelled norms are understood as operator norms. 

\begin{proposition}[Hilbert space inverse bound]
\label{prop:L2-inverse-bound}
The self-adjoint operator $\mathscr{L}\colon Y\subset\mathcal{H}\to\mathcal{H}$ from Lemma~\ref{lem:linearisation-self-adjoint} is invertible and 
$\nm{\mathscr{L}^{-1}}_{\mathcal{H}\to\mathcal{H}}\leq q_0$, 
where $q_0$ is recorded in Appendix~\ref{section-numerical-values}.
\end{proposition}

\begin{proof}
By \eqref{eqn:linearisation-C-bound}, the operator $\mathscr{C}$ is strictly positive and has a bounded inverse.
Consequently,
\begin{align*}
0\leq\mathscr{B}\mathscr{C}^{-1}\mathscr{B}^*
\leq\nm{\mathscr{B}}^2\nm{\mathscr{C}^{-1}}\,\mathrm{Id}
\end{align*}
as quadratic forms on $\mathcal{H}_{\leq N}$.
Subtracting $\mathscr{B}\mathscr{C}^{-1}\mathscr{B}^*$ from $\mathscr{A}$ does not increase its negative eigenvalues, while its positive eigenvalues decrease by at most $\nm{\mathscr{B}}^2\nm{\mathscr{C}^{-1}}$.
Using \eqref{eqn:linearisation-A-bound}, \eqref{eqn:linearisation-A-eigenvalue-bound}, \eqref{eqn:linearisation-B-bound} and \eqref{eqn:linearisation-C-bound}, we verify numerically that the largest eigenvalue of $\mathscr{A}^{-1}$ satisfies 
\begin{align*}
0&<\lambda_{\max}(\mathscr{A}^{-1}),
&
\lambda_{\max}(\mathscr{A}^{-1})
\nm{\mathscr{B}}^2\nm{\mathscr{C}^{-1}}&<1.
\end{align*}
Since the smallest positive eigenvalue of $\mathscr{A}$ is $\lambda_{\max}(\mathscr{A}^{-1})^{-1}$, these inequalities imply that the Schur complement $\mathscr{A}-\mathscr{B}\mathscr{C}^{-1}\mathscr{B}^*$ is invertible and
\begin{align*}
\nm[\Big]{(\mathscr{A}-\mathscr{B}\mathscr{C}^{-1}\mathscr{B}^*)^{-1}}
&\leq
\max\biggl\{
\nm{\mathscr{A}^{-1}},~
\frac{\lambda_{\max}(\mathscr{A}^{-1})}
{1-\lambda_{\max}(\mathscr{A}^{-1})
\nm{\mathscr{B}}^2\nm{\mathscr{C}^{-1}}}
\biggr\}.
\end{align*}
The inverse of $\mathscr{L}$ is therefore given by
\begin{align*}
\mathscr{L}^{-1}
&=
\begin{pmatrix}
\mathrm{Id}\\
-\mathscr{C}^{-1}\mathscr{B}^*
\end{pmatrix}
\bigl(\mathscr{A}-\mathscr{B}\mathscr{C}^{-1}\mathscr{B}^*\bigr)^{-1}
\begin{pmatrix}
\mathrm{Id},&-\mathscr{B}\mathscr{C}^{-1}
\end{pmatrix}
+
\begin{pmatrix}
0&0\\
0&\mathscr{C}^{-1}
\end{pmatrix}.
\end{align*}
Since $\mathscr{C}^{-1}$ is self-adjoint, we have
\begin{align*}
\nm[\Bigg]{
\begin{pmatrix}
\mathrm{Id}\\
-\mathscr{C}^{-1}\mathscr{B}^*
\end{pmatrix}}^2
&=
\nm[\big]{\mathrm{Id}+\mathscr{B}\mathscr{C}^{-2}\mathscr{B}^*}
\leq1+\nm{\mathscr{B}}^2\nm[\big]{\mathscr{C}^{-1}}^2.
\end{align*}
Consequently,
\begin{align*}
\nm{\mathscr{L}^{-1}}_{\mathcal{H}\to\mathcal{H}}
&\leq
\nm{\mathscr{C}^{-1}}
+\Bigl(1+\nm{\mathscr{B}}^2\nm{\mathscr{C}^{-1}}^2\Bigr)
\nm[\Big]{\bigl(\mathscr{A}-\mathscr{B}\mathscr{C}^{-1}\mathscr{B}^*\bigr)^{-1}} 
\leq q_0.
\end{align*}
The bound $q_0$ is computed from \eqref{eqn:linearisation-A-bound}, \eqref{eqn:linearisation-A-eigenvalue-bound}, \eqref{eqn:linearisation-B-bound} and \eqref{eqn:linearisation-C-bound}.
\end{proof}

For the following estimate, we recall from Proposition~\ref{prop:existence} respectively \eqref{eqn:orthogonal_decomposition} the Banach spaces 
\begin{align*}
X&=H^2(\Xi,\R^2)\times\R^{3b-6}\times\R^3,
&
Y&=H^1(\Xi,\R^2)\times\R^{3b-6}\times\R^3, 
&
\mathcal{H}&=L^2(\Xi,\R^2)\times\R^{3b-6}\times\R^3. 
\end{align*}

\begin{lemma}[A priori estimate]
\label{lem:a-priori-estimate}
Let $\mathscr{L}=DF\vert_0\colon X\to Y$ be the linearised operator from Proposition~\ref{prop:linearisation}.  
Then, for every $x\in X$, we have
$\nm{x}_X\leq2\nm{\mathscr{L}x}_Y+\kappa\nm{x}_{\mathcal{H}}$, 
where $\kappa\geq0$ is recorded in Appendix~\ref{section-numerical-values}.
\end{lemma}

\begin{proof}
For $x=(v,\mu,\lambda)\in X$, Proposition~\ref{prop:linearisation} implies
\begin{align*}
\nm{\DN_{\mathrm{dec}}v}_{H^1(\Xi)}
&\leq\nm{\mathscr{L}x}_Y
+\nm{(S_{vv}-\DN_{\mathrm{dec}})v}_{H^1(\Xi)}
+\nm{S_{v\mu}\mu+\rho\lambda}_{H^1(\Xi)}.
\end{align*}
By the definition of the Sobolev norms in \eqref{eqn:Sobolev-norm}, we have
\begin{align*}
\nm{v}_{H^2(\Xi)}^2
&=\nm{\DN_{\mathrm{dec}}v}_{H^1(\Xi)}^2
+\nm{v}_{H^1(\Xi)}^2,
&
\nm{v}_{H^1(\Xi)}^2
&\leq\nm{v}_{L^2(\Xi)}\nm{v}_{H^2(\Xi)}.
\end{align*}
Consequently, the triangle inequality and Young's inequality yield
\begin{align}
\notag
{\nm{x}}_X
&\leq{\nm{v}}_{H^2(\Xi)}
+\sqrt{\abs{\mu}^2+\abs{\lambda}^2}
\\\notag
&\leq{\nm{\mathscr{L}x}}_Y
+\Bigl(1+\nm{S_{vv}-\DN_{\mathrm{dec}}}_{H^1\to H^1}\Bigr)
\nm{v}_{H^1(\Xi)}
\\&\hphantom{{}\leq{}}\notag
+\Bigl(1+
\nm[\big]{\bigl(S_{v\mu},\;\rho\bigr)}_{\R^{3b-3}\to H^1}
\Bigr)
\sqrt{\abs{\mu}^2+\abs{\lambda}^2}
\\\notag
&\leq\nm{\mathscr{L}x}_Y
+\frac{1}{2}\nm{v}_{H^2(\Xi)}
+\frac{1}{2}
\Bigl(1+\nm{S_{vv}-\DN_{\mathrm{dec}}}_{H^1\to H^1}\Bigr)^2
\nm{v}_{L^2(\Xi)}
\\&\hphantom{{}\leq{}}\notag
+\Bigl(1+
\nm[\big]{\bigl(S_{v\mu},\;\rho\bigr)}_{\R^{3b-3}\to H^1}
\Bigr)
\sqrt{\abs{\mu}^2+\abs{\lambda}^2}
\\\notag
&\leq\nm{\mathscr{L}x}_Y
+\frac{1}{2}\nm{x}_X
+\frac{\kappa}{2}\nm{x}_{\mathcal{H}},
\shortintertext{provided that}
\label{eqn:a-priori-kappa}
\kappa
&\geq
\Bigl(1+\nm{S_{vv}-\DN_{\mathrm{dec}}}_{H^1\to H^1}\Bigr)^2
+2\Bigl(1+
\nm[\big]{\bigl(S_{v\mu},\;\rho\bigr)}_{\R^{3b-3}\to H^1}
\Bigr).
\end{align}
Rearranging proves the claimed inequality. 
It remains to control the two operator norms in \eqref{eqn:a-priori-kappa}. 

\emph{Bound on $\nm{S_{vv}-\DN_{\mathrm{dec}}}$.} 
By Proposition~\ref{prop:linearisation} and the identity $\Gamma^{\intercal}\Gamma=\mathrm{Id}$, we have
\begin{align*}
(S_{vv}-\DN_{\mathrm{dec}})v
&=\Gamma^{\intercal}
\bigl(\DN_{\mathrm{dec}}(\Gamma v)-\Gamma\DN_{\mathrm{dec}}v\bigr)
+\Gamma^{\intercal}(\DN_{c_0}-\DN_{\mathrm{dec}})(\Gamma v)
-\sk{\gamma_0,\DN_{c_0}\gamma_0}\,v.
\end{align*}
The proof of \eqref{eqn:C-commutator} does not use the assumption $\Pi_Nv=0$, so the commutator estimate also holds for arbitrary $v\in H^1(\Xi,\R^2)$.
Moreover, for $v\in H^2(\Xi,\R^2)$, Lemma~\ref{lem:decoupled-Fourier} implies
\begin{align*}
\bigl(\DN_{\mathrm{dec}}(\Gamma v)-\Gamma\DN_{\mathrm{dec}}v\bigr)'
&=\DN_{\mathrm{dec}}(\Gamma v')-\Gamma\DN_{\mathrm{dec}}v'
+\DN_{\mathrm{dec}}(\Gamma'v)-\Gamma'\DN_{\mathrm{dec}}v.
\end{align*}
Using $\nm{\DN_{\mathrm{dec}}f}_{L^2}=\nm{f'}_{L^2}$ from Lemma~\ref{lem:decoupled-Fourier} and the product rule, we obtain
\begin{align*}
\nm{\DN_{\mathrm{dec}}(\Gamma'v)-\Gamma'\DN_{\mathrm{dec}}v}_{L^2}
&\leq
\nm{\Gamma''}_{L^2}\nm{v}_{L^\infty}
+2\nm{\Gamma'}_{L^\infty}\nm{v'}_{L^2}
\\
&\leq
\bigl(c_\infty\nm{\Gamma''}_{L^2}
+2\nm{\Gamma'}_{L^\infty}\bigr)\nm{v}_{H^1},
\end{align*}
where $c_\infty=\sqrt{\frac{1}{2}\coth(\pi)}$ as in \eqref{eqn:Sobolev-constant}.
Consequently, applying \eqref{eqn:C-commutator} to $v$ and $v'$, we obtain
\begin{align}
\label{eqn:a-priori-commutator}
\nm[\big]{\DN_{\mathrm{dec}}(\Gamma v)-\Gamma\DN_{\mathrm{dec}}v}_{H^1}
&\leq
\biggl(
\sqrt{\frac{\pi}{6}}\max_{k\in\{1,\ldots,b\}}\nm{\Gamma_k''}_{L^2}
+2\nm{\Gamma'}_{L^\infty}
+c_\infty\nm{\Gamma''}_{L^2}
\biggr)\nm{v}_{H^1}.
\end{align}
By density, \eqref{eqn:a-priori-commutator} extends to every $v\in H^1(\Xi,\R^2)$.

The product rule and the pointwise orthonormality of $\Gamma$ imply
\begin{align}\label{eqn:20261110a}
\nm{\Gamma^{\intercal}f}_{H^1(\Xi)}
&\leq\nm{f}_{H^1(\Xi)}
+\nm{\Gamma'}_{L^\infty(\Xi)}\nm{f}_{L^2(\Xi)}
\end{align}
for every $f\in H^1(\Xi,\R^3)$.
Moreover, since $\abs{\gamma_0}=1$, 
\begin{align}\notag
\nm[\big]{\sk{\gamma_0,\DN_{c_0}\gamma_0}v}_{H^1(\Xi)}
&\leq
\nm{\DN_{c_0}\gamma_0}_{L^\infty(\Xi)}\nm{v}_{H^1(\Xi)}
+\nm[\big]{\sk{\gamma_0,\DN_{c_0}\gamma_0}'}_{L^2(\Xi)}
\nm{v}_{L^\infty(\Xi)}
\\\label{eqn:20261110b}
&\leq
c_\infty
\bigl(2+\nm{\gamma_0'}_{L^\infty(\Xi)}\bigr)
\nm{\DN_{c_0}\gamma_0}_{H^1(\Xi)}
\nm{v}_{H^1(\Xi)}.
\end{align}
Combining \eqref{eqn:20261110a}--\eqref{eqn:20261110b} with \eqref{eqn:C-commutator}, the proof of \eqref{eqn:a-priori-commutator} and Corollary~\ref{cor:decoupled-explicit-bounds}, we obtain
\begin{align}
\notag
\nm{S_{vv}-\DN_{\mathrm{dec}}}_{H^1\to H^1}
&\leq
\Bigl(1+\nm{\Gamma'}_{L^\infty(\Xi)}\Bigr)
\sqrt{\frac{\pi}{6}}
\max_{k\in\{1,\ldots,b\}}\nm{\Gamma_k''}_{L^2}
\\\label{eqn:a-priori-Svv-bound}
&\hphantom{{}\leq{}}
+2\nm{\Gamma'}_{L^\infty(\Xi)}
+\nm{\Gamma''}_{L^\infty(\Xi)}
+C_1(\omega)
+C_1(\omega)\nm{\Gamma'}_{L^\infty(\Xi)}
\\\notag
&\hphantom{{}\leq{}}
+c_\infty
\bigl(2+\nm{\gamma_0'}_{L^\infty(\Xi)}\bigr)
\nm{\DN_{c_0}\gamma_0}_{H^1(\Xi)}.
\end{align}
Here we used $\nm{\Gamma''v}_{L^2}\leq\nm{\Gamma''}_{L^\infty}\nm{v}_{L^2}$ which is obtained as in the proof of \eqref{eqn:a-priori-commutator}, and $C_0(\omega)\leq C_1(\omega)$.
The term $\nm{\Gamma''}_{L^\infty(\Xi)}$ is bounded by interval evaluation.
By density, $S_{vv}-\DN_{\mathrm{dec}}$ extends to a bounded operator on $H^1(\Xi,\R^2)$ satisfying \eqref{eqn:a-priori-Svv-bound}.
The approximation of $\Gamma$ by $G$ and the Sobolev embedding from \eqref{eqn:Sobolev-constant} imply 
\begin{align*}
\nm{\Gamma'}_{L^\infty(\Xi)}
&\leq\nm{G'}_{L^\infty(\Xi)}
+c_\infty\nm{\Gamma-G}_{H^2(\Xi)},
&
\nm{\Gamma''}_{L^2(\Xi)}
&\leq\nm{G''}_{L^2(\Xi)}
+\nm{\Gamma-G}_{H^2(\Xi)}.
\end{align*}
Similarly, the Sobolev embedding and Lemma~\ref{lem:normalization-trace-errors} imply
\begin{align}\notag
\nm{\gamma_0'}_{L^\infty(\Xi)}
&\leq\nm{P'}_{L^\infty(\Xi)}
+c_\infty\nm{\gamma_0-P}_{H^2(\Xi)}.
\shortintertext{Finally,} 
\label{eqn:20261110c}
\nm{\DN_{c_0}\gamma_0}_{H^1(\Xi)}
&\leq\nm{\tilde{h}}_{H^1(\Xi)}
+\nm{\DN_{c_0}\gamma_0-\tilde{h}}_{H^1(\Xi)},
\end{align}
where the last term in \eqref{eqn:20261110c} is bounded by Proposition~\ref{prop:DtN-error}, as in the proof of Proposition~\ref{prop:residual-bound}.
The remaining norms involve Fourier polynomials and are bounded using finite Fourier and interval arithmetic. 

\emph{Bound on $\nm{(S_{v\mu},\;\rho)}$.}
By Proposition~\ref{prop:linearisation} and the Cauchy--Schwarz inequality, we have
\begin{align}
\label{eqn:a-priori-second-operator-norm}
\nm[\big]{(S_{v\mu},\;\rho)}_{\R^{3b-3}\to H^1}^2
&\leq
\sum_{j=1}^{3b-6}\nm{S_{v\mu}e_j}_{H^1(\Xi)}^2
+\sum_{j=1}^3\nm{\rho_j}_{H^1(\Xi)}^2,
\end{align}
where $e_j$ denotes the $j$th standard basis vector of $\R^{3b-6}$, as in \eqref{eqn:low-mode-Fourier-configuration-entries}. 
With $\xi_1,\xi_2,\xi_3$ being the standard basis vectors of $\R^3$, \eqref{eqn:rotation-fields} and \eqref{eqn:frame} imply
\begin{align}
\label{eqn:a-priori-rotation-columns}
\sum_{j=1}^3\nm{\rho_j}_{H^1(\Xi)}^2
&=\nm{\Gamma}_{H^1(\Xi)}^2
=4\pi b+\nm{\Gamma'}_{L^2(\Xi)}^2.
\end{align}
The last norm in \eqref{eqn:a-priori-rotation-columns} is bounded by
\begin{align*}
\nm{\Gamma'}_{L^2(\Xi)}
&\leq\nm{G'}_{L^2(\Xi)}
+\nm{\Gamma-G}_{H^2(\Xi)}.
\end{align*}
For the configuration columns, we recall $h_\mu=\DN_{c_\mu}\gamma_0$ from Lemma~\ref{lem:configuration-linearisation}.
Applying \eqref{eqn:DtN-derivative-remainder} with $\mu=0$ and the normalised boundary velocities $\tilde{\alpha}_j,\tilde{\beta}_j$, we obtain
\begin{align}
\label{eqn:a-priori-configuration-derivative}
\frac{\partial h_\mu}{\partial\mu_j}\Big\vert_{\mu=0}
&=\bigl(\tilde{\beta}_j(\DN_{c_0}-\DN_{\mathrm{dec}})\gamma_0\bigr)'
-\DN_{\mathrm{dec}}\bigl(\tilde{\alpha}_j(\DN_{c_0}-\DN_{\mathrm{dec}})\gamma_0\bigr) 
-(\DN_{c_0}-\DN_{\mathrm{dec}})
\bigl(\tilde{\alpha}_j\DN_{c_0}\gamma_0+\tilde{\beta}_j\gamma_0'\bigr).
\end{align}
By Lemma~\ref{lem:decoupled-Fourier}, we have
\begin{align*}
\nm{\DN_{\mathrm{dec}}f}_{H^1(\Xi)}
&\leq\nm{f}_{H^2(\Xi)},
&
\nm{f'}_{H^1(\Xi)}
&\leq\nm{f}_{H^2(\Xi)}.
\end{align*}
Since $\tilde{\alpha}_j$ and $\tilde{\beta}_j$ have Fourier degree at most one, \eqref{eqn:degree-one-product} and Corollary~\ref{cor:decoupled-explicit-bounds} imply
\begin{align*}
\nm[\bigg]{\frac{\partial h_\mu}{\partial\mu_j}\Big\vert_{\mu=0}}_{H^1(\Xi)}
&\leq
4C_2(\omega)
\Bigl(\nm{\tilde{\alpha}_j}_{L^\infty(\Xi)}
+\nm{\tilde{\beta}_j}_{L^\infty(\Xi)}\Bigr)
\nm{\gamma_0}_{L^2(\Xi)}
\\
&\hphantom{{}\leq{}}
+C_1(\omega)
\nm[\Big]{\sqrt{\tilde{\alpha}_j^2+\tilde{\beta}_j^2}}_{L^\infty(\Xi)}
\Bigl(
\nm{\DN_{c_0}\gamma_0}_{L^2(\Xi)}^2
+\nm{\gamma_0'}_{L^2(\Xi)}^2
\Bigr)^{\frac{1}{2}}.
\end{align*}
Using $\abs{\gamma_0}=1$ and \eqref{eqn:frame-multiplication-bound}, we conclude
\begin{align}
\label{eqn:a-priori-configuration-columns}
\nm{S_{v\mu}e_j}_{H^1(\Xi)}
&\leq
\bigl(1+\nm{\Gamma'}_{L^\infty(\Xi)}\bigr)
\biggl(
4C_2(\omega)\sqrt{2\pi b}
\Bigl(\nm{\tilde{\alpha}_j}_{L^\infty(\Xi)}
+\nm{\tilde{\beta}_j}_{L^\infty(\Xi)}\Bigr)
\\\notag
&\hphantom{{}\leq{}}
+C_1(\omega)
\nm[\Big]{\sqrt{\tilde{\alpha}_j^2+\tilde{\beta}_j^2}}_{L^\infty(\Xi)}
\Bigl(
\nm{\DN_{c_0}\gamma_0}_{L^2(\Xi)}^2
+\nm{\gamma_0'}_{L^2(\Xi)}^2
\Bigr)^{\frac{1}{2}}
\biggr).
\end{align}
The remaining norms involving $\gamma_0$ and $\DN_{c_0}\gamma_0$ are bounded by
\begin{align*}
\nm{\gamma_0'}_{L^2(\Xi)}
&\leq\nm{P'}_{L^2(\Xi)}
+\nm{\gamma_0-P}_{H^2(\Xi)},
&
\nm{\DN_{c_0}\gamma_0}_{L^2(\Xi)}
&\leq\nm{\tilde{h}}_{L^2(\Xi)}
+\nm{\DN_{c_0}\gamma_0-\tilde{h}}_{H^1(\Xi)}.
\end{align*}
These differences are bounded by Lemma~\ref{lem:normalization-trace-errors} and Proposition~\ref{prop:DtN-error}, as in the proof of Proposition~\ref{prop:residual-bound}.
The remaining norms involve Fourier polynomials and are bounded using finite Fourier and interval arithmetic.
Inserting \eqref{eqn:a-priori-rotation-columns} and \eqref{eqn:a-priori-configuration-columns} into \eqref{eqn:a-priori-second-operator-norm} provides the second bound required by \eqref{eqn:a-priori-kappa}.
\end{proof}

Proposition~\ref{prop:L2-inverse-bound} states that the self-adjoint extension $\mathscr{L}\colon Y\subset\mathcal{H}\to\mathcal{H}$ is invertible. 
Combining this with the a priori estimate from Lemma~\ref{lem:a-priori-estimate} allows us to prove that the original linearisation $\mathscr{L}\colon X\to Y$ is an isomorphism and bound its inverse. 

\begin{proposition}[Strong inverse bound]
\label{prop:linear-estimate}
The linearised operator $\mathscr{L}=DF\vert_0\colon X\to Y$ from Proposition~\ref{prop:linearisation} is an isomorphism and 
\(
\nm{\mathscr{L}^{-1}}_{Y\to X}\leq q, 
\)
where $q$ is recorded in Appendix~\ref{section-numerical-values}.
\end{proposition}

\begin{proof}
By Lemma~\ref{lem:a-priori-estimate} and Proposition~\ref{prop:L2-inverse-bound}, we have for every $x\in X$ 
\begin{align*}
\nm{x}_X
&\leq2\nm{\mathscr{L}x}_Y+\kappa\nm{x}_{\mathcal{H}}
\\
&\leq2\nm{\mathscr{L}x}_Y+\kappa q_0\nm{\mathscr{L}x}_{\mathcal{H}}
\leq(2+\kappa q_0)\nm{\mathscr{L}x}_Y. 
\end{align*} 
In particular, $\mathscr{L}\colon X\to Y$ is injective.
To prove surjectivity, let $f=(f_v,f_\mu,f_\lambda)\in Y$.
By Proposition~\ref{prop:L2-inverse-bound}, there exists $x=(v,\mu,\lambda)\in Y$ such that $\mathscr{L}x=f$.
Recalling Proposition~\ref{prop:linearisation}, we obtain
\begin{align*}
\DN_{\mathrm{dec}}v
&=f_v-(S_{vv}-\DN_{\mathrm{dec}})v-S_{v\mu}\mu-\rho\lambda
\in H^1(\Xi,\R^2),
\end{align*} 
where we used \eqref{eqn:a-priori-Svv-bound} and \eqref{eqn:a-priori-second-operator-norm}. 
Since $v\in H^1(\Xi,\R^2)$, the Fourier multiplier formula in Lemma~\ref{lem:decoupled-Fourier} implies that $v\in H^2(\Xi,\R^2)$.
Hence $x\in X$, which proves that $\mathscr{L}\colon X\to Y$ is surjective.
The claimed bound follows for any $q\geq2+\kappa q_0$.
\end{proof}

\subsection{Quadratic estimate}

The third estimate required by Proposition~\ref{prop:existence} controls the nonlinear remainder $\mathscr{N}$ defined as in  Theorem~\ref{thm:contraction}, for which we prove a quantitative Lipschitz bound in this subsection.  
Proposition~\ref{prop:linearisation} computes the derivative of the gradient map $F$ from Definition~\ref{defn:gradient-map} at the origin.
The following statement extends this calculation to arbitrary points in the domain of $F$.

\begin{proposition}
\label{prop:linearisation-general}
The derivative of the gradient map $F$ from Definition~\ref{defn:gradient-map} at a general point $x=(v,\mu,\lambda)\in H^2(\Xi,\R^2)\times B\times\R^3$ is given by 
\begin{align*}
DF\vert_x=
\begin{pmatrix}
S_{vv}(x)&S_{v\mu}(x)&\rho\\
S_{v\mu}(x)^*&S_{\mu\mu}(x)&0\\
\rho^*&0&0
\end{pmatrix},
\end{align*}
where, for $w\in H^2(\Xi,\R^2)$, $h\in\R^{3b-6}$ and $i,j\in\{1,\ldots,3b-6\}$,
\begin{align*}
S_{vv}(x)w
&=(D_v\gamma_v)^{\intercal}\DN_{c_\mu}\bigl((D_v\gamma_v)w\bigr)
-\frac{\sk{v,w}\Gamma^{\intercal}\DN_{c_\mu}\gamma_v
+\sk{\Gamma w,\DN_{c_\mu}\gamma_v}v}{(1+\abs{v}^2)^{3/2}}
\\&\hphantom{{}={}}
-\frac{\sk{\gamma_v,\DN_{c_\mu}\gamma_v}}{1+\abs{v}^2}w
+\frac{3\sk{\gamma_v,\DN_{c_\mu}\gamma_v}\sk{v,w}}{(1+\abs{v}^2)^2}v,
\\
S_{v\mu}(x)h
&=\sum_{j=1}^{3b-6}h_j(D_v\gamma_v)^{\intercal}
\Bigl(\frac{\partial}{\partial\mu_j}\DN_{c_\mu}\Bigr)\gamma_v,
\\
\bigl(S_{v\mu}(x)^*w\bigr)_j
&=\int_\Xi\sk[\Big]{(D_v\gamma_v)w,
\Bigl(\frac{\partial}{\partial\mu_j}\DN_{c_\mu}\Bigr)\gamma_v}\,d\theta,
\\
\bigl(S_{\mu\mu}(x)\bigr)_{ij}
&=\int_\Xi\biggl(
\frac{1}{2}\frac{\partial\tilde{\alpha}_i}{\partial\mu_j}
\bigl(\abs{\gamma_v'}^2-\abs{\DN_{c_\mu}\gamma_v}^2\bigr)
-\frac{\partial\tilde{\beta}_i}{\partial\mu_j}\sk{\gamma_v',\DN_{c_\mu}\gamma_v}
\\&\hphantom{{}={}}
-\sk[\Big]{\tilde{\alpha}_i\DN_{c_\mu}\gamma_v+\tilde{\beta}_i\gamma_v',
\Bigl(\frac{\partial}{\partial\mu_j}\DN_{c_\mu}\Bigr)\gamma_v}
\biggr)\,d\theta. 
\end{align*}
Here, $\tilde{\alpha}_i,\tilde{\beta}_i$ are the normalised boundary velocities from Lemma~\ref{lem:configuration-linearisation}, $D_v\gamma_v$ is given by \eqref{eqn:20261017a}, and $\rho,\rho^*$ are as in Proposition~\ref{prop:linearisation}. 
Moreover, 
$DF|_x$ is symmetric with respect to the product of the $L^2$-inner product on $H^2(\Xi,\R^2)$ and the Euclidean inner products on $\R^{3b-6}$ and $\R^3$.
\end{proposition}

\begin{proof}
Recalling $\gamma_v$ from \eqref{eqn:perturbation}
and $D_v\gamma_v$ from \eqref{eqn:20261019c}, we have for all $w,z\in H^2(\Xi,\R^2)$
\begin{align*}
D_v\gamma_v(w)
&=\frac{\Gamma w}{\sqrt{1+\abs{v}^2}}
-\frac{\sk{v,w}}{1+\abs{v}^2}\gamma_v,
\\
D_v^2\gamma_v(w,z)
&=-\frac{\sk{v,z}\Gamma w+\sk{v,w}\Gamma z}{(1+\abs{v}^2)^{3/2}}
-\frac{\sk{w,z}}{1+\abs{v}^2}\gamma_v
+\frac{3\sk{v,w}\sk{v,z}}{(1+\abs{v}^2)^2}\gamma_v.
\end{align*}
Hence, recalling the Lagrangian $L$ from Remark~\ref{rem:Lagrange_multipliers}, we obtain from \eqref{eqn:reduced-energy-DtN} that
\begin{align*}
D^2_{vv}L\vert_x(w,w)
&=\int_\Xi\sk[\Big]{(D_v\gamma_v)w,\DN_{c_\mu}((D_v\gamma_v)w)}
+\sk[\big]{D_v^2\gamma_v(w,w),\DN_{c_\mu}\gamma_v}\,d\theta
=\int_\Xi\sk[\big]{w,S_{vv}(x)w}\,d\theta.
\end{align*}
Since $\gamma_v$ is independent of $\mu$, \eqref{eqn:reduced-energy-DtN} implies for every $h\in\R^{3b-6}$ that
\begin{align*}
D^2_{v\mu}L\vert_x(w,h)
&=\sum_{j=1}^{3b-6}h_j
\int_\Xi\sk[\Big]{(D_v\gamma_v)w,
\Bigl(\frac{\partial}{\partial\mu_j}\DN_{c_\mu}\Bigr)\gamma_v}\,d\theta
=\int_\Xi\sk{w,S_{v\mu}(x)h}\,d\theta.
\end{align*}
Since $L$ is smooth, its Hessian is symmetric, and therefore
\begin{align*}
D^2_{\mu v}L\vert_x(h,w)
&=D^2_{v\mu}L\vert_x(w,h)
=\sum_{j=1}^{3b-6}h_j\bigl(S_{v\mu}(x)^*w\bigr)_j.
\end{align*}
Since the Lagrange multiplier term in Remark~\ref{rem:Lagrange_multipliers} is independent of $\mu$, the $S_{\mu\mu}(x)$-block is the Hessian of $\mu\mapsto E(\gamma_v,c_\mu)$.
The argument for Lemma~\ref{lem:configuration-linearisation} with $\gamma_0$ replaced by $\gamma_v$ yields the stated formula for $S_{\mu\mu}(x)$.

Since $S_{vv}(x)$ is symmetric, the equality of quadratic forms determines the $vv$-derivative of $L$. 
Finally, recalling Remark~\ref{rem:Lagrange_multipliers}, we have
\begin{align*}
D^2_{v\lambda}L\vert_x(w,\ell)
&=\int_\Xi\sk{w,\rho\,\ell}\,d\theta,
&
D^2_{\mu\lambda}L\vert_x&=0,
&
D^2_{\lambda\lambda}L\vert_x&=0.
\end{align*}
By symmetry of the Hessian, these identities yield the remaining blocks of $DF\vert_x$.
\end{proof}

The constants $r_*,\omega>0$ recorded in Appendix~\ref{section-numerical-values} are chosen such that, for every $\abs{\mu}<r_*$, the configuration $c_\mu$ from \eqref{eqn:perturbation-configuration} admits pairwise disjoint conformal collars of width at least $\omega$ as in \eqref{eqn:conformal_collar_width}.
The following lemma verifies the hypotheses of Proposition~\ref{prop:DtN-differences} on the normalised boundary velocities $\tilde\alpha_{\mu,h},\tilde\beta_{\mu,h}$ 
and provides additional bounds for their derivatives. 
 
\begin{lemma}
\label{lem:configuration-bounds}
For all $h,k\in\R^{3b-6}$, $\abs{\mu},\abs{\nu}<r_*$ and $C_c>0$ as recorded in Appendix~\ref{section-numerical-values} we have 
\begin{align*}
\nm{\tilde{\alpha}_{\mu,h}}_{L^\infty(\Xi)}
+\nm[\big]{\tilde{\beta}_{\mu,h}}_{L^\infty(\Xi)}
&\leq C_c\abs{h},
\\
\nm[\Big]{\tfrac{d}{dt}\big\vert_{t=0}\tilde{\alpha}_{\mu+tk,h}}_{L^\infty(\Xi)}
+\nm[\Big]{\tfrac{d}{dt}\big\vert_{t=0}\tilde{\beta}_{\mu+tk,h}}_{L^\infty(\Xi)}
&\leq C_c\abs{h}\abs{k},
\\[1ex]
\nm{\tilde{\alpha}_{\mu,h}-\tilde{\alpha}_{\nu,h}}_{L^\infty(\Xi)}
+\nm[\big]{\tilde{\beta}_{\mu,h}-\tilde{\beta}_{\nu,h}}_{L^\infty(\Xi)}
&\leq C_c\abs{\mu-\nu}\abs{h},
\\
\nm[\Big]{\tfrac{d}{dt}\big\vert_{t=0}
(\tilde{\alpha}_{\mu+tk,h}-\tilde{\alpha}_{\nu+tk,h})}_{L^\infty(\Xi)}
+\nm[\Big]{\tfrac{d}{dt}\big\vert_{t=0}
\bigl(\tilde{\beta}_{\mu+tk,h}-\tilde{\beta}_{\nu+tk,h}\bigr)}_{L^\infty(\Xi)}
&\leq C_c\abs{\mu-\nu}\abs{h}\abs{k}.
\end{align*}
\end{lemma}

\begin{proof}
By \eqref{eqn:normal-velocity-fourier} and \eqref{eqn:tangential-velocity-fourier}, the boundary velocities in question are Fourier polynomials of degree at most one on each component of $\Xi$.  
The coefficients of these polynomials and their first two derivatives with respect to $\mu$ can be bounded by interval arithmetic, yielding the first two inequalities and a uniform bound for the second derivatives.
The remaining inequalities follow from the fundamental theorem of calculus along the line segment joining $\mu$ and $\nu$.
\end{proof}
 
\begin{lemma}
\label{lem:boundary-normalization}
Let $v\in H^2(\Xi,\R^2)$ satisfy $8c_\infty\nm{v}_{H^2}\leq1$, where $c_\infty=\sqrt{\frac{1}{2}\coth(\pi)}$ as in \eqref{eqn:Sobolev-constant}.
Then, recalling $\gamma_v$ from \eqref{eqn:perturbation}, we have
\begin{align*}
\nm{\gamma_v-\gamma_0}_{H^2}
+\nm{D_v\gamma_v-\Gamma}_{H^2}
&\leq12c_\infty
\Bigl(\nm{\gamma_0}_{H^2}+\nm{\Gamma}_{H^2}\Bigr)\nm{v}_{H^2}.
\end{align*}
Moreover, for $p\in\{1,\frac{3}{2},2\}$, $j\in\{1,2\}$ and $f\in H^j(\Xi,\R^d)$, we have
\begin{align*}
\nm[\Big]{\frac{f}{(1+\abs{v}^2)^p}}_{H^j}
&\leq2\nm{f}_{H^j},
&
\nm[\Big]{\bigl((1+\abs{v}^2)^{-p}-1\bigr)f}_{H^j}
&\leq12c_\infty\nm{v}_{H^2}\nm{f}_{H^j}.
\end{align*}
\end{lemma}
 
\begin{proof}
Let $\phi\in H^2(\Xi,\R^{m\times n})$ and $\psi\in H^j(\Xi,\R^{n\times d})$, where $j\in\{1,2\}$.
By \eqref{eqn:Sobolev-constant}, \eqref{eqn:Sobolev-embedding} and the product rule, we have
\begin{align}\notag
\nm{\phi\psi}_{H^1}
&\leq\nm{\phi}_{L^\infty}\nm{\psi}_{H^1}
+\nm{\phi'}_{L^\infty}\nm{\psi}_{L^2}
\leq2c_\infty\nm{\phi}_{H^2}\nm{\psi}_{H^1},
\\[1ex]\notag
\nm{\phi\psi}_{H^2}
&\leq\nm{\phi}_{L^\infty}\nm{\psi}_{H^2}
+\nm{\phi'}_{L^\infty}
\bigl(2\nm{\psi}_{L^2}^2+4\nm{\psi'}_{L^2}^2\bigr)^{\frac{1}{2}}
+\nm{\phi''}_{L^2}\nm{\psi}_{L^\infty}
\\\notag
&\leq\bigl(2+\sqrt{2}\bigr)c_\infty\nm{\phi}_{H^2}\nm{\psi}_{H^2}.
\shortintertext{Consequently,}
\label{eqn:normalization-product}
\nm{\phi\psi}_{H^j}
&\leq4c_\infty\nm{\phi}_{H^2}\nm{\psi}_{H^j}.
\end{align}
By assumption, $s\vcentcolon=4c_\infty\nm{v}_{H^2}\leq\frac{1}{2}$.
Applying \eqref{eqn:normalization-product} twice, we obtain for every $f\in H^j(\Xi,\R^d)$
\begin{align*}
\nm[\big]{\abs{v}^2f}_{H^j}
&=\nm[\big]{(v^{\intercal}v)f}_{H^j}
\leq4c_\infty\nm[\big]{v^{\intercal}v}_{H^2}\nm{f}_{H^j}
\leq(4c_\infty)^2\nm{v}_{H^2}^2\nm{f}_{H^j}
=s^2\nm{f}_{H^j}.
\end{align*}
By induction in $n\in\N$, we obtain 
\begin{align}
\label{eqn:normalization-powers}
\nm[\big]{(\abs{v}^2)^n f}_{H^j}
&\leq s^{2n}\nm{f}_{H^j}.
\end{align}
Let $a\vcentcolon=(1+\abs{v}^2)^{-\frac{1}{2}}$.
Using the binomial series and \eqref{eqn:normalization-powers}, we obtain
\begin{align*}
\nm{af}_{H^j}
&\leq\sum_{n=0}^\infty\frac{1}{4^n}\binom{2n}{n}s^{2n}\nm{f}_{H^j}
=\frac{1}{\sqrt{1-s^2}}\nm{f}_{H^j}
\leq\frac{4}{3}\nm{f}_{H^j},
\\
\nm{a^3f}_{H^j}
&\leq(1-s^2)^{-\frac{3}{2}}\nm{f}_{H^j}
\leq2\nm{f}_{H^j}.
\end{align*}
Moreover, since $4^{-n}\binom{2n}{n}\leq\frac{1}{2}$ for $n\geq1$, we have
\begin{align}
\label{eqn:normalization-reciprocal-difference}
\nm{(a-1)f}_{H^j}
&\leq\sum_{n=1}^\infty\frac{1}{4^n}\binom{2n}{n}s^{2n}\nm{f}_{H^j}
\leq\frac{s^2}{2(1-s^2)}\nm{f}_{H^j}
\leq\frac{s}{3}\nm{f}_{H^j}.
\end{align}
Recalling $\gamma_v$ from \eqref{eqn:perturbation} and its derivative $D_v\gamma_v$ from \eqref{eqn:20261019c}, we have
\begin{align*}
\gamma_v-\gamma_0
&=(a-1)\gamma_0+a\Gamma v,
&
D_v\gamma_v-\Gamma
&=(a-1)\Gamma-a^3(\gamma_0+\Gamma v)v^{\intercal}.
\end{align*}
Since $\nm{\Gamma v}_{H^2}\leq4c_\infty\nm{\Gamma}_{H^2}\nm{v}_{H^2}=s\nm{\Gamma}_{H^2}$ by \eqref{eqn:normalization-product}, it follows that
\begin{align*}
\nm{\gamma_v-\gamma_0}_{H^2}
&\leq\nm[\big]{(a-1)\gamma_0}_{H^2}+\nm{a\Gamma v}_{H^2}
\leq\frac{s}{3}\nm{\gamma_0}_{H^2}
+\frac{4s}{3}\nm{\Gamma}_{H^2},
\\[.5ex]
\nm{D_v\gamma_v-\Gamma}_{H^2}
&\leq\nm{(a-1)\Gamma}_{H^2}
+\nm{a^3(\gamma_0+\Gamma v)v^{\intercal}}_{H^2}
\\
&\leq\frac{s}{3}\nm{\Gamma}_{H^2}
+2s\Bigl(\nm{\gamma_0}_{H^2}+\nm{\Gamma v}_{H^2}\Bigr)
\leq2s\nm{\gamma_0}_{H^2}
+\bigl(\tfrac{s}{3}+2s^2\bigr)\nm{\Gamma}_{H^2}.
\end{align*}
Since $2s^2\leq s$, adding the last two estimates yields
\begin{align*}
\nm{\gamma_v-\gamma_0}_{H^2}
+\nm{D_v\gamma_v-\Gamma}_{H^2}
&\leq\frac{7s}{3}\nm{\gamma_0}_{H^2}
+\frac{8s}{3}\nm{\Gamma}_{H^2}
\leq\frac{32}{3}c_\infty
\Bigl(\nm{\gamma_0}_{H^2}+\nm{\Gamma}_{H^2}\Bigr)\nm{v}_{H^2}
\end{align*}
which proves the first claim. 
For $p\in\{1,\frac{3}{2},2\}$, the coefficients of the binomial series  
\begin{align*} 
(1+t)^{-p}
&=1+\sum_{n=1}^\infty \frac{p(p+1)\cdots(p+n-1)}{n!}(-t)^n,
&
\abs{t}&<1.
\end{align*}
are bounded by $(n+1)!/n!=n+1$. 
Hence, for $j\in\{1,2\}$ and $f\in H^j(\Xi,\R^d)$,
\begin{align*}
\nm[\Big]{\frac{f}{(1+\abs{v}^2)^p}}_{H^j}
&\leq\sum_{n=0}^\infty(n+1)s^{2n}\nm{f}_{H^j}
=\frac{1}{(1-s^2)^2}\nm{f}_{H^j}
\leq2\nm{f}_{H^j},
\\[1ex]
\nm[\Big]{\bigl((1+\abs{v}^2)^{-p}-1\bigr)f}_{H^j}
&\leq\sum_{n=1}^\infty(n+1)s^{2n}\nm{f}_{H^j}
=\frac{s^2(2-s^2)}{(1-s^2)^2}\nm{f}_{H^j}
\leq\frac{32}{9}s^2\nm{f}_{H^j}. 
\end{align*}
The fact that $s\leq\frac{1}{2}$ and
$\frac{32}{9}s^2
\leq\frac{16}{9}s
\leq3s
=12c_\infty\nm{v}_{H^2}$ concludes the proof. 
\end{proof} 

\begin{proposition}[Nonlinear remainder bound]
\label{prop:nonlinear}
Given $F$ and $\mathscr{L}=DF\vert_0\colon X\to Y$ as in Definition~\ref{defn:gradient-map} and Proposition~\ref{prop:linearisation}, the nonlinear remainder 
$\mathscr{N}(x)=F(x)-F(0)-\mathscr{L}x$ satisfies
\begin{align*}
\nm[\big]{\mathscr{N}(x)-\mathscr{N}(y)}_Y
&\leq K\nm{x-y}_X
\end{align*}
for all $x,y\in X$ with $\nm{x}_X,\nm{y}_X\leq r$ and with the constants $r,K>0$ recorded in Appendix~\ref{section-numerical-values}.
\end{proposition}
 
\begin{proof} 
The radius $r$ recorded in Appendix~\ref{section-numerical-values} satisfies $0<r<r_*$ and $8c_\infty r\leq1$.  
Given $x,y\in X$ with $\nm{x}_X,\nm{y}_X\leq r$, the line segment from $x$ to $y$ is in the domain of $F$, and
\begin{align*}
\mathscr N(x)-\mathscr N(y)
&=\int_0^1\Bigl(DF\vert_{y+s(x-y)}-DF\vert_0\Bigr)(x-y)\,ds.
\end{align*}
It therefore suffices to prove
\begin{align}
\label{eqn:nonlinear-derivative-bound}
\nm[\big]{DF\vert_x-DF\vert_0}_{X\to Y}&\leq K
\end{align}
for all $x=(v,\mu,\lambda)\in X$ with $\nm{x}_X\leq r$.  
Given $z=(w,h,\ell)\in X$, Proposition~\ref{prop:linearisation-general} implies 
\begin{align}
\label{eqn:nonlinear-blocks}
\bigl(DF\vert_x-DF\vert_0\bigr)z
&=
\begin{pmatrix}
\bigl(S_{vv}(x)-S_{vv}(0)\bigr)w
+\bigl(S_{v\mu}(x)-S_{v\mu}(0)\bigr)h
\\
\bigl(S_{v\mu}(x)^*-S_{v\mu}(0)^*\bigr)w
+\bigl(S_{\mu\mu}(x)-S_{\mu\mu}(0)\bigr)h
\\
0
\end{pmatrix}
\end{align}
because $\rho$ and $\rho^*$ are independent of $x$.
Moreover, the $L^2$-adjoint satisfies
\begin{align}\notag
\abs[\Big]{\bigl(S_{v\mu}(x)^*-S_{v\mu}(0)^*\bigr)w}
&=\sup_{\abs{k}=1}
\abs[\bigg]{\int_\Xi
\sk[\Big]{w,\bigl(S_{v\mu}(x)-S_{v\mu}(0)\bigr)k}\,d\theta}
\\\label{eqn:nonlinear-adjoint-bound}
&\leq\nm{w}_{H^2}
\sup_{\abs{k}=1}
\nm[\Big]{\bigl(S_{v\mu}(x)-S_{v\mu}(0)\bigr)k}_{H^1}.
\end{align}
Thus, it suffices to estimate $S_{vv}(x)-S_{vv}(0)$, $S_{v\mu}(x)-S_{v\mu}(0)$ and $S_{\mu\mu}(x)-S_{\mu\mu}(0)$. 
By Proposition~\ref{prop:linearisation-general}, 
\begin{align}\label{eqn:Svv-term1}
\bigl(S_{vv}(x)-S_{vv}(0)\bigr)w
&=(D_v\gamma_v)^{\intercal}\DN_{c_\mu}\bigl((D_v\gamma_v)w\bigr)
-\Gamma^{\intercal}\DN_{c_0}(\Gamma w)
\\[1ex]&\hphantom{{}={}}\label{eqn:Svv-term2}
-\frac{\sk{v,w}\Gamma^{\intercal}\DN_{c_\mu}\gamma_v
+\sk{\Gamma w,\DN_{c_\mu}\gamma_v}v}{(1+\abs{v}^2)^{3/2}}
\\&\hphantom{{}={}}\label{eqn:Svv-term3}
-\biggl(
\frac{\sk{\gamma_v,\DN_{c_\mu}\gamma_v}}{1+\abs{v}^2}
-\sk{\gamma_0,\DN_{c_0}\gamma_0}
\biggr)w
\\&\hphantom{{}={}}\label{eqn:Svv-term4}
+\frac{3\sk{\gamma_v,\DN_{c_\mu}\gamma_v}\sk{v,w}}
{(1+\abs{v}^2)^2}v.
\end{align}
We bound each one of the summands \eqref{eqn:Svv-term1}--\eqref{eqn:Svv-term4} separately. 
By adding and subtracting intermediate terms, \eqref{eqn:Svv-term1} becomes 
\begin{align*}
&(D_v\gamma_v)^{\intercal}\DN_{c_\mu}\bigl((D_v\gamma_v)w\bigr)
-\Gamma^{\intercal}\DN_{c_0}(\Gamma w)
\\
&=(D_v\gamma_v-\Gamma)^{\intercal}\DN_{c_\mu}\bigl((D_v\gamma_v)w\bigr)
+\Gamma^{\intercal}(\DN_{c_\mu}-\DN_{c_0})\bigl((D_v\gamma_v)w\bigr)
+\Gamma^{\intercal}\DN_{c_0}\bigl((D_v\gamma_v-\Gamma)w\bigr).
\end{align*}
Estimates \eqref{eqn:normalization-product}, \eqref{eqn:DtN-differences-1} and \eqref{eqn:DtN-differences-3} imply
\begin{align*}
\nm[\big]{(D_v\gamma_v-\Gamma)^{\intercal}
\DN_{c_\mu}\bigl((D_v\gamma_v)w\bigr)}_{H^1}
&\leq(4c_\infty)^2\bigl(1+C_1(\omega)\bigr)
\nm[\big]{D_v\gamma_v-\Gamma}_{H^2}
\nm{D_v\gamma_v}_{H^2}\nm{w}_{H^2},
\\[1ex]
\nm[\big]{\Gamma^{\intercal}(\DN_{c_\mu}-\DN_{c_0})
\bigl((D_v\gamma_v)w\bigr)}_{H^1}
&\leq(4c_\infty)^2C_\Lambda\abs{\mu}
\nm{\Gamma}_{H^2}\nm{D_v\gamma_v}_{H^2}\nm{w}_{H^2},
\\[1ex]
\nm[\big]{\Gamma^{\intercal}\DN_{c_0}
\bigl((D_v\gamma_v-\Gamma)w\bigr)}_{H^1}
&\leq(4c_\infty)^2\bigl(1+C_1(\omega)\bigr)
\nm{\Gamma}_{H^2}
\nm[\big]{D_v\gamma_v-\Gamma}_{H^2}\nm{w}_{H^2}.
\end{align*}
Consequently, the term \eqref{eqn:Svv-term1} is bounded by 
\begin{align}
\label{eqn:Svv-estimate1}
&\nm[\Big]{
(D_v\gamma_v)^{\intercal}\DN_{c_\mu}\bigl((D_v\gamma_v)w\bigr)
-\Gamma^{\intercal}\DN_{c_0}(\Gamma w)
}_{H^1}
\\\notag
&\leq(4c_\infty)^2\biggl(
\bigl(1+C_1(\omega)\bigr)
\Bigl(\nm{D_v\gamma_v}_{H^2}+\nm{\Gamma}_{H^2}\Bigr)
\nm[\big]{D_v\gamma_v-\Gamma}_{H^2}
+C_\Lambda\abs{\mu}
\nm{\Gamma}_{H^2}\nm{D_v\gamma_v}_{H^2}
\biggr)\nm{w}_{H^2}.
\end{align} 
Since $\nm{v}_{H^2}\leq\nm{x}_X\leq r$, Lemma~\ref{lem:boundary-normalization} implies
\begin{align*}
\nm[\big]{D_v\gamma_v-\Gamma}_{H^2}
&\leq12c_\infty
\bigl(\nm{\gamma_0}_{H^2}+\nm{\Gamma}_{H^2}\bigr)\nm{x}_X,
\\[1ex]
\nm{D_v\gamma_v}_{H^2}
&\leq\nm{\Gamma}_{H^2}
+\nm[\big]{D_v\gamma_v-\Gamma}_{H^2}
\leq\nm{\Gamma}_{H^2}
+12c_\infty r
\bigl(\nm{\gamma_0}_{H^2}+\nm{\Gamma}_{H^2}\bigr).
\end{align*}
Since $\abs{\mu}\leq\nm{x}_X$, each term on the right-hand side of inequality \eqref{eqn:Svv-estimate1} is therefore bounded by an explicitly computable constant times $\nm{x}_X\nm{w}_{H^2}$.

In the second summand \eqref{eqn:Svv-term2}, Lemma~\ref{lem:boundary-normalization}, \eqref{eqn:normalization-product} and \eqref{eqn:DtN-differences-1} imply
\begin{align}\notag
&\nm[\bigg]{
\frac{\sk{v,w}\Gamma^{\intercal}\DN_{c_\mu}\gamma_v
+\sk{\Gamma w,\DN_{c_\mu}\gamma_v}v}{(1+\abs{v}^2)^{3/2}}
}_{H^1}
\\\notag
&\leq2(4c_\infty)^3
\Bigl(
\nm{v}_{H^2}\nm{w}_{H^2}\nm{\Gamma}_{H^2}
+\nm{\Gamma}_{H^2}\nm{w}_{H^2}\nm{v}_{H^2}
\Bigr)\nm{\DN_{c_\mu}\gamma_v}_{H^1}
\\\label{eqn:Svv-estimate2}
&\leq4(4c_\infty)^3\bigl(1+C_1(\omega)\bigr)
\nm{\Gamma}_{H^2}\nm{\gamma_v}_{H^2}
\nm{x}_X\nm{w}_{H^2}.
\end{align}
By the triangle inequality and Lemma~\ref{lem:boundary-normalization}, we have
\begin{align}
\label{eqn:quadratic-gamma-bound}
\nm{\gamma_v}_{H^2}
&\leq\nm{\gamma_0}_{H^2}+\nm{\gamma_v-\gamma_0}_{H^2}
\leq\nm{\gamma_0}_{H^2}
+12c_\infty r\bigl(\nm{\gamma_0}_{H^2}+\nm{\Gamma}_{H^2}\bigr).
\end{align}
Hence, \eqref{eqn:Svv-estimate2} is again bounded by an explicitly computable constant times $\nm{x}_X\nm{w}_{H^2}$. 

In the third summand \eqref{eqn:Svv-term3}, addition and subtraction of intermediate terms yields
\begin{align*}
&\frac{\sk{\gamma_v,\DN_{c_\mu}\gamma_v}}{1+\abs{v}^2}
-\sk{\gamma_0,\DN_{c_0}\gamma_0}
\\
&=\bigl((1+\abs{v}^2)^{-1}-1\bigr)
\sk{\gamma_0,\DN_{c_0}\gamma_0}
+(1+\abs{v}^2)^{-1}
\bigl(
\sk{\gamma_v-\gamma_0,\DN_{c_\mu}\gamma_v}
+\sk{\gamma_0,\DN_{c_\mu}\gamma_v-\DN_{c_0}\gamma_0}
\bigr).
\end{align*}
By Lemma~\ref{lem:boundary-normalization} and \eqref{eqn:normalization-product}, we obtain
\begin{align}
\notag
&\nm[\bigg]{
\biggl(
\frac{\sk{\gamma_v,\DN_{c_\mu}\gamma_v}}{1+\abs{v}^2}
-\sk{\gamma_0,\DN_{c_0}\gamma_0}
\biggr)w
}_{H^1}
\\\label{eqn:Svv-estimate3}
&\leq(4c_\infty)^2
\biggl(
12c_\infty\nm{v}_{H^2}
\nm{\gamma_0}_{H^2}\nm{\DN_{c_0}\gamma_0}_{H^1}
+2\nm{\gamma_v-\gamma_0}_{H^2}\nm{\DN_{c_\mu}\gamma_v}_{H^1}
\\[-.5ex]\notag&\hphantom{{}={}}
+2\nm{\gamma_0}_{H^2}
\nm{\DN_{c_\mu}\gamma_v-\DN_{c_0}\gamma_0}_{H^1}
\biggr)\nm{w}_{H^2}.
\end{align}
By \eqref{eqn:DtN-differences-1}, \eqref{eqn:DtN-differences-3} and Lemma~\ref{lem:boundary-normalization}, using $\abs{\mu}\leq\nm{x}_X$, we obtain
\begin{align*}
\nm{\DN_{c_\mu}\gamma_v}_{H^1}
&\leq\bigl(1+C_1(\omega)\bigr)\nm{\gamma_v}_{H^2},
\\[1ex]
\nm{\DN_{c_\mu}\gamma_v-\DN_{c_0}\gamma_0}_{H^1}
&\leq\nm{\DN_{c_\mu}(\gamma_v-\gamma_0)}_{H^1}
+\nm{(\DN_{c_\mu}-\DN_{c_0})\gamma_0}_{H^1}
\\
&\leq\bigl(1+C_1(\omega)\bigr)\nm{\gamma_v-\gamma_0}_{H^2}
+C_\Lambda\abs{\mu}\nm{\gamma_0}_{H^2}
\\
&\leq
\Bigl(
12c_\infty\bigl(1+C_1(\omega)\bigr)
\bigl(\nm{\gamma_0}_{H^2}+\nm{\Gamma}_{H^2}\bigr)
+C_\Lambda\nm{\gamma_0}_{H^2}
\Bigr)\nm{x}_X.
\end{align*}
Together with \eqref{eqn:quadratic-gamma-bound} and Lemma~\ref{lem:boundary-normalization}, the bounds for $\DN_{c_\mu}\gamma_v$ and $\DN_{c_\mu}\gamma_v-\DN_{c_0}\gamma_0$ imply that the right-hand side of \eqref{eqn:Svv-estimate3} is bounded by an explicitly computable constant times $\nm{x}_X\nm{w}_{H^2}$.

For the last summand \eqref{eqn:Svv-term4}, estimates \eqref{eqn:normalization-product} and \eqref{eqn:DtN-differences-1} and Lemma~\ref{lem:boundary-normalization} with $p=2$ imply
\begin{align}\notag
\nm[\bigg]{
\frac{3\sk{\gamma_v,\DN_{c_\mu}\gamma_v}\sk{v,w}}{(1+\abs{v}^2)^2}v
}_{H^1}
&\leq6(4c_\infty)^4
\nm{\gamma_v}_{H^2}\nm{\DN_{c_\mu}\gamma_v}_{H^1}
\nm{v}_{H^2}^2\nm{w}_{H^2}
\\\label{eqn:Svv-estimate4}
&\leq6(4c_\infty)^4\bigl(1+C_1(\omega)\bigr)
r\nm{\gamma_v}_{H^2}^2
\nm{x}_X\nm{w}_{H^2}.
\end{align}
By estimate \eqref{eqn:quadratic-gamma-bound}, the coefficient in \eqref{eqn:Svv-estimate4} is uniformly bounded.
Combining 
\eqref{eqn:Svv-estimate1},
\eqref{eqn:Svv-estimate2},
\eqref{eqn:Svv-estimate3} and 
\eqref{eqn:Svv-estimate4}, 
we conclude that $\nm{(S_{vv}(x)-S_{vv}(0))w}_{H^1}$ is bounded by an explicitly computable constant times $\nm{x}_X\nm{w}_{H^2}$.

We next estimate $S_{v\mu}(x)-S_{v\mu}(0)$.  
Recalling the directional derivative $\dot\DN_{\mu,h}$ of the Dirichlet-to-Neumann operator from \eqref{eqn:dotDN},  Proposition~\ref{prop:linearisation-general} implies 
\begin{align}
\label{eqn:Svmu-splitting}
\bigl(S_{v\mu}(x)-S_{v\mu}(0)\bigr)h
&=(D_v\gamma_v-\Gamma)^{\intercal}\dot\DN_{\mu,h}\gamma_v
+\Gamma^{\intercal}\dot\DN_{\mu,h}(\gamma_v-\gamma_0)
+\Gamma^{\intercal}(\dot\DN_{\mu,h}-\dot\DN_{0,h})\gamma_0.
\end{align}
By estimates \eqref{eqn:normalization-product}, \eqref{eqn:DtN-differences-2} and \eqref{eqn:DtN-differences-4}, we have
\begin{align*}
\nm[\big]{(D_v\gamma_v-\Gamma)^{\intercal}\dot\DN_{\mu,h}\gamma_v}_{H^1}
&\leq4c_\infty C_\Lambda
\nm{D_v\gamma_v-\Gamma}_{H^2}
\nm{\gamma_v}_{H^1}\abs{h},
\\[1ex]
\nm[\big]{\Gamma^{\intercal}\dot\DN_{\mu,h}(\gamma_v-\gamma_0)}_{H^1}
&\leq4c_\infty C_\Lambda
\nm{\Gamma}_{H^2}
\nm{\gamma_v-\gamma_0}_{H^1}\abs{h},
\\[1ex]
\nm[\big]{\Gamma^{\intercal}(\dot\DN_{\mu,h}-\dot\DN_{0,h})\gamma_0}_{H^1}
&\leq4c_\infty C_\Lambda
\bigl(1+12C_c(1+C_3(\omega))\bigr)
\abs{\mu}\nm{\Gamma}_{H^2}\nm{\gamma_0}_{H^2}\abs{h}.
\end{align*} 
By Lemma~\ref{lem:boundary-normalization}, both $\nm{D_v\gamma_v-\Gamma}_{H^2}$ and $\nm{\gamma_v-\gamma_0}_{H^1}$ are bounded by
$12c_\infty\bigl(\nm{\gamma_0}_{H^2}+\nm{\Gamma}_{H^2}\bigr)\nm{x}_X$.
Moreover, \eqref{eqn:quadratic-gamma-bound} implies
\begin{align*}
\nm{\gamma_v}_{H^1}+\nm{\Gamma}_{H^2}
&\leq(1+12c_\infty r)
\bigl(\nm{\gamma_0}_{H^2}+\nm{\Gamma}_{H^2}\bigr).
\end{align*}
Consequently, the estimates for the three summands of \eqref{eqn:Svmu-splitting}, together with $\abs{\mu}\leq\nm{x}_X$, yield
\begin{align}
\label{eqn:Svmu-estimate}
\nm[\big]{(S_{v\mu}(x)-S_{v\mu}(0))h}_{H^1}
&\leq4c_\infty C_\Lambda\biggl(
12c_\infty(1+12c_\infty r)
\bigl(\nm{\gamma_0}_{H^2}+\nm{\Gamma}_{H^2}\bigr)^2
\\\notag
&\hphantom{{}={}}
+\bigl(1+12C_c(1+C_3(\omega))\bigr)
\nm{\Gamma}_{H^2}\nm{\gamma_0}_{H^2}
\biggr)\nm{x}_X\abs{h}.
\end{align}
It remains to estimate $S_{\mu\mu}(x)-S_{\mu\mu}(0)$.
For $h,k\in\R^{3b-6}$, the terms in
$\sk{(S_{\mu\mu}(x)-S_{\mu\mu}(0))h,k}$
involving derivatives of $\tilde\alpha$ are, by Proposition~\ref{prop:linearisation-general},
\begin{align}
\label{eqn:20261114c}
\frac{1}{2}\int_\Xi\biggl(
&\Bigl(\frac{d}{dt}\Big\vert_{t=0}\tilde\alpha_{\mu+th,k}\Bigr)
\bigl(\abs{\gamma_v'}^2-\abs{\DN_{c_\mu}\gamma_v}^2\bigr)
-\Bigl(\frac{d}{dt}\Big\vert_{t=0}\tilde\alpha_{th,k}\Bigr)
\bigl(\abs{\gamma_0'}^2-\abs{\DN_{c_0}\gamma_0}^2\bigr)
\biggr)\,d\theta.
\end{align}
We estimate the terms in \eqref{eqn:20261114c} involving $\abs{\gamma_v'}^2$ and $\abs{\DN_{c_\mu}\gamma_v}^2$ separately.
\begin{align}\notag
&\frac{1}{2}\abs[\bigg]{\int_\Xi\biggl(
\Bigl(\frac{d}{dt}\Big\vert_{t=0}\tilde\alpha_{\mu+th,k}\Bigr)\abs{\gamma_v'}^2
-\Bigl(\frac{d}{dt}\Big\vert_{t=0}\tilde\alpha_{th,k}\Bigr)\abs{\gamma_0'}^2
\biggr)\,d\theta}
\\\notag
&=\frac{1}{2}\abs[\bigg]{\int_\Xi\biggl(
\Bigl(\frac{d}{dt}\Big\vert_{t=0}
\bigl(\tilde\alpha_{\mu+th,k}-\tilde\alpha_{th,k}\bigr)\Bigr)\abs{\gamma_v'}^2
+\Bigl(\frac{d}{dt}\Big\vert_{t=0}\tilde\alpha_{th,k}\Bigr)
\sk{\gamma_v'-\gamma_0',\gamma_v'+\gamma_0'}
\biggr)\,d\theta}
\\\notag
&\leq\frac{C_c}{2}\biggl(
\abs{\mu}\nm{\gamma_v'}_{L^2}^2
+\nm{\gamma_v'-\gamma_0'}_{L^2}
\bigl(\nm{\gamma_v'}_{L^2}+\nm{\gamma_0'}_{L^2}\bigr)
\biggr)\abs{h}\abs{k}
\\\label{eqn:20261114a}
&\leq\frac{C_c}{2}\biggl(
\nm{\gamma_v}_{H^2}^2
+12c_\infty\bigl(\nm{\gamma_0}_{H^2}+\nm{\Gamma}_{H^2}\bigr)
\bigl(\nm{\gamma_v}_{H^2}+\nm{\gamma_0}_{H^2}\bigr)
\biggr)\nm{x}_X\abs{h}\abs{k}.
\end{align}
The first inequality follows from Lemma~\ref{lem:configuration-bounds} and Cauchy--Schwarz, while the second uses Lemma~\ref{lem:boundary-normalization} and $\abs{\mu}\leq\nm{x}_X$.
By \eqref{eqn:quadratic-gamma-bound}, the coefficient in \eqref{eqn:20261114a} is uniformly bounded in terms of the recorded constants.
For the terms in \eqref{eqn:20261114c} involving $\abs{\DN_{c_\mu}\gamma_v}^2$, we similarly obtain
\begin{align}\notag
&\frac{1}{2}\abs[\bigg]{\int_\Xi\biggl(
\Bigl(\frac{d}{dt}\Big\vert_{t=0}\tilde\alpha_{\mu+th,k}\Bigr)
\abs{\DN_{c_\mu}\gamma_v}^2
-\Bigl(\frac{d}{dt}\Big\vert_{t=0}\tilde\alpha_{th,k}\Bigr)
\abs{\DN_{c_0}\gamma_0}^2
\biggr)\,d\theta}
\\\notag
&\leq\frac{1}{2}\abs[\bigg]{\int_\Xi
\Bigl(\frac{d}{dt}\Big\vert_{t=0}
\bigl(\tilde\alpha_{\mu+th,k}-\tilde\alpha_{th,k}\bigr)\Bigr)
\abs{\DN_{c_\mu}\gamma_v}^2\,d\theta}
\\\notag
&\hphantom{{}={}}
+\frac{1}{2}\abs[\bigg]{\int_\Xi
\Bigl(\frac{d}{dt}\Big\vert_{t=0}\tilde\alpha_{th,k}\Bigr)
\sk{\DN_{c_\mu}\gamma_v-\DN_{c_0}\gamma_0,
\DN_{c_\mu}\gamma_v+\DN_{c_0}\gamma_0}\,d\theta}
\\\label{eqn:20261114b}
&\leq\frac{C_c}{2}\biggl(
\abs{\mu}\nm{\DN_{c_\mu}\gamma_v}_{L^2}^2
+\nm{\DN_{c_\mu}\gamma_v-\DN_{c_0}\gamma_0}_{L^2}
\bigl(\nm{\DN_{c_\mu}\gamma_v}_{L^2}
+\nm{\DN_{c_0}\gamma_0}_{L^2}\bigr)
\biggr)\abs{h}\abs{k}.
\end{align}
The first inequality follows by adding and subtracting an intermediate term, while the second uses Lemma~\ref{lem:configuration-bounds}.
By \eqref{eqn:DtN-differences-1}, \eqref{eqn:DtN-differences-3} and Lemma~\ref{lem:boundary-normalization}, we have
\begin{align}\notag
\nm{\DN_{c_\mu}\gamma_v}_{L^2}
&\leq\bigl(1+C_1(\omega)\bigr)\nm{\gamma_v}_{H^2},
\\\notag
\nm{\DN_{c_0}\gamma_0}_{L^2}
&\leq\bigl(1+C_1(\omega)\bigr)\nm{\gamma_0}_{H^2},
\\\label{eqn:quadratic-DtN-difference}
\nm{\DN_{c_\mu}\gamma_v-\DN_{c_0}\gamma_0}_{L^2}
&\leq
\Bigl(
12c_\infty\bigl(1+C_1(\omega)\bigr)
\bigl(\nm{\gamma_0}_{H^2}+\nm{\Gamma}_{H^2}\bigr)
+C_\Lambda\nm{\gamma_0}_{H^2}
\Bigr)\nm{x}_X.
\end{align}
Substituting these bounds into \eqref{eqn:20261114b}, using $\abs{\mu}\leq\nm{x}_X$ and \eqref{eqn:quadratic-gamma-bound}, bounds its right-hand side by an explicitly computable constant times $\nm{x}_X\abs{h}\abs{k}$.

The terms in $\sk{(S_{\mu\mu}(x)-S_{\mu\mu}(0))h,k}$ involving derivatives of $\tilde\beta$ satisfy
\begin{align}\notag
&\abs[\bigg]{\int_\Xi\biggl(
\Bigl(\frac{d}{dt}\Big\vert_{t=0}\tilde\beta_{\mu+th,k}\Bigr)
\sk{\gamma_v',\DN_{c_\mu}\gamma_v}
-\Bigl(\frac{d}{dt}\Big\vert_{t=0}\tilde\beta_{th,k}\Bigr)
\sk{\gamma_0',\DN_{c_0}\gamma_0}
\biggr)\,d\theta}
\\\notag
&=\abs[\bigg]{\int_\Xi\biggl(
\Bigl(\frac{d}{dt}\Big\vert_{t=0}
\bigl(\tilde\beta_{\mu+th,k}-\tilde\beta_{th,k}\bigr)\Bigr)
\sk{\gamma_v',\DN_{c_\mu}\gamma_v}
\\\notag
&\hphantom{{}={}}
+\Bigl(\frac{d}{dt}\Big\vert_{t=0}\tilde\beta_{th,k}\Bigr)
\bigl(
\sk{\gamma_v'-\gamma_0',\DN_{c_\mu}\gamma_v}
+\sk{\gamma_0',\DN_{c_\mu}\gamma_v-\DN_{c_0}\gamma_0}
\bigr)
\biggr)\,d\theta}
\\\label{eqn:20261114d}
&\leq C_c\biggl(
\abs{\mu}\nm{\gamma_v'}_{L^2}\nm{\DN_{c_\mu}\gamma_v}_{L^2}
+\nm{\gamma_v'-\gamma_0'}_{L^2}\nm{\DN_{c_\mu}\gamma_v}_{L^2}
+\nm{\gamma_0'}_{L^2}\nm{\DN_{c_\mu}\gamma_v-\DN_{c_0}\gamma_0}_{L^2}
\biggr)\abs{h}\abs{k}
\\\label{eqn:20261114e}
&\leq
C_c\biggl(
\bigl(1+C_1(\omega)\bigr)\nm{\gamma_v}_{H^2}^2
+12c_\infty\bigl(1+C_1(\omega)\bigr)
\bigl(\nm{\gamma_0}_{H^2}+\nm{\Gamma}_{H^2}\bigr)\nm{\gamma_v}_{H^2}
\\[-1ex]&\hphantom{{}={}}\notag
+\nm{\gamma_0}_{H^2}
\Bigl(
12c_\infty\bigl(1+C_1(\omega)\bigr)
\bigl(\nm{\gamma_0}_{H^2}+\nm{\Gamma}_{H^2}\bigr)
+C_\Lambda\nm{\gamma_0}_{H^2}
\Bigr)\biggr)
\nm{x}_X\abs{h}\abs{k}.
\end{align} 
Inequality \eqref{eqn:20261114d} follows from Lemma~\ref{lem:configuration-bounds} and Cauchy--Schwarz, while \eqref{eqn:20261114e} follows from Lemma~\ref{lem:boundary-normalization}, \eqref{eqn:DtN-differences-1} and \eqref{eqn:quadratic-DtN-difference}.
By \eqref{eqn:quadratic-gamma-bound}, the coefficient in \eqref{eqn:20261114e} is uniformly bounded in terms of the recorded constants. 

For the remaining terms in $\sk{(S_{\mu\mu}(x)-S_{\mu\mu}(0))h,k}$, we first estimate the directional derivatives of the Dirichlet-to-Neumann operator.
By \eqref{eqn:DtN-differences-2}, \eqref{eqn:DtN-differences-4} and Lemma~\ref{lem:boundary-normalization}, we have
\begin{align}
\notag
\nm[\big]{\dot\DN_{\mu,h}\gamma_v}_{L^2}
&\leq C_\Lambda\nm{\gamma_v}_{H^1}\abs{h}
\leq C_\Lambda\nm{\gamma_v}_{H^2}\abs{h},
\\[1ex]\notag
\nm[\big]{\dot\DN_{\mu,h}\gamma_v-\dot\DN_{0,h}\gamma_0}_{L^2}
&\leq\nm[\big]{\dot\DN_{\mu,h}(\gamma_v-\gamma_0)}_{L^2}
+\nm[\big]{(\dot\DN_{\mu,h}-\dot\DN_{0,h})\gamma_0}_{L^2}
\\ 
&\leq C_\Lambda\biggl(
\nm{\gamma_v-\gamma_0}_{H^1}
+\Bigl(1+12C_c\bigl(1+C_3(\omega)\bigr)\Bigr)
\abs{\mu}\nm{\gamma_0}_{H^2}
\biggr)\abs{h}
\label{eqn:quadratic-DtN-derivative-bounds}\\\notag
&\leq C_\Lambda\biggl(
12c_\infty\bigl(\nm{\gamma_0}_{H^2}+\nm{\Gamma}_{H^2}\bigr)
+\Bigl(1+12C_c\bigl(1+C_3(\omega)\bigr)\Bigr)
\nm{\gamma_0}_{H^2}
\biggr)\nm{x}_X\abs{h}.
\end{align}
By Proposition~\ref{prop:linearisation-general}, the remaining terms in $\sk{(S_{\mu\mu}(x)-S_{\mu\mu}(0))h,k}$ can be estimated by adding and subtracting intermediate terms in each scalar product, and by applying Lemma~\ref{lem:configuration-bounds} and the Cauchy--Schwarz inequality:
\begin{align}
\notag
&\abs[\bigg]{\int_\Xi\biggl(
\tilde\alpha_{\mu,k}\sk{\DN_{c_\mu}\gamma_v,\dot\DN_{\mu,h}\gamma_v}
+\tilde\beta_{\mu,k}\sk{\gamma_v',\dot\DN_{\mu,h}\gamma_v}
-\tilde\alpha_{0,k}\sk{\DN_{c_0}\gamma_0,\dot\DN_{0,h}\gamma_0}
-\tilde\beta_{0,k}\sk{\gamma_0',\dot\DN_{0,h}\gamma_0}
\biggr)\,d\theta}
\notag\\
&\leq C_c\abs{k}\biggl(
\Bigl(
\abs{\mu}\bigl(\nm{\DN_{c_\mu}\gamma_v}_{L^2}+\nm{\gamma_v'}_{L^2}\bigr)
+\nm{\DN_{c_\mu}\gamma_v-\DN_{c_0}\gamma_0}_{L^2}
+\nm{\gamma_v'-\gamma_0'}_{L^2}
\Bigr)
\nm[\big]{\dot\DN_{\mu,h}\gamma_v}_{L^2}
\label{eqn:quadratic-Smm-remaining}
\\\notag
&\hphantom{{}={}}
+\Bigl(
\nm{\DN_{c_0}\gamma_0}_{L^2}+\nm{\gamma_0'}_{L^2}
\Bigr)
\nm[\big]{\dot\DN_{\mu,h}\gamma_v-\dot\DN_{0,h}\gamma_0}_{L^2}
\biggr).
\end{align}
By \eqref{eqn:quadratic-gamma-bound}, \eqref{eqn:DtN-differences-1} and \eqref{eqn:quadratic-DtN-derivative-bounds}, the factors in \eqref{eqn:quadratic-Smm-remaining} involving $\gamma_v$, $\DN_{c_\mu}\gamma_v$ and $\dot\DN_{\mu,h}\gamma_v$ are uniformly bounded, the last one proportionally to $\abs{h}$.
The differences involving these quantities are controlled by Lemma~\ref{lem:boundary-normalization}, \eqref{eqn:quadratic-DtN-difference} and \eqref{eqn:quadratic-DtN-derivative-bounds}, each with a factor $\nm{x}_X$.
Hence, \eqref{eqn:quadratic-Smm-remaining} is bounded by an explicitly computable constant times $\nm{x}_X\abs{h}\abs{k}$.
Combining this with \eqref{eqn:20261114a}, \eqref{eqn:20261114b} and \eqref{eqn:20261114e}, and taking the supremum over $\abs{h}=\abs{k}=1$, we conclude that the operator norm $\nm{S_{\mu\mu}(x)-S_{\mu\mu}(0)}$ is bounded by an explicitly computable constant times $\nm{x}_X$.

By the block decomposition \eqref{eqn:nonlinear-blocks} and estimate \eqref{eqn:nonlinear-adjoint-bound} for the $L^2$-adjoint, we have
\begin{align*}
\nm[\big]{DF\vert_x-DF\vert_0}_{X\to Y}
&\leq
\nm[\big]{S_{v\mu}(x)-S_{v\mu}(0)}_{\R^{3b-6}\to H^1}
\\&\hphantom{{}={}}
+\max\Bigl\{
\nm[\big]{S_{vv}(x)-S_{vv}(0)}_{H^2\to H^1},
\nm{S_{\mu\mu}(x)-S_{\mu\mu}(0)}
\Bigr\}.
\end{align*}
Indeed, the norm of a symmetric $2\times2$ matrix with nonnegative entries is bounded by its maximum row sum.
The bounds \eqref{eqn:Svv-estimate1}--\eqref{eqn:Svv-estimate4}, \eqref{eqn:Svmu-estimate}, \eqref{eqn:20261114a}, \eqref{eqn:20261114b}, \eqref{eqn:20261114e} and \eqref{eqn:quadratic-Smm-remaining} imply that each operator norm on the right-hand side is bounded by $\nm{x}_X$ times an explicitly computable coefficient.
The companion code computes these coefficients and records as $K$ an upper bound for $r$ times the mixed-block coefficient plus $r$ times the larger of the two diagonal-block coefficients.
Since $\nm{x}_X\leq r$, we conclude that
\begin{align*}
\nm[\big]{D\mathscr N\vert_x}_{X\to Y}
&=\nm[\big]{DF\vert_x-DF\vert_0}_{X\to Y}
\leq K,
\end{align*}
which proves the claim.
\end{proof}

\subsection{Fix point argument}
 
We apply Proposition~\ref{prop:existence} to the approximate solution from Section~\ref{sec:initial}.

\begin{corollary}
\label{cor:existence}
Let $r>0$ be as in Proposition~\ref{prop:nonlinear}.
Then there exist $v\in H^2(\Xi,\R^2)$ and $\mu\in B$ such that
\begin{align*}
F(v,\mu,0)&=0,
&
\nm{v}_{H^2}^2+\abs{\mu}^2&\leq r^2.
\end{align*}
The harmonic extension $\Phi_{\gamma_v,c_\mu}\colon\Omega_{c_\mu}\to\B^3$ of the boundary map $\gamma_v$ is a branched free boundary minimal immersion which is smooth up to the boundary and has no branch points on $\partial\Omega_{c_\mu}$.
\end{corollary}

\begin{proof}
By Propositions~\ref{prop:residual-bound}, \ref{prop:linear-estimate} and~\ref{prop:nonlinear}, the linearised operator $\mathscr{L}=DF\vert_0\colon X\to Y$ is an isomorphism and \eqref{eqn:assumptions} holds for the respective constants $\varepsilon,q,K$. 
The numerical values recorded in Appendix~\ref{section-numerical-values} satisfy
\begin{align}
\label{eqn:numerical_conditions}
0<r&<r_*, &
q\varepsilon+qKr&\leq r, &
qK&<1, &
q_0r&<1.
\end{align}
Hence, $\{\mu\in\R^{3b-6}\mid\abs{\mu}\leq r\}\subset B$ and \eqref{eqn:radii-conditions} holds.
In particular, assumption~\ref{prop:existence-a} of Proposition~\ref{prop:existence} is satisfied.

To verify assumption~\ref{prop:existence-b}, let $v\in H^2(\Xi,\R^2)$ with $\nm{v}_{H^2}\leq r$.
By Proposition~\ref{prop:linearisation}, we have $\mathscr{L}(0,0,a)=(\rho a,0,0)$ for every $a\in\R^3$.
Hence, Proposition~\ref{prop:L2-inverse-bound} implies $\abs{a}\leq q_0\nm{\rho a}_{L^2}$.
Recalling the definition of $\mathcal{R}^v$ from Lemma~\ref{lem:zeros-F} and using \eqref{eqn:rotation-fields-v}, we obtain
\begin{align*}
a^{\intercal}\mathcal{R}^v a
&=\sum_{i,j=1}^3 a_i a_j
\int_\Xi\sk{\rho_j,\rho_i^v}\,d\theta
\\
&=\int_\Xi
\sk[\Big]{\rho a,
\rho a+\Gamma^{\intercal}(a\times\Gamma v)
-\sk{\gamma_0,a\times\Gamma v}\,v}\,d\theta
\\
&=\nm{\rho a}_{L^2}^2
+\int_\Xi
\Bigl(
\sk{\Gamma\rho a,a\times\Gamma v}
+\sk{\rho a,v}^2
\Bigr)\,d\theta
\\
&\geq\nm{\rho a}_{L^2}^2
-\abs{a}\nm{\rho a}_{L^2}\nm{v}_{L^2}
\geq(1-q_0r)\nm{\rho a}_{L^2}^2.
\end{align*}
Here we used the identity $\sk{\gamma_0,a\times\Gamma v}=-\sk{\rho a,v}$, the fact that $\Gamma^{\intercal}\Gamma=\mathrm{Id}$ and the Cauchy--Schwarz inequality. 
As $q_0r<1$ and $\nm{\rho a}_{L^2}>0$ for all $a\in\R^3\setminus\{0\}$, we have $a^{\intercal}\mathcal R^v a>0$ for all $a\in\R^3\setminus\{0\}$.
Hence, $\mathcal R^v$ is invertible and assumption~\ref{prop:existence-b} is satisfied. 
The claim then follows by applying Proposition~\ref{prop:existence}. 
\end{proof}

\section{Embeddedness and Morse index}
\label{sec:embeddedness} 

Let $v\in H^2(\Xi,\R^2)$ and $\mu\in B$ be as in Corollary~\ref{cor:existence}.
Then $(\gamma,c)\vcentcolon=(\gamma_v,c_\mu)\in H^2(\Xi,\Sp^2)\times\mathcal{C}$ is a critical point of the reduced energy \eqref{eqn:reduced-energy-DtN} and
\begin{align*}
\Phi=\Phi_{\gamma,c}\colon\Omega_c\subset\Sp^2\to\B^3
\end{align*}
is a branched free boundary minimal immersion which is smooth up to the boundary and has no branch points on $\partial\Omega_c$.
The goal of this section is to prove that $\Phi$ is actually an embedding and to determine the Morse index and nullity of the resulting free boundary minimal surface, as stated in Theorem~\ref{thm:main}.
The idea is to verify that the spectral index of $\Phi$ is equal to $1$, which implies embeddedness by \cite[Proposition 8.1]{FraserSchoen2016} and allows us to derive sharp upper bounds for the Morse index. 
We then derive the matching lower bound by analysing the Hessian of the Lagrangian $L$ from Remark~\ref{rem:Lagrange_multipliers}.
Throughout this section, we reuse the Fourier approximations and quantitative estimates employed in the proof of Corollary~\ref{cor:existence}.

\subsection{Spectral index and nullity}
Following \cite[Definition~2.5]{KarpukhinSchulz}, we consider the spectral quadratic form
\begin{align}
\label{eqn:spectral-form}
Q_\spec(f,f)&\vcentcolon=\int_{\Omega_c}\abs{df}^2
-\int_{\partial\Omega_c}\abs{\partial_\eta\Phi}f^2
\end{align}
for $f\in H^1(\Omega_c)$, where $\eta$ denotes the outward unit conormal along $\partial\Omega_c$ and all geometric quantities are taken with respect to the round metric on $\Omega_c\subset\Sp^2$.
The \emph{spectral index} $\ind_\spec(\Phi)$ is defined as the maximal dimension of a subspace of $H^1(\Omega_c)$ on which $Q_\spec$ is negative definite, and the \emph{spectral nullity} $\nul_\spec(\Phi)$ is the dimension of its kernel, i.\,e.~the space of functions $f\in H^1(\Omega_c)$ satisfying $Q_\spec(f,h)=0$ for all $h\in H^1(\Omega_c)$.

To compute the spectral index and nullity of $\Phi$, we reduce the problem to a quadratic form on $H^{\frac{1}{2}}(\Xi)$.
Recalling the diffeomorphism $\psi_c\colon\Xi\to\partial\Omega_c$ given by \eqref{eqn:20261009b} on each component of $\Xi$ and the Dirichlet-to-Neumann operator $\DN_c$ defined (for scalar functions) as in \eqref{eqn:DtN}, we consider the quadratic form
\begin{align}
\label{eqn:boundary_spectral_form}
Q_\spec^\partial(u,u)&\vcentcolon=\int_\Xi u\mathscr{S}u\,d\theta,
&
\mathscr{S}&\vcentcolon=\DN_c-\abs{\gamma'},
\end{align}
defined initially for smooth $u\colon\Xi\to\R$ and then extended continuously to $H^{\frac{1}{2}}(\Xi)$.

\begin{lemma}
\label{lem:spectral-index-boundary}
The spectral index $\ind_\spec(\Phi)$ is equal to the maximal dimension of a subspace of $H^{\frac{1}{2}}(\Xi)$ on which the quadratic form $Q_\spec^\partial$ defined in \eqref{eqn:boundary_spectral_form} is negative definite, and the spectral nullity $\nul_\spec(\Phi)$ is equal to the dimension of the kernel of $Q_\spec^\partial$.
\end{lemma}
\begin{proof}
Every $f\in H^1(\Omega_c)$ decomposes uniquely as $f=h+f_0$, where $h\in H^1(\Omega_c)$ is harmonic and $f_0\in H^1_0(\Omega_c)$.
By integration by parts and since $f_0$ has zero trace, we have
\begin{align*}
Q_\spec(h,f_0)&=0,
&
Q_\spec(f_0,f_0)&=\int_{\Omega_c}\abs{df_0}^2.
\end{align*}
Thus, the space of harmonic functions and $H^1_0(\Omega_c)$ are orthogonal with respect to $Q_\spec$, and $Q_\spec$ is positive definite on $H^1_0(\Omega_c)$.
Consequently, restricting $Q_\spec$ to the space of harmonic functions preserves both the spectral index and the spectral nullity.

The map $h\mapsto h\circ\psi_c$ defines an isomorphism from the space of harmonic functions in $H^1(\Omega_c)$ to $H^{\frac{1}{2}}(\Xi)$.
Since $\Phi$ is conformal, we have
\begin{align*}
\abs{\gamma_k'}&=\sin(r_k)\abs{\partial_\eta\Phi}\circ\psi_{c,k}
\end{align*}
for each $k\in\{1,\ldots,b\}$.
Hence, recalling Green's formula \eqref{eqn:DtN-Green}, we obtain for every harmonic function $h\in H^1(\Omega_c)$ with smooth boundary values $u=h\circ\psi_c$  
\begin{align*}
Q_\spec(h,h)
&=\int_\Xi u\DN_c u\,d\theta-\int_\Xi\abs{\gamma'}u^2\,d\theta
=Q_\spec^\partial(u,u).
\end{align*}
By density and continuity, this identity holds for all harmonic functions $h\in H^1(\Omega_c)$, which completes the proof.
\end{proof}

We necessarily have $\ind_\spec(\Phi)\geq1$ since $\DN_c1=0$ and 
\begin{align*}
Q_\spec^\partial(1,1)&=-\int_\Xi\abs{\gamma'}\,d\theta<0.
\end{align*}
By the free boundary condition \eqref{eqn:fbhm_eq}, conformality of $\Phi$ and the definition \eqref{eqn:DtN} of the Dirichlet-to-Neumann operator, we have
$\DN_c\gamma=\abs{\gamma'}\gamma$. 
Denoting the Euclidean components of $\gamma$ by $\gamma^1,\gamma^2,\gamma^3$, we thus have
$\mathscr{S}\gamma^i=0$ for each $i\in\{1,2,3\}$. 
To conclude $\nul_\spec(\Phi)\geq3$, we still need to verify that $\gamma^1,\gamma^2,\gamma^3$ are linearly independent. 
(At this stage, we have not ruled out the possibility that $\Phi$ is a branched covering of the equatorial disc, despite the absence of branch points along $\partial\Omega_c$.)

\begin{lemma}
\label{lem:spectral-nullity-lower-bound}
Let $a_{k,0},a_{k,n},b_{k,n}$ denote the coefficients of the Fourier polynomial $P$ from Section~\ref{sec:initial}, and let $\delta$ be as in Lemma~\ref{lem:normalization-trace-errors} with $\varphi=P$.
If the $3\times3$-matrix
\begin{align*}
\sum_{k=1}^b\biggl(
a_{k,0}a_{k,0}^{\intercal}
+\sum_{n=1}^M
\Bigl(a_{k,n}a_{k,n}^{\intercal}
+b_{k,n}b_{k,n}^{\intercal}\Bigr)
\biggr)-r^2(1+\delta)\mathrm{Id}
\end{align*}
is positive definite, then $\gamma^1,\gamma^2,\gamma^3$ are linearly independent and $\nul_\spec(\Phi)\geq3$.
\end{lemma}
\begin{proof}
Recalling that $\gamma_0=P/\abs{P}$, by \eqref{eqn:perturbation} we have
\begin{align*}
\abs{\gamma-\gamma_0}^2
&=2-\frac{2}{\sqrt{1+\abs{v}^2}}
\leq\abs{v}^2.
\end{align*}
Hence, by Corollary~\ref{cor:existence}, for every $a\in\R^3$ we have
\begin{align*}
\nm{\sk{a,\gamma}}_{L^2}
&\geq\nm{\sk{a,\gamma_0}}_{L^2}
-\abs{a}\nm{\gamma-\gamma_0}_{L^2}
\geq\frac{\nm{\sk{a,P}}_{L^2}}{\nm{P}_{L^\infty}}-r\abs{a}.
\end{align*}
By orthonormality of the Fourier basis \eqref{eq:realFourierBasis}, we have
\begin{align*}
\int_\Xi PP^{\intercal}\,d\theta
&=\sum_{k=1}^b\biggl(
a_{k,0}a_{k,0}^{\intercal}
+\sum_{n=1}^M
\Bigl(a_{k,n}a_{k,n}^{\intercal}
+b_{k,n}b_{k,n}^{\intercal}\Bigr)
\biggr).
\end{align*}
Moreover, as in the proof of Lemma~\ref{lem:normalization-trace-errors}, we have $\nm{P}_{L^\infty}^2\leq1+\delta$.
Hence, by assumption,
\begin{align*}
\nm{\sk{a,P}}_{L^2}^2
=\sk[\Big]{a,\Bigl(\int_\Xi PP^{\intercal}\,d\theta\Bigr)a}
>r^2\nm{P}_{L^\infty}^2\abs{a}^2
\end{align*}
for every $a\in\R^3\setminus\{0\}$. 
Consequently, $\nm{\sk{a,\gamma}}_{L^2}>0$ for every $a\in\R^3\setminus\{0\}$, which proves that $\gamma^1,\gamma^2,\gamma^3$ are linearly independent.
Since $\mathscr{S}\gamma^i=0$ for each $i\in\{1,2,3\}$ the claim $\nul_\spec(\Phi)\geq3$ follows from Lemma~\ref{lem:spectral-index-boundary}. 
\end{proof} 

To obtain an upper bound for $\ind_\spec(\Phi)+\nul_\spec(\Phi)$, we recall the Fourier truncation $\Pi_N$ from Definition~\ref{def:fourier-truncation}.
Given $N\in\N$, we consider the symmetric operator
\begin{align*}
\mathscr{S}_N&\vcentcolon=\Pi_N\mathscr{S}\Pi_N\colon\Pi_N L^2(\Xi)\to\Pi_N L^2(\Xi)
\end{align*}
and denote its eigenvalues in increasing order and counted with multiplicity by $\lambda_j(\mathscr{S}_N)$.

\begin{lemma}
\label{lem:spectral-index-upper-bound}
Suppose that there exist constants $a_N,b_N>0$ and $m_N\geq0$ such that
\begin{align*}
\lambda_5(\mathscr{S}_N)&\geq b_N,
&
Q_\spec^\partial(w,w)&\geq a_N\nm{w}_{L^2}^2,
&
\abs[\big]{Q_\spec^\partial(u,w)}&\leq m_N\nm{u}_{L^2}\nm{w}_{L^2}
\end{align*}
for all $u\in\Pi_N L^2(\Xi)$ and $w\in H^{\frac{1}{2}}(\Xi)$ with $\Pi_Nw=0$.
If $a_Nb_N>m_N^2$, then
\begin{align*}
\ind_\spec(\Phi)+\nul_\spec(\Phi)&\leq4.
\end{align*}
\end{lemma}
\begin{proof}
Let $W\subset\Pi_N L^2(\Xi)$ be the span of four orthonormal eigenfunctions of $\mathscr{S}_N$ corresponding to its four smallest eigenvalues.
By assumption, for every $u\in\Pi_N L^2(\Xi)\cap W^\perp$ and $w\in H^{\frac{1}{2}}(\Xi)$ with $\Pi_Nw=0$, we have
\begin{align}
\label{eqn:20261124-a}
Q_\spec^\partial(u+w,u+w)
&\geq b_N\nm{u}_{L^2}^2
-2m_N\nm{u}_{L^2}\nm{w}_{L^2}
+a_N\nm{w}_{L^2}^2.
\end{align}
Since $a_N,b_N>0$ and $a_Nb_N>m_N^2$, the right-hand side of \eqref{eqn:20261124-a} is positive unless $u=w=0$.
Thus, $Q_\spec^\partial$ is positive definite on the subspace
\begin{align*}
\bigl\{f\in H^{\frac{1}{2}}(\Xi)\mid\Pi_Nf\perp W\bigr\},
\end{align*}
which has codimension four in $H^{\frac{1}{2}}(\Xi)$.
Since the sum of the spectral index and nullity is the maximal dimension of a subspace on which $Q_\spec^\partial$ is nonpositive, the claim follows from Lemma~\ref{lem:spectral-index-boundary}.
\end{proof}

\paragraph{High-mode estimate.}
Recalling the common conformal collar width $\omega$ recorded in Appendix~\ref{section-numerical-values}, the definition \eqref{eqn:boundary_spectral_form} and the estimates \eqref{eqn:20261029d} and \eqref{eqn:20261029h} from the proof of Lemma~\ref{lem:decoupled-collar} yield
\begin{align}
\label{eqn:spectral-high-mode}
Q_\spec^\partial(w,w)
&\geq\bigl((N+1)\tanh(\omega(N+1))-\nm{\gamma'}_{L^\infty}\bigr)\nm{w}_{L^2}^2
\end{align}
for every $w\in H^{\frac{1}{2}}(\Xi)$ with $\Pi_Nw=0$. 
Indeed, since the function $x\mapsto x\tanh(\omega x)$ is increasing for $x\geq0$, estimate \eqref{eqn:spectral-high-mode} follows, first for Fourier polynomials, and then in general by density and continuity of $Q_\spec^\partial$. 
Recall from the proof of Proposition~\ref{prop:nonlinear} that $8c_\infty r\leq1$, so that Lemma~\ref{lem:boundary-normalization} applies by Corollary~\ref{cor:existence}.
Combining this lemma with \eqref{eqn:Sobolev-embedding}, we obtain
\begin{align}
\notag
\nm{\gamma'-P'}_{L^\infty}
&\leq
c_\infty\nm{\gamma_0-P}_{H^2}
+c_\infty\nm{\gamma-\gamma_0}_{H^2}
\\\label{eqn:spectral-weight-error}
&\leq
c_\infty\nm{\gamma_0-P}_{H^2}
+12c_\infty^2r\bigl(\nm{\gamma_0}_{H^2}+\nm{\Gamma}_{H^2}\bigr)
=\vcentcolon\varepsilon_\gamma.
\end{align}
The normalisation error $\nm{\gamma_0-P}_{H^2}$ in \eqref{eqn:spectral-weight-error} is bounded by Lemma~\ref{lem:normalization-trace-errors}, and we have
\begin{align*}
\nm{\gamma_0}_{H^2}&\leq\nm{P}_{H^2}+\nm{\gamma_0-P}_{H^2},
&
\nm{\Gamma}_{H^2}&\leq\nm{G}_{H^2}+\nm{\Gamma-G}_{H^2},
\end{align*}
where $G$ is the Fourier polynomial approximation of $\Gamma$ constructed using \eqref{eqn:frame-Fourier}.
Combining \eqref{eqn:spectral-high-mode} and \eqref{eqn:spectral-weight-error}, we may take
\begin{align}
\label{eqn:spectral-aN}
a_N&=(N+1)\tanh(\omega(N+1))-\nm{P'}_{L^\infty}-\varepsilon_\gamma
\end{align}
in Lemma~\ref{lem:spectral-index-upper-bound}, provided that \eqref{eqn:spectral-aN} is positive.

\paragraph{Low-mode estimate.}
Recalling the Fourier polynomial $\tilde{h}$ from \eqref{eqn:tilde_h}, we approximate $\abs{\gamma'}$ by the Fourier polynomial $\sk{P,\tilde{h}}$ already employed in \eqref{eqn:Fourier-Fourier-approximation}.
Since $\abs[\big]{\abs{\gamma'}-\abs{P'}}\leq\abs{\gamma'-P'}$, estimate \eqref{eqn:spectral-weight-error} and the triangle inequality yield
\begin{align}
\label{eqn:spectral-weight-approximation-error}
\nm[\big]{\abs{\gamma'}-\sk{P,\tilde{h}}}_{L^\infty}
&\leq\varepsilon_\gamma+\nm[\big]{\abs{P'}-\sk{P,\tilde{h}}}_{L^\infty}
=\vcentcolon\varepsilon_\spec.
\end{align}
The second summand in the definition \eqref{eqn:spectral-weight-approximation-error} of $\varepsilon_\spec$ can be bounded by interval arithmetic using the Fourier polynomials $P$ and $\tilde{h}$.

We approximate $\mathscr{S}_N$ by the symmetric operator
\begin{align}
\label{eqn:spectral-low-mode-approximation}
\mathscr{S}_N^0&\vcentcolon=\Pi_N(\DN_{c_0}-\sk{P,\tilde{h}})\Pi_N.
\end{align}
By \eqref{eqn:DtN-differences-3} and Corollary~\ref{cor:existence}, we have
$\nm{\Pi_N(\DN_c-\DN_{c_0})u}_{H^1}\leq C_\Lambda r\nm{u}_{H^1}$
for every $u\in\Pi_N L^2(\Xi)$.
Applying this estimate to the eigenfunctions of the operator $\Pi_N(\DN_c-\DN_{c_0})\Pi_N$, we conclude that the absolute values of all its eigenvalues are bounded by $C_\Lambda r$.
Since this operator is symmetric with respect to the $L^2$-scalar product, it follows that
\begin{align}
\label{eqn:spectral-DtN-error}
\nm[\big]{\Pi_N(\DN_c-\DN_{c_0})\Pi_N}_{L^2\to L^2}
&\leq C_\Lambda r.
\end{align}
By the definition \eqref{eqn:spectral-low-mode-approximation} of $\mathscr{S}_N^0$, estimates \eqref{eqn:spectral-weight-approximation-error} and \eqref{eqn:spectral-DtN-error} imply
$\nm{\mathscr{S}_N-\mathscr{S}_N^0}_{L^2\to L^2}\leq C_\Lambda r+\varepsilon_\spec$.
Consequently,  
\begin{align}
\label{eqn:spectral-low-mode-perturbation}
\lambda_5(\mathscr{S}_N)
&\geq\lambda_5(\mathscr{S}_N^0)-C_\Lambda r-\varepsilon_\spec.
\end{align}
It remains to obtain a lower bound for $\lambda_5(\mathscr{S}_N^0)$.
For $k\in\{1,\ldots,b\}$ and $s\in\{0,\ldots,2N\}$, let $u_{(2N+1)(k-1)+s+1}\colon\Xi\to\R$ be the function equal to the Fourier basis element $f_s$ defined in \eqref{eq:realFourierBasis} on the $k$-th component of $\Xi$ and zero elsewhere. 
Then $u_1,\ldots,u_{(2N+1)b}$ form an orthonormal basis of $\Pi_N L^2(\Xi)$.
For each $i\in\{1,\ldots,(2N+1)b\}$, let $\tilde{f}_i$ be the boundary trace of a harmonic function on $\Omega_{c_0}$ of the form \eqref{eqn:mfs-field}, chosen to approximate $u_i$.
We then define the Fourier polynomial
\begin{align}\label{eqn:spectral-DtN-approximation}
\tilde{h}_i
&\vcentcolon=\DN_{\mathrm{dec}}u_i
+\Pi_N(\DN_{c_0}-\DN_{\mathrm{dec}})\tilde{f}_i.
\end{align}
The Fourier coefficients of $\tilde{h}_i$ in \eqref{eqn:spectral-DtN-approximation} can be computed explicitly using Lemma~\ref{lem:decoupled-Fourier} and Proposition~\ref{prop:mfs-fourier-coefficients}.
Since $\Pi_Nu_i=u_i$ and $\DN_{\mathrm{dec}}$ commutes with $\Pi_N$, Corollary~\ref{cor:decoupled-explicit-bounds} yields
\begin{align}\label{eqn:spectral-DtN-column-error}
\nm{\Pi_N\DN_{c_0}u_i-\tilde{h}_i}_{L^2}
&=\nm[\Big]{\Pi_N(\DN_{c_0}-\DN_{\mathrm{dec}})(u_i-\tilde{f}_i)}_{L^2}
\leq C_0(\omega)\nm{u_i-\tilde{f}_i}_{L^2},
\end{align}
where the last norm is bounded by Proposition~\ref{prop:mfs-truncation-bounds}.
We define the symmetric matrix $\tilde{\mathscr{S}}_N$ with entries
\begin{align}\label{eqn:spectral-low-mode-matrix}
\tilde{s}_{ij}
&\vcentcolon=\int_\Xi
\frac{1}{2}\bigl(\tilde{h}_iu_j+\tilde{h}_ju_i\bigr)
-\sk{P,\tilde{h}}u_iu_j
\,d\theta.
\end{align}
The matrix entries \eqref{eqn:spectral-low-mode-matrix} can be computed using finite Fourier arithmetic. 
Using interval arithmetic, we enclose the entries of $\tilde{\mathscr{S}}_N$ and bound its fifth eigenvalue from below by computing a rigorous eigenvalue bound for the midpoint matrix and accounting for the interval error.  
Identifying $\tilde{\mathscr{S}}_N$ with the corresponding operator on $\Pi_N L^2(\Xi)$, we obtain for every $u=\sum_{i=1}^{(2N+1)b}t_i u_i\in\Pi_N L^2(\Xi)$
\begin{align*}
\abs[\Big]{\sk[\big]{(\mathscr{S}_N^0-\tilde{\mathscr{S}}_N)u,u}_{L^2}}
&\leq\nm[\bigg]{\sum_{i=1}^{(2N+1)b}t_i(\Pi_N\DN_{c_0}u_i-\tilde{h}_i)}_{L^2}\nm{u}_{L^2}.
\end{align*}
Since $\sum_{i=1}^{(2N+1)b}t_i^2=\nm{u}_{L^2}^2$, the Cauchy--Schwarz inequality and estimate \eqref{eqn:spectral-DtN-column-error} yield
\begin{align}\label{eqn:spectral-low-mode-error}
\nm[\big]{\mathscr{S}_N^0-\tilde{\mathscr{S}}_N}_{L^2\to L^2}
&\leq C_0(\omega)
\biggl(\sum_{i=1}^{(2N+1)b}\nm{u_i-\tilde{f}_i}_{L^2}^2\biggr)^{\frac{1}{2}}
=\vcentcolon\varepsilon_N,
\end{align}
where we also used the symmetry of the operators.
The constant $\varepsilon_N$ in \eqref{eqn:spectral-low-mode-error} is bounded using Proposition~\ref{prop:mfs-truncation-bounds}.
Combining \eqref{eqn:spectral-low-mode-perturbation} and \eqref{eqn:spectral-low-mode-error}, we obtain
\begin{align}
\label{eqn:spectral-low-mode-bound}
\lambda_5(\mathscr{S}_N)
&\geq\lambda_5(\tilde{\mathscr{S}}_N)-\varepsilon_N-C_\Lambda r-\varepsilon_\spec.
\end{align}
Thus, a lower bound for $\lambda_5(\tilde{\mathscr{S}}_N)$ and upper bounds for the error terms in \eqref{eqn:spectral-low-mode-bound} yield a lower bound $b_N$ for $\lambda_5(\mathscr{S}_N)$ as required in Lemma~\ref{lem:spectral-index-upper-bound}, provided that $b_N>0$.

\paragraph{Coupling estimate.}
By the Fourier multiplier formula in Lemma~\ref{lem:decoupled-Fourier}, we have
$(\mathrm{Id}-\Pi_N)\DN_{\mathrm{dec}}\Pi_N=0$.
Hence, recalling the factorisation \eqref{eqn:decoupled-factorization} with $T=e^{-\omega\DN_{\mathrm{dec}}/2}$ and $\nm{S_c}_{L^2\to L^2}\leq1$, we obtain
\begin{align}
\label{eqn:spectral-DtN-coupling}
\nm[\big]{(\mathrm{Id}-\Pi_N)\DN_c\Pi_N}_{L^2\to L^2}
&=\frac{1}{\omega}\nm[\big]{(\mathrm{Id}-\Pi_N)TS_cT\Pi_N}_{L^2\to L^2}
\leq\frac{1}{\omega}e^{-\omega(N+1)/2}.
\end{align}
Combining \eqref{eqn:spectral-DtN-coupling} and \eqref{eqn:spectral-weight-approximation-error}, we obtain
\begin{align}
\label{eqn:spectral-mN}
\nm[\big]{(\mathrm{Id}-\Pi_N)\mathscr{S}\Pi_N}_{L^2\to L^2}
&\leq\frac{1}{\omega}e^{-\omega(N+1)/2}
+\nm[\big]{(\mathrm{Id}-\Pi_N)\sk{P,\tilde{h}}\Pi_N}_{L^2\to L^2}
+\varepsilon_\spec.
\end{align}
Since $P$ and $\tilde{h}$ are Fourier polynomials, the operator $(\mathrm{Id}-\Pi_N)\sk{P,\tilde{h}}\Pi_N$ is represented by a finite matrix whose columns are given by the coefficients of the Fourier polynomials
$(\mathrm{Id}-\Pi_N)(\sk{P,\tilde{h}}u_i)$
with respect to the orthonormal Fourier basis \eqref{eq:realFourierBasis} on each component of $\Xi$.
Since
\begin{align*}
Q_\spec^\partial(u,w)
&=\sk[\big]{(\mathrm{Id}-\Pi_N)\mathscr{S}\Pi_Nu,w}_{L^2}
\end{align*}
for $u\in\Pi_N L^2(\Xi)$ and $w\in H^{\frac{1}{2}}(\Xi)$ with $\Pi_Nw=0$, an upper bound for the right-hand side of \eqref{eqn:spectral-mN} yields a constant $m_N$ as required in Lemma~\ref{lem:spectral-index-upper-bound}.

In summary, if the assumption of Lemma~\ref{lem:spectral-nullity-lower-bound} is satisfied and the estimates \eqref{eqn:spectral-aN}, \eqref{eqn:spectral-low-mode-bound} and \eqref{eqn:spectral-mN} yield constants $a_N,b_N>0$ and $m_N\geq0$ satisfying $a_Nb_N>m_N^2$, then Lemmata~\ref{lem:spectral-nullity-lower-bound} and~\ref{lem:spectral-index-upper-bound} imply
\begin{align}
\label{eqn:spectral-index-nullity}
\ind_\spec(\Phi)&=1,
&
\nul_\spec(\Phi)&=3.
\end{align}
To deduce the absence of branch points from \eqref{eqn:spectral-index-nullity}, we prove the following statement.

\begin{lemma}[Embeddedness]
\label{lem:embeddedness}
Let $(\Omega,g_0)$ be a compact, connected, oriented Riemannian surface of genus zero with $2\leq b\in\N$ boundary components, and let $\Phi\colon\Omega\to\B^3$ be a branched conformal free boundary minimal immersion which is smooth up to the boundary and has no branch points along $\partial\Omega$.
If $\ind_\spec(\Phi)=1$, then $\nul_\spec(\Phi)=3$ and $\Phi$ has no branch points and is an embedding.
\end{lemma}

\begin{proof}
Since $\Phi$ has no branch points along $\partial\Omega$, we may choose a smooth function $\rho\colon\Omega\to\interval{0,\infty}$ satisfying $\rho=\abs{\partial_\eta\Phi}$ along $\partial\Omega$ and consider the conformal metric $g=\rho^2g_{0}$ on $\Omega$.
By conformal invariance of harmonicity and the Dirichlet energy on surfaces, recalling the free boundary condition \eqref{eqn:fbhm_eq}, we have
\begin{align*}
&\left\{\begin{aligned}
\Delta_g\Phi&=0 &&\text{ in }\Omega,
\\
\partial_{\eta_g}\Phi&=\Phi &&\text{ on }\partial\Omega 
\end{aligned}\right.
&
Q_\spec(f,f)&=\int_{(\Omega,g)}\abs{df}_g^2
-\int_{(\partial\Omega,g)}f^2,
\end{align*}
where $\eta_g$ denotes the outward unit conormal with respect to $g$.
Since $Q_\spec(1,1)<0$ and
\begin{align*}
Q_\spec(1,f)&=-\int_{(\partial\Omega,g)}f,
\end{align*}
any function $f\in H^1(\Omega)$ with zero boundary mean and $Q_\spec(f,f)<0$ would give rise to a two-dimensional subspace spanned by $1$ and $f$ on which $Q_\spec$ is negative definite.
Since $\ind_\spec(\Phi)=1$, we thus have
\begin{align*}
\int_{(\Omega,g)}\abs{df}_g^2&\geq\int_{(\partial\Omega,g)}f^2
\end{align*}
for every $f\in H^1(\Omega)$ with zero boundary mean with respect to $g$.
Hence, by the variational characterisation of the first nonzero Steklov eigenvalue, the latter is at least $1$, and hence equal to $1$, since the components of $\Phi$ are Steklov eigenfunctions with eigenvalue $1$. 

Suppose that $p\in\Omega$ is an interior branch point of $\Phi$.
Since the coordinate functions of $\Phi$ span a space of dimension at least two, we may choose $a\in\R^3$ such that $f\vcentcolon=\sk{a,\Phi}$ does not vanish identically but satisfies $f(p)=0$.
Moreover, since $d\Phi(p)=0$, we have $df(p)=0$, so that $f$ vanishes to order at least two at $p$.
Since $\Omega$ has genus zero, it follows from \cite[Lemma~3.2]{KarpukhinKokarevPolterovich2014} that $f$ has at least three nodal domains, contradicting Courant's nodal domain theorem for first Steklov eigenfunctions (see \cite[§\,2.1]{KarpukhinKokarevPolterovich2014}).

If the image of $\Phi$ were contained in a plane through the origin, then $\Phi$ would be an unbranched covering of an equatorial disc and hence a diffeomorphism, contradicting $b\geq2$.
Thus, the coordinate functions of $\Phi$ are linearly independent, so $\nul_\spec(\Phi)=3$ by \cite[Corollary~1.3]{KarpukhinKokarevPolterovich2014}.
Moreover, embeddedness follows from \cite[Proposition~8.1]{FraserSchoen2016}.
\end{proof}

\subsection{Sharp Morse index bounds}
\label{sec:upper_index_bound}

Given a compact, oriented Riemannian surface $(\Omega,g)$ with boundary and a conformal free boundary minimal immersion $\Phi\colon\Omega\to\B^3$, we denote by $\ind_\area(\Phi)$ and $\nul_\area(\Phi)$ its Morse index and nullity for the area functional with respect to normal variations. 
Following \cite[§\,2.1]{KarpukhinSchulz}, the second variation of the Dirichlet energy \eqref{eqn:energy} at $\Phi$ is given by the quadratic form
\begin{align}
\label{eqn:energy-spectral-form}
Q_\ener(V,V)&=\sum_{i=1}^3Q_\spec(V^i,V^i),
\end{align}
on the space $\mathcal{V}_\Phi$ of smooth maps $V\colon\Omega\to\R^3$ satisfying the constraint $\sk{V,\Phi}=0$ along $\partial\Omega$, where $Q_\spec$ is defined as in \eqref{eqn:spectral-form} with respect to the metric $g$.
The \emph{energy index} $\ind_\ener(\Phi)$ and the \emph{energy nullity} $\nul_\ener(\Phi)$ are defined as the maximal dimension of a subspace of $\mathcal{V}_\Phi$ on which $Q_\ener$ is negative definite and the dimension of the kernel of $Q_\ener$, respectively.
  
\begin{proposition}
\label{prop:Morse-index-upper-bound}
Let $(\Omega,g)$ be a compact, connected, oriented Riemannian surface of genus zero with $3\leq b\in\N$ boundary components, and let $\Phi\colon\Omega\to\B^3$ be a conformal free boundary minimal immersion.
If $\ind_\spec(\Phi)=1$, then
\begin{align*}
\ind_\area(\Phi)&\leq3b-3,
&
\ind_\area(\Phi)+\nul_\area(\Phi)&\leq3b.
\end{align*}
Moreover, if $\ind_\area(\Phi)=3b-3$, then $\nul_\area(\Phi)=3$.
\end{proposition}

\begin{proof}
By conformal invariance of the quadratic forms \eqref{eqn:spectral-form} and \eqref{eqn:energy-spectral-form}, we may assume that $g$ is the metric induced by $\Phi$.
Since $\ind_\spec(\Phi)=1$ by assumption, we have $\ind_\ener(\Phi)\leq3$ by \cite[Corollary~2.17]{KarpukhinSchulz}, applied with the trivial symmetry group.
By \cite[Theorem~1]{Lima2022} (see also \cite[Theorems~3.8--9 and Lemma~A.5]{KarpukhinSchulz}), we then have
\begin{align*}
\ind_\area(\Phi)&\leq3b-6+\ind_\ener(\Phi)\leq3b-3.
\end{align*}
By Lemma~\ref{lem:embeddedness}, we have $\nul_\spec(\Phi)=3$ and $\Phi$ is an embedding. 
Its components are quadratically independent on $\partial\Omega$, since the zero set on $\Sp^2$ of a nonzero homogeneous quadratic polynomial cannot contain three pairwise disjoint smooth embedded closed curves; cf.~\cite[proof of Proposition~5.17]{Tran}.
The argument leading to \cite[(2.14)]{KarpukhinSchulz} uses only the spectral index, spectral nullity and quadratic independence of the components of $\Phi$.
Applying estimate \cite[(2.14)]{KarpukhinSchulz} (with $n=2$ therein), we obtain
\begin{align}
\label{eqn:energy_index+nullity}
\ind_\ener(\Phi)+\nul_\ener(\Phi)&\leq(n+1)+\tfrac{1}{2}n(n+1)=6.
\end{align}
(The additional assumption $\ind_\ener(\Phi)=n+1$ in \cite[Proposition~2.23]{KarpukhinSchulz} is used only afterwards to determine the energy nullity.)
Combining \eqref{eqn:energy_index+nullity} with \cite[Theorems~3.8--9 and Lemma~A.5]{KarpukhinSchulz} yields
\begin{align*}
\ind_\area(\Phi)+\nul_\area(\Phi)
&\leq\ind_\ener(\Phi)+\nul_\ener(\Phi)+3b-6
\leq3b.
\end{align*}
Finally, since its Euler characteristic is $\chi(\Omega)=2-b<0$, the surface $(\Omega,g)$ does not admit any nonzero Killing field tangent to $\partial\Omega$ by \cite[Remark~3.2]{KarpukhinSchulz}, because Killing fields are conformal.
Hence, the ambient rotations induce three linearly independent normal Jacobi fields along $\Phi$, so $\nul_\area(\Phi)\geq3$.
If $\ind_\area(\Phi)=3b-3$, the second inequality in the statement of the proposition implies $\nul_\area(\Phi)\leq3$, which proves the equality case.
\end{proof}

\begin{proof}[Proof of Theorem~\ref{thm:index-upper-bound}]
Let $\Phi\colon\Sigma\to\B^3$ denote the inclusion map.
By \cite[Theorem~1.1]{ChodoshGianocca}, the first Steklov eigenvalue of $\Sigma$ is equal to $1$, which implies $\ind_\spec(\Phi)=1$ (see \cite[Remark~2.7]{KarpukhinSchulz}).
The claim thus follows from Proposition~\ref{prop:Morse-index-upper-bound}.
\end{proof}

To obtain a matching lower bound for the Morse index of our surface, we consider the reduced energy from Definition~\ref{defn:reduced-energy} in the parametrisation \eqref{eqn:perturbation} as a functional
\begin{align}
\label{eqn:mathcalE}
\mathcal{E}\colon H^2(\Xi,\R^2)\times B&\to\R 
\\\notag
(v,\mu)&\mapsto E(\gamma_v,c_\mu)
=\widetilde{E}(\Phi_{\gamma_v,c_\mu},g_{\Sp^2}), 
\end{align}
where $\widetilde{E}$ is the two-variable energy defined in \eqref{eqn:2variable_energy}. 
By Corollary~\ref{cor:critical-slice} and Proposition~\ref{prop:critical-points}, every critical point $(v_*,\mu_*)$ of $\mathcal{E}$ with $\gamma_{v_*}$ nonconstant yields a branched free boundary minimal immersion $\Phi_*\vcentcolon=\Phi_{\gamma_{v_*},c_{\mu_*}}$.
If $\Phi_*$ has no branch points, its Morse index $\ind_\area(\Phi_*)$ can be bounded from below using the Hessian of $\mathcal{E}$.    
Given a quadratic form $Q$, we denote by $\ind(Q)$ the maximal dimension of a subspace on which $Q$ is negative definite.

\begin{proposition}
\label{prop:reduced-energy-index}
Let $(v_*,\mu_*)$ be a critical point of \eqref{eqn:mathcalE} and suppose that $\Phi_*\vcentcolon=\Phi_{\gamma_{v_*},c_{\mu_*}}$ has no branch points.
Denoting the Hessian of $\mathcal{E}$ at $(v_*,\mu_*)$ by
$Q_\mathcal{E}\vcentcolon=D^2\mathcal{E}\vert_{(v_*,\mu_*)}$, 
we have
\[
\ind_\area(\Phi_*)\geq\ind(Q_\mathcal{E}).
\]
\end{proposition}
 
\begin{proof}
By Proposition~\ref{prop:critical-points}, the map $\Phi_*\vcentcolon=\Phi_{\gamma_{v_*},c_{\mu_*}}\colon\Omega_{c_{\mu_*}}\to\B^3$ is a conformal free boundary minimal immersion which is smooth up to the boundary.
To compare the second variations of energy and area, we first pull back the maps and metrics to the fixed domain $\Omega_{c_*}$, where $c_*\vcentcolon=c_{\mu_*}$. 
Recalling the diffeomorphisms from \eqref{eqn:diffeos_f_c}, let 
\begin{align*}
\tilde{f}_c&\vcentcolon=f_c\circ f_{c_*}^{-1}\colon\Omega_{c_*}\to\Omega_c.
\end{align*}
For $(v,\mu)\in H^2(\Xi,\R^2)\times B$, we define the pullbacks
\begin{align*}
\tilde{\Phi}_{v,\mu}&\vcentcolon=\Phi_{\gamma_v,c_\mu}\circ\tilde{f}_{c_\mu},
&
\tilde{g}_\mu&\vcentcolon=\tilde{f}_{c_\mu}^*g_{\Sp^2}.
\end{align*}
Note that $\tilde{\Phi}_{v,\mu}\circ\psi_{c_*}=\gamma_v$ by \eqref{eqn:diffeos_compatibility}.
By invariance of the Dirichlet energy under pullback, we have
\begin{align*}
\mathcal{E}(v,\mu)
&=\widetilde{E}(\Phi_{\gamma_v,c_\mu},g_{\Sp^2})
=\widetilde{E}(\tilde{\Phi}_{v,\mu},\tilde{g}_\mu).
\end{align*}
Given $(w,h)\in H^2(\Xi,\R^2)\times\R^{3b-6}$, we consider the one-parameter families
\begin{align*}
v(t)&\vcentcolon=v_*+tw,
&
\mu(t)&\vcentcolon=\mu_*+th,
&
\Psi_t&\vcentcolon=\tilde{\Phi}_{v(t),\mu(t)},
&
g_t&\vcentcolon=\tilde{g}_{\mu(t)},
\end{align*}
so that at $t=0$ we have $\Psi_0=\Phi_*$, $g_0=g_{\Sp^2}$ and $\Psi_t(\partial\Omega_{c_*})\subset\Sp^2$.
By \eqref{eqn:20261020b} in the proof of Lemma~\ref{lem:smoothness}, the family $\Psi_t$ depends smoothly on $t$ in $H^{\frac{5}{2}}(\Omega_{c_*},\R^3)$ and hence in $C^1(\Omega_{c_*},\R^3)$ by Sobolev embedding.
Thus, $\Psi_t$ remains an immersion for sufficiently small $\abs{t}$, and $\operatorname{area}(\Psi_t)$ depends twice differentiably on $t$.
By definition of $Q_\mathcal{E}$ and invariance of the Dirichlet energy under pullback, we have
\begin{align}
\label{eqn:energyhessianrelation}
Q_\mathcal{E}((w,h),(w,h))
&=\frac{d^2}{dt^2}\Big\vert_{t=0}\mathcal{E}(v(t),\mu(t))
=\frac{d^2}{dt^2}\Big\vert_{t=0}\widetilde{E}(\Psi_t,g_t).
\end{align} 
By the energy--area inequality \cite[equation~(8)]{Lima2022}, we have
\begin{align}
\label{eqn:energy-area-inequality}
\widetilde{E}(\Psi_t,g_t)&\geq\operatorname{area}(\Psi_t),
\end{align}
with equality at $t=0$ since $\Phi_*$ is conformal with respect to $g_{\Sp^2}$.
Therefore, $\widetilde{E}(\Psi_t,g_t)-\operatorname{area}(\Psi_t)$ has a minimum at $t=0$.
Taking second derivatives in \eqref{eqn:energy-area-inequality} at $t=0$ and recalling \eqref{eqn:energyhessianrelation}, we obtain
\begin{align}
\label{eqn:energymorserelation}
Q_\mathcal{E}\bigl((w,h),(w,h)\bigr)
&\geq\frac{d^2}{dt^2}\Big\vert_{t=0}\operatorname{area}(\Psi_t).
\end{align}
Finally, the second variation formula for area (cf.~\cite[§\,3]{FraserSchoen2016}) extends by smooth approximation to the family $\Psi_t$ and yields
\begin{align*}
\frac{d^2}{dt^2}\Big\vert_{t=0}\operatorname{area}(\Psi_t)
&=Q_\area(u,u)
\vcentcolon=\int_{(\Omega_{c_*},g_*)}
\Bigl(\abs{du}_{g_*}^2-\abs{\II}_{g_*}^2u^2\Bigr)
-\int_{(\partial\Omega_{c_*},g_*)}u^2,
\end{align*} 
where $g_*\vcentcolon=(\Phi_*)^*g_{\R^3}$, where $\II$ denotes the second fundamental form of $\Phi_*$ and $u\vcentcolon=\sk{\dot{\Psi}_0,\nu}$ is the normal component of the variation field.
By \eqref{eqn:energymorserelation}, the linear map $(w,h)\mapsto u(w,h)$ is injective on any subspace $W$ on which $Q_\mathcal{E}$ is negative definite, and $Q_\area$ is negative definite on its image.
Hence, $\ind_\area(\Phi_*)\geq\ind(Q_\mathcal{E})$.
\end{proof}

We now relate the Hessian of $\mathcal{E}$ to that of the Lagrangian
$L\colon H^2(\Xi,\R^2)\times B\times\R^3\to\R$ from Remark~\ref{rem:Lagrange_multipliers}.
We also recall the map $F\colon H^2(\Xi,\R^2)\times B\times\R^3\to H^1(\Xi,\R^2)\times\R^{3b-6}\times\R^3$ from Definition~\ref{defn:gradient-map} and the matrix $\mathcal{R}^v\in\R^{3\times3}$ from Lemma~\ref{lem:zeros-F}.

\begin{proposition}
\label{prop:Lagrangian-index}
Let $(v_*,\mu_*,\lambda_*)$ be a zero of $F$ such that the matrix $\mathcal{R}^{v_*}$ is invertible and $\Phi_*\vcentcolon=\Phi_{\gamma_{v_*},c_{\mu_*}}$ has no branch points.
Denoting the Hessian of $L$ at $(v_*,\mu_*,\lambda_*)$ by
$Q_L\vcentcolon=D^2L\vert_{(v_*,\mu_*,\lambda_*)}$, 
and recalling the notation $Q_\mathcal{E}$ from Proposition~\ref{prop:reduced-energy-index}, we have
\begin{align*}
\ind(Q_L)-3&=\ind(Q_\mathcal{E})\leq\ind_\area(\Phi_*).
\end{align*}
\end{proposition}

\begin{proof}
By Lemma~\ref{lem:zeros-F}, we have $\lambda_*=0$ and $(v_*,\mu_*)$ is a critical point of $\mathcal{E}$.
The inequality thus follows from Proposition~\ref{prop:reduced-energy-index}.
To prove the equality, we set $Z\vcentcolon=H^2(\Xi,\R^2)\times\R^{3b-6}$ and recall the rotation fields $\rho_i$ of \eqref{eqn:rotation-fields}, the functions $\rho, \rho^*$ of Proposition~\ref{prop:linearisation}, and the perturbed rotation fields $\rho_i^{v}$ of \eqref{eqn:rotation-fields-v}. 
We define the perturbation directions corresponding to the rotation fields $\rho_i^{v_*}$, and the linear map $R\colon\R^3 \to Z$, by
\begin{align*}
r_i&\vcentcolon=(\rho_i^{v_*},0)\in Z,
&
R\colon\R^3&\to Z,
&
Ra&\vcentcolon=\sum_{i=1}^3a_ir_i.
\end{align*}
By rotational invariance of $\mathcal{E}$, recalling \eqref{eqn:energy-symmetries}, we have
$D\mathcal{E}\vert_{(v,\mu)}(\rho_i^v,0)=0$ 
for all $(v,\mu)\in H^2(\Xi,\R^2)\times B$ and $i\in\{1,2,3\}$.
Differentiating at the critical point $(v_*,\mu_*)$ yields
\begin{align}
\label{eq:qvanishesonrotations}
Q_\mathcal{E}(z,r_i)&=0
\end{align}
for every $z\in Z$ and $i\in\{1,2,3\}$.
We consider the linear map $C\colon Z\to\R^3$
given by $C(w,h)=\rho^*w$, and set $T\vcentcolon=\ker C$.
Since $CR=(\mathcal{R}^{v_*})^{\intercal}$ is invertible by assumption, every $z\in Z$ decomposes uniquely as $z=s+Ra$ with $s\in T$ and $a\in\R^3$.
By \eqref{eq:qvanishesonrotations}, we have $Q_\mathcal{E}(z,z)=Q_\mathcal{E}(s,s)$ and hence
\begin{align*}
\ind(Q_\mathcal{E})&=\ind(Q_\mathcal{E}|_T).
\end{align*}
Recalling the definition of the Lagrangian $L$ from Remark~\ref{rem:Lagrange_multipliers}, we have
\begin{align*}
L(v,\mu,\lambda)&=\mathcal{E}(v,\mu)+\sk{\lambda,\rho^*v}.
\end{align*}
Consequently, for $z=s+Ra\in Z$ and $\ell\in\R^3$, we obtain
\begin{align*}
Q_L((z,\ell),(z,\ell))
&=Q_\mathcal{E}(z,z)+2\sk{\ell,Cz}
=Q_\mathcal{E}(s,s)+2\sk{\ell,CRa}.
\end{align*}
Since $CR$ is invertible, the change of variables $(z,\ell)\mapsto(s,CRa,\ell)$ is a linear isomorphism $Z\times\R^3\to T\times\R^3\times\R^3$.
Writing $\xi\vcentcolon=CRa$, we obtain the decomposition
\begin{align*}
Q_L&=Q_\mathcal{E}|_T\oplus\bigl((\xi,\ell)\mapsto2\sk{\ell,\xi}\bigr).
\end{align*}
Since
$2\sk{\ell,\xi}=\tfrac{1}{2}\abs{\xi+\ell}^2-\tfrac{1}{2}\abs{\xi-\ell}^2$,
the second summand has index $3$.
By invariance of the index under linear isomorphisms and additivity under direct sums, we conclude that
\begin{align*}
\ind(Q_L)&=\ind(Q_\mathcal{E}|_T)+3=\ind(Q_\mathcal{E})+3.
\qedhere
\end{align*}
\end{proof}

It remains to obtain a lower bound for the index of the Hessian of $L$ at the zero $x\vcentcolon=(v,\mu,0)$ of $F$ obtained in Corollary~\ref{cor:existence}.
For this purpose, we use the low-mode block $\mathscr{A}$ of the linearised operator $\mathscr{L}$ from \eqref{eqn:LABC}.

\begin{lemma}
\label{lem:Lagrangian-index-lower-bound}
Let the constant $K$ be as in Proposition~\ref{prop:nonlinear} and $N\in\N$ as in the decomposition \eqref{eqn:orthogonal_decomposition}.
Let $\lambda_1(\mathscr{A})\leq\lambda_2(\mathscr{A})\leq\ldots$ denote the eigenvalues of $\mathscr{A}$, counted with multiplicity. 
If
\begin{align}
\label{eqn:Lagrangian-index-criterion}
\lambda_j(\mathscr{A})+K\sqrt{1+N^2}<0
\end{align}
for some $j\in\{1,\ldots,\dim\mathcal{H}_{\leq N}\}$, then $\ind(D^2L\vert_x)\geq j$.
\end{lemma}

\begin{proof}
Recalling Remark~\ref{rem:Lagrange_multipliers}, we have
\begin{align}
\label{eqn:Lagrangian-Hessians}
D^2L\vert_x(z,w)&=\sk{DF\vert_x z,w}_{\mathcal{H}},
&
D^2L\vert_0(z,w)&=\sk{\mathscr{L}z,w}_{\mathcal{H}}
\end{align}
for all $z,w\in X$.
In analogy with \eqref{eqn:LABC}, we define the symmetric operator
\begin{align*}
\mathscr{A}_x&\vcentcolon=\Pi_N DF\vert_x\Pi_N\colon\mathcal{H}_{\leq N}\to\mathcal{H}_{\leq N}.
\end{align*}
By Corollary~\ref{cor:existence} and estimate \eqref{eqn:nonlinear-derivative-bound}, we have
$\nm{\mathscr{A}_x-\mathscr{A}}_{X\to Y}\leq K$, 
where we used that $\Pi_N$ is an orthogonal projection with respect to the inner products in both $X$ and $Y$.
By the definition of the Sobolev norms in \eqref{eqn:Sobolev-norm}, we have
$\nm{z}_X\leq\sqrt{1+N^2}\nm{z}_Y$ 
for all $z\in\mathcal{H}_{\leq N}$.
Applying these estimates to the eigenvectors of $\mathscr{A}_x-\mathscr{A}$ and using its symmetry with respect to the $\mathcal{H}$-inner product, we obtain
\begin{align*}
\nm{\mathscr{A}_x-\mathscr{A}}_2&\leq K\sqrt{1+N^2}.
\end{align*}
Consequently, by \eqref{eqn:Lagrangian-index-criterion}, we have
$\lambda_j(\mathscr{A}_x)
\leq\lambda_j(\mathscr{A})+\nm{\mathscr{A}_x-\mathscr{A}}_2
<0$. 
Since $\mathcal{H}_{\leq N}\subset X$, equation \eqref{eqn:Lagrangian-Hessians} implies that $D^2L\vert_x$ is negative definite on the span of the first $j$ eigenvectors of $\mathscr{A}_x$.
\end{proof}

In summary, the embeddedness of $\Phi$ follows from \eqref{eqn:spectral-index-nullity} by Lemma~\ref{lem:embeddedness}.
The upper bounds for its Morse index and nullity follow from Proposition~\ref{prop:Morse-index-upper-bound}, while a matching lower bound for the Morse index follows from Proposition~\ref{prop:Lagrangian-index} and Lemma~\ref{lem:Lagrangian-index-lower-bound}, provided that \eqref{eqn:Lagrangian-index-criterion} holds for $j=3b$.

\section{Asymmetry and area}
\label{sec:asymmetry}

We recall from Section~\ref{sec:initial} that the approximate solution $(\gamma_0,c_0)$ is given by a circle configuration $c_0=(p_1,\ldots,p_b;r_1,\ldots,r_b)\in\mathcal{C}$ with $b=61$ and a Fourier polynomial $P\in\mathcal P_M$ defining $\gamma_0=P/\abs{P}$.
Here $p_k\in\Sp^2$ and $r_k\in\interval{0,\pi}$ are the centre and geodesic radius of the $k$-th disc defining $\Omega_{c_0}\subset\Sp^2$.
The restriction $P_k$ of $P$ to the $k$-th component of $\Xi$ is an $\R^3$-valued Fourier polynomial of degree at most $M$ given by \eqref{eq:PkFourierCoeffs}.
Given $v\in H^2(\Xi,\R^2)$ and $\mu\in B$ as in Corollary~\ref{cor:existence} we recall from \eqref{eqn:perturbation}
and \eqref{eqn:perturbation-configuration} that $\gamma_v\colon\Xi\to\Sp^2$ is the boundary map and $c_\mu$ is the circle configuration of the free boundary minimal immersion $\Phi_{\gamma_v,c_\mu}\colon\Omega_{c_\mu}\to\B^3$.
Lemma~\ref{lem:embeddedness} implies that $\Phi_{\gamma_v, c_\mu}$ is an embedding and that $\gamma_v$ parametrises $b$ pairwise disjoint simple closed curves in $\Sp^2$, provided that \eqref{eqn:spectral-index-nullity} holds.  

We aim to prove that $\gamma_v(\Xi)$ has no nontrivial symmetries in $\Ogroup(3)$, which implies the same for $\Phi_{\gamma_v,c_\mu}(\Omega_{c_\mu})$.
To do this, we associate with each boundary component $\gamma_{v,k}$ the area vector of the enclosed spherical disc.
Any symmetry of $\gamma_v(\Xi)$ permutes these vectors.
We then show that three of these vectors are linearly independent and that the norm of each of them differs from the norms of all the other vectors, to conclude that the only symmetry of $\gamma_v(\Xi)$ is the identity.

\paragraph{The area vector.}
Given $\gamma\in H^1(\Sp^1,\R^3)$, we define its \emph{area vector} by 
\begin{align}\label{eqn:area-vector}
A(\gamma)\vcentcolon=\frac{1}{2}\int_0^{2\pi}\gamma(\theta)\times\gamma'(\theta)\,d\theta.
\end{align}
If $\gamma$ parametrises a simple closed curve in $\Sp^2$ contained in an open hemisphere, we denote by $D_\gamma$ the spherical disc bounded by $\gamma$ and contained in that hemisphere.
By Stokes' theorem, we have 
$A(\gamma)=\pm\int_{D_\gamma}x$, where the sign depends on the orientation of $\gamma$. 
For every $R\in\Ogroup(3)$, we have $D_{R\circ\gamma}=R(D_\gamma)$ and hence
\begin{align*}
\int_{D_{R\circ\gamma}}x=R\int_{D_\gamma}x.
\end{align*}
Using \eqref{eq:PkFourierCoeffs}, we see that 
\begin{align*}
A(P_k)
&=\frac{1}{2}\int_0^{2\pi}P_k\times P_k'\,d\theta
\\
&=\frac{1}{2}\int_0^{2\pi}
\biggl(\frac{a_{k,0}}{\sqrt{2\pi}}+\frac{1}{\sqrt{\pi}}\sum_{n=1}^M\bigl(a_{k,n}\cos(n\theta)+b_{k,n}\sin(n\theta)\bigr)\biggr)
\\
&\phantom{{}={}}\times
\biggl(\frac{1}{\sqrt{\pi}}\sum_{n=1}^M\bigl(-na_{k,n}\sin(n\theta)+nb_{k,n}\cos(n\theta)\bigr)\biggr)\,d\theta
\\
&=\frac{1}{2\pi}\sum_{n=1}^M n\int_0^{2\pi}a_{k,n}\times b_{k,n}\,d\theta
=\sum_{n=1}^M n\,a_{k,n}\times b_{k,n}.
\end{align*}
This yields precise values for $A(P_k)$ in terms of the Fourier coefficients $a_{k,n}$, $b_{k,n}$ without having to perform any numerical integration.

\paragraph{Norm uniqueness check.} 
We first bound the difference of $A(\gamma_{v,k})$ and $A(P_k)$, and then use this estimate to show that three of the vectors $A(\gamma_{v,k})$ have unique norm. 

\begin{lemma}\label{lem:AsymNormDifferenceBound-B}
For all $k\in\{1,\ldots,b\}$, we have
\[
\abs[\big]{A(\gamma_{v,k})-A(P_k)}
\leq \pi\bigl(1+\nm{P'}_{L^\infty(\Xi)}\bigr)\varepsilon_\gamma
< c_{d},  
\] 
where $\varepsilon_\gamma$ is as in \eqref{eqn:spectral-weight-error} and $c_d$ is recorded in Appendix~\ref{section-numerical-values}.
\end{lemma}

\begin{proof}
We have
\begin{align*}
\abs{A(\gamma_{v,k})-A(P_k)}
&\leq \frac12\int_0^{2\pi}
\abs[\Big]{\gamma_{v,k}\times(\gamma_{v,k}'-P_k')
      +(\gamma_{v,k}-P_k)\times P_k'}\,d\theta\\
&\leq \frac12\int_0^{2\pi}
\bigl(\abs{\gamma_{v,k}'-P_k'}
 +\abs{\gamma_{v,k}-P_k}\,\abs{P_k'}\bigr)\,d\theta\\
&\leq \pi\bigl(1+\nm{P'}_{L^\infty(\Xi)}\bigr)\varepsilon_\gamma.
\end{align*} 
In the first step we used the definition of $A$;
in the second step we used the triangle inequality and $\abs{\gamma_{v,k}}=1$;
in the third step we used the estimates yielding \eqref{eqn:spectral-weight-error} together with the Sobolev embedding \eqref{eqn:Sobolev-embedding} to bound $\nm{\gamma_v'-P'}_{L^\infty(\Xi)}$ and $\nm{\gamma_v-P}_{L^\infty(\Xi)}$ by $\varepsilon_\gamma$. 
Inserting the previously computed bounds for $\nm{P'}_{L^\infty(\Xi)}$ and $\varepsilon_\gamma$ proves the claim. 
\end{proof}

\begin{lemma}\label{lem:AsymNormUniqueness}
For each $i \in \{42, 36, 53\}$, for all $j \in \{1, \ldots, b\}$ with $j \neq i$ we have 
\[
\abs{A(\gamma_{v,i})}\neq\abs{A(\gamma_{v,j})}.
\]
\end{lemma}

\begin{proof}
We check
$\abs[\big]{\abs{A(P_i)}-\abs{A(P_j)}}>2c_d$
using interval arithmetic for $i,j$ as in the statement (see Appendix~\ref{section-numerical-values}).
Lemma~\ref{lem:AsymNormDifferenceBound-B} then yields 
$\abs[\big]{\abs{A(\gamma_{v,i})}-\abs{A(\gamma_{v,j})}}>0$.
\end{proof}

\paragraph{Linear independence check.} 
To prove the asymmetry of our free boundary minimal surface, we combine Lemma~\ref{lem:AsymNormUniqueness} with the following linear independence check. 

\begin{lemma}\label{lem:AsymLinIndep}
    The vectors $\{A(\gamma_{v,i}) \st i \in \{36, 42, 53\}\}$ are linearly independent in $\R^3$.
\end{lemma}

\begin{proof}
Lemma~\ref{lem:AsymNormDifferenceBound-B} gives interval enclosures for $A(\gamma_{v,i})$ for $i \in \{36,42,53\}$.
We check using interval arithmetic that 
$0\neq\det\bigl(A(\gamma_{v,36})\big\vert A(\gamma_{v,42})\big\vert A(\gamma_{v,53})\bigr)$ (see Appendix \ref{section-numerical-values}), which proves the claim.
\end{proof}
 
\begin{proposition}[Asymmetry]
\label{prop:Asym}
Suppose that \eqref{eqn:spectral-index-nullity} holds and let $\surf\vcentcolon=\Phi_{\gamma_v,c_\mu}(\Omega_{c_\mu})$.
Then the only $R\in\Ogroup(3)$ satisfying $R\surf=\surf$ is $R=\mathrm{Id}$.
\end{proposition} 

\begin{proof}
By the strong maximum principle, $\Phi_{\gamma_v,c_\mu}$ maps the interior of $\Omega_{c_\mu}$ into the open unit ball.
Therefore, $\surf\cap\Sp^2=\gamma_v(\Xi)$.
By \eqref{eqn:spectral-weight-error} and \eqref{eqn:spectral-aN}, we have
\begin{align*}
\nm{\gamma_v'}_{L^\infty(\Xi)}
&\leq\nm{P'}_{L^\infty(\Xi)}+\varepsilon_\gamma
=2\tanh(2\omega)-a_1
<\tfrac{1}{2},
\end{align*}
using the bound for $a_1$ recorded in Appendix~\ref{section-numerical-values}.
Consequently, each curve $\gamma_{v,k}$ has length less than $\pi$ and is thus contained in the open hemisphere centred at $\gamma_{v,k}(0)$.

Let $R\in\Ogroup(3)$ satisfy $R\surf=\surf$, which implies that $R\gamma_v(\Xi)=\gamma_v(\Xi)$. 
In particular, $R$ permutes the area vectors defined in \eqref{eqn:area-vector} up to sign.
By Lemma~\ref{lem:AsymNormUniqueness}, $R$ therefore preserves the curves $\gamma_{v,i}(\Sp^1)$ for $i\in\{36,42,53\}$ and hence the corresponding enclosed spherical discs $D_{\gamma_{v,i}}\subset\Sp^2$.
Since $A(\gamma_{v,i})=\pm\int_{D_{\gamma_{v,i}}}x$, it follows that $R$ fixes the vectors $A(\gamma_{v,36})$, $A(\gamma_{v,42})$ and $A(\gamma_{v,53})$.
Since by Lemma~\ref{lem:AsymLinIndep} these three vectors are linearly independent, this implies that $R=\mathrm{Id}$.
\end{proof}

\paragraph{Area bounds.}
By weak conformality of the free boundary minimal immersion $\Phi_{\gamma_v,c_\mu}$, we have
\begin{align*}
\operatorname{area}(\Phi_{\gamma_v,c_\mu})
&\vcentcolon=\operatorname{area}\bigl(\Omega_{c_\mu},\Phi_{\gamma_v,c_\mu}^*g_{\R^3}\bigr) 
=\frac{1}{2}\int_{\Omega_{c_\mu}}\abs[\big]{d\Phi_{\gamma_v,c_\mu}}_{g_{\Sp^2}}^2
=E(\gamma_v,c_\mu).
\end{align*}
We derive a rigorous enclosure for this quantity by estimating the boundary length.

\begin{proposition}
\label{proposition:area-bounds}
The free boundary minimal immersion $\Phi_{\gamma_v,c_\mu}$ from Corollary~\ref{cor:existence} satisfies
\[
\operatorname{area}(\Phi_{\gamma_v,c_\mu})\in[a_{\mathrm{low}},a_{\mathrm{high}}]
\]
for $a_{\mathrm{low}},a_{\mathrm{high}}$ recorded in Appendix~\ref{section-numerical-values}.
\end{proposition}

\begin{proof}
We first enclose the boundary length of the approximate solution $(\gamma_0,c_0)$ from Section~\ref{sec:initial}.
Writing
\[
\sigma_k\vcentcolon=\abs{\gamma_{0,k}'}
=\frac{\abs{P_k\times P_k'}}{\abs{P_k}^2},
\]
we first compute a Fourier polynomial approximation $s\colon\Xi\to\R$ as for the expression in equation \eqref{eqn:frame-Fourier} with error bounded in $L^2(\Xi)$.
As in the construction following \eqref{eqn:frame-Fourier}, we use the Lipschitz method to bound the error: 
we evaluate $\sigma-s$ at many points and use a coarse bound for the first derivative of it.
This gives a $C^0$-estimate on each interval and by weighing the $C^0$-bounds appropriately, we obtain an $L^2$-bound for $\sigma-s$, and therefore
\begin{align*}
\abs[\Big]{\operatorname{length}(\gamma_0)-2\pi\sum_{k=1}^b\widehat{s}_k(0)}
&=\abs[\bigg]{\int_\Xi(\sigma-s)\,d\theta}
\leq\sqrt{2\pi b}\,\nm{\sigma-s}_{L^2(\Xi)},
\end{align*}
which gives an enclosure for $\operatorname{length}(\gamma_0)$.
It remains to obtain an enclosure for $\operatorname{length}(\gamma_v)$ from that for $\operatorname{length}(\gamma_0)$. 
We have
\begin{align*}
\abs[\big]{\operatorname{length}(\gamma_v)-\operatorname{length}(\gamma_0)}
&\leq \int_\Xi \abs{\gamma_v'-\gamma_0'}\,d\theta\\
&\leq \sqrt{2\pi b}\,\nm{\gamma_v'-\gamma_0'}_{L^2(\Xi)}\\
&\leq \frac{\sqrt{2\pi b}}{2}\,\nm{\gamma_v-\gamma_0}_{H^2(\Xi)}\\
&\leq 6c_\infty\sqrt{2\pi b}\,\bigl(\nm{\gamma_0}_{H^2(\Xi)}+\nm{\Gamma}_{H^2(\Xi)}\bigr)\nm{v}_{H^2(\Xi)}\\[.5ex]
&\leq 6c_\infty\sqrt{2\pi b}\,\bigl(\nm{\gamma_0}_{H^2(\Xi)}+\nm{\Gamma}_{H^2(\Xi)}\bigr)r,
\end{align*}
where in the first step we used the reverse triangle inequality;
in the second step we used Cauchy--Schwarz;
in the third step we used $2\abs{n}\leq 1+n^2$ and \eqref{eqn:Sobolev-norm};
in the fourth step we used Lemma~\ref{lem:boundary-normalization};
and in the fifth step we used Corollary~\ref{cor:existence}.
The right hand side of the chain of inequalities has been bounded before, and we compute $\operatorname{length}(\gamma_0)$. 
This gives an enclosure for $\operatorname{length}(\gamma_v)$. 
By \eqref{eqn:reduced-energy-DtN}, the free boundary condition and conformality, we have (cf.~\cite[Proposition~2.4]{Li2020})
\begin{align*}
2\operatorname{area}(\Phi_{\gamma_v,c_\mu})
=2E(\gamma_v,c_\mu)
&=\int_\Xi\sk{\gamma_v,\DN_{c_\mu}\gamma_v}\,d\theta
=\int_\Xi\abs{\gamma_v'}\,d\theta
=\operatorname{length}(\gamma_v)
\end{align*}
which proves the claim. 
\end{proof}

\paragraph{Proof of Theorem~\ref{thm:main}.}
Let $v\in H^2(\Xi,\R^2)$ and $\mu\in B$ be as in Corollary~\ref{cor:existence}.
Then 
\[
\Phi_{\gamma_v,c_\mu}\colon\Omega_{c_\mu}\to\B^3
\] 
is a branched free boundary minimal immersion which is smooth up to the boundary and we set 
\(
\surf\vcentcolon=\Phi_{\gamma_v,c_\mu}(\Omega_{c_\mu}).
\)
We verify numerically that the assumption of Lemma~\ref{lem:spectral-nullity-lower-bound} is satisfied and that the estimates \eqref{eqn:spectral-aN}, \eqref{eqn:spectral-low-mode-bound} and \eqref{eqn:spectral-mN} yield constants $a_1,b_1>0$ and $m_1\geq0$ satisfying $a_1b_1>m_1^2$ (see Appendix~\ref{section-numerical-values}).
Hence, \eqref{eqn:spectral-index-nullity} holds, and Lemma~\ref{lem:embeddedness} implies that $\Phi_{\gamma_v,c_\mu}$ is an embedding.
Consequently, $\surf$ has the same topology as the domain $\Omega_{c_\mu}\subset\Sp^2$, proving~\ref{thm:main-topo}.
Claim~\ref{thm:main-sym} about the trivial symmetry group follows from Proposition~\ref{prop:Asym}. 

To prove statement \ref{thm:main-index} about the Morse index, we recall the low-mode block $\mathscr{A}$ of the linearised operator from \eqref{eqn:LABC} and claim that
\begin{align}
\label{eqn:Morse-index-numerical-criterion}
\lambda_{3b}(\mathscr{A})+K\sqrt{1+N^2}&<0,
\end{align}
where $N$ is as in \eqref{eqn:orthogonal_decomposition} and $K$ is as in Proposition~\ref{prop:nonlinear}.
Using the enclosure of $\mathscr{A}$ obtained in the derivation of \eqref{eqn:linearisation-A-bound}, we construct a matrix $W$ with $3b$ columns and check
$W^{\intercal}(\mathscr{A}+K\sqrt{17}\,\mathrm{Id})W<0$.
Since $N=4$, this gives \eqref{eqn:Morse-index-numerical-criterion}.
Since $\mathcal{R}^v$ is invertible by the proof of Corollary~\ref{cor:existence}, Proposition~\ref{prop:Lagrangian-index} and Lemma~\ref{lem:Lagrangian-index-lower-bound} imply
\begin{align*}
\ind_\area(\Phi_{\gamma_v,c_\mu})&\geq\ind(D^2L\vert_x)-3\geq3b-3.
\end{align*}
Combining this with Proposition~\ref{prop:Morse-index-upper-bound}, we conclude that
$\ind_\area(\Phi_{\gamma_v,c_\mu})=3b-3$ and $\nul_\area(\Phi_{\gamma_v,c_\mu})=3$.


\appendix

\section{Values of numerical constants}
\label{section-numerical-values}

{
\renewcommand{\arraystretch}{1.2} 
\setlength{\heavyrulewidth}{1.2pt}
\setlength{\lightrulewidth}{0.3pt}
\setlength{\cmidrulewidth}{\lightrulewidth}
\setlength{\aboverulesep}{1.5pt}
\setlength{\belowrulesep}{1.5pt}
\noindent
\begin{tabularx}{\linewidth}{@{}l l >{\raggedright\arraybackslash}X@{}}
\toprule
\textbf{Name} & \textbf{Value} & \textbf{Reference}\\
\midrule[0.8pt]
$b,M$ & $61,\ 64$ & Section~\ref{sec:initial}\\
\midrule
    & $96$ & Eqn.~\eqref{eqn:tilde_h}\\
\cmidrule(l){2-3}
$N$ & $4$ & Eqn.~\eqref{eqn:orthogonal_decomposition}\\
\cmidrule(l){2-3}
    & $1$ & Eqn.~\eqref{eqn:spectral-low-mode-approximation}\\
\midrule
$\widetilde N$ & $97$ & Eqns.~\eqref{eqn:Fourier-Fourier-DtN-approximation}, \eqref{eqn:configuration-DtN-approximation}\\
\midrule
Fixed centres & $18,29,33$ & Eqn.~\eqref{eqn:infinitesimal-mobius-on-fixed-parameters}\\
\midrule
$r_*,\epsilon_2$ & $10^{-4},\ 0.46423$ & Lemma~\ref{lem:Mobius-transversality}\\
\midrule
$\omega,C_c$ & $1.8,\ 335.171$ & Eqn.~\eqref{eqn:conformal_collar_width}, Lemma~\ref{lem:configuration-bounds}\\
\midrule
$\varepsilon$ & $9.670\times10^{-42}$ & Proposition~\ref{prop:residual-bound}\\
\midrule
$q_0,\kappa,q$ & $8621,\ 356086,\ 3.070\times10^9$ & Proposition~\ref{prop:L2-inverse-bound}, Lemma~\ref{lem:a-priori-estimate}, Proposition~\ref{prop:linear-estimate}\\
\midrule
$r,K$ & $10^{-24},\ 4.279\times10^{-11}$ & Proposition~\ref{prop:nonlinear}\\
\midrule
$\inf_\Xi\abs{\gamma_0'}$ & $\geq 0.05628$ & Eqn.~\eqref{eqn:frame}\\
\midrule
$\nm{\gamma_0}_{H^2},\ \nm{\Gamma}_{H^2}$ & $\leq 19.684,\ \leq 55.638$ & Eqn.~\eqref{eqn:spectral-weight-error}\\
\midrule
$\lambda_{\min}$ of the shifted Gram matrix & $\geq 127.109$ & Lemma~\ref{lem:spectral-nullity-lower-bound}\\
\midrule
$\varepsilon_\gamma$ & $\leq 4.537\times10^{-22}$ & Eqn.~\eqref{eqn:spectral-weight-error}, Lemma \ref{lem:AsymNormDifferenceBound-B}\\
\midrule
$\varepsilon_\spec$ & $\leq 4.537\times10^{-22}$ & Eqn.~\eqref{eqn:spectral-weight-approximation-error}\\
\midrule
$\varepsilon_N$ ($N=1$) & $\leq 0.001648955$ & Eqn.~\eqref{eqn:spectral-low-mode-error}\\
\midrule
$\lambda_5(\tilde{\mathscr S}_1)$ & $\geq 0.09125909$ & Eqn.~\eqref{eqn:spectral-low-mode-bound}\\
\midrule
$a_1,b_1,m_1$ & $1.93704,\ 0.0896101,\ 0.0921361$ & Lemma~\ref{lem:spectral-index-upper-bound}\\
\midrule
$c_d$ & $1.511\times10^{-21}$ & Lemma~\ref{lem:AsymNormDifferenceBound-B}\\
\midrule
$\min_{i\in\{36,42,53\},\,j\ne i}\abs[\big]{\abs{A(P_i)}-\abs{A(P_j)}}$ & $\geq 2.93213\times10^{-5}$ & Lemma~\ref{lem:AsymNormUniqueness}\\
\midrule
$\det\bigl(A(\gamma_{v,36}),A(\gamma_{v,42}),A(\gamma_{v,53})\bigr)$ & $\leq -9.51471\times10^{-7}$ & Lemma~\ref{lem:AsymLinIndep}\\
\midrule
$a_{\mathrm{low}}$ & $11.201930735916036032529$ & Proposition~\ref{proposition:area-bounds}\\
\midrule
$a_{\mathrm{high}}$ & $11.201930735916036032537$ & Proposition~\ref{proposition:area-bounds}\\
\bottomrule
\end{tabularx}
\par}


\clearpage
\bibliography{fbms-bibtex}

\providecommand{\bysame}{\leavevmode\hbox to3em{\hrulefill}\thinspace}
\providecommand{\MR}{\relax\ifhmode\unskip\space\fi MR }
\providecommand{\MRhref}[2]{%
  \href{http://www.ams.org/mathscinet-getitem?mr=#1}{#2}
}
\providecommand{\href}[2]{#2}
\begin{thebibliography}{10}

\bibitem{Almgren1965}
Frederick Almgren, \emph{The theory of varifolds}, Mimeographed notes,
  Princeton, 1965.

\bibitem{BarrientosGaticaMaischak2002}
Mauricio~A. Barrientos, Gabriel~N. Gatica, and Matthias Maischak, \emph{A
  posteriori error estimates for linear exterior problems via mixed-{FEM} and
  {D}t{N} mappings}, M2AN Math. Model. Numer. Anal. \textbf{36} (2002), no.~2,
  241--272.

\bibitem{Brakke1992}
Kenneth~A. Brakke, \emph{The surface evolver}, Experiment. Math. \textbf{1}
  (1992), no.~2, 141--165.

\bibitem{ButtsworthHodgkinson2026}
Timothy Buttsworth and Liam Hodgkinson, \emph{Computationally assisted proof of
  a novel {$\mathsf{O}(3)\times\mathsf{O}(10)$}-invariant {E}instein metric on
  {$S^{12}$}}, J. Lond. Math. Soc. (2) \textbf{113} (2026), no.~2, Paper No.
  e70477, 44.

\bibitem{CarlottoFranzSchulz2022}
Alessandro Carlotto, Giada Franz, and Mario Schulz, \emph{Free boundary minimal
  surfaces with connected boundary and arbitrary genus}, Camb. J. Math.
  \textbf{10} (2022), no.~4, 835--857.

\bibitem{CarlottoSchulzWiygul2025Stackings}
Alessandro Carlotto, Mario~B. Schulz, and David Wiygul, \emph{Disc stackings
  and their {M}orse index}, Adv. Nonlinear Stud. \textbf{25} (2025), no.~3,
  756--806.

\bibitem{CarlottoSchulzWiygul2025Nonuniqueness}
\bysame, \emph{Infinitely many pairs of free boundary minimal surfaces with the
  same topology and symmetry group}, Mem. Amer. Math. Soc. \textbf{313} (2025),
  no.~1591.

\bibitem{CharltonHellerHellerTraizet2024}
Steven Charlton, Lynn Heller, Sebastian Heller, and Martin Traizet,
  \emph{Minimal surfaces and alternating multiple zetas}, Preprint,
  arXiv:2407.07130.

\bibitem{ChenOsting2026}
Haoyu Chen and Braxton Osting, \emph{Computation of free boundary minimal tori
  via extremal {Steklov} eigenvalues}, Preprint, arXiv:2610.25787.

\bibitem{ChenOsting2026Zoo}
\bysame, \emph{A spectral zoo of genus-zero free boundary minimal surfaces with
  many boundary components}, Preprint, arXiv:2611.31936.

\bibitem{ChodoshGianocca}
Otis Chodosh and Matilde Gianocca, \emph{Free boundary minimal surfaces of
  genus zero}, Preprint, arXiv:2609.13054.

\bibitem{ColdingDeLellis2003}
Tobias~H. Colding and Camillo De~Lellis, \emph{The min-max construction of
  minimal surfaces}, Surveys in differential geometry, {V}ol. {VIII} ({B}oston,
  {MA}, 2002), Surv. Differ. Geom., vol.~8, Int. Press, Somerville, MA, 2003,
  pp.~75--107.

\bibitem{DaLioPigati2020}
Francesca Da~Lio and Alessandro Pigati, \emph{Free boundary minimal surfaces: a
  nonlocal approach}, Ann. Sc. Norm. Super. Pisa Cl. Sci. (5) \textbf{20}
  (2020), no.~2, 437--489.

\bibitem{DeLellisPellandini2010}
Camillo De~Lellis and Filippo Pellandini, \emph{Genus bounds for minimal
  surfaces arising from min-max constructions}, J. Reine Angew. Math.
  \textbf{644} (2010), 47--99.

\bibitem{Douglas1931}
Jesse Douglas, \emph{Solution of the problem of {P}lateau}, Trans. Amer. Math.
  Soc. \textbf{33} (1931), no.~1, 263--321.

\bibitem{DziukElliott2013}
Gerhard Dziuk and Charles~M. Elliott, \emph{Finite element methods for surface
  {PDE}s}, Acta Numer. \textbf{22} (2013), 289--396.

\bibitem{DziukHutchinson1996}
Gerhard Dziuk and John~E. Hutchinson, \emph{On the approximation of unstable
  parametric minimal surfaces}, Calc. Var. Partial Differential Equations
  \textbf{4} (1996), no.~1, 27--58. \MR{1379192}

\bibitem{Evans2010PDE}
Lawrence~C. Evans, \emph{Partial differential equations}, second ed., Graduate
  Studies in Mathematics, vol.~19, American Mathematical Society, Providence,
  RI, 2010.

\bibitem{FolhaPacardZolotareva2017}
Abigail Folha, Frank Pacard, and Tatiana Zolotareva, \emph{Free boundary
  minimal surfaces in the unit 3-ball}, Manuscripta Math. \textbf{154} (2017),
  no.~3-4, 359--409.

\bibitem{FranzKetoverSchulz}
Giada Franz, Daniel Ketover, and Mario~B. Schulz, \emph{Genus one critical
  catenoid}, J. Reine Angew. Math. (to appear), Preprint, arXiv:2409.12588.

\bibitem{FranzSchulz2026}
Giada Franz and Mario~B. Schulz, \emph{Topological control for min-max free
  boundary minimal surfaces}, J. Eur. Math. Soc. (2026), published online
  first.

\bibitem{FraserSchoen2011}
Ailana Fraser and Richard Schoen, \emph{The first {S}teklov eigenvalue,
  conformal geometry, and minimal surfaces}, Adv. Math. \textbf{226} (2011),
  no.~5, 4011--4030.

\bibitem{FraserSchoen2016}
\bysame, \emph{Sharp eigenvalue bounds and minimal surfaces in the ball},
  Invent. Math. \textbf{203} (2016), no.~3, 823--890.

\bibitem{GirouardLagace2021}
Alexandre Girouard and Jean Lagac\'e, \emph{Large {S}teklov eigenvalues via
  homogenisation on manifolds}, Invent. Math. \textbf{226} (2021), no.~3,
  1011--1056.

\bibitem{Goluzin1969}
G.~M. Goluzin, \emph{Geometric theory of functions of a complex variable},
  Translations of Mathematical Monographs, vol.~26, American Mathematical
  Society, Providence, RI, 1969.

\bibitem{Harbrecht2009}
Helmut Harbrecht, \emph{On the numerical solution of {P}lateau's problem},
  Appl. Numer. Math. \textbf{59} (2009), no.~11, 2785--2800.

\bibitem{HoangIzsakMaros2024}
Hieu~T. Hoang, Ferenc Izs\'ak, and G\'abor Maros, \emph{Method of fundamental
  solutions: new approximation results and applications}, J. Comput. Appl.
  Math. \textbf{448} (2024), Paper No. 115934, 10.

\bibitem{KapouleasLi2021}
Nikolaos Kapouleas and Martin Man-chun Li, \emph{Free boundary minimal surfaces
  in the unit three-ball via desingularization of the critical catenoid and the
  equatorial disc}, J. Reine Angew. Math. \textbf{776} (2021), 201--254.

\bibitem{KapouleasWiygul2023}
Nikolaos Kapouleas and David Wiygul, \emph{Free-boundary minimal surfaces with
  connected boundary in the $3$-ball by tripling the equatorial disc}, J.
  Differential Geom. \textbf{123} (2023), no.~2, 311--362.

\bibitem{KapouleasZou2021}
Nikolaos Kapouleas and Jiahua Zou, \emph{Free boundary minimal surfaces in the
  {E}uclidean three-ball close to the boundary}, Preprint, arXiv:2111.11308.

\bibitem{KarpukhinKokarevPolterovich2014}
Mikhail Karpukhin, Gerasim Kokarev, and Iosif Polterovich, \emph{Multiplicity
  bounds for {S}teklov eigenvalues on {R}iemannian surfaces}, Ann. Inst.
  Fourier (Grenoble) \textbf{64} (2014), no.~6, 2481--2502.

\bibitem{KarpukhinKusnerMcGrathStern}
Mikhail Karpukhin, Robert Kusner, Peter McGrath, and Daniel Stern,
  \emph{Embedded minimal surfaces in $\mathbb{S}^3$ and $\mathbb{B}^3$ via
  equivariant eigenvalue optimization}, Preprint, arXiv:2402.13121.

\bibitem{KarpukhinSchulz}
Mikhail Karpukhin and Mario~B. Schulz, \emph{Index and nullity of minimal
  surfaces via representation theory}, Preprint, arXiv:2610.03023.

\bibitem{KarpukhinStern2024}
Mikhail Karpukhin and Daniel Stern, \emph{From {S}teklov to {L}aplace: free
  boundary minimal surfaces with many boundary components}, Duke Math. J.
  \textbf{173} (2024), no.~8, 1557--1629.

\bibitem{KetoverEquivariant}
Daniel Ketover, \emph{Equivariant min-max theory}, Preprint, arXiv:1612.08692.

\bibitem{KetoverFBMS}
\bysame, \emph{Free boundary minimal surfaces of unbounded genus}, Preprint,
  arXiv:1612.08691.

\bibitem{Ketover2019}
\bysame, \emph{Genus bounds for min-max minimal surfaces}, J. Differential
  Geom. \textbf{112} (2019), no.~3, 555--590.

\bibitem{Lee2003}
John~M. Lee, \emph{Introduction to smooth manifolds}, Graduate Texts in
  Mathematics, vol. 218, Springer-Verlag, New York, 2003.

\bibitem{Li2015}
Martin Man-chun Li, \emph{A general existence theorem for embedded minimal
  surfaces with free boundary}, Comm. Pure Appl. Math. \textbf{68} (2015),
  no.~2, 286--331.

\bibitem{Li2020}
\bysame, \emph{Free boundary minimal surfaces in the unit ball: recent advances
  and open questions}, Proceedings of the {I}nternational {C}onsortium of
  {C}hinese {M}athematicians 2017, Int. Press, Boston, MA, 2020, pp.~401--435.

\bibitem{Lima2022}
Vanderson Lima, \emph{Bounds for the {M}orse index of free boundary minimal
  surfaces}, Asian J. Math. \textbf{26} (2022), no.~2, 227--252.

\bibitem{McGrathZou2024}
Peter McGrath and Jiahua Zou, \emph{On the areas of genus zero free boundary
  minimal surfaces embedded in the unit 3-ball}, J. Geom. Anal. \textbf{34}
  (2024), no.~9, Paper No. 274, 14.

\bibitem{Moore2007}
J.~D. Moore, \emph{Second variation of energy for minimal surfaces in
  {R}iemannian manifolds}, Mat. Contemp. \textbf{33} (2007), 215--243, Workshop
  on Differential Geometry.

\bibitem{OishiTanabeOgitaRump2007}
Shin'ichi Oishi, Kunio Tanabe, Takeshi Ogita, and Siegfried~M. Rump,
  \emph{Convergence of {R}ump's method for inverting arbitrarily
  ill-conditioned matrices}, J. Comput. Appl. Math. \textbf{205} (2007), no.~1,
  533--544.

\bibitem{OudetKaoOsting2021}
\'{E}douard Oudet, Chiu-Yen Kao, and Braxton Osting, \emph{Computation of free
  boundary minimal surfaces via extremal {S}teklov eigenvalue problems}, ESAIM
  Control Optim. Calc. Var. \textbf{27} (2021), Paper No. 34, 30.

\bibitem{PinkallPolthier1993}
Ulrich Pinkall and Konrad Polthier, \emph{Computing discrete minimal surfaces
  and their conjugates}, Experiment. Math. \textbf{2} (1993), no.~1, 15--36.

\bibitem{Pitts1981}
Jon~T. Pitts, \emph{Existence and regularity of minimal surfaces on
  {R}iemannian manifolds}, Mathematical Notes, vol.~27, Princeton University
  Press, Princeton, N.J.; University of Tokyo Press, Tokyo, 1981.

\bibitem{Platt2026}
Daniel Platt, \emph{Non-uniqueness and symmetries for the {Nirenberg} problem
  using computer assistance}, Preprint, arXiv:2603.29544.

\bibitem{PlattSchulzWoodCode}
Daniel Platt, Mario~B. Schulz, and Albert Wood, \emph{Free boundary minimal
  surface in the unit ball with trivial symmetry group}, 2026, URL:
  \texttt{github.com/danielplatt/certified-fbms}.

\bibitem{Polthier2005}
Konrad Polthier, \emph{Computational aspects of discrete minimal surfaces},
  Global theory of minimal surfaces, Clay Math. Proc., vol.~2, Amer. Math.
  Soc., Providence, RI, 2005, pp.~65--111. \MR{2167256}

\bibitem{ReitererTrubowitz2019}
Michael Reiterer and Eugene Trubowitz, \emph{Choptuik's critical spacetime
  exists}, Comm. Math. Phys. \textbf{368} (2019), no.~1, 143--186.

\bibitem{SakakibaraShimizu2024}
Koya Sakakibara and Yuuki Shimizu, \emph{Numerical analysis of the {P}lateau
  problem by the method of fundamental solutions}, J. Sci. Comput. \textbf{100}
  (2024), no.~1, Paper No. 2, 20.

\bibitem{SapagovasEtAl2024}
Mifodijus Sapagovas, Vytautas B{\=u}da, Saulius Maskeli{\=u}nas, Olga
  {\v{S}}tikonien{\.{e}}, and Art{\=u}ras {\v{S}}tikonas, \emph{Minimal
  surfaces and the plateau problem: numerical methods and applications},
  Informatica (Vilnius) \textbf{35} (2024), no.~2, 401--420.

\bibitem{SchulzGallery}
Mario~B. Schulz, \emph{Geometric {A}nalysis {G}allery}, 2020, URL:
  \texttt{mbschulz.github.io}.

\bibitem{Smith1982}
Francis Smith, \emph{On the existence of embedded minimal $2$-spheres in the
  $3$-sphere, endowed with an arbitrary {R}iemannian metric}, Ph.D. thesis,
  1982.

\bibitem{Struwe1988}
Michael Struwe, \emph{Plateau's problem and the calculus of variations},
  Mathematical Notes, vol.~35, Princeton University Press, Princeton, NJ, 1988.

\bibitem{Taylor2011}
Michael~E. Taylor, \emph{Partial differential equations {I}. {B}asic theory},
  second ed., Applied Mathematical Sciences, vol. 115, Springer, New York,
  2011.

\bibitem{Thomson1904}
J.~J. Thomson, \emph{{XXIV. On the structure of the atom: an investigation of
  the stability and periods of oscillation of a number of corpuscles arranged
  at equal intervals around the circumference of a circle; with application of
  the results to the theory of atomic structure}}, The {L}ondon, {E}dinburgh,
  and {D}ublin {P}hilosophical {M}agazine and {J}ournal of {S}cience \textbf{7}
  (1904), no.~39, 237--265.

\bibitem{Tran}
Hung Tran, \emph{The {M}orse index of prismatic free boundary minimal surfaces
  in the unit ball}, Preprint, arXiv:2610.02078.

\end{thebibliography}

\printaddress

\end{document}